\documentclass[12pt]{article}
\usepackage{amsmath}
\usepackage{amsfonts}
\usepackage{graphics}
\usepackage{amsthm}
\usepackage{pictex}
\usepackage{makeidx}
\numberwithin{equation}{section}
\numberwithin{figure}{section}
\newcommand{\HP}{Her\-mite-Pad\'e}
\newcommand{\RH}{Riemann-Hilbert}
\newcommand{\PHI}{\mathbf{\Phi}}
\newcommand{\Cap}{\mathop\textrm{cap}}
\newcommand{\supp}{\mathop\textup{supp}}
\newcommand{\BVP}{boundary value problem}
\newcommand{\res}{\mathop\textup{res}}
\renewcommand{\O}{\mathcal{O}}
\renewcommand{\Re}{\mathop\textrm{Re }}
\newcommand{\Ai}{\textup{Ai}}

\newtheorem{theorem}{Theorem}[section]
\newtheorem{proposition}[theorem]{Proposition}

\newtheorem{Definition}[theorem]{Definition}
\newenvironment{definition}{\begin{Definition}\rm}{\end{Definition}}
\newtheorem{Remark}[theorem]{Remark}
\newenvironment{remark}{\begin{Remark}\rm}{\end{Remark}}

\makeindex

\begin{document}
\title{Asymptotics of Hermite-Pad\'e rational approximants for two
    analytic functions with separated pairs of branch points
    (case of genus $0$)}
\author{A.~I.~Aptekarev \\ {\normalsize Keldysh Institute of Applied Mathematics,} \\
          {\normalsize Russian Academy of Sciences, Russia} \and
       A.~B.~J.~Kuijlaars, W.~Van Assche \\
        {\normalsize Department of Mathematics} \\
        {\normalsize Katholieke Universiteit Leuven, Belgium}}

\maketitle

\begin{abstract}
We investigate the asymptotic behavior for type II Hermite-Pad\'e
approximation to two functions, where each function has two branch
points and the pairs of branch points are separated. We give a
classification of the cases such that the limiting counting
measures for the poles of the \HP\ approximants are described by an
algebraic function $h$ of order $3$ and genus $0$. This situation
gives rise to a vector-potential equilibrium problem for measures
$\lambda$, $\mu_1$, and $\mu_2$, and the poles of the common
denominator are asymptotically distributed like $\lambda/2$. We
also work out the strong asymptotics for the corresponding \HP\
approximants by using a $3\times 3$ Riemann-Hilbert problem that
characterizes this Hermite-Pad\'e approximation problem.
\end{abstract}
\tableofcontents


\section{Introduction} \label{sec:1}
\subsection{Definition of \HP{} approximants and general statement of the problem}
\label{sec:1.1}
Let $\vec{f} = (f_1,\ldots,f_p)$ be a vector of Laurent series near infinity
\begin{equation}  \label{eq:1.1}
   f_j(z) = \sum_{k=0}^\infty \frac{f_{j,k}}{z^k},\qquad j=1,\ldots,p.
\end{equation}
The \textit{\HP{} rational approximants} (of type II)
\[  \pi_{\vec{n}} = \left( \frac{Q_{\vec{n}}^{(1)}}{P_{\vec{n}}}, \ldots,
    \frac{Q_{\vec{n}}^{(p)}}{P_{\vec{n}}} \right) \]
for the vector $\vec{f}$ and multi-index $\vec{n} = (n_1, \ldots, n_p) \in \mathbb{N}^p$
are defined by
\[  \deg P_{\vec{n}} \leq |\vec{n}| = n_1 + \cdots + n_p, \]
\begin{equation}   \label{eq:1.2}
    f_j(z)P_{\vec{n}}(z) - Q_{\vec{n}}^{(j)}(z) =:
    R_{\vec{n}}^{(j)}(z) = \O\left(\frac{1}{z^{n_j+1}}\right), \qquad z \to \infty,
\end{equation}
where the $Q_{\vec{n}}^{(j)}$ are polynomials,
for $j=1,\ldots, p$.
This definition is equivalent to a homogeneous linear system of equations
for the coefficients of the polynomial $P_{\vec{n}}$.
This system always has a solution, but the solution is not necessarily unique.
In the case of uniqueness (up to a multiplicative constant) and in
case any non-trivial solution has full degree $\vec{n}$,
the multi-index $\vec{n}$ is called \textit{normal} and the polynomial $P_{\vec{n}}$ can
be normalized as monic
\[  P_{\vec{n}}(z) = \prod_{k=1}^{|\vec{n}|} (z-z_{k,\vec{n}}). \]
The \HP{} approximants $\pi_{\vec{n}}$ provide the best local
(near infinity) simultaneous rational approximation of the vector
$(f_1,\ldots,f_p)$ of Laurent series (\ref{eq:1.1}). The
construction (\ref{eq:1.2}) was introduced by Hermite \cite{1} in
connection with his proof of the transcendence of $e$. See the
papers \cite{2,3,67,5,6,66} and the monograph \cite{4} for more details.

In this paper we study  the asymptotic behavior of the \textit{diagonal}
\HP{} approximants ($\vec{n}=(n,n)$) for two functions $f_1$ and $f_2$ with
branch points at the points $A_1=\{a_1,b_1\}$ and $A_2=\{a_2,b_2\}$ respectively.
We say that
\begin{equation} \label{eq:1.3}
   f_j \in \mathcal{A}(\overline{\mathbb{C}} \setminus A_j), \quad  A_j = \{ a_j,b_j\}
\end{equation}
if the Laurent expansion (\ref{eq:1.1}) is convergent in a neighborhood of
infinity and has an analytic continuation along any path in $\mathbb C \setminus A_j$.
A typical example is the function
\index{Functions!f1f2@$f_1$, $f_2$}
\index{Parameters!a1a2@$a_1$, $a_2$}
\index{Parameters!b1b2@$b_1$, $b_2$}
\index{Point sets!A1A2@$A_1$, $A_2$}
\[ f_j(z) = \log((z-a_j)/(z-b_j)). \]
Actually the class of functions we allow is larger and we give a more precise
definition in Section 2.
We shall assume that the pairs of branch points do not coincide, i.e., $A_1 \neq A_2$, although
they might have non-empty intersection. Our goal is
\begin{itemize}
\item To determine the limiting distribution of the zeros of the common denominator $P_{\vec{n}}$,
which are the poles of the \HP{} approximants:
\begin{equation}  \label{eq:1.4}
   \nu_{P_{\vec{n}}} = \frac{1}{2n} \sum_{k=1}^{2n} \delta(z-z_{k,\vec{n}}) \stackrel{*}{\to} \ ? \qquad n \to \infty.
\end{equation}
\index{Measures!nuP@$\nu_{P_{\vec{n}}}$}
\item To obtain asymptotic formulas for the \HP{} polynomials $P_{\vec{n}}$ and the
functions $R_{\vec{n}}^{(j)}$.
\item To prove convergence theorems for the \HP{} approximants
\[  \lim_{n \to \infty} \frac{Q_{\vec{n}}^{(j)}(z)}{P_{\vec{n}}(z)} = f_j(z), \quad z \in \Omega_j , \ j=1,2, \]
and to describe the domains of convergence $\Omega_1$ and $\Omega_2$, depending on the location of the points
$a_1,a_2,b_1,b_2 \in \mathbb{C}$.
\end{itemize}
For the situation under consideration ($p=2$ and $\vec{n}=(n,n)$) we will use the notation
\[   \pi_{\vec{n}} = \pi_n, \quad P_{\vec{n}} = P_n, \quad Q_{\vec{n}}^{(j)} = Q_n^{(j)}, \quad R_{\vec{n}}^{(j)} =R_n^{(j)}, \qquad j=1,2, \]
and we assume that $P_{\vec{n}}$ is monic. \index{Functions!pin@$\pi_n$}\index{Functions!Pn@$P_n$}\index{Functions!Qnj@$Q_n^{(1)}$, $Q_n^{(2)}$}\index{Functions!Rnj@$R_n^{(1)}$, $R_n^{(2)}$}
The rigorous definition of the classes of functions under consideration and the
statement of the results of this paper will be presented in the next section. In the
following subsections of this introduction we give a brief historical review of the
analytic aspects of the Pad\'e and \HP{} approximants in order to introduce some basic
notions and problems. Then we conclude the introduction with a general description of
the results in this paper.

\subsection{Pad\'e approximants (analytic aspect). Motivation for \HP{} analysis} \label{sec:1.2}
The special case of \HP{} approximation (\ref{eq:1.2}) for $p=1$ (i.e., the best local rational
approximation to one function near infinity) corresponds to \textit{Pad\'e approximation}.
When the coefficients of the Laurent series (\ref{eq:1.1}) are the moments of a positive measure $\mu$ supported on
$\mathbb{R}$
\begin{equation}   \label{eq:1.5}
   f(z) = \sum_{k=0}^\infty \frac{c_k}{z^{k+1}} = \int_{\mathbb{R}} \frac{d\mu(x)}{z-x}, \qquad
   c_k = \int_{\mathbb{R}} x^k\, d\mu(x),
\end{equation}
then the denominators $P_n$ of the Pad\'e approximants are polynomials \textit{orthogonal} to the powers $\leq n-1$
with respect to $\mu$:
\begin{equation}  \label{eq:1.6}
   \int_{\mathbb{R}} P_n(x) x^k\, d\mu(x) = 0, \qquad k=0,1,\ldots,n-1.
\end{equation}
In the general case when we are dealing with general coefficients in (\ref{eq:1.1}) the orthogonality
relations become \textit{non-Hermitian} or \textit{complex}.

The analytic theory of Pad\'e approximants for the \textit{real case} (\ref{eq:1.5}) is based on remarkable classical results.
These include Markov's theorem \cite{7} on the locally uniform convergence outside the convex hull $\Delta$ of
the compact support of $\mu$:
\[   \lim_{n \to \infty} \pi_n(z) = \int \frac{d\mu(x)}{z-x},
\qquad z \in  \overline{\mathbb{C}} \setminus \Delta,  \]
where $K$ is compact, and the theory of Bernstein and Szeg\H{o} \cite{8,9} on the strong
asymptotics of orthogonal polynomials satisfying (\ref{eq:1.6}).

The analytic theory for the \textit{complex case}, particularly for
functions $f \in \mathcal{A}(\overline{\mathbb{C}} \setminus A)$, where $A$ is a finite set
of points in $\mathbb{C}$, started to be developed not so long ago. Nuttall \cite{10} has put
forward the important relation between the
\textit{maximal} domain of analyticity for the analytic function $f$ and the domain of
convergence of the diagonal Pad\'e approximants. The Pad\'e approximants, which are single
valued rational functions, approximate a holomorphic branch of the analytic function in the
domain of their convergence. At the same time most of the poles of the rational approximants
tend to the boundary of the domain of convergence and the support of their limiting distribution
models the cuts which make the function $f$ single valued.
Nuttall has conjectured (and proved for some important special cases \cite{10,11,12,13}) that
these cuts have a minimal logarithmic capacity among all cuts converting the function to a
single valued branch. Thus the domain of convergence corresponds to the \textit{maximal}
(in the sense of \textit{minimal} boundary)
domain of holomorphicity for the analytic function $f \in \mathcal{A}(\overline{\mathbb{C}} \setminus A)$.
The complete proof of Nuttall's conjecture (even in a more general setting where the set $A$
has capacity $0$) was obtained by Stahl. In a series of papers \cite{14,15} he proved:
\begin{itemize}
\item Existence of a domain $\Omega^*$ such that $f$ is holomorphic in $\Omega^*$
($f \in H(\Omega^*)$) and the boundary
$\Delta = \partial \Omega^*$ has the property
\[     \Cap \Delta = \min_{\partial \Omega: f \in H(\Omega)} \Cap \partial \Omega . \]
\item The weak limit of the pole counting measure (\ref{eq:1.4})
\[   \nu_{P_n} \stackrel{*}{\to} \lambda, \qquad \supp \lambda = \Delta, \]
and weak asymptotics for the denominators of the Pad\'e approximants
\[    \lim_{n \to \infty} \log |P_n(z)| = - V^\lambda(z), \qquad z \in \Omega^*, \]
where $V^\lambda$ is the logarithmic potential \index{Functions!Vlambda@$V^{\lambda}$}
\[   V^\lambda(z) = \int \log \frac{1}{|z-t|}\, d\lambda(t)  \]
of an extremal measure $\lambda$ minimizing the energy functional
\[  I(\lambda) = \iint \log \frac{1}{|x-t|}\, d\lambda(t)\,d\lambda(x) =
\min_{\supp \mu \subset \Delta,\ \mu(\Delta)=1} I(\mu) \]
among all probability measures on $\Delta$.
\end{itemize}
The extremal measure possesses the equilibrium properties
\begin{equation}   \label{eq:1.8}
  \begin{cases}  V^\lambda = \textrm{const.} \\[10pt]
    \displaystyle \frac{\partial V^\lambda}{\partial n_+}
    = \frac{\partial V^\lambda}{\partial n_-}
  \end{cases}
            \textrm{a.e. on } \Delta,
\end{equation}
where $\frac{\partial}{\partial n_{\pm}}$ denotes the normal derivatives
on $\Delta$ (which is a finite union of analytic arcs). The conditions (\ref{eq:1.8})
characterize the measure $\lambda$ and its support $\Delta$.
These results lead to Stahl's main convergence theorem
\[   \lim_{n \to \infty} \pi_n(z) = f(z), \qquad z \in \Omega^* \ \textrm{in capacity}, \]
where the convergence in capacity is defined in the same manner as convergence in measure.

The generalization of these notions and results from Pad\'e approximation to
\HP{} approximation is a very difficult and challenging problem.
As we will see, the geometry of the domains of convergence and the extremal compact sets
where the poles of the \HP{} approximants accumulate is much more diverse and complicated.
The analytic techniques capable of proving the asymptotics and convergence results require a
significant development in comparison with the methods appropriate for Pad\'e approximation and
orthogonal polynomials. These circumstances give a good motivation and direction for the development
of the analytic aspect of \HP{} approximation.

\subsection{Short survey of asymptotic results for general classes of \HP{} approximants}
\label{sec:1.3}
Perhaps one of the first results on the asymptotics of \HP{} polynomials was obtained
by Kalyagin \cite{16}. For a special class of multiple orthogonal polynomials
$P_{n,n}^{\alpha,\beta,\gamma}(x)=x^{2n}+\cdots$ generalizing the Jacobi polynomials:
\[   \int_{\Delta_j} P_{n,n}^{\alpha,\beta,\gamma}(x) x^k w(x)\, dx = 0,
\qquad k=0,1,\ldots,n-1,\ j=1,2, \]
with $\Delta_1=[-1,0]$, $\Delta_2=[0,1]$ and weight function
\[   w(x) = |x+1|^\alpha |x|^\gamma |x-1|^\beta, \qquad x \in \Delta_j,\ j=1,2, \]
and for the corresponding functions of the second kind
\begin{equation}  \label{eq:1.9}
   R_{n}^{(j)}(z) = \frac{1}{2\pi i} \int_{\Delta_j} \frac{P_{n,n}^{\alpha,\beta,\gamma}(x)}{x-z} w(x)\, dx, \qquad j=1,2,
\end{equation}
he proved the strong asymptotics (Szeg\H{o} type asymptotics) as $n \to \infty$:
\begin{eqnarray*}
    P_{n,n}^{\alpha,\beta,\gamma}(z) & = & C_0^{-n} \Phi_0^{-n}(z) [ F_0(z) + o(1)],
       \qquad z \in \Omega_0=\overline{\mathbb{C}}\setminus (\Delta_1 \cup \Delta_2), \\
    R_{n}^{(j)}(z) & = & C_j^{-n} \Phi_j^{-n}(z) [ F_j(z) + o(1)],
       \qquad z \in \Omega_j=\overline{\mathbb{C}}\setminus \Delta_j,\ \ j=1,2,
\end{eqnarray*}
and the corresponding formulas on the intervals $\Delta_j$ $(j=1,2)$, where
the convergence is uniform on compact subsets of the indicated domains. Here
\[    F_\ell, 1/F_\ell \in H(\Omega_\ell), \qquad \ell=0,1,2,   \]
depend on $\alpha,\beta,\gamma$ and they are analogs of the Szeg\H{o} function.
The main terms $\Phi_0,\Phi_1,\Phi_2$
of the asymptotics are the single valued branches of an \textit{algebraic function}:
\[     \Phi^3(z) + 3\Phi^2(z) + \left(3- \frac{27}{4} z^2\right) \Phi(z) -1 = 0, \]
where the branch points are at $A=\{-1,0,1\}$. We note that $\Phi$ is independent of $\alpha,\beta,\gamma$
and depends on the supports $\Delta_1$, $\Delta_2$ of the weight, i.e., on the set $A$.
The function $\Phi$ is a rational function on the \textit{three-sheeted Riemann surface}
\begin{equation}  \label{eq:1.10}
    \mathfrak{R}(A) = \overline{ \mathfrak{R}_0 \cup \mathfrak{R}_1 \cup \mathfrak{R}_2}
\end{equation}
obtained by glueing the sheets $\mathfrak{R}_j = \overline{\mathbb{C}}\setminus \Delta_j$ $(j=1,2)$ to the sheet
$\mathfrak{R}_0 = \overline{\mathbb{C}}\setminus (\Delta_1 \cup \Delta_2)$ so that the upper and lower sides of
the cuts on two neighboring sheets are identified. The function $\Phi$ is defined (up to a multiplicative constant)
by its \textit{divisor} (set of poles and zeros),
\begin{equation} \label{eq:1.11}
    \Phi(z) = \begin{cases}
              \displaystyle \frac{1}{C_0z^2} + \cdots, & z \to \infty^{(0)}, \\[15pt]
              \displaystyle \frac{z}{C_j} + \cdots, & z \to \infty^{(j)},\ j=1,2,
              \end{cases}
\end{equation}
and the normalization is chosen so that $C_0C_1C_2=1$. The Riemann surface $\mathfrak{R}$ and rational functions on it play an
important role for the asymptotic analysis of \HP{} approximants.

Another important notion of \textit{vector potential equilibrium} was introduced by Gonchar and
Rakhmanov in \cite{17} (see also \cite{18,19}). Let $\{\Delta_1,\ldots,\Delta_p\}$ be a
collection of compact sets in $\mathbb{C}$ and let $D=(d_{i,j})_{i,j=1}^p$ be a real symmetric
nonsingular positive definite matrix. An additional condition on $D$ to be compatible with
$(\Delta_1,\ldots,\Delta_p)$ is that $d_{i,j} \geq 0$ whenever
$\Delta_i \cap \Delta_j \neq \emptyset$. For a vector of measures
\[  \vec{\mu} = (\mu_1,\ldots,\mu_p), \quad \supp \mu_j \subset \Delta_j,\quad j=1,\ldots,p, \]
the energy functional $I(\vec{\mu})$ is defined as
\begin{equation}  \label{eq:1.12}
  I(\vec{\mu}) = \sum_{i=1}^p\sum_{j=1}^p d_{i,j} I(\mu_i,\mu_j),
\end{equation}
where $I(\mu_i,\mu_j)$ is the mutual energy of two scalar measures
\[    I(\mu_i,\mu_j) = \int_{\Delta_i} \int_{\Delta_j} \log \frac{1}{|x-t|} \,
    d\mu_i(x)\,d\mu_j(t) . \]
The extremal vector measure $\vec{\lambda}$,
minimizing the energy functional (\ref{eq:1.12}) among all $\vec{\mu}$
where all $\mu_j$ are probability measures 
possesses the equilibrium properties
\begin{equation}  \label{eq:1.13}
    U_j^{\vec{\lambda}}(x) = \sum_{i=1}^p d_{ij} V^{\lambda_i}(x)
    \begin{cases}
     = \kappa_j, & x \in \supp \lambda_j = \Delta_j^*, \\
     \geq \kappa_j, & x \in \Delta_j \setminus \Delta_j^*.
    \end{cases}
\end{equation}
Here the vector $\vec{U}^{\vec{\lambda}} = (U_1^{\vec{\lambda}},\ldots,U_p^{\vec{\lambda}})$
is called the \textit{vector potential of the vector valued measure $\vec{\lambda}$ with
respect to the interaction matrix $D$}.

In the paper \cite{17} Gonchar and Rakhmanov investigated the \HP{} approximants
(\ref{eq:1.2}) for the system of Markov-type functions (\ref{eq:1.5})
\begin{equation}  \label{eq:1.14}
  f_j(z) = \hat{\mu}_j(z) = \int_{\Delta_j} \frac{d\mu_j(x)}{z-x}, \qquad \mu_j' > 0
  \textrm{ on } \Delta_j \subset \mathbb{R}, \quad j=1,\ldots,p,
\end{equation}
where $\Delta_j$ $(j=1,\ldots,p)$ are non-overlapping intervals
\begin{equation}  \label{eq:1.15}
   \stackrel{\circ}{\Delta}_i \cap \stackrel{\circ}{\Delta}_j = \emptyset,\qquad i \neq j,
\end{equation}
where $\stackrel{\circ}{\Delta}$ denotes the interior of the interval $\Delta$.
They proved the weak asymptotics for the common denominator $P_{\vec{n}}$ where
$\vec{n} = (n, \ldots, n)$ of the \HP{} approximants as $n \to \infty$:
\[    \nu_{P_{\vec{n}}} \stackrel{*}{\to} \frac{1}{p} \sum_{j=1}^p \lambda_j, \]
\[   \frac1n \log |P_{\vec{n}}(z)| \to - \sum_{j=1}^p V^{\lambda_j}(z),
     \qquad z \in \overline{\mathbb{C}} \setminus \sum_{j=1}^p \Delta_j^*, \]
where $\lambda_j$ $(j=1,\ldots,p)$ are the components of the extremal (equilibrium) vector
measure with matrix of interaction
\begin{equation}  \label{eq:1.16}
       d_{j,j} = 2, \quad d_{i,j} = 1, \qquad 1 \leq i\neq j \leq p.
\end{equation}
The potentials of the components of the equilibrium measure (after normalization) can be
harmonically continued through the intervals $\Delta_j^*$ forming the Riemann surface
(as in (\ref{eq:1.10}) where the index $j$ runs from $1$ to $p$ and the cuts on the sheets
join the endpoints of the supports $\Delta_j^*$). This fact had been noticed in \cite{20}.
Thus the notion of rational function (\ref{eq:1.11}) on $\mathfrak{R}$ and vector
equilibrium problem (\ref{eq:1.13}) are equivalent and they are related by
\[  \exp \left( - V^{\lambda_j}(z) \right) = \left| \frac{\Phi_j(z)}{C_j} \right|,
    \qquad j=1,\ldots,p,
    \quad \exp \left( \sum_{j=1}^p V^{\lambda_j}(z) \right) =
    \left| \frac{\Phi_0(z)}{C_0}\right|, \]
where the normalization constants $C_1,\ldots,C_p$ and the equilibrium constants
$\kappa_1,\ldots,\kappa_p$ are connected by a linear system of equations.
The following convergence theorem was proved in \cite{17}:
\[ \lim_{n \to \infty} \frac{Q_{\vec{n}}^{(j)}(z)}{P_{\vec{n}}(z)}
     = \begin{cases}   \displaystyle \int_{\Delta_j} \frac{d\mu_j(x)}{z-x},
                    & z \in \Omega_j^*, \\[12pt]
                       \infty, & z \in  (\overline{\mathbb{C}} \setminus \Delta))
                       \setminus \Omega_j^*,
       \end{cases} \qquad j=1,\ldots,p,  \]
where
\[   \Omega_j^* = \{ z : |\Phi_j(z)| > |\Phi_0(z)|\}. \]
Note that, in view of (\ref{eq:1.11}), $\Omega_j^* \neq \emptyset$. The existence of
the non-empty domain $\overline{\mathbb{C}} \setminus \overline{\Omega_j^*}$ depends
on the input geometry, i.e., on the size and the
location of the $\Delta_j$ $(j=1,\ldots,p)$.

Thus the results in \cite{17} show that there are two new phenomena for the asymptotic
behavior of \HP{} approximants as compared to Pad\'e approximants:
\begin{enumerate}
\item The components $\Delta_j^*$  of the support of the pole counting measure
do not correspond to cuts making
the functions $f_j$ $(j=1,\ldots,p)$ holomorphic. We call this the \textit{pushing effect}:
it might happen that $\Delta_j \setminus \Delta_j^* \neq \emptyset$ for some
$j \in \{1,\ldots,p\}$.
\item The appearance of domains of divergence inside the domain of holomorphicity of $f_j$:
it might happen that $(\overline{\mathbb{C}} \setminus \Delta_j) \setminus \Omega_j^* \neq
\emptyset$ for some $j \in \{1,\ldots,p\}$.
\end{enumerate}
These two phenomena are related by
\[  \Delta_j \setminus \Delta_j^* \neq \emptyset \begin{array}{c} \Rightarrow \\
    \not \Leftarrow \end{array}
                  (\overline{\mathbb{C}} \setminus \Delta_j) \setminus \Omega_j^* \neq \emptyset . \]
The system (\ref{eq:1.14})--(\ref{eq:1.15}) was introduced in 1919 by Angelesco \cite{21} as
a system for which all the multi-indices of the \HP{} approximants are normal, and this system
was later rediscovered in \cite{22}. Another system of Markov-type functions (\ref{eq:1.5})
with normal diagonal multi-indices for the \HP{} approximants was introduced by Nikishin in
\cite{23}. A system (\ref{eq:1.14}) is a Nikishin system of order $p$ if $\Delta_j = \Delta$
for $j=1,\ldots,p$ and $d\mu_j/d\mu_1$ $(j=2,\ldots,p)$ have analytic continuation from
$\Delta$ and form a Nikishin system of order $p-1$ with respect to another interval $F$ for
which $F \cap \Delta = \emptyset$. The asymptotic behavior of the denominators of the \HP{}
approximants for a Nikishin system is similar to the behavior of Pad\'e approximants in
the sense that the two phenomena for Angelesco systems do not appear (see \cite{23,24,25,26,27}). However, for a Nikishin system a new effect appears for the functions of the second kind (\ref{eq:1.9})
\begin{equation}  \label{eq:1.17}
   R_{\vec{n}}^{(j)}(z) = P_{\vec{n}}(z) f_j(z) - Q_{\vec{n}}^{(j)}(z), \qquad j=2,\ldots,p.
\end{equation}
They have extra zeros which accumulate on the interval $F$ and are dense on this interval
as $n \to \infty$. These are extra interpolation points for the \HP{} approximants, apart
from the interpolation condition at $\infty$.

To conclude this survey of the results for the \textit{real case}, i.e.,
\HP{} approximants for a vector of
Markov type functions, we mention the recent paper \cite{19} on mixed
Angelesco-Nikishin systems defined by a
graph-tree, and the papers \cite{28,29} where the strong Szeg\H{o}-type asymptotics
of the \HP{} polynomials
for Angelesco systems and Nikishin systems was obtained.

The analytic theory of \HP{} approximants for the \textit{complex case} has been initiated
by Nuttall. In the two pioneering papers \cite{30,31} of 1981 he obtained some asymptotic
formulas for \HP{} approximants to functions with
separated complex branch points \cite{30} (a complex analog of an Angelesco system) and to
functions meromorphic on the
same Riemann surface \cite{31} (i.e., functions with the same set of branch points, like a
Nikishin system for the real case). The results of \cite{30} were verified by some heuristic
considerations and numerical experiments, and the paper \cite{31} contains rigorous theorems.
In his fundamental work \cite{3}  of 1984, Nuttall made an attempt to formulate a general
conjecture about the asymptotic behavior, as $n \to \infty$, of the diagonal
$\vec{n} = (n,n,\ldots,n)$ \HP{} polynomials. On the basis of his conjecture lies a
$(p+1)$-sheeted Riemann surface like (\ref{eq:1.10})) which
depends on the set of $p$ functions which are being approximated. He showed how to
determine this Riemann surface for some special classes of functions, but the general
case was left as an open problem. Nevertheless, assuming the
existence of the appropriate Riemann surface $\mathfrak{R}$, he conjectured that the strong
asymptotics can be described by solutions of some boundary value problem on $\mathfrak{R}$.
For the main term of the asymptotics one would have
\begin{equation}  \label{eq:1.18}
  |P_{\vec{n}}(z)|^{1/n} \to  |\Phi^{-1}(z)| = \exp \left( -\Re G(z) \right),
  \qquad  n \to \infty,
\end{equation}
where $G$ is an Abelian integral of the third kind with logarithmic poles at
$\infty^{(\ell)}$, $(\ell=0,1,\ldots,p)$ with residues
\begin{equation}  \label{eq:1.19}
   G(z) = \begin{cases}
     -p \log z + \O(1), & z \to \infty^{(0)}, \\
     \log z + \O(1), & z \to \infty^{(j)},\qquad j=1,\ldots,p,
         \end{cases}
\end{equation}
and elsewhere $G$ is analytic in the local variable. If the genus of $\mathfrak{R}$ is
greater than zero, then an additional condition on $G$ is imposed: all the periods of $G$
are purely imaginary. Such a function is unique up to an additive constant and $\Re G(z)$
is a single valued function on $\mathfrak{R}$ (see, for example, \cite{27}).
We note that the condition $\textrm{genus}(\mathfrak{R})=0$ implies single-valuedness of
$\Phi$ in (\ref{eq:1.18}), therefore $\Phi$ is a rational function on $\mathfrak{R}$,
uniquely (up to a multiplicative normalization) defined on $\mathfrak{R}$ by its
divisor (\ref{eq:1.11}).

After Nuttall's results there were practically no other rigorous results for the
complex case of \HP{}
approximants\footnote{except for the special case of \HP\ approximation to the $e^z$
\cite{55,56,57,51,52,49,58,59,60,61,62}.}.
One of the reasons was the absence of suitable techniques for the analysis of strong
asymptotics of non-Hermitian (complex) orthogonal polynomials which could be adopted
to \HP{} approximation. There was recently substantial progress in proving new results
for the strong asymptotics of orthogonal polynomials
by means of a matrix-valued Riemann-Hilbert method. The method is based on the reformulation of
the definition
(\ref{eq:1.6}) of the orthogonal polynomials in terms of a $2\times 2$ matrix-valued \RH{} problem
(due to Fokas, Its, and Kitaev \cite{34,33}) and the steepest descent analysis of this
\RH{} problem for $n\to \infty$ (due to Deift and Zhou \cite{35}). This method was
initially designed to study asymptotics for integrable PDEs and was later applied to prove
asymptotic results for polynomials orthogonal on the real axis with respect to real valued analytic
weights, including varying weights (depending on $n$) \cite{68,38,39,36,37,40,63,41}
and related questions from random matrix theory. It has later been noticed \cite{43,42,44}
that the method also works
for the non-Hermitian orthogonality in the complex plane with respect to complex valued
weights. In \cite{45}
multiple orthogonality for \HP{} polynomials was reformulated in terms of a
$(p+1)\times (p+1)$ matrix-valued \RH{} problem.

In this survey subsection we introduced the main players for the asymptotics
of the complex case of \HP{} approximation
in the historical order of their appearance on the scene. They are
\begin{itemize}
\item An appropriate $(p+1)$-sheeted Riemann surface $\mathfrak{R}$.
\item Vector potential equilibrium problems.
\item A standard Abelian integral (\ref{eq:1.19}) and rational functions on $\mathfrak{R}$.
\item Analysis of a $(p+1)\times (p+1)$ matrix-valued \RH{} problem as a method for proving the asymptotic results.
\end{itemize}

Concluding this survey of results which have important and direct influence on
this paper, we would like to mention the papers \cite{43,48,46,47,64,65,49,53,54}
which have been written during the work
on the present paper and in which some of the ideas and methods elaborated
here have already been implemented.

\subsection{General description of the limiting behavior of \HP{} approximants} \label{sec:1.4}
The analysis of numerical computations of the zeros of the polynomials $P_{\vec{n}}$,
which are the common denominators for the \HP{} approximants to functions of the
class (\ref{eq:1.3}), shows that the support for the limiting zero distribution has a
different geometry, depending on the position of the branch points.
In Figures \ref{fig:1}--\ref{fig:5} we present the results of the computations for the
functions
\[    f_j(z) = \log \left( \frac{z-a_j}{z-b_j} \right), \qquad j=1,2.  \]
There are several typical patterns for the support. We describe those that are related
to a Riemann surface of genus zero.
\begin{description}
\item[\boldmath Case I.]
When the pairs $\{a_1,b_1\}$ and $\{a_2,b_2\}$ are
`far away' from each other, the zeros of $P_{\vec{n}}$ accumulate on two disjoint
arcs $\Delta_1$ and $\Delta_2$, which are joining the branch points $a_1$ and $b_1$
for the function $f_1$ and the branch points $a_2$ and $b_2$ of the function $f_2$,
respectively (Figure \ref{fig:1}). Each of the arcs $\Delta_1$
and $\Delta_2$  accumulates half of the zeros of $P_{\vec{n}}$.
\begin{figure}[ht]
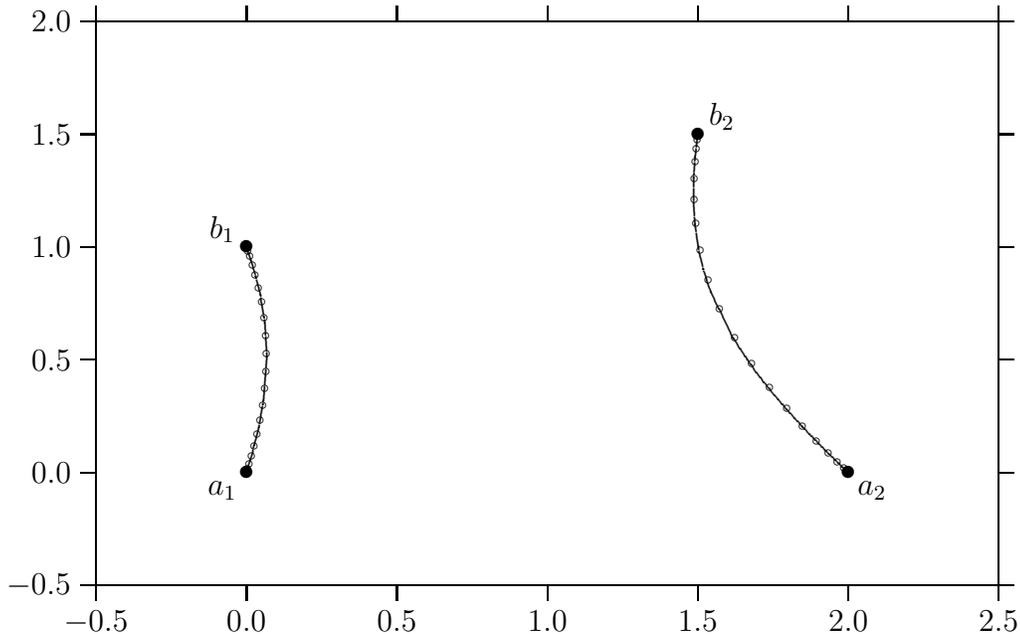

\centering
\mbox{
\input trajectoriesI.tex }
\caption{$a_1=0$, $b_1=i$, $a_2=2$, $b_2=(3+3i)/2$ (case I). The
zeros of $P_{n,n}$ accumulate on two disjoint arcs $\Delta_1$ and $\Delta_2$ that
connect the branch points. The zeros of $P_{20,20}$ are indicated by $\circ$.}
\label{fig:1}
\end{figure}
\item[\boldmath Case II.] When the pairs $\{a_1,b_1\}$ and $\{a_2,b_2\}$ are
closer to each other, then it may happen that the zeros of $P_{\vec{n}}$
accumulate on a set $\Delta$
that connects all four branch points $a_1$, $b_1$, $a_2$ and $b_2$ as
shown in Figure~\ref{fig:2}. In this case the support can not be split
into two separate pieces of equal mass.
\begin{figure}[ht]
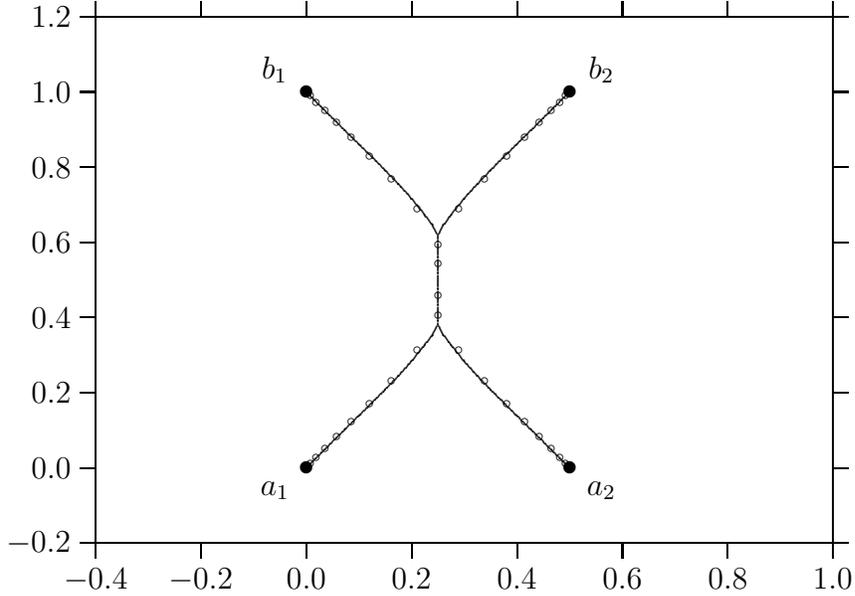

\centering
\mbox{
\input trajectoriesII.tex }
\caption{$a_1=i$, $b_1=0$, $a_2=1/2+i$, $b_2=1/2$,  (case II).
The zeros of $P_{n,n}$ accumulate on a contour $\Delta$ that connects
all four branch points. The zeros of $P_{20,20}$ are indicated by $\circ$.}
\label{fig:2}
\end{figure}

\item[\boldmath Case III.] For certain configurations of pairs
$\{a_1,b_1\}$ and $\{a_2,b_2\}$,
the zeros of $P_{\vec{n}}$ accumulate on a set $\Delta_0$
that consists of two disjoint arcs $\Delta_1^*$ and $\Delta_2$ (as in case I),
but contrary to the case I, the arcs are not joining all branch points,
see Figure \ref{fig:3}. One of the branch points ($b_1$ in Figure \ref{fig:3})
does not belong to $\Delta_0$. The arc $\Delta_1^*$ connects $a_1$ with a point
$b^*$ which is different from $b_1$. The other arc $\Delta_2$ connects
$a_2$ with $b_2$. Both arcs accumulate half of the zeros of $P_{\vec{n}}$.
\begin{figure}[ht]
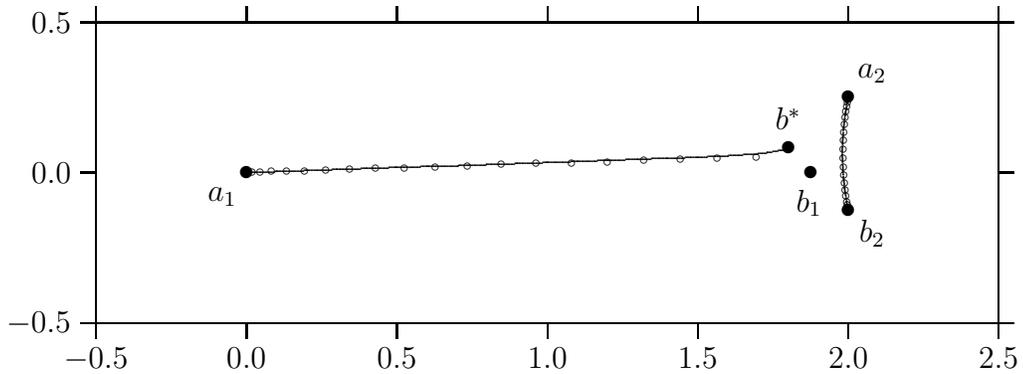

\centering
\mbox{
\input trajectoriesIII.tex }
\caption{$a_1=0$, $b_1=15/8$, $a_2=2+i/4$, $b_2=2-i/8$ (case III). The zeros
of $P_{n,n}$ accumulate on two disjoint arcs $\Delta_1^*$ and $\Delta_2$. The
branch point $b_1$ is not contained in $\Delta_1^*$. The zeros of
$P_{20,20}$ are indicated by $\circ$.}
\label{fig:3}
\end{figure}
\end{description}

The above three cases concern situations where the
four branch points $a_j$, $b_j$ are all distinct. The final two cases deal
with situations where $\{a_1, b_1\}$ and $\{a_2, b_2\}$ have a common point,
and we take it so that
\[ b_1 = b_2 = b. \]

\begin{description}

\item[\boldmath Case IV.] The pairs $\{a_1, b\}$ and $\{a_2, b\}$ are such that
the zeros of $P_{\vec{n}}$ accumulate on a set $\Delta$ that connnects all
three branch points.

\begin{figure}[ht]
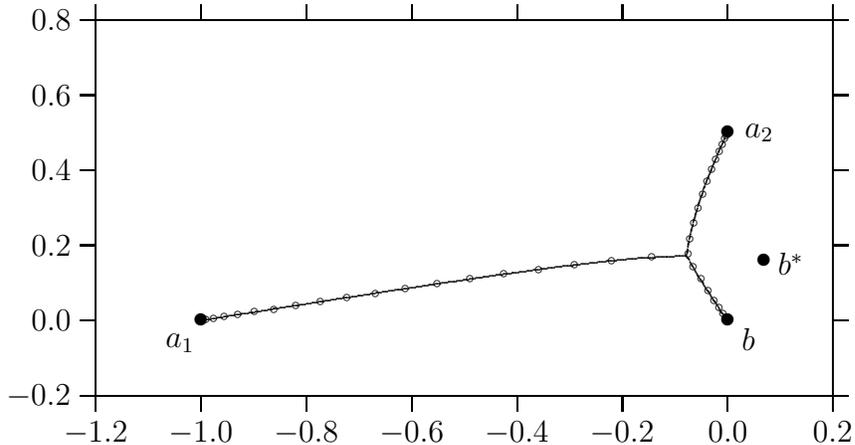

\centering
\mbox{
\input trajectoriesIV.tex }
\caption{$a_1=-1$, $b=0$, $a_2=i/2$ (case $\textup{IV}$)}
\label{fig:4}
\end{figure}

\item[\boldmath Case V.] Certain positions of the pairs $\{a_1, b\}$ and
$\{a_2, b\}$ are such that
the zeros of $P_{\vec{n}}$ accumulate on a set $\Delta_0$ that consists of
two disjoint arcs $\Delta_1^*$ and $\Delta_2$ as in case III. Thus, $\Delta_2$
connects $a_2$ with $b$, but the other arc $\Delta_1^*$ connects $a_1$
with a point $b^*$ which is different from $b$, and which does not belong
to $\Delta_2$. Some of this situations can be realized as limiting
cases of case III where the point $b_1$ tends to $b_2$. These cases
are not part of the case V, but rather belong to case III.
The case V contains the situations that cannot be realized as
limiting cases of case III.

\begin{figure}[ht]
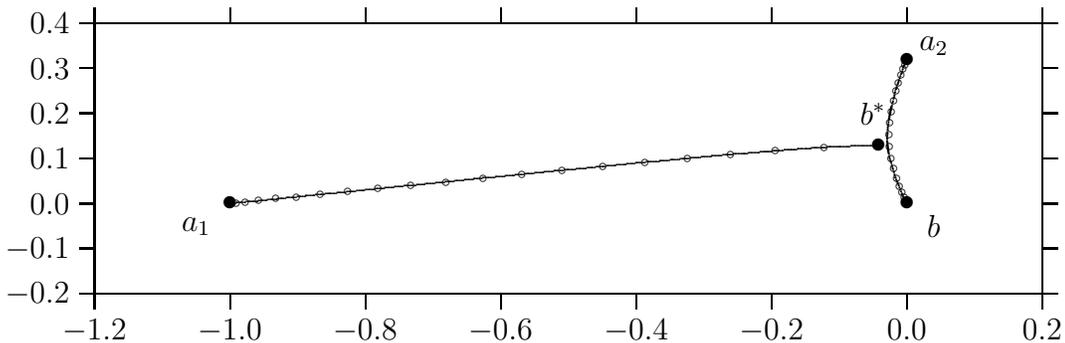

\centering
\mbox{
\input trajectoriesV.tex }
\caption{$a_1=-1$, $b=0$, $a_2=0.32i$ (case $\textup{V}$)}
\label{fig:5}
\end{figure}
\end{description}

\begin{figure}[ht]
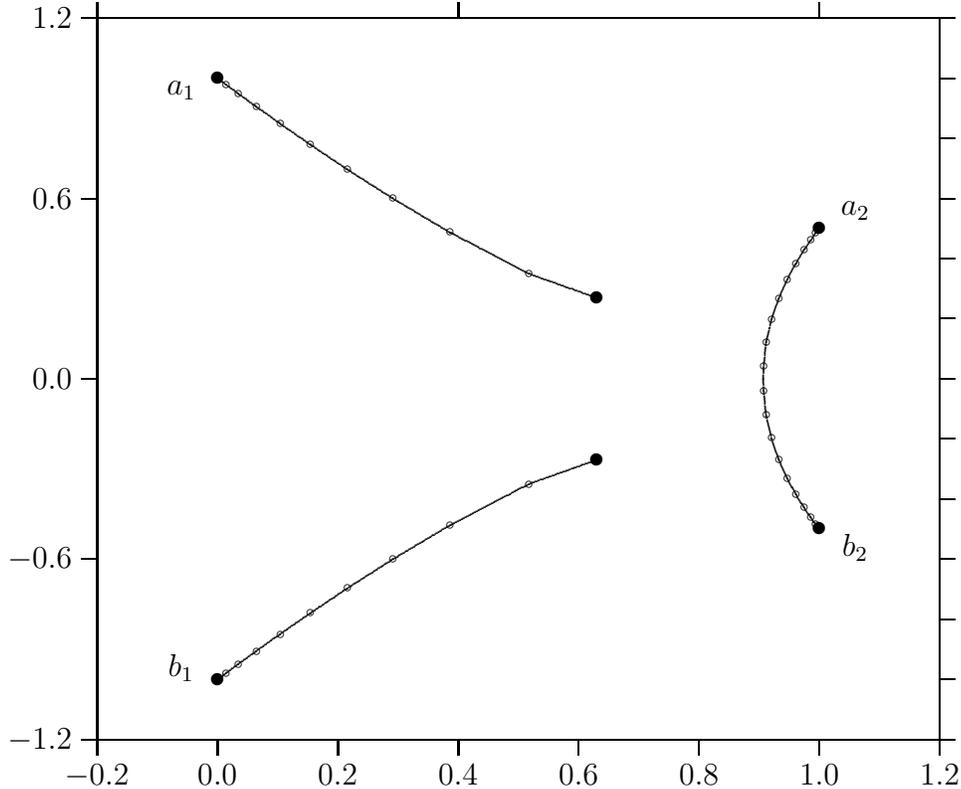

\centering
\mbox{
\input trajectoriesVI.tex }
\caption{$a_1=i$, $b_1=-1$, $a_2=1+i/2$, $b_2=1-i/2$ (higher genus case).}
\label{fig:5b}
\end{figure}

There are also `higher genus cases' such as the one
shown in Figure \ref{fig:5b}. Here the zeros of $P_{\vec{n}}$ accumulate on a set $\Delta$
consisting of three disjoint arcs.
Each of the branch points is contained in $\Delta$, but only two
of them (in the figure it is $a_2$ and $b_2$) are on the same arc $\Delta_2$.
This arc accumulates half of the zeros of $P_{\vec{n}}$. The other
half are on the two remaining arcs. We will not treat the higher genus
cases in this paper.

There are also critical cases when there is a transition from one case to another case,
e.g., the critical case where the two arcs $\Delta_1$ and $\Delta_2$ have
one common point, which is the transition from case I to case II.

The cases I and III  present a complex generalization of an Angelesco system
(\ref{eq:1.14})--(\ref{eq:1.15}).
As we will show, the limiting distribution of the poles of the
\HP{} approximants for these cases is the
equilibrium measure of the vector potential extremal problem (\ref{eq:1.12})--(\ref{eq:1.13})
with the Angelesco interaction matrix (\ref{eq:1.16}), and the support
$\Delta:=\Delta_1 \cup \Delta_2$  of this measure is characterized
by a vector analog of Stahl's symmetry property
\[   \frac{\partial U_j^{\vec{\lambda}}}{\partial n_+} =
    \frac{\partial U_j^{\vec{\lambda}}}{\partial n_-}
    \textrm{ on  $\Delta_j$},\quad j=1,2. \]
The cases II, IV and V have features of a Nikishin system because of the presence of the
common piece and we cannot separate $\Delta$ into two disjoint pieces.
For these cases a new curve $E$ appears  on which the finite
zeros of the function of the second kind
(i.e., extra interpolation points) accumulate and an extremal vector potential problem
describes both the limiting
distribution of the poles of the \HP{} approximants and the extra interpolation points.

The three-sheeted Riemann surface which is appropriate for the cases I, III and V has
four branch points at the end points of the components of the support of the
vector equilibrium measure (for case I they are the branch points of the functions
which are being approximated). Thus this Riemann surface has genus zero.
The three-sheeted Riemann surfaces for the cases II and IV have a more complicated sheet structure
but their genus is also zero.

\subsection{Objectives, structure, main results, and tools of the paper}  \label{sec:1.5}
In this paper we prove the asymptotic formulas and
the convergence theorem of the \HP{} approximants for the cases when the appropriate
Riemann surface is of genus zero. In a forthcoming paper  we plan to describe the
cases when the appropriate Riemann surfaces are of genus
one and two. In this paper we will handle cases I, II, III, IV, V that give rise to a
genus zero Riemann surface as described above.

Rigorous definitions of the class of functions which we are approximating corresponding to
the geometrical cases  and
statements of the results will be presented
in the next section. The rest of the paper is devoted to the proofs of the theorems.
Here we present \textit{a general description of the obtained results}.

\subsubsection{Geometry of the problem}  \label{sec:1.5.1}
Our approach to the geometry of the problem,
i.e., finding the appropriate Riemann surface, classification of the geometrical cases, etc.,
is based on the algebraic function $h$ defined by \index{Functions!h@$h$}
\begin{equation}   \label{eq:1.20}
    h^3(z) - 3\ \frac{P_2(z)}{\Pi_4(z)}\ h(z) + 2\ \frac{P_1(z)}{\Pi_4(z)} = 0,
\end{equation}
where the polynomial $\Pi_4$ is defined by the input parameters (the branch points of the
functions we are approximating) \index{Functions!Pi4@$\Pi_4$}
\[    \Pi_4(z) = (z-a_1)(z-b_1)(z-a_2)(z-b_2), \]
and the polynomials $P_2$ and $P_1$ contain two parameters \index{Functions!P1@$P_1$} \index{Functions!P2@$P_2$}
\begin{eqnarray} \label{eq:1.21}
     P_2(z) & = & z^2 - d_1z + d_2, \nonumber \\
     P_1(z) & = & z-c , \\
     d_1 & = & \frac{2c+a_1+b_1+a_2+b_2}{3},  \nonumber
\end{eqnarray} \index{Parameters!d1d2@$d_1$, $d_2$} \index{Parameters!c@$c$}
which can be determined by the input data $a_1,a_2,b_1,b_2$, using the information about
the geometrical case to which the input data belongs. The function $h$ has three branches
$h_0,h_1,h_2$ which we fix at infinity as \index{Functions!h0h1h2@$h_0$, $h_1$, $h_2$}
\begin{equation}  \label{eq:1.22}
   \begin{cases}
   \displaystyle h_0(z) = - \frac{2}{z} + \cdots & \\[10pt]
   \displaystyle h_j(z) = \frac{1}{z} + \cdots, & j=1,2.
   \end{cases}
  \qquad z \to \infty.
\end{equation}
The function $h$ is a rational function on its Riemann surface $\mathfrak{R}$. \index{Surfaces!R@$\mathfrak{R}$}
It has
\begin{itemize}
\item \textit{four poles} at the points on $\mathfrak{R}$ for which the projection
$\pi : \mathfrak{R} \to
 \overline{\mathbb{C}}$ are the points $\{a_1,a_2,b_1,b_2\}$:
\begin{equation}  \label{eq:1.23}
    h(\xi) = \infty \quad \Rightarrow \quad \xi \in \mathfrak{R}:
    \pi(\xi) \in \{a_1,a_2,b_1,b_2\};
\end{equation}
\item \textit{three zeros} at infinity on each sheet and \textit{one zero} at the point
for which the projection is the parameter $c$ from (\ref{eq:1.21}):
\begin{equation}  \label{eq:1.24}
    h(\xi) = 0 \quad \Rightarrow \quad
    \begin{cases}
      \xi = \infty^{(i)}, & i=0,1,2, \\
      \xi \in \mathfrak{R}: \pi(\xi) = c.
    \end{cases}
\end{equation}
\end{itemize}
If the zero at $c$ cancels a pole from $\{a_1,a_2,b_1,b_2\}$, for example when $c=b_1$,
then the equation (\ref{eq:1.20}) for the algebraic function $h$ reduces to
\begin{equation}  \label{eq:1.25}
    h^3(z) -  \frac{3(z-d)}{(z-a_1)(z-a_2)(z-b_2)}\ h(z) + \frac{2}{(z-a_1)(z-a_2)(z-b_2)} = 0,
\end{equation}
where \index{Parameters!d@$d$}
\[    d = \frac{a_1+a_2+b_2}{3}. \]
The discriminant $\mathcal{D}$ of the equation (\ref{eq:1.25}) is \index{Functions!D@$\mathcal{D}$}
\[   \mathcal{D} = \frac{(z-d)^3-(z-a_1)(z-a_2)(z-b_2)}{(z-a_1)(z-a_2)(z-b_2)}, \]
hence the function $h$ from (\ref{eq:1.25}) has branch points at $a_1,a_2,b_2$
(the poles of $h$) and at the point \index{Parameters!bstar@$b^*$}
\begin{equation}  \label{eq:1.26}
    b^* = \frac{a_1a_2b_2-d^3}{a_1a_2+a_1b_2+a_2b_2-3d^2} .
\end{equation}
The Riemann surface corresponding to equation (\ref{eq:1.25}) has genus zero.
This Riemann surface was considered for the first time by Nuttall \cite{3,30} by means
of another equation, and independently by Kalyagin \cite{16,20}. In our work this Riemann
surface will be an appropriate $\mathfrak{R}$ for the geometrical cases III and V
(for Case V the functions we are approximating have a joint branch point,
i.e., $b_1=b_2=b$ in (\ref{eq:1.3})).

For the geometrical Case I the appropriate Riemann surface is defined by the general
equation (\ref{eq:1.20}) where the two unknown parameters in (\ref{eq:1.21}) are
determined by the two conditions
\begin{enumerate}
\item the genus of $\mathfrak{R}$ is zero,
\item the monodromy condition
\[   h_0 \in H(\overline{\mathbb{C}} \setminus (\widetilde{\Delta}_1 \cup
        \widetilde{\Delta}_2)),\quad
     h_j \in H(\overline{\mathbb{C}} \setminus \widetilde{\Delta}_j),\ j=1,2,\
     \widetilde{\Delta}_1 \cap \widetilde{\Delta}_2 = \emptyset. \]
\end{enumerate}
Here $\widetilde{\Delta}_j$ is an arbitrary Jordan arc joining the points $a_j$ and $b_j$
($j=1,2$). \index{Contours!Delta1tildeDelta2tilde@$\widetilde{\Delta}_1$, $\widetilde{\Delta}_2$}
The discriminant $\mathcal{D}$ of the equation (\ref{eq:1.20}) is
\[    \mathcal{D} = \frac{\widetilde{\mathcal{D}}}{\Pi_4^3},
    \qquad \widetilde{\mathcal{D}} = P_2^3-\Pi_4 P_1^2, \]
where the polynomial $\widetilde{\mathcal{D}}$ has degree $4$ due to (\ref{eq:1.21}).
Condition 1 above implies that \index{Functions!Dtilde@$\widetilde{\mathcal{D}}$}
\begin{equation} \label{eq:1.27}
   \widetilde{\mathcal{D}}(z) = \textrm{const } (z-z_1)^2(z-z_2)^2,
\end{equation}
which gives a system of two algebraic equations (of high order) for the determination
of the two unknown parameters in (\ref{eq:1.21}). This system of algebraic equations
has several solutions with different monodromy properties, for example,
$h_0 \in H(\overline{\mathbb{C}} \setminus \widetilde{\Delta}_1)$,
$h_1 \in H(\overline{\mathbb{C}} \setminus (\widetilde{\Delta}_1 \cup \widetilde{\Delta}_2))$,
$h_2 \in H(\overline{\mathbb{C}} \setminus \widetilde{\Delta}_1)$, hence condition 2
above chooses the right solution.
We will not work with this system of algebraic equations because it is too cumbersome.
Instead we introduce in the
next section a substitution for (\ref{eq:1.21}) which automatically fulfills (\ref{eq:1.27}).

Thus the geometrical part of the problem consists of the determination of the algebraic
function $h$
starting from the input branch points $A_j = \{a_j,b_j\}$, ($j=1,2$). As we already
mentioned above, we
characterize in this paper the position of $A_1$ and $A_2$ which guarantee that the
algebraic curve (\ref{eq:1.20})
is of genus zero.

\subsubsection{Standard functions for the asymptotics}  \label{sec:1.5.2}
Since the geometrical analysis gives us the algebraic function $h$ and its Riemann
surface $\mathfrak{R}$,
we can define the functions which allow us to state the asymptotical results. The Abelian
integral, see (\ref{eq:1.19}), is defined as \index{Functions!G@$G$}
\begin{equation}  \label{eq:1.28'}
    G(\xi) = \int_{a_1}^\xi h(z)\, dz, \qquad \xi \in \mathfrak{R},
\end{equation}
and the function \index{Functions!Phi@$\Phi$}
\begin{equation}  \label{eq:1.28}
   \Phi(\xi) = \exp G(\xi), \qquad \xi \in \mathfrak{R}
\end{equation}
is a single valued (rational) function on $\mathfrak{R}$, which is a consequence of the
fact that the genus of
$\mathfrak{R}$ is zero. The local selection of the branches of the algebraic function $h$
in (\ref{eq:1.22})
gives us the local definition of the branches for the algebraic function $\Phi$
(see (\ref{eq:1.11})): \index{Functions!Phi0Phi1Phi2@$\Phi_0$, $\Phi_1$, $\Phi_2$} \index{Parameters!C0C1C2@$C_0$, $C_1$, $C_2$}
\[  \Phi_0(z) = \frac{1}{C_0z^2} + \cdots, \quad \Phi_j(z) = \frac{z}{C_j} + \cdots,
        \quad j=1,2, \qquad z \to \infty, \]
\[   C_1C_2C_0=1, \quad C_1 > 0. \]
Using the union of the analytic curves defined by \index{Contours!Gamma@$\Gamma$}
\[   \Gamma = \{ z \in \mathbb{C} : |\Phi_\ell(z)| = | \Phi_k(z)|,\
        \textrm{for some } 0 \leq \ell < k \leq 2\} , \]
we define the holomorphic branches of $\Phi$ (and respectively $h$)
globally: \index{Contours!gamma1gamma2gamma3@$\gamma_0$, $\gamma_1$, $\gamma_2$}
\[   \Phi_\ell \in H (\mathbb{C} \setminus \gamma_\ell),
        \quad \gamma_\ell \subset \Gamma,\qquad \ell=0,1,2. \]
Details will be given in the next section.

Thus the appropriate Riemann surface $\mathfrak{R} =
\overline{\mathfrak{R}_0 \cup \mathfrak{R}_1 \cup \mathfrak{R}_2}$ \index{Surfaces!R0R1R2@$\mathfrak{R}_0$, $\mathfrak{R}_1$, $\mathfrak{R}_2$}
can be realized as three sheets of the extended complex plane cut along
the contours $\gamma_\ell$:
\[   \mathfrak{R}_\ell = \overline{\mathbb{C}} \setminus \gamma_\ell, \]
and pasted through \index{Contours!gamma01gamma02gamma12@$\gamma_{0,1}$, $\gamma_{0,2}$, $\gamma_{1,2}$}
\[  \overline{\mathfrak{R}_\ell} \cap \overline{\mathfrak{R}_k} = \gamma_{\ell,k},
    \qquad \ell \neq k,\ \ell,k=0,1,2, \]
so that the upper and the lower sides of the cuts on two neighboring sheets are identified.

\subsubsection{Asymptotic results and convergence}  \label{sec:1.5.3}
Using the global definition of the branches of the standard algebraic functions
$h$ and $\Phi$ (depending on the geometrical case) we can sketch our asymptotic results.
The limiting distribution (\ref{eq:1.4}) of the
zeros of the diagonal $\vec{n}=(n,n)$ \HP{} denominator $P_{\vec{n}}$ and a strong
asymptotic formula can
be written \textit{in a unique way for all geometrical cases} under consideration.
In the present paper we prove that \index{Measures!lambda@$\lambda$}
\[   \nu_{P_{\vec{n}}} \stackrel{*}{\to} \frac{1}{2} \lambda, \qquad n \to \infty, \]
where the real valued measure $\lambda$, which is of total mass two, is given by
\begin{equation}  \label{eq:1.29}
   d\lambda(z) = \frac{1}{2\pi i} \left[( h_0)_+ - (h_0)_- \right] \, dz,
   \qquad z \in \gamma_0,
\end{equation}
where the subscript $+$ or $-$ as usual denotes the limiting value of the function
taken from the left or right, respectively, when traversing the contour according
to the orientation determined by the complex line element $dz$. The strong
asymptotics for $P_{\vec{n}}$ is \index{Functions!F0@$F_0$}
\[   P_{\vec{n}}(z) = C_0^{-n} \Phi_0^{-n}(z) \left( F_0(z) + \O(1/n) \right),
    \qquad z \in  \overline{\mathbb{C}} \setminus \gamma_0, \]
where $F_0 \in H(\overline{\mathbb{C}}\setminus \gamma_0)$ is an analog of the
Szeg\H{o} function defined by mean of a certain boundary value problem
(see details in the next section).

The asymptotic formulas for the functions of the second kind (\ref{eq:1.9}) and
(\ref{eq:1.17}) essentially depend on the geometrical case. For the cases I and III
we prove that \index{Functions!F1F2@$F_1$, $F_2$}
\begin{equation}  \label{eq:1.30}
  R_{\vec{n}}^{(j)}(z) = C_j^{-n} \Phi_j^{-n} \left( F_j(z) + \O(1/n) \right),
  \qquad z \in  \overline{\mathbb{C}}
   \setminus \gamma_j, \quad j=1,2,
\end{equation}
for certain functions $F_j \in H(\overline{\mathbb{C}} \setminus \gamma_j)$ $(j=1,2)$.
For the case II, IV  and V
the answer for the functions of the second kind is more involved.
A general feature for these cases is that the asymptotic formula (\ref{eq:1.30})
remains valid in a neighborhood of $\infty$ which is smaller than
$\overline{\mathbb{C}} \setminus  \gamma_j$, $j=1,2$,
and around the $\gamma_j$ $(j=1,2)$ there may appear a domain where one of the main
terms of the asymptotics of $\{\Phi_1,\Phi_2\}$ changes from one to the other, so as
a result the functions of the second kind on the boundary of this domain have an oscillatory
asymptotic behavior which leads to an accumulation of its zeros there.

 From (\ref{eq:1.17}) we have that
\[    f_j(z) - \frac{Q_{\vec{n}}^{(j)}(z)}{P_{\vec{n}}(z)} =
        \frac{R_{\vec{n}}^{(j)}(z)}{P_{\vec{n}}(z)}, \]
hence our asymptotic results give a complete picture of the convergence of the
\HP{} approximants with a description of the possible regions of divergence and the sets
of the accumulation points of the extra interpolation points.

\subsubsection{Tools}  \label{sec:1.5.4}
To prove the asymptotic results we start with a $3\times 3$ matrix-valued \RH{}
boundary value problem  characterizing
\HP{} approximants for two functions, which was proposed in \cite{45}.
Then, using the information about the geometry of the problem,
we develop the steepest descent method of Deift and Zhou \cite{35} which was already
successfully used  for $2\times 2$ matrix-valued \RH{} problems.

We can say that for the geometrical cases I and III the jump matrices in the
initial \RH{} problem have a block structure
and therefore most of the steps of the asymptotic analysis of the solution of the
$3\times 3$ matrix-valued \RH{} problem
can be reduced to the $2 \times 2$ problem, which has  already been developed.
As a new feature we like to mention the procedure of finding the explicit solutions
for the matrix \RH{} problem with non-varying jumps (i.e., jumps which do not depend on $n$),
reducing the matrix problem to a boundary value problem on the corresponding Riemann surface.

For the geometrical cases II, IV and V the jump matrices of the \RH{} problem do
not possess this block structure. We introduce a new decomposition of the jump matrix to
block structure jump matrices (which can be treated by the traditional local decomposition)
and a jump matrix with exponentially growing non-diagonal terms. Nevertheless, due to the
analyticity of the solution of the \RH{} problem, the contour on which this growing jump
occurs can be moved in the domain where the growing terms of the jump become exponentially
decaying. We call this new decomposition a \textit{global opening of the lenses}.
This new procedure brings about new curves which initially have not been present
in the statement of the problem.\footnote{The phenomenon of global opening of
lenses was first discovered during the research leading to this paper.
It turns out to be a general feature for higher order RH problems and
was also used in the papers \cite{48,65}.} For example, we will discover  the analytic curves on
which the extra interpolation points accumulate. We also like to mention a wonderful picture:
when the contours of the jumps make their global movement they may meet each other or one
may pass through another and the corresponding jumps interact: we can see interference,
transparent penetration or even annihilation.

\section{Rigorous definitions and statements of the results}  \label{sec:2}
\subsection{Class of functions and reformulation of the \HP{} approximation problem}
\label{sec:2.1}
Now we will be more precise in the definition of the class of functions (\ref{eq:1.3})
we are approximating and
we state a matrix \RH{} problem as a reformulation of the \HP{} approximation problem.

Let $a$ and $b$ be points in the complex plane and let $\Delta$ be a Jordan rectifiable
arc joining $a$ and $b$: \index{Contours!Delta@$\Delta$}
\begin{equation} \label{eq:2.1}
    a,b \in \mathbb{C}, \quad \Delta = \{\Delta(t): t \in [0,1]\},
        \quad  \Delta(0)=a, \Delta(1)=b.
\end{equation}
We will consider functions $f$ of the form \index{Functions!f@$f$}
\begin{equation}  \label{eq:2.2}
    f(z) = \frac{1}{2\pi i} \int_\Delta \frac{w(\xi)}{\xi-z}\, d\xi,
\end{equation}
where $w$ is some ``nice'' function on $\Delta$ as specialized in
the following definition. \index{Functions!w@$w$}
Note that
\begin{equation} \label{eq:2.3}
     w(\xi) = f_+(\xi) - f_-(\xi), \qquad \xi \in \Delta,
\end{equation}
where $f_{\pm}$ denote the boundary values of $f$ from the left
and right using the orientation on $\Delta$ from $a$ and $b$.

\begin{definition} \label{def:2.1}
Let $\alpha, \beta > -1$ and let $\Omega$ be a domain containing \index{Domains!Omega@$\Omega$}
the arc $\Delta$ joining $a$ and $b$. \index{Parameters!alphabeta@$\alpha$, $\beta$}
Then we say that the function $f$ given in (\ref{eq:2.2})
belongs to the class $\mathcal{A}(a,\alpha;b,\beta;\Omega)$, \index{Function classes!AaalphabbetaOmega@$\mathcal{A}(a,\alpha;b,\beta;\Omega)$}
\begin{equation}  \label{eq:2.4}
     f \in \mathcal{A}(a,\alpha;b,\beta;\Omega),
\end{equation}
if the function $w$ given in (\ref{eq:2.3}) satisfies \index{Functions!w0@$w_0$}
\begin{equation}  \label{eq:2.5}
   \begin{cases}
      w(\xi) = w_0(\xi) (a -\xi)^\alpha (\xi-b)^\beta, & \alpha, \beta > -1, \\
      w_0, 1/w_0 \in H(\Omega).
   \end{cases}
\end{equation}
\end{definition}

Here we fix a holomorphic branch on $\Delta \subset \Omega$ of $(a-z)^\alpha(z-b)^\beta$.
Thus a function $f \in \mathcal{A}((a,\alpha;b,\beta;\Omega)$ has analytic continuation
across the arc $\Delta$ into the domain $\Omega$  situated on the next sheet of the
Riemann surface of $f$. Moreover in the
representation (\ref{eq:2.2}) of $f(z)$ we can choose any arc $\Delta$ in $\Omega$
that is homotopic in $\Omega$ to the original arc, as
long as we keep the endpoints fixed.

For functions $f_1$ and $f_2$ satisfying (\ref{eq:2.2})--(\ref{eq:2.3}), the definition
of the \HP{} approximants
(\ref{eq:1.2}) is equivalent (by Cauchy's theorem) to the determination of the \textit{multiple orthogonal
polynomial} $P_{n_1,n_2}$ which satisfies \index{Functions!Pn1n2@$P_{n_1,n_2}$}
\begin{equation}  \label{eq:2.6}
   \int_{\Delta_j} P_{n_1,n_2}(\xi) \xi^k w_j(\xi)\, d\xi = 0, \qquad k=0,1,\ldots,n_j-1,\ j=1,2,
\end{equation}
and the \textit{functions of the second kind} are given by \index{Functions!Rn1n2@$R_{n_1,n_2}^{(1)}$, $R_{n_1,n_2}^{(2)}$}
\begin{equation}  \label{eq:2.7}
   R_{n_1,n_2}^{(j)}(z) = \frac{1}{2\pi i} \int_{\Delta_j} \frac{P_{n_1,n_2}(\xi)w_j(\xi)}{\xi-z} \, d\xi, \qquad j=1,2.
\end{equation}
Therefore the error of the approximation is \index{Functions!Qn1n2@$Q_{n_1,n_2}^{(1)}$, $Q_{n_1,n_2}^{(2)}$}
\begin{equation}  \label{eq:2.8}
   f_j(z) - \frac{Q_{n_1,n_2}^{(j)}(z)}{P_{n_1,n_2}(z)} = \frac{R_{n_1,n_2}^{(j)}(z)}{P_{n_1,n_2}(z)}.
\end{equation}
We note that the right hand side of (\ref{eq:2.8}) is independent of the normalization of $P_{n_1,n_2}$. We will fix
the normalization of $P_{n_1,n_2}$ by choosing it to be a monic polynomial of degree $n_1+n_2$. Then
\begin{equation}  \label{eq:2.9}
   P_{n_1,n_2}(z) = z^{n_1+n_2} + \cdots,
\end{equation}
and
\begin{equation}  \label{eq:2.10}
   \begin{cases}
  \displaystyle  R_{n_1,n_2}^{(j)}(z) = \frac{1+\O(1/z)}{m_{n_1,n_2}^{(j)} z^{n_j+1}}, & \quad z \to \infty \\[12pt]
  \displaystyle  \frac{1}{m_{n_1,n_2}^{(j)}} = \frac{-1}{2\pi i} \int_{\Delta_j} P_{n_1,n_2}(\xi)\xi^{n_j} w_j(\xi)\, d\xi .
  \end{cases}
\end{equation}
However, the existence of $P_{n_1,n_2}$ satisfying (\ref{eq:2.6}) and (\ref{eq:2.9}) is
not guaranteed: we recall from the introduction that this is only guaranteed for normal
multi-indices. \index{Parameters!mn1n2@$m_{n_1,n_2}^{(1)}$, $m_{n_1,n_2}^{(2)}$}

For normal multi-indices the \HP{} polynomials $P_{n_1,n_2}$ and the functions of the second
kind (\ref{eq:2.7}) can also be defined by means of a matrix-valued \RH{} problem \cite{45}.
We state the problem for the diagonal multi-indices (the case of our interest)
\begin{equation}  \label{eq:2.11}
       \vec{n} = (n,n), \quad \vec{n}_1 = (n-1,n), \quad \vec{n}_2= (n,n-1) .
\end{equation}
If $\vec{n}, \vec{n}_1, \vec{n}_2$ are normal indices, then the
$3\times 3$ matrix-valued analytic function \index{Matrices!Y@$Y$}
\begin{equation}  \label{eq:2.12}
    Y(z) = \begin{pmatrix}
     P_{\vec{n}}(z) & R_{\vec{n}}^{(1)}(z) & R_{\vec{n}}^{(2)}(z) \\
     m_1P_{\vec{n}_1}(z) & m_1R_{\vec{n}_1}^{(1)}(z) & m_1R_{\vec{n}_2}^{(2)}(z) \\
     m_2P_{\vec{n}_2}(z) & m_2R_{\vec{n}_2}^{(1)}(z) & m_2R_{\vec{n}_2}^{(2)}(z)
            \end{pmatrix}, \qquad z \in \mathbb{C} \setminus \Delta_0, \quad \Delta_0 = \Delta_1 \cup \Delta_2,
\end{equation}
with $m_1=m_{\vec{n}_1}^{(1)}$ and $m_2 = m_{\vec{n}_2}^{(2)}$ (see (\ref{eq:2.10})),
is a solution of the following \RH{} problem:  \index{Contours!Delta0@$\Delta_0$}
\begin{equation} \label{eq:2.13}
\left. \parbox{5in}{
\begin{enumerate}
\item $Y \in H^{3\times 3}(\mathbb{C} \setminus \Delta_0)$ with boundary values
$Y_{\pm}$ on $\stackrel{\circ}{\Delta}_0$.
\item Jump relation on $\Delta_0$
\[   Y_+(\xi) = Y_-(\xi) W(\xi), \qquad \xi \in \stackrel{\circ}{\Delta}_0. \]
\item Normalization near infinity:
\[   Y(z) = \left[ I + \O\left( \frac{1}{z} \right)  \right]
    \begin{pmatrix} z^{2n} & 0 & 0 \\
                    0 & z^{-n} & 0 \\
                    0 & 0 & z^{-n}
     \end{pmatrix}, \qquad z \to \infty.  \]
\end{enumerate}}
\right\}
\end{equation}
Here the jump matrix is \index{Matrices!W@$W$} \index{Functions!w1w2@$w_1$, $w_2$}
\begin{equation}  \label{eq:2.14}
   W(\xi) = \begin{pmatrix}
              1 & w_1(\xi) & w_2(\xi) \\
              0 & 1 & 0 \\
              0 & 0 & 1
             \end{pmatrix},
\end{equation}
where we assume that
\begin{eqnarray*}
      w_1 = 0 & & \textrm{on } \Delta_2\setminus \Delta_1 \\
      w_2 = 0 & & \textrm{on } \Delta_1 \setminus \Delta_2,
\end{eqnarray*}
and $I$ is the $3 \times 3$ identity matrix. \index{Matrices!I@$I$}
If we suppose in addition that the entries of the jump matrix (\ref{eq:2.14})
belong to the class (\ref{eq:2.5}), then the matrix \RH{} boundary value problem
(\ref{eq:2.13}), supplemented with conditions
at the endpoints of $\Delta_0 = \Delta_1  \cup \Delta_2$, has a unique solution.
The endpoint conditions when $z \to a_1$ are
\begin{equation}  \label{eq:2.15}
   Y(z) = \begin{cases}
    \displaystyle \begin{pmatrix}
                  \O(1) & \O(|z-a_1|^{\alpha_1}) & \O(1) \\
                  \O(1) & \O(|z-a_1|^{\alpha_1}) & \O(1) \\
                  \O(1) & \O(|z-a_1|^{\alpha_1}) & \O(1)
                  \end{pmatrix}, & -1 < \alpha_1 < 0, \\[25pt]
      \displaystyle \begin{pmatrix}
                  \O(1) & \O(\log |z-a_1|) & \O(1) \\
                  \O(1) & \O(\log|z-a_1|) & \O(1) \\
                  \O(1) & \O(\log|z-a_1|) & \O(1)
                  \end{pmatrix}, &  \alpha_1 = 0, \\[25pt]
         \displaystyle \begin{pmatrix}
                  \O(1) & \O(1) & \O(1) \\
                  \O(1) & \O(1) & \O(1) \\
                  \O(1) & \O(1) & \O(1)
                  \end{pmatrix}, & \alpha_1 > 0,
  \end{cases}
\end{equation}
Here $a_1$ is assumed to be an endpoint of $\Delta_1$ and not an
endpoint of $\Delta_2$. If $a_1$ coincides with an
endpoint of $\Delta_2$ (e.g., $a_1=a_2$), then we will assume that
\begin{equation}  \label{eq:2.16}
    a_1 = a_2 \Rightarrow \alpha_1=\alpha_2,
\end{equation}
and the endpoint conditions as $z \to a_1$ become
\begin{equation}  \label{eq:2.17}
   Y(z) = \begin{cases}
    \displaystyle \begin{pmatrix}
                  \O(1) & \O(|z-a_1|^{\alpha_1}) & \O(|z-a_1|^{\alpha_1}) \\
                  \O(1) & \O(|z-a_1|^{\alpha_1}) & \O(|z-a_1|^{\alpha_1}) \\
                  \O(1) & \O(|z-a_1|^{\alpha_1}) & \O(|z-a_1|^{\alpha_1})
                  \end{pmatrix}, & -1 < \alpha_1 = \alpha_2 < 0, \\[25pt]
      \displaystyle \begin{pmatrix}
                  \O(1) & \O(\log |z-a_1|) & \O(\log |z-a_1|) \\
                  \O(1) & \O(\log|z-a_1|) & \O(\log |z-a_1|) \\
                  \O(1) & \O(\log|z-a_1|) & \O(\log |z-a_1|)
                  \end{pmatrix}, &  \alpha_1 = \alpha_2 = 0, \\[25pt]
         \displaystyle \begin{pmatrix}
                  \O(1) & \O(1) & \O(1) \\
                  \O(1) & \O(1) & \O(1) \\
                  \O(1) & \O(1) & \O(1)
                  \end{pmatrix}, & \alpha_1 = \alpha_2 > 0,
  \end{cases}
\end{equation}
When $z$ approaches an endpoint other than $a_1$ then the
endpoint conditions are similar to (\ref{eq:2.15})--(\ref{eq:2.17}).
We then have

\begin{proposition} \label{prop:2.1}
Suppose the \RH{} boundary value problem \eqref{eq:2.13} has the matrix jump $W$ given by
\eqref{eq:2.14} with entries satisfying \eqref{eq:2.5} and endpoint conditions
\eqref{eq:2.15} or \eqref{eq:2.16}--\eqref{eq:2.17} (and similar conditions at the other
endpoints). If this \RH{} problem has a solution, then this solution is unique
and its entries \eqref{eq:2.12} are the \HP{} polynomials (multiple orthogonal polynomials)
\eqref{eq:2.6} and the functions of the second kind \eqref{eq:2.7} corresponding to
the normal index $\vec{n} = (n,n)$ of the \HP{} approximants for the functions
\begin{equation}  \label{2.18}
    f_j \in \mathcal{A}(a_j,\alpha_j;b_j,\beta_j;\Omega_j), \qquad j=1,2.
\end{equation}
\end{proposition}
\begin{proof}
See \cite{45} for the relation between this $3\times 3$ matrix-valued \RH{} problem and multiple
orthogonal polynomials. The uniqueness when the endpoints conditions are imposed follows in
a similar way as in \cite{41}.
\end{proof}

Our goal is to find a solution of the problem (\ref{eq:2.13})--(\ref{eq:2.17}) for large $n$.
This would give the normality of the \HP{} approximants for large $n$ and their asymptotics
as $n \to \infty$. For this purpose we apply a steepest descent method for asymptotic analysis
(as $n \to \infty$) of the matrix \RH{} problem. One of the main ingredients of this method
consists in finding a location for the arcs $\Delta_1$ and $\Delta_2$ within the domains
$\Omega_1$ and $\Omega_2$ such that the jump matrices in the \RH{}
problem admit a factorization as a product of matrices which tend exponentially fast
(as $n \to \infty$) to the identity matrix outside $\Delta = \Delta_1 \cup \Delta_2$ and a
matrix independent of $n$. The location of these arcs depends on the position
of the points $\{a_1,a_2,b_1,b_2\}$.
Thus the form of the asymptotics of the solution of the \RH{} boundary value problem and
the method of obtaining it depend strongly on the position of the points $\{a_1,b_1,a_2,b_2\}$.

We will assume throughout that the domains $\Omega_1$ and $\Omega_2$ are such that
they contain the ``optimal'' arcs $\Delta_1$ and $\Delta_2$. We also assume that
the domains are large enough so that the deformations we need to do in the course
of the steepest descent analysis can be performed within $\Omega_1$ and $\Omega_2$.
This is in particular important for the cases II, IV, and V where one of the steps involves
global deformation of contours.

\subsection{Geometry of the problem. Cases I and II}  \label{sec:2.2}
In this section we define a class of positions of the points \index{Point sets!A@$A$}
\begin{equation}  \label{eq:2.18}
   A:= \{a_1,b_1;a_2,b_2\}
\end{equation}
which we call the \textbf{geometrical cases I and II}
\begin{equation*}
     A \in \textrm{I} \cup \textrm{II}.
\end{equation*}
This class characterizes analytic functions (\ref{eq:2.4})
\index{Functions!f1f2@$f_1$, $f_2$} \index{Function classes!AajalphajbjbetajOmegaj@$\mathcal{A}(a_j,\alpha_j;b_j,\beta_j;\Omega_j)$}
\begin{equation}  \label{eq:2.19}
f_j \in \mathcal{A}(a_j,\alpha_j;b_j,\beta_j,\Omega_j), \qquad j=1,2,
\end{equation}
whose \HP\ asymptotics are described by means of algebraic functions of the third order
with the only branch points $a_1$, $b_1$, $a_2$, $b_2$. \index{Parameters!a1a2@$a_1$, $a_2$} \index{Parameters!b1b2@$b_1$, $b_2$}
We define the two classes I and II
and we fix for the  input data (\ref{eq:2.18}) the global branches of the algebraic
functions needed for the presentation of the asymptotic results and we introduce
measures for the description of the limiting behavior of the zeros of the \HP\ polynomials.

\subsubsection{New coordinates for the input data}  \label{sec:2.2.1}
The algebraic function $h$ (see (\ref{eq:1.20})) associated with the input data
(\ref{eq:2.18}) satisfies the equation \index{Functions!h@$h$}
\begin{equation}  \label{eq:2.20}
     h^3(z) - 3\ \frac{P_2(z)}{\Pi_4(z)}\ h(z) + 2\ \frac{P_1(z)}{\Pi_4(z)} = 0,
\end{equation}
and plays a key role for the classification of the geometrical cases. We recall from the
introduction that $\Pi_4(z) = (z-a_1)(z-b_1)(z-a_2)(z-b_2)$ and the two unknown parameters of the monic
polynomials $P_1$ and $P_2$, of degree one and two respectively (see (\ref{eq:1.21})), have to
satisfy --- for the cases under consideration in this paper --- the condition that the genus of $h$ is zero.
It follows from this that the discriminant $\mathcal{D}$ of the equation (\ref{eq:2.20}) \index{Functions!D@$\mathcal{D}$} \index{Functions!Dtilde@$\widetilde{\mathcal{D}}$}
\[   \mathcal{D} = \frac{\widetilde{\mathcal{D}}}{\Pi_4^3(z)}, \quad \widetilde{\mathcal{D}} = P_2^3-\Pi_4 P_1^2, \]
has zeros of even multiplicity, i.e., (\ref{eq:2.18}) are the only branch points of $h$. As we already
mentioned in the introduction, this gives us two algebraic equations (of rather high order) for the determination
of the unknown parameters of $P_1$ and $P_2$. In practice it is better to consider the algebraic function $h$
\textit{with a special class of coefficients} in (\ref{eq:2.20}) satisfying the condition that
the zeros of the discriminant of $h$ have multiplicity two.
In other words, starting from known coefficients in (\ref{eq:2.20}), we will obtain the set $A$ in
(\ref{eq:2.18}) as the set of zeros for the polynomial
\begin{equation}  \label{eq:2.21}
   \Pi_4(z) = (z-a_1)(z-b_1)(z-a_2)(z-b_2) .
\end{equation}
To realize this, we introduce new coordinates $(k,p,s,c)$ \index{Parameters!k@$k$} \index{Parameters!p@$p$} \index{Parameters!s@$s$} \index{Parameters!c@$c$}
in $\mathbb{C}^4$ for the input data.
The coefficients in (\ref{eq:2.20}) are given by the following
\begin{proposition}  \label{prop:2.2}
Suppose we are given \index{Point sets!B@$B$}
\begin{equation} \label{eq:2.23}
    B = (k,p,s,c) \in \mathbb{C}^4,
\end{equation}
and we define \index{Functions!P1@$P_1$} \index{Functions!P2@$P_2$}
\begin{equation}   \label{eq:2.23'}
    P_1(z) = z-c, \quad P_2(z) = P_1^2(z)+ 2p P_1(z) + s^2.
\end{equation}
Then  \index{Functions!Pi4@$\Pi_4$}
\[   \Pi_4(z) = \frac{P_2^3(z)-\widetilde{\mathcal{D}}(z)}{P_1^2(z)}, \quad \textrm{with }\
    \widetilde{\mathcal{D}}(z) = (kP_1^2(z)+3psP_1(z)+s^3)^2, \]
is a polynomial of degree $4$ in the variable $z$:
\begin{equation}  \label{eq:2.24}
   \Pi_4 = P_1^4 + 6pP_1^3 + (3s^2+12p^2-k^2) P_1^2 + (8p^3+12s^2p-6psk)P_1+3s^4+3p^2s^2-2s^3k.
\end{equation}
\end{proposition}
The proof of Proposition \ref{prop:2.2} is given in Subsection \ref{sec:3.1.1}.

\medskip

So given $k,p,s,c$ we construct the polynomials $P_1$, $P_2$, and $\Pi_4$
as in (\ref{eq:2.23'}) and (\ref{eq:2.24}) and we use these polynomials in
the definition (\ref{eq:2.20}) of the algebraic function $h$.
Then it follows from Proposition \ref{prop:2.2} that $\widetilde{D}$ is a perfect square,
and  (\ref{eq:2.20}) has  branch points at the zeros of $\Pi_4$
only, which we denote by $a_1$, $b_1$, $a_2$, $b_2$ (taken in some order)
and we assume that they are all distinct. This defines a mapping
\begin{equation} \label{eq:2.22}
    B  = (k, p,s,c) \mapsto  A = \{a_1, b_1, a_2, b_2\}
\end{equation}
Thus, if we use the coordinates
$B$ in (\ref{eq:2.23}) then we can determine both the set $A$ of four points in (\ref{eq:2.18}) and all
the coefficients of the equation (\ref{eq:2.20}) of the algebraic function $h$, which has
its branch points in $A$.
Generically, this function $h$ has genus zero, that is, with the exception of a lower
dimensional subset of the four dimensional $(k,p,s,c)$-space. However, this class
also contains functions of higher genus
(for example, functions with branch points of the third order). These degenerate examples are not
connected with our task and we exclude them from our considerations.

\subsubsection{The function $\Phi$ and the contour $\Gamma$}  \label{sec:2.2.2}
Since the algebraic function $h$ has genus zero, the exponential of the Abelian integral
\begin{equation}  \label{eq:2.25}
     \PHI = \exp \left( \int h(z)\, dz  \right)
\end{equation}
is also an algebraic function with the same Riemann surface as $h$. \index{Functions!PHI@$\PHI$}
The equation (\ref{eq:2.25}) defines $\PHI$ up to a multiplicative constant.

\begin{proposition}  \label{prop:2.3}
For a suitable normalization we have that the function $\PHI$ satisfies the equation
\[  \PHI^3(z) + q_1(z)\PHI^2(z) + q_2(z) \PHI(z) + q_0 = 0, \]
where $q_j$ are polynomials of degree $\leq j$ $(j=0,1,2)$ given by \index{Functions!q0q1q2@$q_0$, $q_1$, $q_2$}
\begin{eqnarray}  \label{eq:2.26}
    q_1(z) &=& -6\sqrt{3} p P_1(z) + \frac{1}{\sqrt{3}} \left[ (k-2\sqrt{3}p) \kappa_+ + (k+2\sqrt{3} p)\kappa_- \right],
   \nonumber \\
    q_2(z) &=& -\kappa_+ \kappa_- \left( P_1^2(z) + 4p P_1(z) + \frac{9s^2-4k^2+36p^2}{9} \right), \\
    q_0 & = & \frac{2\sqrt3}{243} (3s-2k)^2 \kappa_+^2\kappa_-^2,  \nonumber
\end{eqnarray}
where $P_1(z)=z-c$ and \index{Parameters!kappapm@$\kappa_+$, $\kappa_-$}
\begin{equation} \label{eq:2.26'}
    \kappa_{\pm} = 3s+k\pm 3\sqrt{3} p.
    \end{equation}
\end{proposition}
The proof of Proposition is in Subsection \ref{sec:3.1.1}.

\medskip
We assume from now on that $k$, $p$, and $s$ are such that \index{Parameters!q0@$q_0$}
\[ q_0 \neq 0 \]
This also implies that $\kappa_+ \neq 0$ and $\kappa_- \neq 0$.

The function $\PHI$ is defined in (\ref{eq:2.25}) up to a multiplicative constant. In what follows we often choose
a normalization such that the product of all three branches of $\Phi$ equals one: \index{Functions!Phi@$\Phi$}
\begin{equation}  \label{eq:2.27}
    \Phi_0\Phi_1\Phi_2 = 1.
\end{equation}
Therefore the function $\Phi$ with normalization (\ref{eq:2.27}) is related to the function $\PHI$ from
Proposition \ref{prop:2.3} by \index{Functions!Phi0Phi1Phi2@$\Phi_0$, $\Phi_1$, $\Phi_2$}
\[   \Phi(z) = -\frac{\PHI(z)}{q_0^{1/3}},  \]
where we recall our assumption that $q_0 \neq 0$.
Taking into account the explicit expression (\ref{eq:2.26}) for the equation above, we see that the
three branches of the normalized function $\Phi$ behave at infinity as
\begin{equation} \label{eq:2.28}
    \begin{cases}
   \displaystyle   \Phi_j(z) = \frac{z}{C_j} + \cdots, & j=1,2, \\[12pt]
   \displaystyle   \Phi_0(z) = \frac{1}{C_0z^2} + \cdots
   \end{cases} \qquad z \to \infty,
\end{equation}
where \index{Parameters!C0C1C2@$C_0$, $C_1$, $C_2$}
\[   C_1 = - \frac{q_0^{1/3}}{\kappa_+}, \quad C_2 = \frac{q_0^{1/3}}{\kappa_-}, \quad C_0= \frac{1}{C_1C_2}.  \]
 From the equation for $\Phi$ we can obtain a parametrization of the curve
\begin{equation}  \label{eq:2.29}
   \Gamma = \{ z : \ |\Phi_j(z)| = |\Phi_k(z)| \textup{ for some }  0 \leq j < k \leq 2 \}
\end{equation}
in terms of the function \index{Contours!Gamma@$\Gamma$}
\begin{equation} \label{eq:2.30}
J(\nu,z) = \nu^3 + A(z)\nu^2 + B(z) \nu + C(z)
\end{equation}
where \index{Functions!Jnuz@$J(\nu,z)$}
\begin{eqnarray}  \label{eq:2.31}
   A(z) & = & \frac{3q_0-q_1(z)q_2(z)}{q_0} \nonumber \\
   B(z) & = & \frac{q_0q_1^3(z)+q_2^3(z)-5q_0q_1(z)q_2(z) + 3q_0^2}{q_0^2} \\
   C(z) & = &  \frac{2q_0q_1^3(z) - q_1^2(z)q_2^2(z) + 2q_2^3(z) - 4q_0q_1(z)q_2(z)+q_0^2}{q_0^2} \nonumber
\end{eqnarray}
and $q_0,q_1,q_2$ are the coefficients (\ref{eq:2.26}) of the equation for $\PHI$. \index{Functions!A@$A$} \index{Functions!B@$B$} \index{Functions!C@$C$}

\begin{proposition}  \label{prop:2.4}
The set $\Gamma$ given in \eqref{eq:2.29} for the function $\Phi$ given in \eqref{eq:2.25}
satisfies
\begin{equation} \label{eq:2.30'}
  \Gamma = \{ z: \  J(\nu,z) = \nu^3 + A(z)\nu^2 + B(z) \nu + C(z) = 0
    \textup{ for some } \nu \in [-2,2] \},
\end{equation}
where $A$, $B$, and $C$ are given by \eqref{eq:2.31}.
\end{proposition}
The proof of Proposition \ref{prop:2.4} is also given in Subsection \ref{sec:3.1.1}.

\subsubsection{Structure of $\Gamma$. Definition of cases I and II} \label{sec:2.2.3}
Here we shall use the global structure of $\Gamma$ to define the geometrical cases I and II. Since the polynomial
$J(\nu,z)$ in (\ref{eq:2.30}) is of degree $6$ in the variable $z$, the contour $\Gamma$ consists of six
trajectories $z_j(\nu)$ parametrized by $\nu \in [-2,2]$. When $\nu=2$ we have
\[   J(2,z) = - \frac{\Pi_4(z)[P_1(z)(3s+k)+9sp]^2}{q_0^2} .  \]
These trajectories start from the points $a_1$, $b_1$, $a_2$, $b_2$ and two trajectories start from the point \index{Parameters!alpha@$\alpha$}
\[   \alpha = c + \frac{9sp}{3s+k}.  \]
Here we assume that $\alpha$ is different from $a_j,b_j$, $j=1,2$. The case when $\alpha$
coincides with one of the $a_j,b_j$, $j=1,2$ will be treated in the next subsection
(see the case III).
We denote by \index{Contours!gammaajgammabjgammaalphaj@$\gamma_{a_j}$, $\gamma_{b_j}$, $\gamma_{\alpha_j}$}
\begin{equation}  \label{eq:2.32}
   \gamma_{a_1}, \gamma_{b_1}, \gamma_{a_2}, \gamma_{b_2}, \gamma_{\alpha_1}, \gamma_{\alpha_2}
\end{equation}
these trajectories which are then continuously extended as $\nu$ decreases from $2$, see Figure \ref{fig:6}.
When $\nu=-2$ we have \index{Parameters!beta1beta2beta3@$\beta_1$, $\beta_2$, $\beta_3$}
\begin{equation}  \label{eq:2.33}
   J(-2,z) = \frac{[q_0-q_1(z)q_2(z)]^2}{q_0^2} = \textrm{const } (z-\beta_1)^2(z-\beta_2)^2(z-\beta_3)^2.
\end{equation}
These six trajectories therefore meet pairwise at the points $\beta_1, \beta_2, \beta_3$, which are the zeros
of $q_0-q_1q_2$ (see Figure \ref{fig:7}).
\begin{figure}[ht]
\parbox{3in}{\centering
\framebox{
\unitlength 1pt
\linethickness{0.5pt}
\begin{picture}(185,160)(0,0)
\put(20,40){\circle*{4}}
\put(108,12){\circle*{4}}
\put(120,135){\circle*{4}}
\put(158,80){\circle*{4}}
\put(100,85){\circle*{4}}
\put(20,40){\vector(4,1){14}}
\put(108,12){\vector(-1,1){10}}
\put(100,85){\vector(-2,1){12}}
\put(100,85){\vector(2,-1){12}}
\put(120,135){\vector(0,-1){14}}
\put(158,80){\vector(-2,1){12}}
\put(20,30){\makebox(0,0)[cc]{\footnotesize $a_1$}}
\put(108,2){\makebox(0,0)[cc]{\footnotesize $b_1$}}
\put(100,75){\makebox(0,0)[cc]{\footnotesize $\alpha$}}
\put(158,70){\makebox(0,0)[cc]{\footnotesize $b_2$}}
\put(120,145){\makebox(0,0)[cc]{\footnotesize $a_2$}}
\put(42,38){\makebox(0,0)[cc]{\footnotesize $\gamma_{a_1}$}}
\end{picture}}
\caption{Start ($\nu=2$) of the trajectories of $\Gamma$}
\label{fig:6}} \quad
\parbox{3in}{\centering
\framebox{
\unitlength 1pt
\linethickness{0.5pt}
\begin{picture}(185,160)(0,0)
\put(20,40){\circle*{4}}
\put(108,12){\circle*{4}}
\put(120,135){\circle*{4}}
\put(158,80){\circle*{4}}
\put(70,45){\circle*{4}}
\put(125,95){\circle*{4}}
\put(170,145){\circle*{4}}
\put(58,53){\vector(3,-2){10}}
\put(82,37){\vector(-3,2){10}}
\put(113,107){\vector(1,-1){10}}
\put(137,83){\vector(-1,1){10}}
\put(158,157){\vector(1,-1){10}}
\put(182,133){\vector(-1,1){10}}
\put(20,30){\makebox(0,0)[cc]{\footnotesize $a_1$}}
\put(108,2){\makebox(0,0)[cc]{\footnotesize $b_1$}}
\put(158,70){\makebox(0,0)[cc]{\footnotesize $b_2$}}
\put(120,145){\makebox(0,0)[cc]{\footnotesize $a_2$}}
\put(70,55){\makebox(0,0)[cc]{\footnotesize $\beta_1$}}
\put(125,85){\makebox(0,0)[cc]{\footnotesize $\beta_2$}}
\put(180,150){\makebox(0,0)[cc]{\footnotesize $\beta_3$}}
\end{picture}}
\caption{Finish ($\nu=-2$) of the trajectories of $\Gamma$}
\label{fig:7}}
\end{figure}

\begin{definition} \label{def:2.2}
We say that the set of points $\{a_1,b_1,a_2,b_2\}$
belongs to the \textup{geometrical cases I or II}
\begin{equation}  \label{eq:2.34}
      A = \{a_1,b_1,a_2,b_2\} \in \textup{I} \cup \textup{II}
\end{equation}
if there exist coordinates $B=(k,p,s,c)$ which are
mapped by \eqref{eq:2.22} to $A$ such that
\begin{enumerate}
\item the algebraic function $z(\nu)$ defined by the equation $J(\nu,z)=0$, see \eqref{eq:2.30}, has
no branch points on $(-2,2)$, i.e., that by analytic continuation \index{Functions!znu@$z(\nu)$}
the trajectories \eqref{eq:2.32} are defined globally for $\nu \in [-2,2]$
\[   \gamma_{a_1}(\nu), \ \gamma_{b_1}(\nu), \
    \gamma_{a_2}(\nu), \ \gamma_{b_2}(\nu), \
    \gamma_{\alpha_1}(\nu), \ \gamma_{\alpha_2}(\nu),
    \quad \nu \in [-2,2]. \]
\item When $\nu = -2$ we have
\[   \gamma_{a_1}(-2) = \gamma_{b_1}(-2), \quad \gamma_{a_2}(-2) = \gamma_{b_2}(-2).  \]
\end{enumerate}
\end{definition}


Now, for $A \in \textrm{I} \cup \textrm{II}$ we can
define two arcs in $\mathbb{C}$ \index{Contours!gamma1gamma2@$\gamma_1$, $\gamma_2$} \index{Contours!gammaa1b1gammaa2b2@$\gamma_{a_1,b_1}$, $\gamma_{a_2,b_2}$}
\begin{equation}  \label{eq:2.35}
    \gamma_j = \gamma_{a_j,b_j} = \gamma_{a_j} \cup \gamma_{b_j}, \qquad j=1,2,
\end{equation}
each connecting $a_j$ and $b_j$ $j=1,2$,
and a closed analytic curve \index{Contours!gammaalpha@$\gamma_{\alpha}$}
\[   \gamma_\alpha := \gamma_{\alpha_1} \cup \gamma_{\alpha_2}.  \]

\begin{definition}  \label{def:2.2'}
Given $A \in \textup{I} \cup \textup{II}$, we say that (see Figures~{\rm \ref{fig:8}} and {\rm \ref{fig:9}})
\begin{enumerate}
\item $A \in \textup{I}$ if $\gamma_1 \cap \gamma_2 = \emptyset$.
\item $A \in \textup{II}$ if $\gamma_1$ and $\gamma_2$ have two points of intersection
\[   \gamma_{1} \cap \gamma_2  = \{ c_1,c_2\}. \]
\end{enumerate}
\end{definition}

In the case II we assume that the branch points $a_j, b_j$ and
the intersection points $c_j$ are labeled as shown in Figure~\ref{fig:9}.
That is, if we follow the curve $\gamma_j$ starting at $a_j$ ($j=1,2$) then
we first meet $c_1$ and then $c_2$. \index{Parameters!c1c2@$c_1$, $c_2$}

\begin{figure}[ht]
\parbox{3in}{\centering
\framebox{
\unitlength 1pt
\linethickness{0.5pt}
\begin{picture}(185,165)(0,0)
\put(20,80){\circle*{4}}
\put(60,20){\circle*{4}}
\put(160,20){\circle*{4}}
\put(160,140){\circle*{4}}
\put(80,15){\circle*{4}}
\put(20,85){$a_1$}
\put(50,15){$b_1$}
\put(165,145){$a_2$}
\put(165,15){$b_2$}
\put(85,10){$\alpha$}
\qbezier(20,80)(40,70)(60,20)
\qbezier(160,20)(130,80)(160,140)
\put(25,60){$\gamma_1$}
\put(150,80){$\gamma_2$}
\put(70,60){$\gamma_\alpha$}
\thinlines
\qbezier(51,40.5)(80,37)(80,15)
\qbezier(80,15)(78,-8)(53,0)
\qbezier(51,40.5)(33,6)(53,0)
\qbezier(51,40.5)(65,70)(56,84)
\qbezier(56,84)(26,124)(3,98)
\qbezier(3,98)(-12,80)(10,52)
\qbezier(10,52)(16,45)(51,40.5)
\end{picture}}
\caption{Global trajectories, case I}
\label{fig:8}} \quad
\parbox{3in}{\centering
\framebox{
\unitlength 1pt
\linethickness{0.5pt}
\begin{picture}(180,160)(0,0)
\put(80,20){\circle*{4}}
\put(70,140){\circle*{4}}
\put(150,15){\circle*{4}}
\put(160,150){\circle*{4}}
\put(120,40){\circle*{4}}
\put(115,120){\circle*{4}}
\put(124,80){\circle*{4}}
\put(55,145){$a_1$}
\put(65,15){$b_1$}
\put(165,155){$a_2$}
\put(155,10){$b_2$}
\put(122,120){$c_1$}
\put(127,40){$c_2$}
\put(113,80){$\alpha$}
\put(90,125){$\gamma_1$}
\put(123,145){$\gamma_2$}
\put(15,80){$\gamma_\alpha$}
\qbezier(160,150)(124,144)(115,120)
\qbezier(115,120)(100,80)(120,40)
\qbezier(120,40)(130,20)(150,15)
\qbezier(70,140)(100,140)(115,120)
\qbezier(115,120)(145,80)(120,40)
\qbezier(120,40)(110,24)(80,20)
\thinlines
\qbezier(30,80)(30,160)(70,160)
\qbezier(70,160)(100,160)(115,120)
\qbezier(115,120)(130,80)(120,40)
\qbezier(120,40)(110,0)(70,0)
\qbezier(70,0)(30,0)(30,80)
\end{picture}}
\caption{Global trajectories, case II}
\label{fig:9}}
\end{figure}

\subsubsection{Riemann surface for case I.
Definition of the global branches for the algebraic functions $h$ and $\Phi$}  \label{sec:2.2.4}
We denote for the case I (see Figure~\ref{fig:8}) \index{Contours!Delta1Delta2@$\Delta_1$, $\Delta_2$} \index{Contours!delta1delta2@$\delta_1$, $\delta_2$} \index{Contours!Delta0@$\Delta_0$}
\begin{equation}  \label{eq:2.35'}
   \Delta_j := \delta_j := \gamma_j, \quad j=1,2, \quad \Delta_0:=\Delta_1\cup\Delta_2.
\end{equation}

\begin{figure}[ht]
\unitlength 2pt
\linethickness{0.4pt}
\centering
\begin{picture}(150,80)(0,55)
\put(50,140){\line(1,0){100}}
\put(150,140){\line(-5,-4){25}}
\put(125,120){\line(-1,0){100}}
\put(25,120){\line(5,4){25}}
\put(65,135){\circle*{2}}
\put(55,125){\circle*{2}}
\put(110,135){\circle*{2}}
\put(95,125){\circle*{2}}
\qbezier(110,135)(90,135)(95,125)
\qbezier(65,135)(70,130)(55,125)
\put(50,110){\line(1,0){100}}
\put(150,110){\line(-5,-4){25}}
\put(125,90){\line(-1,0){100}}
\put(25,90){\line(5,4){25}}
\put(65,105){\circle*{2}}
\put(55,95){\circle*{2}}
\qbezier(65,105)(70,100)(55,95)
\put(50,80){\line(1,0){100}}
\put(150,80){\line(-5,-4){25}}
\put(125,60){\line(-1,0){100}}
\put(25,60){\line(5,4){25}}
\put(110,75){\circle*{2}}
\put(95,65){\circle*{2}}
\qbezier(110,75)(90,75)(95,65)
\multiput(55,125)(0,-2){15}{\line(0,-1){1}}
\multiput(65,135)(0,-2){15}{\line(0,-1){1}}
\multiput(95,125)(0,-2){30}{\line(0,-1){1}}
\multiput(110,135)(0,-2){30}{\line(0,-1){1}}
\put(20,130){\makebox(0,0)[cc]{$\mathfrak{R}_0$}}
\put(20,100){\makebox(0,0)[cc]{$\mathfrak{R}_1$}}
\put(20,70){\makebox(0,0)[cc]{$\mathfrak{R}_2$}}
\put(60,137){\makebox(0,0)[cc]{\small $a_1$}}
\put(50,125){\makebox(0,0)[cc]{\small $b_1$}}
\put(115,135){\makebox(0,0)[cc]{\small $a_2$}}
\put(100,125){\makebox(0,0)[cc]{\small $b_2$}}
\put(70,130){\makebox(0,0)[cc]{\small $\Delta_1$}}
\put(62,125){\footnotesize $+$}
\put(58,129){\footnotesize $-$}
\put(98,136){\makebox(0,0)[cc]{\small $\Delta_2$}}
\put(96,129){\footnotesize $+$}
\put(90,129){\footnotesize $-$}
\put(70,100){\makebox(0,0)[cc]{\small $\Delta_1$}}
\put(62,95){\footnotesize $+$}
\put(58,99){\footnotesize $-$}
\put(98,76){\makebox(0,0)[cc]{\small $\Delta_2$}}
\put(96,69){\footnotesize $+$}
\put(90,69){\footnotesize $-$}
\end{picture}
\caption{Riemann surface for case I}
\label{fig:10}
\end{figure}
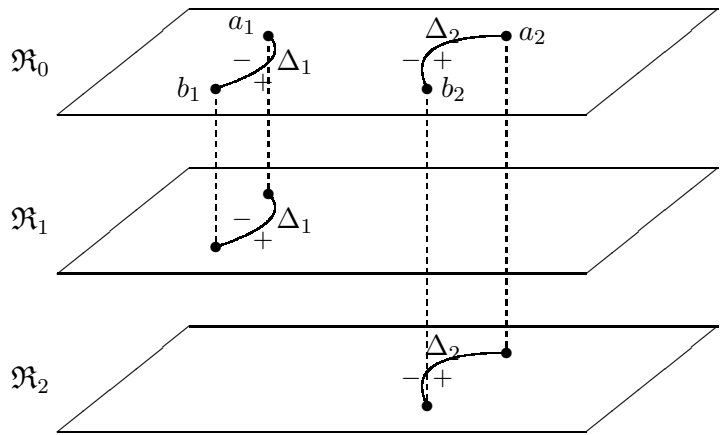

\begin{definition}  \label{def:2.3}
For a set of points $A\in \textup{I}$ the corresponding Riemann surface \index{Surfaces!R@$\mathfrak{R}$} \index{Surfaces!R0R1R2@$\mathfrak{R}_0$, $\mathfrak{R}_1$, $\mathfrak{R}_2$}
\begin{equation}  \label{eq:2.36}
       \mathfrak{R}(A) = \overline{\mathfrak{R}_0 \cup \mathfrak{R}_1 \cup \mathfrak{R}_2}
\end{equation}
is formed by glueing the sheets of the complex plane, cut along the arcs $\Delta_j$ in \eqref{eq:2.35}
\[   \mathfrak{R}_j = \overline{\mathbb{C}} \setminus \Delta_j, \qquad j=1,2, \]
to the sheet
\[   \mathfrak{R}_0 = \overline{\mathbb{C}} \setminus (\Delta_1 \cup \Delta_2), \]
so that the positive (negative) side of the cut on one sheet is identified with the negative (positive) side on the
neighboring sheet as in Figure~{\rm \ref{fig:10}}.
The curves are oriented from $a_j$ to $b_j$ and the positive side is on the left, the negative side on the right
of the oriented curve.
\end{definition}

The function $h$ in (\ref{eq:2.20}) is therefore a rational function on the Riemann surface
$\mathfrak{R}(A)$ in (\ref{eq:2.36})
with divisor (\ref{eq:1.23})--(\ref{eq:1.24}), that is, $h$ has simple poles at the branch
points $a_j, b_j$, simple zeros at the points at infinity on each of the three sheets,
and a fourth simple pole at a value whose projection is the parameter $c$.
The \textit{local definition of the branches of $h$ at infinity}
(\ref{eq:1.22}) can be extended globally as \index{Functions!h0h1h2@$h_0$, $h_1$, $h_2$}
\[   h_0 \in H(\overline{\mathbb{C}} \setminus (\Delta_1 \cup \Delta_2)), \]
\[    h_j \in H(\overline{\mathbb{C}} \setminus \Delta_j), \qquad j=1,2.  \]
The function $\Phi$ in (\ref{eq:2.25}) is also a rational function on $\mathfrak{R}(A)$ with
a double zero at the point at infinity on the sheet $\mathfrak{R}_0$ and a simple
pole at the points at infinity on the other two sheets, see (\ref{eq:2.28}),
and normalization (\ref{eq:2.27}). The local definition of the branches of $\Phi$ at infinity
(\ref{eq:2.28}) can be extended globally as
\begin{equation}  \label{eq:2.36'}
   \Phi_j \in H(\overline{\mathbb{C}} \setminus \Delta_j), \qquad j=0,1,2.
\end{equation}

The contour $\Gamma$ can now be written as a union \index{Contours!Gamma01Gamma02Gamma12@$\Gamma_{0,1}$, $\Gamma_{0,2}$, $\Gamma_{1,2}$}
$\Gamma = \Gamma_{0,1} \cup \Gamma_{0,2} \cup \Gamma_{1,2}$ where
\begin{equation}  \label{eq:2.37''}
   \Gamma_{j,k} = \{ z \in \mathbb C \mid |\Phi_j(z)| = |\Phi_k(z)| \},
    \qquad 0 \leq j < k \leq 2.
\end{equation}
In (\ref{eq:2.37''}) it is understood that for $z$ on one of the cuts $\Delta_1$
and $\Delta_2$ we should take
limiting values. If equality holds on one side of the cut then equality holds
on the other side as well.

The three branches $\Phi_j$, $j=0,1,2$, also determine a number
of regions in the complex plane \index{Domains!Omegajkl@$\Omega_{j,k,\ell}$}
\begin{equation}  \label{eq:2.37'}
   \Omega_{j,k,\ell} = \{ z \in \mathbb{C}: |\Phi_j(z)| < |\Phi_k(z)| < |\Phi_\ell(z)| \},
   \qquad  j,k,\ell=0,1,2.
\end{equation}
We also define \index{Domains!Omegajk@$\Omega_{j,k}$}
\begin{equation}  \label{eq:2.37}
  \Omega_{j,k} = \{ z \in \mathbb{C}: |\Phi_j(z)| < |\Phi_k(z)|\}, \qquad j,k=0,1,2.
\end{equation}
The contours $\Gamma_{j,k}$ and the regions $\Omega_{j,k,\ell}$ give a partition
of the complex plane.

\begin{figure}[ht]
\centering
\framebox{
\unitlength 1pt
\linethickness{0.5pt}
\begin{picture}(120,100)
\put(110,10){\circle*{4}}
\put(90,90){\circle*{4}}
\put(115,5){$b_1$}
\put(95,95){$a_1$}
\qbezier(110,10)(80,40)(90,90)
\put(10,60){\circle*{4}}
\put(50,20){\circle*{4}}
\put(10,65){$a_2$}
\put(55,20){$b_2$}
\qbezier(10,60)(35,50)(50,20)
\qbezier(0,70)(20,90)(50,60)
\qbezier(50,60)(80,30)(60,10)
\qbezier(60,10)(45,-5)(15,20)
\qbezier(15,20)(-15,50)(0,70)
\put(0,0){$\Omega_{0,1,2}$}
\put(10,30){$\Omega_{0,2,1}$}
\put(95,50){$\Gamma_{0,1}$}
\put(40,40){$\Gamma_{0,2}$}
\put(20,80){$\Gamma_{1,2}$}
\end{picture}}
\framebox{
\unitlength 1pt
\linethickness{0.5pt}
\begin{picture}(120,100)
\put(110,10){\circle*{4}}
\put(90,90){\circle*{4}}
\put(115,5){$b_1$}
\put(95,95){$a_1$}
\qbezier(110,10)(80,40)(90,90)
\put(20,60){\circle*{4}}
\put(60,20){\circle*{4}}
\put(23,63){$a_2$}
\put(60,23){$b_2$}
\qbezier(20,60)(40,50)(60,20)
\qbezier(14,66)(26,76)(38,64)
\qbezier(38,64)(46,56)(46,40)
\qbezier(46,40)(24,40)(16,46)
\qbezier(16,46)(6,58)(14,66)
\qbezier(46,40)(70,40)(70,24)
\qbezier(70,24)(70,14)(62,14)
\qbezier(62,14)(46,14)(46,40)
\put(0,0){$\Omega_{0,1,2}$}
\put(95,50){$\Gamma_{0,1}$}
\put(35,70){$\Gamma_{1,2}$}
\put(65,5){$\Gamma_{0,1}$}
\put(-10,40){$\Omega_{0,2,1}$}
\put(25,10){$\Omega_{1,0,2}$}
\put(0,45){\vector(2,1){20}}
\put(35,15){\vector(2,1){20}}
\end{picture}}
\framebox{\unitlength 1pt
\linethickness{0.5pt}
\begin{picture}(120,100)
\put(110,10){\circle*{4}}
\put(90,90){\circle*{4}}
\put(115,5){$b_1$}
\put(95,95){$a_1$}
\qbezier(110,10)(80,40)(90,90)
\put(40,50){\circle*{4}}
\put(55,30){\circle*{4}}
\put(35,56){$a_2$}
\put(57,22){$b_2$}
\qbezier(40,50)(50,40)(55,30)
\qbezier(30,60)(40,70)(55,50)
\qbezier(55,50)(70,30)(65,20)
\qbezier(65,20)(60,10)(40,30)
\qbezier(40,30)(20,50)(30,60)
\put(0,0){$\Omega_{0,1,2}$}
\put(95,50){$\Gamma_{0,1}$}
\put(55,52){$\Gamma_{0,1}$}
\put(10,30){$\Omega_{1,0,2}$}
\put(20,35){\vector(2,1){18}}
\end{picture}}
\caption{Possible partitions of $\mathbb{C}$ by $\Gamma$ for case I}
\label{fig:11}
\end{figure}
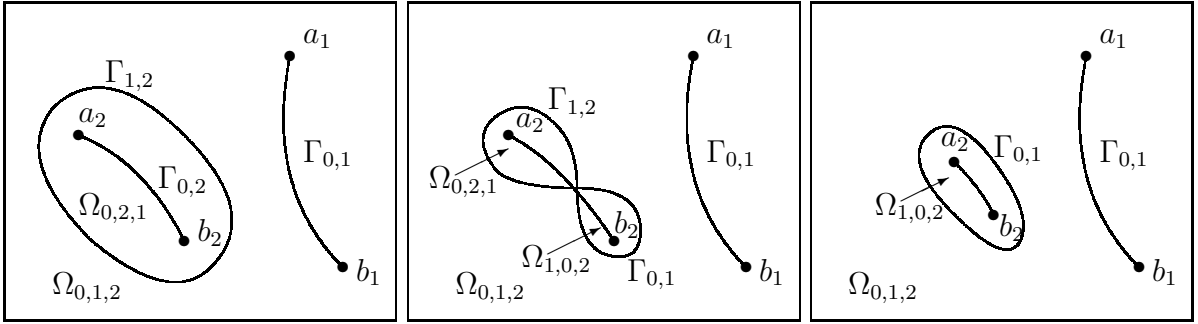

\begin{proposition}  \label{prop:2.4A}
For $A \in \textup{I}$ we have
\begin{enumerate}
\item There are open sets $U_j$ such that
\[ \Delta_j \subset U_j, \qquad U_j \setminus \Delta_j \subset \Omega_{0,j}. \]
\item There exist sets $A \in \textup{I}$ such that (see Figure~{\rm \ref{fig:11})}
\[   \Omega_{1,0} \neq \emptyset \quad \textrm{or} \quad \Omega_{2,0} \neq \emptyset.  \]
\end{enumerate}
\end{proposition}

Figure~{\rm \ref{fig:11}} illustrates possible partitions of the
complex plane by the regions $\Omega_{j,k,\ell}$ and the contours $\Gamma_{j,k}$.
In the left-most figure we have that $\Omega_{1,0}$ and $\Omega_{2,0}$ are both empty.
In the other two figures the set $\Omega_{1,0}$ is non-empty and it contains
part of $\Delta_2$ (as in the middle figure) or all of $\Delta_2$
(as in right-most figure) in its interior.

Now we introduce measures which will describe the limiting behavior of the \HP\ approximants.

\begin{theorem}  \label{thm:2.1}
For $A$ in case $\textup{I}$ we have that the \index{Measures!lambda@$\lambda$}
jump of the function $h_0$ on $\Delta_0=\Delta_1 \cup \Delta_2$ produces a positive measure $\lambda$
      of total mass $2$
\begin{equation}  \label{eq:2.370}
   d\lambda(\xi) = \frac{1}{2\pi i} \Bigl( (h_0)_+(\xi) - (h_0)_-(\xi)\Bigr) \, d\xi, \qquad \xi \in \Delta_0.
\end{equation}
The measure $\lambda$ consists of two measures of mass 1 with densities (with respect
to complex line element $d\xi$)
\begin{equation}  \label{eq:2:371}
  \lambda_j'(\xi) = \frac{1}{2\pi i} \Bigl( h_{0+}(\xi)-h_{0-}(\xi) \Bigr)\Bigr|_{\Delta_j},
\end{equation}
each supported on $\Delta_j$, \index{Measures!lambda1lambda2@$\lambda_1$, $\lambda_2$}
\[  \lambda = \begin{cases}
     \lambda_1 & \textrm{on $\Delta_1$}, \\
     \lambda_2 & \textrm{on $\Delta_2$},
    \end{cases}  \]
and
\[   \lambda_j'(\xi) = \frac{m_j(\xi)}{\sqrt{(\xi-a_j)(\xi-b_j)}}, \qquad j=1,2, \]
where $m_j$ is analytic on $\Delta_j$ for $j=1,2$. \index{Functions!m1m2@$m_1$, $m_2$}
\end{theorem}

The proof of Proposition \ref{prop:2.4A} and Theorem \ref{thm:2.1} is in Subsection \ref{sec:3.1.2}.

\subsubsection{Riemann surface for the case II. Definition of the global branches for the
algebraic functions $h$ and $\Phi$}   \label{sec:2.2.5}
In the case II we have $\gamma_1 \cap \gamma_2 = \{c_1,c_2\}$.
We assume that the points $a_j$, $b_j$ and $c_j$ are labelled
as shown in Figure~\ref{fig:9}. We use $\gamma_{a_j, c_1}$ to denote \index{Contours!gammaa1c1gammaa2c2@$\gamma_{a_1,c_1}$, $\gamma_{a_2,c_2}$}
the part of the arc $\gamma_j$ between $a_j$ and $c_1$ and similarly
for $\gamma_{b_j,c_2}$. \index{Contours!gammab1c1gammab2c2@$\gamma_{b_1,c_1}$, $\gamma_{b_2,c_2}$}
Then we denote for the case II \index{Contours!Delta1tildeDelta2tilde@$\widetilde{\Delta}_1$, $\widetilde{\Delta}_2$}
\begin{equation}  \label{eq:2.372}
  \widetilde{\Delta}_j := \gamma_{a_j,c_1} \cup \gamma_{b_j,c_2}, \quad j=1,2,\quad
  \begin{cases}
   E_1 := \gamma_2 \setminus \widetilde{\Delta}_2, \\
   E_2 := \gamma_1 \setminus \widetilde{\Delta}_1.
  \end{cases}
\end{equation}
The arcs $E_1$ and $E_2$ form a boundary of a lens-shaped domain $G$ \index{Contours!E1E2@$E_1$, $E_2$}
\index{Domains!G@$G$}
\begin{equation}  \label{eq:2.373}
   \partial G := E_1 \cup E_2.
\end{equation}
Note that the analytic curve $\gamma_\alpha$ has to pass through the points $c_1, c_2$ and that
it divides the domain $G$ into two parts (otherwise it would give a contradiction with the
maximum principle for harmonic functions). We set \index{Contours!Delta12@$\Delta_{1,2}$}
\begin{equation}  \label{eq:2.374}
  \Delta_{1,2} := \gamma_\alpha \cap G.
\end{equation}
Finally, we denote \index{Contours!Delta1Delta2@$\Delta_1$, $\Delta_2$} \index{Contours!Delta0@$\Delta_0$}
\begin{equation}  \label{eq:2.375}
  \Delta_j:= \widetilde{\Delta}_j \cup \Delta_{1,2}, \quad
  \Delta_0 := \Delta_1 \cup \Delta_2 = \widetilde{\Delta}_1 \cup \Delta_{1,2} \cup \widetilde{\Delta}_2.
\end{equation}

\begin{definition}  \label{def:2.3B}
For a set of points $A \in \textup{II}$ the corresponding
Riemann surface \index{Surfaces!R@$\mathfrak{R}$} \index{Surfaces!R0R1R2@$\mathfrak{R}_0$, $\mathfrak{R}_1$, $\mathfrak{R}_2$}
\[   \mathfrak{R}(A) := \overline{\mathfrak{R}_0 \cup \mathfrak{R}_1 \cup \mathfrak{R}_2}  \]
is formed by glueing the sheets of the cut complex plane
\begin{equation}  \label{eq:2.36B}
   \begin{cases}
    \mathfrak{R}_1 := \overline{\mathbb{C}} \setminus ( \Delta_1 \cup E_1), \\
    \mathfrak{R}_2 := \overline{\mathbb{C}} \setminus (\widetilde{\Delta}_2 \cup E_1) =
                  \overline{\mathbb{C}} \setminus \gamma_2
   \end{cases}
\end{equation}
to the sheet
\[   \mathfrak{R}_0 := \overline{\mathbb{C}} \setminus \Delta_0
                 = \overline{\mathbb{C}} \setminus (\Delta_1 \cup \widetilde{\Delta}_2), \]
and along $E_1$ the sheets $\mathfrak{R}_1$ and $\mathfrak{R}_2$ are glued to each other (see Figure~{\rm \ref{fig:10B}}).
\end{definition}

\medskip

\begin{figure}[ht]
\unitlength 2pt
\linethickness{0.4pt}
\centering
\begin{picture}(150,80)(0,55)
\put(50,140){\line(1,0){100}}
\put(150,140){\line(-5,-4){25}}
\put(125,120){\line(-1,0){100}}
\put(25,120){\line(5,4){25}}
\put(75,135){\circle*{2}}
\put(65,125){\circle*{2}}
\put(110,135){\circle*{2}}
\put(100,125){\circle*{2}}
\put(90,132){\circle*{2}}
\put(85,128){\circle*{2}}
\qbezier(110,135)(100,134)(90,132)
\qbezier(75,135)(85,134)(90,132)
\qbezier(110,75)(100,74)(90,72)
\qbezier(75,105)(85,104)(90,102)
\qbezier(65,125)(75,126)(85,128)
\qbezier(65,95)(75,96)(85,98)
\qbezier(100,125)(90,126)(85,128)
\qbezier(100,65)(90,66)(85,68)
\qbezier(85,68)(80,70)(90,72)
\qbezier(85,98)(80,100)(90,102)
\qbezier(90,132)(87,130)(85,128)
\qbezier(90,102)(87,100)(85,98)
\put(50,110){\line(1,0){100}}
\put(150,110){\line(-5,-4){25}}
\put(125,90){\line(-1,0){100}}
\put(25,90){\line(5,4){25}}
\put(75,105){\circle*{2}}
\put(65,95){\circle*{2}}
\put(50,80){\line(1,0){100}}
\put(150,80){\line(-5,-4){25}}
\put(125,60){\line(-1,0){100}}
\put(25,60){\line(5,4){25}}
\put(110,75){\circle*{2}}
\put(100,65){\circle*{2}}
\put(90,102){\circle*{2}}
\put(85,98){\circle*{2}}
\multiput(65,125)(0,-2){15}{\line(0,-1){1}}
\multiput(75,135)(0,-2){15}{\line(0,-1){1}}
\multiput(100,125)(0,-2){30}{\line(0,-1){1}}
\multiput(110,135)(0,-2){30}{\line(0,-1){1}}
\multiput(90,132)(0,-2){15}{\line(0,-1){1}}
\multiput(85,128)(0,-2){15}{\line(0,-1){1}}
\put(20,130){\makebox(0,0)[cc]{$\mathfrak{R}_0$}}
\put(20,100){\makebox(0,0)[cc]{$\mathfrak{R}_1$}}
\put(20,70){\makebox(0,0)[cc]{$\mathfrak{R}_2$}}
\put(70,137){\makebox(0,0)[cc]{\small $a_1$}}
\put(60,125){\makebox(0,0)[cc]{\small $b_1$}}
\put(70,107){\makebox(0,0)[cc]{\small $a_1$}}
\put(60,95){\makebox(0,0)[cc]{\small $b_1$}}
\put(115,135){\makebox(0,0)[cc]{\small $a_2$}}
\put(105,125){\makebox(0,0)[cc]{\small $b_2$}}
\put(115,75){\makebox(0,0)[cc]{\small $a_2$}}
\put(105,65){\makebox(0,0)[cc]{\small $b_2$}}
\put(92,137){\makebox(0,0)[cc]{\small $c_1$}}
\put(82,124){\makebox(0,0)[cc]{\small $c_2$}}
\put(92,107){\makebox(0,0)[cc]{\small $c_1$}}
\put(82,94){\makebox(0,0)[cc]{\small $c_2$}}
\end{picture}
\caption{Riemann surface for case II}
\label{fig:10B}
\end{figure}

\begin{remark}  \label{rem:2.2} \
\begin{enumerate}
\item The defined Riemann surface possesses a certain non-symmetry with respect to the pairs
$\{a_1,b_1\}$ and $\{a_2,b_2\}$. We also can use a dual $\mathfrak{R}$ given by
\[   \mathfrak{R}_0:= \overline{\mathbb{C}} \setminus \Delta_0, \quad
     \mathfrak{R}_2:= \overline{\mathbb{C}} \setminus (\Delta_2 \cup E_2), \quad
     \mathfrak{R}_1:= \overline{\mathbb{C}} \setminus \gamma_1.  \]
\item Although all three sheets are glued together at the points $c_1,c_2$, it can easily be checked
that these points are not branch points of $\mathfrak{R}$.
\item Note that the $\mathfrak{R}_1$ sheet is a disconnected set. It consists of two
components: a domain $G_1$ \index{Domains!G1@$G_1$} bounded by $E_1$ and $\Delta_{1,2}$
\[   \partial G_1 := E_1 \cup \Delta_{1,2}, \]
and the domain $\overline{\mathbb{C}} \setminus (\overline{G_1} \cup \widetilde{\Delta}_1)$.
\end{enumerate}
\end{remark}

\begin{figure}[ht]
\centering
\unitlength 4pt
\linethickness{0.4pt}
\begin{picture}(60,60)
\put(18,52){\circle*{1}}
\put(14,18){\circle*{1}}
\put(40,48){\circle*{1}}
\put(32,18){\circle*{1}}
\put(28,42){\circle*{1}}
\put(24,26){\circle*{1}}
\put(17,54){$a_1$}
\put(12,16){$b_1$}
\put(41,49){$a_2$}
\put(33,17){$b_2$}
\qbezier(28,42)(24,34)(24,26)
\qbezier(28,42)(17,35)(24,26)
\qbezier(28,42)(33,35)(24,26)
\qbezier(24,26)(27.5,21.5)(32,18)
\qbezier(24,26)(19.5,21.5)(14,18)
\qbezier(28,42)(33.5,45.5)(40,48)
\qbezier(28,42)(25.5,45.5)(18,52)
\qbezier(28,42)(34,58)(50,54)
\qbezier(50,54)(66,50)(54,24)
\qbezier(54,24)(44,2)(32,8)
\qbezier(32,8)(24,12)(24,26)
\put(0,40){$\Omega_{0,1,2}$}
\put(55,23){$\gamma_\alpha$}
\put(21.5,32){$\Omega_{0,1,2}$}
\put(26,35){$\Omega_{0,1,2}$}
\put(44,42){$\Omega_{0,2,1}$}
\put(23,48){$\gamma_1$}
\put(33,48){$\gamma_2$}
\end{picture}
\caption{Partition of $\mathbb{C}$ by $\Gamma$ into $\Omega_{j,k,\ell}$. The case II}
\label{fig:11B}
\end{figure}
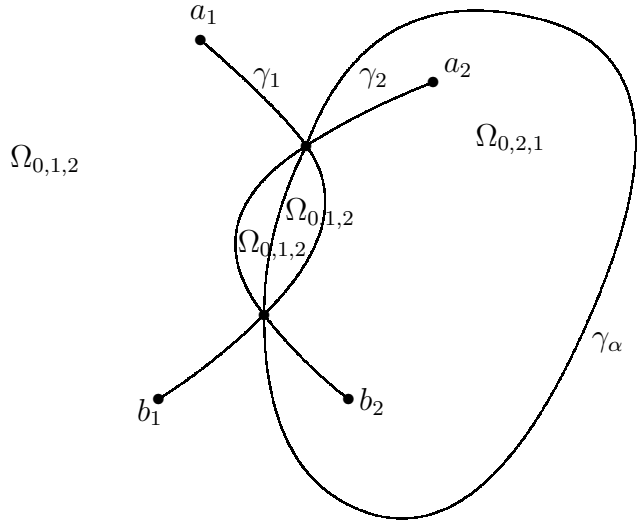

The structure of the sheets (\ref{eq:2.36B}) defines the global branches of the functions
$h$ and $\Phi$: \index{Functions!h0h1h2@$h_0$, $h_1$, $h_2$} \index{Functions!Phi0Phi1Phi2@$\Phi_0$, $\Phi_1$, $\Phi_2$}
\begin{eqnarray}  \label{eq:2.376}
    &&h_0, \Phi_0 \in H(\overline{\mathbb{C}} \setminus \Delta_0), \nonumber \\
    &&h_1,\Phi_1 \in H(\mathbb{C} \setminus( \widetilde{\Delta}_1 \cup \overline{G_1})) \cup H(G_1), \\
    &&h_2, \Phi_2 \in H(\mathbb{C} \setminus (\widetilde{\Delta}_2 \cup E_1)).  \nonumber
\end{eqnarray}
More precisely, in the domains
\[   \overline{\mathbb{C}} \setminus \Delta_0, \quad
     \mathbb{C} \setminus (\widetilde{\Delta}_1 \cup \overline{G_1}), \quad
     \mathbb{C} \setminus (\widetilde{\Delta}_2 \cup E_1),  \]
the branches
\[   (h_0,\Phi_0), \quad (h_1,\Phi_1),\quad (h_2,\Phi_2)  \]
are respectively the result of analytic continuation of (\ref{eq:1.22}) and (\ref{eq:2.28}) from
point at infinity and the branches $(h_1,\Phi_1)$ in $G_1$ are the result of analytic continuation
of $(h_0,\Phi_0)$ through $\Delta_{1,2}$ or $(h_2,\Phi_2)$ through $E_1$.
Using the continuity of the global branches of $\Phi$ along $\gamma_1, \gamma_2, \gamma_\alpha$ and
the maximum principle, we obtain a partition of $\mathbb{C}$ by $\Gamma$ into domains
$\Omega_{j,k,\ell}$ (see \eqref{eq:2.37'} for the definition), as is shown in Figure~\ref{fig:11B}.

\begin{proposition}  \label{prop:2.4B}
If $A \in \textup{II}$
then the contour $\Gamma = \Gamma_{0,1} \cup \Gamma_{0,2} \cup \Gamma_{1,2}$
defined in \eqref{eq:2.29}--\eqref{eq:2.30} has the following structure: \index{Contours!Gamma01Gamma02Gamma12@$\Gamma_{0,1}$, $\Gamma_{0,2}$, $\Gamma_{1,2}$}
\begin{eqnarray*}
   \Gamma_{0,1} & := & \{ z: |\Phi_0(z)| = |\Phi_1(z)| \} = \Delta_1, \\
   \Gamma_{0,2} & := & \{ z: |\Phi_0(z)| = |\Phi_2(z)| \} = \widetilde{\Delta}_2 , \\
   \Gamma_{1,2} & := & \{ z: |\Phi_1(z)| = |\Phi_2(z)| \} = E_1 \cup E_2 \cup (\gamma_\alpha \setminus \Delta_{1,2}),
\end{eqnarray*}
and for the domains $\Omega_{j,k,\ell}$ we have (see Figure~{\rm \ref{fig:11B}})
\begin{equation}  \label{eq:2.377}
   \begin{cases}
    \partial \Omega_{0,1,2} = \gamma_\alpha \cup \widetilde{\Delta}_1 \cup E_1 \cup E_2, \\
    \partial \Omega_{0,2,1} = (\gamma_\alpha \setminus \Delta_{1,2}) \cup \widetilde{\Delta}_2 \cup E_2.
    \end{cases}
\end{equation}
\end{proposition}

The proof of Proposition \ref{prop:2.4B} is in Subsection \ref{sec:3.1.3},
as is the proof of the following Theorem \ref{thm:2.1B} which
introduces the measures.

\begin{theorem} \label{thm:2.1B}
For $A \in \textup{II}$ we have:
\begin{enumerate}
\item The jump of $h_0$ on $\Delta_0$ produces a positive measure $\lambda$ of total mass $2$ \index{Measures!lambda@$\lambda$}
\[   \frac{1}{2\pi i} \Bigl( h_{0+}(\xi)-h_{0-}(\xi) \Bigr)\, d\xi =: d\lambda(\xi), \qquad \xi \in \Delta_0. \]
The measure $\lambda$ consists of two measures $\lambda_1$ and $\tilde{\lambda}_2$ supported on $\Delta_1$
and $\widetilde{\Delta}_2$ \index{Measures!lambda1lambda2tilde@$\lambda_1$, $\tilde{\lambda_2}$}
\[   \lambda = \begin{cases}
               \lambda_1 & \textrm{on $\Delta_1$}, \\
               \tilde{\lambda}_2 & \textrm{on $\widetilde{\Delta}_2$},
                \end{cases}   \]
with respective densities \index{Functions!m1m2@$m_1$, $m_2$}
\[   \lambda_1'(\xi) = \frac{m_1(\xi)}{\sqrt{(\xi-a_1)(\xi-b_1)}}, \qquad m_1 \in H(\widetilde{\Delta}_1) \cap
                                  H(\Delta_{1,2}),  \]
\[   \tilde{\lambda}_2'(\xi) = \frac{m_2(\xi)}{\sqrt{(\xi-a_2)(\xi-b_2)}}, \qquad m_2 \in H(\widetilde{\Delta}_2).  \]
\item The jump of $h_1$ on $E_1$ produces a positive measure $\mu_1$ \index{Measures!mu1@$\mu_1$}
\[  \frac{1}{2\pi i} \Bigl( h_{1+}(\xi)-h_{1-}(\xi) \Bigr)\, d\xi =: d\mu_1(\xi), \qquad \xi \in E_1 \]
and $\mu_1' \in H(E_1)$.
\item There are connections among the total masses of these measures:
\[   |\lambda_1| + |\tilde{\lambda}_2| = 2, \quad |\lambda_1|-|\mu_1|=1.  \]
\end{enumerate}
\end{theorem}

Thus the Riemann surface for the case II produces a system of three positive measures
\[   \{\lambda_1,\tilde{\lambda}_2,\mu_1 \}, \quad
    \begin{cases}
      \supp(\lambda_1) = \Delta_1, \\
      \supp(\tilde{\lambda}_2) = \widetilde{\Delta}_2, \\
      \supp(\mu_1) = E_1,
     \end{cases}
     \quad \begin{cases}
           |\lambda_1| + |\tilde{\lambda}_2| = 2, \\
           |\lambda_1|-|\mu_1|=1.
           \end{cases}  \]
If we consider the dual Riemann surface (see Remark~\ref{rem:2.2}, item 1), then we arrive at
a dual system of three positive measures \index{Measures!lambda1tildelambda2@$\tilde{\lambda}_1$, $\lambda_2$} \index{Measures!mu2@$\mu_2$}
 \begin{equation}  \label{eq:2.40B}
   \{\lambda_2,\tilde{\lambda}_1,\mu_2 \}, \quad
    \begin{cases}
      \supp(\lambda_2) = \Delta_2, \\
      \supp(\tilde{\lambda}_1) = \widetilde{\Delta}_1, \\
      \supp(\mu_2) = E_2,
     \end{cases}
     \quad \begin{cases}
           |\lambda_2| + |\tilde{\lambda}_1| = 2, \\
           |\lambda_2|-|\mu_2|=1,
           \end{cases}
\end{equation}
and we have
\[    \lambda_1 + \tilde{\lambda}_2 = \lambda_2 + \tilde{\lambda}_1 = \lambda.  \]

\subsection{Geometry of the problem. Case III}  \label{sec:2.3}
Recall that the geometrical case III is such that the zeros of $P_{\vec{n}}$
accumulate on two disjoint arcs which do not contain all the branch points.
We assume that the branch points are numbered such that one arc connects
$a_2$ and $b_2$ and such that the other arc connects $a_1$ and $b^*$ with
$b^* \neq b_1$, see Figure \ref{fig:3}.
We again start from the algebraic function $h$ in (\ref{eq:2.20}), however we now use its reduced form
(\ref{eq:1.25}). Without loss of generality for the set \index{Point sets!A@$A$}
$A = \{a_1,b_1;a_2,b_2\}$
we assume that
\begin{equation}  \label{eq:2.38}
  |a_1-b_1| \geq |a_2-b_2|, \quad
   \textrm{dist}(a_1,[a_2,b_2]) \geq \textrm{dist}(b_1,[a_2,b_2]).
\end{equation}
For a set $A$ with conditions (\ref{eq:2.38}) we set \index{Point sets!Aprime@$A'$}
\begin{equation}  \label{eq:2.39}
    A' = \{ a_1;a_2,b_2\}
\end{equation}
and associate with the triple $A'$ the algebraic function (\ref{eq:1.25}): \index{Functions!h@$h$}
\begin{equation}  \label{eq:2.40}
   h^3 - 3 \frac{z- \displaystyle \frac{a_1+a_2+b_2}{3}}{(z-a_1)(z-a_2)(z-b_2)}\ h + \frac{2}{(z-a_1)(z-a_2)(z-b_2)} = 0.
\end{equation}
We recall from (\ref{eq:1.25}) that $h$ in (\ref{eq:2.40}) has four branch points at
\index{Parameters!a1a2@$a_1$, $a_2$} \index{Parameters!b1b2@$b_1$, $b_2$} \index{Parameters!bstar@$b^*$}
\begin{equation}  \label{eq:2.41}
   a_1,\ a_2,\ b_2,\ b^* = \frac{a_1a_2b_2-\left( \displaystyle \frac{a_1+a_2+b_2}{3} \right)^3}{a_1a_2+a_1b_2+a_2b_2-
  3\left( \displaystyle \frac{a_1+a_2+b_2}{3} \right)^3} .
\end{equation}
We do not need to use the coordinates (\ref{eq:2.23}) now since we already have the explicit expressions
in terms of the input data $A'$ for the coefficients of the equation (\ref{eq:2.40}). However, in order to avoid
too cumbersome expressions we set, without loss of generality,
\begin{equation}  \label{eq:2.42}
   a_1:=-1, \quad b_2:=0, \quad a_2:= a,
\end{equation}
and $a \neq -1$, $a \neq 0$. \index{Parameters!a@$a$}

\begin{proposition}  \label{prop:2.5}
The exponential function of the Abelian integral \index{Functions!Phi@$\Phi$}
\[   \Phi = \exp \left( \int h(z)\, dz \right) \]
of the function $h$ in \eqref{eq:2.40} is, up to a multiplicative constant, an algebraic function satisfying
the equation
\begin{equation}  \label{eq:2.43}
  \Phi^3 + q_1(z) \Phi^2 + q_2(z) \Phi + q_0 = 0,
\end{equation}
where  $\textup{deg } q_j \leq j$ for $j=0,1,2$ and the $q_j$ are rational functions of $a_1,a_2,b_2$. If we take
into account \eqref{eq:2.42} then we have \index{Functions!q0q1q2@$q_0$, $q_1$, $q_2$}
\begin{eqnarray} \label{eq:2.44}
   q_1(z) & = & z(a-1)(a^2+\frac{5a}{2} +1) + \frac{a(a^2+4a+1)}{2} \nonumber \\
   q_2(z) & = & - \kappa\left( \frac{27z^2}{4} - \frac92 (a-1)z - \frac14 (a^2+10a+1) \right) \\
   q_0 & = & \kappa^2  \nonumber
\end{eqnarray}
where $\kappa = a^2(a+1)^2/4$. \index{Parameters!kappa@$\kappa$}
\end{proposition}

The proof of Proposition \ref{prop:2.5} is in Subsection \ref{sec:3.1.4}.

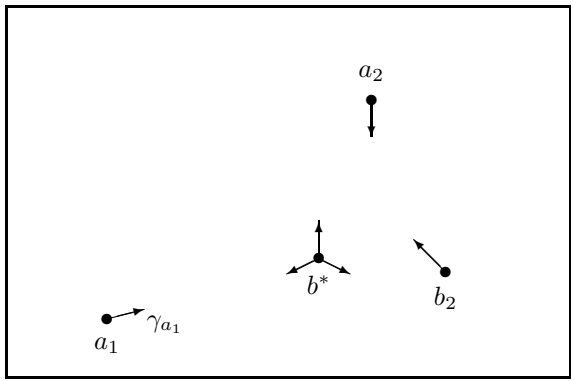
\begin{figure}[ht]
\centering
\framebox{
\unitlength 1pt
\linethickness{0.5pt}
\begin{picture}(185,120)(0,0)
\put(20,12){\circle*{4}}
\put(120,95){\circle*{4}}
\put(148,30){\circle*{4}}
\put(100,35){\circle*{4}}
\put(20,12){\vector(4,1){14}}
\put(100,35){\vector(0,1){14}}
\put(100,35){\vector(2,-1){12}}
\put(100,35){\vector(-2,-1){12}}
\put(120,95){\vector(0,-1){14}}
\put(148,30){\vector(-1,1){12}}
\put(20,2){\makebox(0,0)[cc]{\footnotesize $a_1$}}
\put(100,25){\makebox(0,0)[cc]{\footnotesize $b^*$}}
\put(148,20){\makebox(0,0)[cc]{\footnotesize $b_2$}}
\put(120,105){\makebox(0,0)[cc]{\footnotesize $a_2$}}
\put(42,10){\makebox(0,0)[cc]{\footnotesize $\gamma_{a_1}$}}
\end{picture}}
\caption{Start ($\nu=2$) of the trajectories of $\Gamma$ for Case III}
\label{fig:12}
\end{figure}

Now substitute the $q_0,q_1,q_2$ from (\ref{eq:2.44}) into the expressions (\ref{eq:2.31}) for the
coefficients $A,B,C$ (see Proposition \ref{prop:2.4}) of \index{Functions!Jnuz@$J(\nu,z)$} \index{Functions!A@$A$} \index{Functions!B@$B$} \index{Functions!C@$C$}
\[ J(\nu,z) = \nu^3+A(z)\nu^2+B(z)\nu+C(z). \]
Then as in Proposition \ref{prop:2.4},  we obtain a parametrization of the contour
$\Gamma$ \index{Contours!Gamma@$\Gamma$}
\begin{eqnarray} \nonumber
  \Gamma & := & \{ z\ : \ |\Phi_j(z)| = |\Phi_k(z)|\ \textup{for some } 0 \leq j < k \leq 2\} \\
  & = &  \{ z\ :  \ J(\nu,z) =0 \, \textup{for some } \nu \in[-2,2] \}
 \label{eq:2.44'}
\end{eqnarray}
as the union of six trajectories. When $\nu=2$ the three trajectories $\gamma_{a_1}$, $\gamma_{a_2}$,
$\gamma_{b_2}$ \index{Contours!gammaa1gammaa2gammab2@$\gamma_{a_1}$, $\gamma_{a_2}$, $\gamma_{b_2}$}
start from the points $a_1, a_2, b_2$ and the other three trajectories $\gamma_{b^*}^{(j)}$ $(j=1,2,3)$
start from the point $b^*$ (see Figure \ref{fig:12}). \index{Contours!gammabstarj@$\gamma_{b^*}^{(1)}$, $\gamma_{b^*}^{(2)}$, $\gamma_{b^*}^{(3)}$}
When $\nu=-2$ these trajectories meet pairwise at the points $\beta_1,\beta_2,\beta_2$ given by
(\ref{eq:2.33}) (as in Figure \ref{fig:7}). \index{Parameters!beta1beta2beta3@$\beta_1$, $\beta_2$, $\beta_3$}

\begin{definition}  \label{def:2.4}
We call the triple $A'$ in \eqref{eq:2.39} \textup{acceptable} for case III if
\begin{enumerate}
\item The algebraic function $z(\nu)$ defined by \eqref{eq:2.44'} has no branch points on $(-2,2)$; \index{Functions!znu@$z(\nu)$}
\item When $\nu=-2$ we have
\[   \gamma_{a_2}(-2) = \gamma_{b_2}(-2), \quad \gamma_{a_1}(-2) = \gamma_{b^*}^{(j)}(-2)
    \mbox{ for some } j \in \{1,2,3\}. \]
\end{enumerate}
\end{definition}

If part 2 of the definition is satisfied we assume without loss of generality
that the trajectories starting from $b^*$ are numbered such that
\[ \gamma_{a_1}(-2) = \gamma_{b^*}^{(1)}(-2). \]
It follows from Definition \ref{def:2.4} that  the trajectories
$\gamma_{a_2}$, $\gamma_{b_2}$, $\gamma_{a_1}$, $\gamma_{b^*}^{(1)}$
for an acceptable triple $A'$ are defined globally
for $\nu \in [-2,2]$ and we can define the arcs joining the branch points (\ref{eq:2.41})
as \index{Contours!Delta1starDelta2@$\Delta_1^*$, $\Delta_2$} \index{Contours!Delta0@$\Delta_0$}
\begin{equation}  \label{eq:2.44*}
  \Delta_1^* := \gamma_{a_1} \cup \gamma_{b^*}^{(1)}, \quad \Delta_2 := \gamma_{a_2} \cup \gamma_{b_2},
   \quad \Delta_0 := \Delta_1^* \cup \Delta_2.
\end{equation}
Then we define for an acceptable triple $A'$ the Riemann surface as in subsection \ref{sec:2.2.4} (see (\ref{eq:2.36}))
\index{Surfaces!Rstar@$\mathfrak{R}^*$}
\begin{equation}  \label{eq:2.45}
  \mathfrak{R}^*(A') = \overline{\mathfrak{R}_0^* \cup \mathfrak{R}_1^* \cup \mathfrak{R}_2^*},
\end{equation}
with three sheets \index{Surfaces!R0starR1starR2star@$\mathfrak{R}_0^*$, $\mathfrak{R}_1^*$, $\mathfrak{R}_2^*$}
\[  \mathfrak{R}_0^* = \overline{\mathbb{C}} \setminus (\Delta_1^* \cup \Delta_2), \quad
    \mathfrak{R}_1^* = \overline{\mathbb{C}} \setminus \Delta_1^*, \quad
    \mathfrak{R}_2^* = \overline{\mathbb{C}} \setminus \Delta_2, \]
and the global branches of the algebraic functions $h$ and $\Phi$, defined by (\ref{eq:2.40}) and (\ref{eq:2.43}),
are \index{Functions!h0h1h2@$h_0$, $h_1$, $h_2$} \index{Functions!Phi0Phi1Phi2@$\Phi_0$, $\Phi_1$, $\Phi_2$}
\begin{equation} \label{eq:2.46}
     h_0,\Phi_0 \in H(\mathbb{C}\setminus(\Delta_1^* \cup \Delta_2)), \quad
     h_1, \Phi_1 \in H(\mathbb{C} \setminus \Delta_1^*), \quad
     h_2,\Phi_2 \in H(\mathbb{C} \setminus \Delta_2).
\end{equation}

\begin{figure}[ht]
\centering
\framebox{
\unitlength 1.5pt
\linethickness{0.5pt}
\begin{picture}(110,60)
\put(10,10){\circle*{3}}
\put(60,35){\circle*{3}}
\put(80,30){\circle*{3}}
\put(70,50){\circle*{3}}
\put(5,15){$a_1$}
\put(58,25){$b^*$}
\put(82,32){$b_2$}
\put(73,53){$a_2$}
\qbezier(10,10)(35,25)(60,35)
\qbezier(80,30)(60,40)(70,50)
\qbezier(60,35)(60,60)(70,60)
\qbezier(60,35)(75,20)(80,20)
\qbezier(70,60)(90,60)(90,30)
\qbezier(80,20)(90,20)(90,30)
\put(32,30){\small $\Gamma_{0,1}$}
\put(35,0){$\Omega_{0,1,2}$}
\put(70,40){\small $\Gamma_{0,2}$}
\put(91,45){$\Omega_{1,0,2}$}
\put(90,47){\vector(-1,0){10}}
\end{picture}}
\framebox{
\unitlength 1.5pt
\linethickness{0.5pt}
\begin{picture}(110,60)
\put(10,30){\circle*{3}}
\put(60,30){\circle*{3}}
\put(70,30){\circle*{3}}
\put(95,30){\circle*{3}}
\put(7,20){$a_1$}
\put(56,20){$b^*$}
\put(68,34){$a_2$}
\put(95,35){$b_2$}
\put(10,30){\line(1,0){50}}
\put(70,30){\line(1,0){25}}
\qbezier(60,30)(65,40)(70,40)
\qbezier(70,40)(75,40)(80,30)
\qbezier(60,30)(65,20)(70,20)
\qbezier(70,20)(75,20)(80,30)
\qbezier(80,30)(87,45)(95,45)
\qbezier(95,45)(105,45)(105,30)
\qbezier(80,30)(87,15)(95,15)
\qbezier(95,15)(105,15)(105,30)
\put(35,0){$\Omega_{0,1,2}$}
\put(35,35){\small $\Gamma_{0,1}$}
\put(63,45){$\Omega_{1,0,2}$}
\put(88,50){$\Omega_{0,2,1}$}
\put(103,13){\small $\Gamma_{1,2}$}
\put(83,24){\small $\Gamma_{0,2}$}
\put(65,43){\vector(0,-1){10}}
\put(90,48){\vector(0,-1){10}}
\end{picture}}
\caption{Partition of $\mathbb{C}$ by $\Gamma$. The case III (acceptable triple)}
\label{fig:13}
\end{figure}
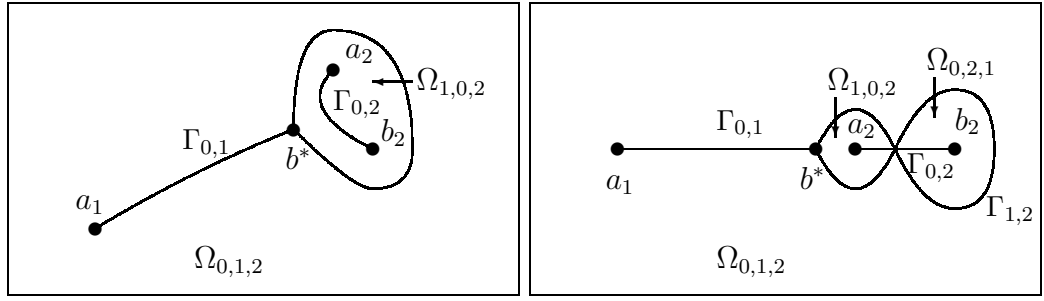

Possible partitions of $\mathbb{C}$ by $\Gamma$ into the domains
$\Omega_{j,k,\ell}$ are shown in Figure~\ref{fig:13}. We have

\begin{proposition}  \label{prop:2.4II}
\
\begin{enumerate}
\item There exists an open set $U_2$ such that
\[ \Delta_2 \subset U_2, \qquad U_2 \setminus \Delta_2 \subset \Omega_{0,2} \]
and there exists an open set $U_1$ such that
\[ \Delta_1^* \setminus \{b^*\} \subset U_1, \qquad U_1 \setminus \Delta_1^* \subset \Omega_{0,1}. \]
\item For any acceptable triple $A'$ we have
\[ \Omega_{1,0} \neq \emptyset \qquad \textup{ and } \qquad b^* \in \partial \Omega_{0,1}.  \]
In addition we have that $\Omega_{1,0}$ is connected.
\end{enumerate}
Moreover the trajectories $\gamma_{b^*}^{(j)}$ $(j=1,2,3)$ start from $b^*$ at an angle $2\pi/3$ and split
a neighborhood of $b^*$ into two domains belonging to $\Omega_{0,1}$ (the boundaries of the domain contain
$\Delta_1^*$) and one domain belonging to $\Omega_{1,0}$ (the boundary
contains $\gamma_{b^*}^{(2)}$ and $\gamma_{b^*}^{(3)}$). See Figure {\rm \ref{fig:13}}.
\end{proposition}

Now we introduce the measures.

\begin{theorem}  \label{thm:2.2}
For an acceptable triple $A'$ we have that \index{Measures!lambda@$\lambda$}
the jump of the branch $h_0$ on $\Delta = \Delta_1^* \cup \Delta_2$ produces a positive measure $\lambda$
on $\gamma$ of total mass $2$. The measure consists of two measures of mass one:
\[   \lambda'(\xi) = \begin{cases}
       \lambda_1'(\xi) = \frac{1}{2\pi i} [ (h_0)_+(\xi) - (h_0)_-(\xi) ], & \xi \in \Delta_1^*, \\[10pt]
       \lambda_2'(\xi) = \frac{1}{2\pi i} [ (h_0)_+(\xi) - (h_0)_-(\xi) ], & \xi \in \Delta_2,
     \end{cases} \]
and \index{Measures!lambda1lambda2@$\lambda_1$, $\lambda_2$}
\[   \lambda_1'(\xi) = \sqrt{\frac{\xi-b^*}{\xi-a_1}}\ m_1(\xi), \quad
     \lambda_2'(\xi) = \frac{m_2(\xi)}{\sqrt{(\xi-a_2)(\xi-b_2)}}, \]
with $m_1$ analytic in a neighborhood of $\Delta_1^*$ and \index{Functions!m1m2@$m_1$, $m_2$}
$m_2$ analytic in a neighborhood of $\Delta_2$.
\end{theorem}

The proofs of Proposition \ref{prop:2.4II} and Theorem \ref{thm:2.2}
are in Subsection \ref{sec:3.1.4},

We can now define the geometrical Case III.
\begin{definition}  \label{def:2.5}
We say that a set of points $A:=\{a_1,b_1;a_2,b_2\}$ satisfying \eqref{eq:2.38} belongs to the
\textup{geometrical Case III}
\[    A \in \textup{III} \]
if
\begin{enumerate}
\item The triple $A' = \{ a_1;a_2,b_2\}$ is acceptable.
\item $b_1 \in \Omega_{1,0}$.
\end{enumerate}
\end{definition}

For $A \in \textup{III}$ we assign to $A$ the Riemann surface $\mathfrak{R}^*$ and the algebraic functions
$h$ and $\Phi$ corresponding to its triple $A'$, i.e., we use (\ref{eq:2.40}), (\ref{eq:2.43}), and
(\ref{eq:2.45}).

\subsection{Geometry of the problem. Common branch point: cases IV and V} \label{sec:2.4}
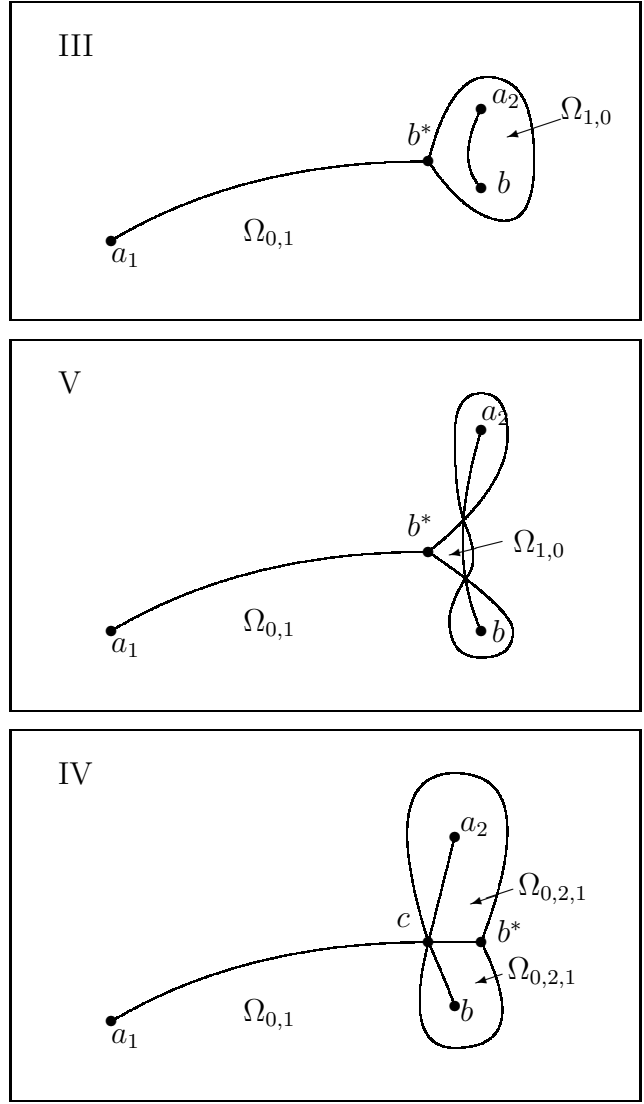
\begin{figure}[ht]
\centering
\framebox{
\unitlength 2pt
\linethickness{0.4pt}
\begin{picture}(105,50)
\put(0,45){$\textup{III}$}
\put(10,10){\circle*{2}}
\put(80,35){\circle*{2}}
\put(70,25){\circle*{2}}
\put(80,20){\circle*{2}}
\qbezier(10,10)(35,25)(70,25)
\qbezier(80,35)(75,25)(80,20)
\qbezier(70,25)(75,45)(85,40)
\qbezier(85,40)(90,37.5)(90,25)
\qbezier(90,25)(90,10)(80,15)
\qbezier(80,15)(75,17.5)(70,25)
\put(10,6){$a_1$}
\put(66,28){$b^*$}
\put(82,36){$a_2$}
\put(83,19){$b$}
\put(35,10){$\Omega_{0,1}$}
\put(95,33){$\Omega_{1,0}$}
\put(95,33){\vector(-3,-1){10}}
\end{picture}}\\ \medskip
\framebox{
\unitlength 2pt
\linethickness{0.4pt}
\begin{picture}(105,60)
\put(0,55){$\textup{V}$}
\put(10,10){\circle*{2}}
\put(80,48){\circle*{2}}
\put(70,25){\circle*{2}}
\put(80,10){\circle*{2}}
\qbezier(10,10)(35,25)(70,25)
\qbezier(80,48)(73,25)(80,10)
\qbezier(70,25)(85,37)(85,47)
\qbezier(85,47)(85,55)(80,55)
\qbezier(80,55)(75,55)(75,47)
\qbezier(80,48)(73,25)(80,10)
\qbezier(75,47)(75,35)(77,30)
\qbezier(70,25)(86,14)(86,10)
\qbezier(86,10)(86,5)(80,5)
\qbezier(80,5)(74,5)(74,12)
\qbezier(74,12)(74,15)(77,20)
\qbezier(77,30)(80,23)(77,20)
\put(10,6){$a_1$}
\put(66,28){$b^*$}
\put(80,50){$a_2$}
\put(82,8){$b$}
\put(35,10){$\Omega_{0,1}$}
\put(86,25){$\Omega_{1,0}$}
\put(84,27){\vector(-4,-1){10}}
\end{picture}}\\ \medskip
\framebox{
\unitlength 2pt
\linethickness{0.4pt}
\begin{picture}(105,60)
\put(0,55){IV}
\put(10,10){\circle*{2}}
\put(75,45){\circle*{2}}
\put(70,25){\circle*{2}}
\put(75,13){\circle*{2}}
\put(80,25){\circle*{2}}
\qbezier(10,10)(35,25)(70,25)
\qbezier(70,25)(75,25)(80,25)
\qbezier(75,13)(72.5,19.5)(70,25)
\qbezier(70,25)(60,57)(75,57)
\qbezier(75,57)(92,57)(80,25)
\qbezier(75,5)(90,5)(80,25)
\qbezier(70,25)(65,5)(75,5)
\qbezier(75,45)(72.5,35)(70,25)
\put(10,6){$a_1$}
\put(83,25){$b^*$}
\put(76,46){$a_2$}
\put(76,10){$b$}
\put(64,28){$c$}
\put(35,10){$\Omega_{0,1}$}
\put(85,18){$\Omega_{0,2,1}$}
\put(87,34){$\Omega_{0,2,1}$}
\put(84,19){\vector(-3,-1){5}}
\put(86,35){\vector(-3,-1){8}}
\end{picture}}
\caption{The subclasses $\textup{III}$, $\textup{V}$ and IV for a common branch point}
\label{fig:14}
\end{figure}

In this section we present a complete classification of the geometry of the \HP{} approximation problem
for the two functions (\ref{eq:1.3}) with a common branch point, i.e.,
\[     f_j \in \mathcal{A}(\overline{\mathbb{C}} \setminus \{ a_j, b\}), \qquad a_1 \neq a_2 . \]
i.e., \index{Point sets!A@$A$}
\[ A = \{a_1,b;a_2,b\}.  \]
Again, as in the case $\textup{III}$ we associate with the triple \index{Point sets!Aprime@$A'$}
\begin{equation}  \label{eq:2.47}
   A' = \{a_1;a_2,b\}
\end{equation}
the algebraic function $h$ defined in (\ref{eq:2.40}), \index{Functions!h@$h$}
and the Abelian exponential $\Phi$ given by (\ref{eq:2.43}). \index{Functions!Phi@$\Phi$}
There are three possibilities for the position of the points  $A$:
\begin{enumerate}
\item The triple $A'$ in (\ref{eq:2.47}) is acceptable (see Definition \ref{def:2.4}) and
$b \in \overline{\Omega}_{1,0}$. Then $A$ can be realized as a limiting situation of case III
and we say
\[     A \in \textup{III} ,  \]
(see top figure of Figure \ref{fig:14}).
\item The triple $A'$ is acceptable but $b \notin \overline{\Omega}_{1,0}$ (see Figure \ref{fig:14} in the middle).
\item The triple $A'$ is not acceptable (see bottom figure of Figure \ref{fig:14}).
\end{enumerate}

\begin{definition}  \label{def:2.5*}
We call the subclasses {\rm 2} and {\rm 3} above the cases V and IV respectively.
\end{definition}

As in the previous subsection (see (\ref{eq:2.42})) we set, without loss of generality, \index{Parameters!a@$a$}
\[     a_1=-1, \quad b=0, \quad a_2:=a, \qquad |a| < 1. \]
In order to give the explicit description of these three subclasses we define
three regions $D_j$ of the disc $D = \{a:\ |a| < 1 \}$ as follows \index{Domains!D@$D$}
\begin{equation}  \label{eq:2.47'}
    \begin{cases}
     a \in D_1 \Longleftrightarrow A = \{-1,0;a,0\} \in \textup{III}, \\
     a \in D_2 \Longleftrightarrow A \in \textup{V}, \\
     a \in D_3 \Longleftrightarrow A \in \textup{IV}.
         \end{cases}
\end{equation}
We have \index{Domains!D1D2D3@$D_1$, $D_2$, $D_3$}
\begin{theorem}  \label{thm:2.4}
\
\begin{enumerate}
\item Two branches of the algebraic curve
\[   a^4 + (4-4\nu)a^3+(22-8\nu)a^2+(4-4\nu)a+1=0 \]
lying in $D$ for $\nu \in [-2,2]$ form the boundary $\partial D_1$ of the domain $D_1$.
\item The algebraic curve \index{Functions!Phat@$\widehat{P}(a,\nu)$}
\begin{multline*}
  \widehat{P}(a,\nu):= 16a^{12} + 96 a^{11} + (336-108\nu)a^{10} + (800-540\nu)a^9 + (2169-1404\nu)a^8 \\
   +\ (4932-2376\nu)a^7 + (6630-2808\nu)a^6 + (4932-2376\nu)a^5 \\
   +\ (2169-1404\nu)a^4 + (800-540\nu)a^3 + (336-108\nu)a^2 + 96a + 16 = 0,
\end{multline*}
which for $\nu=2$ is factorized as
\[   (a-1)^4(2a+1)^4(a+2)^4=0, \]
has six branches in $D$ for $\nu \in [-2,2]$. These six branches start (for $\nu=2$) from the points $1$
(two branches) and $-1/2$ (four branches) and give  the outer boundary $\partial D_2^{\textup{out}}$ of the
region $D_2$. The boundary $\partial D_2$ therefore consists of
\[    \partial D_2 = \partial D_1 \cup \partial D_2^{\textup{out}} .  \]
\item The region $D_3$ is the open set bounded by $\partial D_2^{\textup{out}}$ and $\{|a| =1\}$.
\end{enumerate}
\end{theorem}

The proof of Theorem \ref{thm:2.4} is in section \ref{sec:3.1.5}.

\begin{figure}[ht]
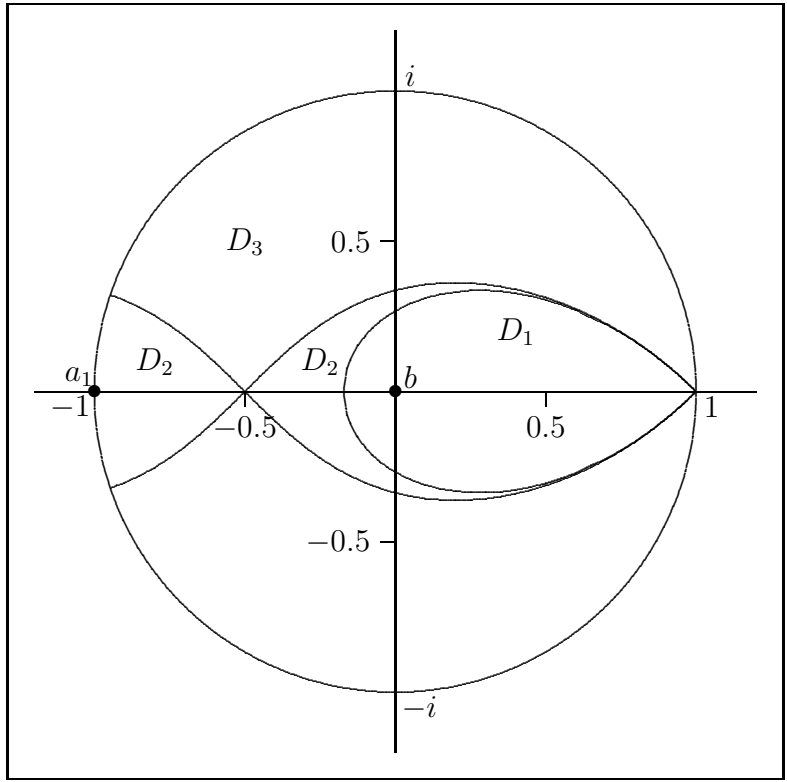

\centering
\framebox{\input fig15 }
\caption{The regions $D_1$, $D_2$, $D_3$ and their boundaries}
\label{fig:15}
\end{figure}

Now we have to assign the Riemann surface to these subclasses and we have
to define the global branches of the algebraic functions $h$ and $\Phi$. The subclass $\textup{III}$ is already
considered (see subsection~\ref{sec:2.3}). For the subclass $\textup{V}$ we again
use (as in $\textup{III}$) all the attributes of the acceptable triple $A'$, i.e., the Riemann surface
$\mathfrak{R}^*(A')$ in (\ref{eq:2.45}) and the global selection of the branches (\ref{eq:2.46}). For
the subclass $\textup{V}$ we can also make item 2 of Proposition~\ref{prop:2.4II} more precise.

\begin{proposition} \label{prop:2.6}
For $A \in \textup{V}$ we have (see Figure {\rm \ref{fig:14}}, second picture) that $b \in \overline{\Omega}_{0,2,1}$.
\end{proposition}

The proof of Proposition \ref{prop:2.6} is in section \ref{sec:3.1.5}.

We denote the connected component of $\Omega_{0,2,1}$ containing $b$ by $G_1$ \index{Domains!G1@$G_1$}
and its boundary by $E_1$.  \index{Contours!E1@$E_1$}
Also we set (see (\ref{eq:2.44*})) We set \index{Contours!Delta12@$\Delta_{1,2}$}
\[   \Delta_{1,2} := \Delta_2 \cap G_1. \]
For the subclass $\textup{V}$ we use the measures $\lambda_1$, $\lambda_2$ introduced in Theorem~\ref{thm:2.2}.
We set \index{Measures!lambda1lambda2@$\lambda_1$, $\lambda_2$}
\index{Measures!lambda1tildelambda12@$\tilde{\lambda}_1$, $\lambda_{1,2}$}
\[  \tilde{\lambda}_1:=\lambda_1, \quad \lambda_{1,2} := \lambda_2\Bigr|_{\Delta_{1,2}},
    \quad \lambda_2:= \lambda_2, \]
where $\lambda_1$ and $\lambda_2$ at the right hand sides are from (\ref{eq:2.46}). In
addition to these measures we shall use the balayage of the measure $\lambda_{1,2}$ on the
contour $E_1$. We denote this measure by $\mu_1$: \index{Measures!mu1@$\mu_1$}
\begin{equation}  \label{eq:2.48}
   V^{\mu_1} = V^{\lambda_{1,2}}, \qquad \textrm{on $E_1$}.
\end{equation}

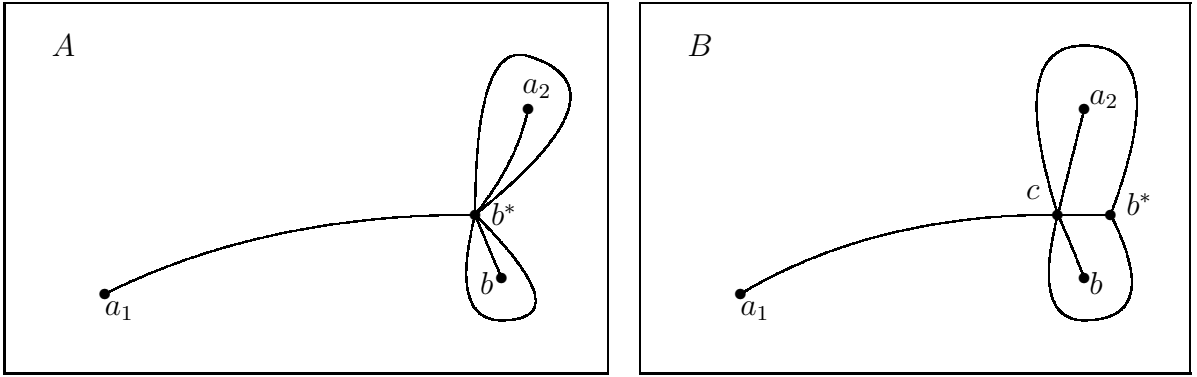
\begin{figure}[ht]
\centering
\framebox{
\unitlength 2pt
\linethickness{0.4pt}
\begin{picture}(100,60)
\put(0,55){$A$}
\put(10,10){\circle*{2}}
\put(90,45){\circle*{2}}
\put(85,13){\circle*{2}}
\put(80,25){\circle*{2}}
\qbezier(10,10)(40,25)(80,25)
\qbezier(80,25)(88,35)(90,45)
\qbezier(80,25)(82,20)(85,13)
\qbezier(80,25)(80,58)(90,55)
\qbezier(90,55)(110,49)(80,25)
\qbezier(80,25)(75,5)(85,5)
\qbezier(85,5)(100,5)(80,25)
\put(10,6){$a_1$}
\put(83,23){$b^*$}
\put(89,48){$a_2$}
\put(81,10){$b$}
\end{picture}}\quad
\framebox{
\unitlength 2pt
\linethickness{0.4pt}
\begin{picture}(90,60)
\put(0,55){$B$}
\put(10,10){\circle*{2}}
\put(75,45){\circle*{2}}
\put(70,25){\circle*{2}}
\put(75,13){\circle*{2}}
\put(80,25){\circle*{2}}
\qbezier(10,10)(35,25)(70,25)
\qbezier(70,25)(75,25)(80,25)
\qbezier(75,13)(72.5,19.5)(70,25)
\qbezier(70,25)(60,57)(75,57)
\qbezier(75,57)(92,57)(80,25)
\qbezier(75,5)(90,5)(80,25)
\qbezier(70,25)(65,5)(75,5)
\qbezier(75,45)(72.5,35)(70,25)
\put(10,6){$a_1$}
\put(83,25){$b^*$}
\put(76,46){$a_2$}
\put(76,10){$b$}
\put(64,28){$c$}
\end{picture}}
\caption{Non-acceptable triples: left is $A$-critical, right is $B$-strictly non-acceptable}
\label{fig:16}
\end{figure}

The last case IV requires special treatment. We start with a characterization of the
contour $\Gamma=\{ z: |\Phi_j(z)|=|\Phi_k(z)|, j\neq k, j,k=0,1,2\}$ for a \textit{non-acceptable} triple $A'$.
Consider the alternative to Definition \ref{def:2.4} of an acceptable triple (see Figure~\ref{fig:16})
\begin{definition}  \label{def:2.4*}
\
\begin{description}
\item[A.] If the algebraic function $z(\nu)$ has a branch point on $(-2,2)$ (see (\ref{eq:2.44'}))
then the triple $A'$ is called \textup{critical}. \index{Functions!znu@$z(\nu)$}
\item[B.]`The triple $A'$ is called \textup{strictly non-acceptable} when
\begin{enumerate}
\item The first condition in Definition \ref{def:2.4} still holds.
\item For $\nu=2$ we have (maybe after renumbering of trajectories $\gamma_{b^*}^{(j)}$)
\[    \gamma_{a_1}(-2)=\gamma_{b^*}^{(1)}(-2), \quad \gamma_{a_2}(-2) = \gamma_{b^*}^{(2)}(-2), \quad
      \gamma_b(-2) = \gamma_{b^*}^{(3)}(-2). \]
\end{enumerate}
\end{description}
\end{definition}
 From the definition we derive (see Figure~\ref{fig:16})

\begin{proposition}  \label{prop:2.7}
\
\begin{description}
\item[A.] For a critical triple one of the trajectories $\gamma_\alpha$, $\alpha \in \{a_2,b\}$ meets
the trajectory $\gamma_{b^*}^{(2)}$ or $\gamma_{b^*}^{(3)}$, i.e.,
\[ \exists \nu_m \in [-2,2), \exists \alpha \in \{a_2,b_2\}, \exists j \in \{2,3\}:
   \gamma_\alpha(\nu_m) = \gamma_{b^*}^{(j)}(\nu_m). \]
For a critical triple the meeting point of the trajectories coincides with $b^*$
\[   \gamma_\alpha(\nu_m) = \gamma_{b^*}^{(j)}(\nu_m) = b^*. \]
\item[B.] For a strictly non-acceptable triple the arcs \index{Contours!gamma1gamma2gamma3@$\gamma_0$, $\gamma_1$, $\gamma_2$}
\begin{equation}  \label{eq:2.49}
  \gamma_1 = \gamma_{a_1} \cup \gamma_{b^*}^{(1)}, \quad \gamma_2 = \gamma_{a_2} \cup \gamma_{b^*}^{(2)},
    \quad  \gamma_3 = \gamma_b \cup \gamma_{b^*}^{(3)}
\end{equation}
intersect at one point $c$ \index{Parameters!c@$c$}
\[   \{c \} = \gamma_1  \cap \gamma_2 \cap \gamma_3 \neq \emptyset.  \]
\end{description}
\end{proposition}

The proof of this proposition is also in Subsection \ref{sec:3.1.5}.

\begin{figure}[ht]
\unitlength 2pt
\linethickness{0.4pt}
\centering
\framebox{\begin{picture}(90,60)
\put(0,55){$\mathfrak{R}_0^*$}
\put(10,6){\small $a_1$}
\put(76,46){\small $a_2$}
\put(76,10){\small $b$}
\put(83,25){\small $b^*$}
\put(10,10){\circle*{2}}
\put(75,45){\circle*{2}}
\put(70,25){\circle{2}}
\put(75,13){\circle*{2}}
\put(80,25){\circle*{2}}
\qbezier(10,10)(35,25)(70,25)
\qbezier(75,13)(72.5,19.5)(70,25)
\qbezier(75,45)(72.5,35)(70,25)
\end{picture}}
\framebox{\begin{picture}(90,60)
\put(0,55){$\mathfrak{R}_1^*$}
\put(10,6){\small $a_1$}
\put(76,46){\small $a_2$}
\put(76,10){\small $b$}
\put(83,25){\small $b^*$}
\put(10,10){\circle*{2}}
\put(75,45){\circle*{2}}
\put(70,25){\circle{2}}
\put(75,13){\circle*{2}}
\put(80,25){\circle*{2}}
\qbezier(10,10)(35,25)(70,25)
\qbezier(70,25)(75,25)(80,25)
\end{picture}}
\vskip10pt
\framebox{\begin{picture}(90,60)
\put(0,55){$\mathfrak{R}_2^*$}
\put(10,6){\small $a_1$}
\put(76,46){\small $a_2$}
\put(76,10){\small $b$}
\put(83,25){\small $b^*$}
\put(10,10){\circle*{2}}
\put(75,45){\circle*{2}}
\put(70,25){\circle{2}}
\put(75,13){\circle*{2}}
\put(80,25){\circle*{2}}
 \qbezier(70,25)(75,25)(80,25)
 \qbezier(75,13)(72.5,19.5)(70,25)
 \qbezier(75,45)(72.5,35)(70,25)
\end{picture}}
\caption{Sheet structure of the Riemann surface for a strictly non-acceptable triple (case $\textup{IV}$)}
\label{fig:17}
\end{figure}
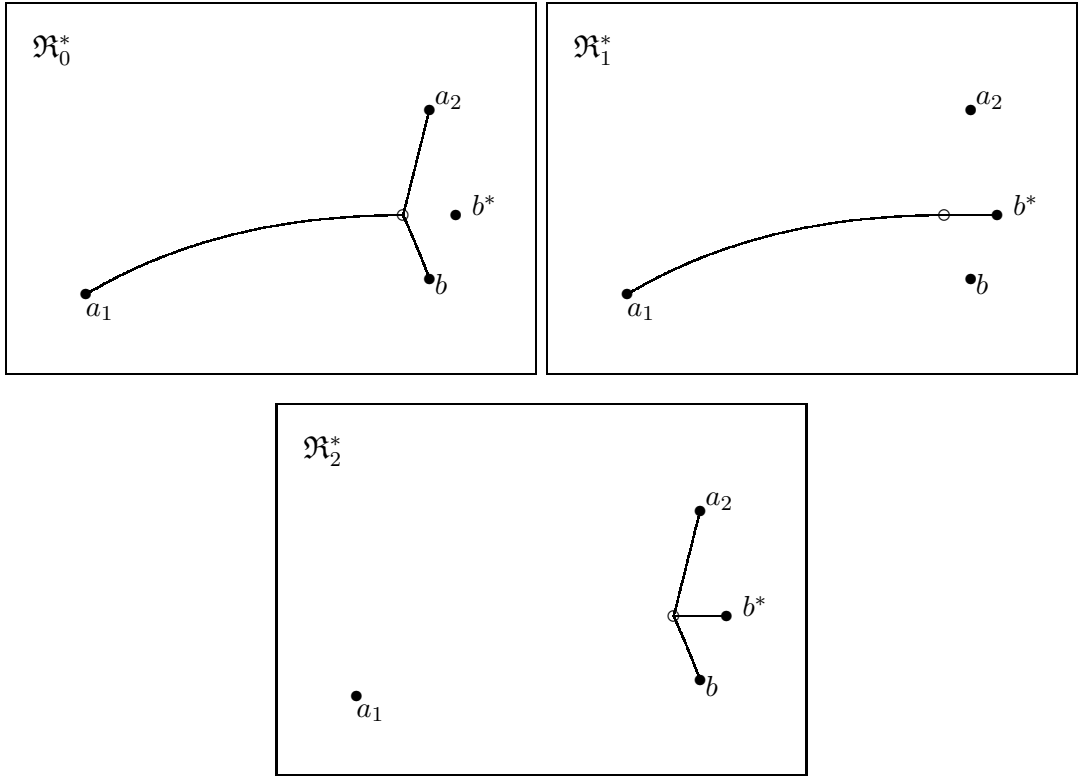

For strictly non-acceptable triples we now assign a Riemann surface and define the global branches for the
algebraic functions $h$ and $\Phi$. We have that each arc from (\ref{eq:2.49}) is divided into two parts by the
point $c$: \index{Contours!gammaa1cgammaa2c@$\gamma_{a_1,c}$, $\gamma_{a_2,c}$, $\gamma_{b,c}$}
\index{Contours!gammacbstar1gammacbstar2@$\gamma_{c,b^*}^{(1)}$, $\gamma_{c,b^*}^{(2)}$, $\gamma_{c,b^*}^{(3)}$}
\index{Contours!gamma1gamma2gamma3@$\gamma_0$, $\gamma_1$, $\gamma_2$}
\begin{equation}  \label{eq:2.50}
   \gamma_1 = \gamma_{a_1,c} \cup \gamma_{c,b^*}^{(1)}, \quad
     \gamma_2 = \gamma_{a_2,c} \cup \gamma_{c,b^*}^{(2)}, \quad
     \gamma_3 = \gamma_{b,c} \cup \gamma_{c,b^*}^{(3)}.
\end{equation}
We take three sheets of the complex plane cut as indicated (see Figure~\ref{fig:17}): \index{Surfaces!R0starR1starR2star@$\mathfrak{R}_0^*$, $\mathfrak{R}_1^*$, $\mathfrak{R}_2^*$}
\begin{eqnarray}  \label{eq:2.51}
  \mathfrak{R}^*_0 & = & \overline{\mathbb{C}} \setminus (\gamma_{a_1,c} \cup \gamma_{a_2,c} \cup \gamma_{b,c}) \nonumber \\
  \mathfrak{R}^*_1 & = & \overline{\mathbb{C}} \setminus \gamma_1  \\
  \mathfrak{R}^*_2 & = & \overline{\mathbb{C}} \setminus (\gamma_{c,b^*}^{(1)} \cup \gamma_{a_2,c} \cup \gamma_{b,c}), \nonumber
\end{eqnarray}
and then we glue them together so that they form a Riemann surface of genus zero: \index{Surfaces!Rstar@$\mathfrak{R}^*$}
\[   \mathfrak{R}^*(A') = \overline{\mathfrak{R}^*_0 \cup \mathfrak{R}^*_1 \cup \mathfrak{R}^*_2}, \]
which we assign to the strictly non-acceptable triple $A'$. The algebraic functions $h$ and $\Phi$ for the
set of parameters $A'=\{a_1;a_2,b\}$ are single valued rational functions on $\mathfrak{R}^*(A')$. The structure
of the three sheets (\ref{eq:2.51}) defines the global branches of $h$ and
$\Phi$:  \index{Functions!h0h1h2@$h_0$, $h_1$, $h_2$} \index{Functions!Phi0Phi1Phi2@$\Phi_0$, $\Phi_1$, $\Phi_2$}
\begin{eqnarray}  \label{eq:2.52}
   h_0,\Phi_0 &\in& H(\mathbb{C} \setminus  (\gamma_{a_1,c} \cup \gamma_{a_2,c} \cup \gamma_{b,c})), \nonumber \\
   h_1,\Phi_1 &\in& H(\mathbb{C} \setminus  \gamma_1) , \\
   h_2, \Phi_2 &\in& H(\mathbb{C} \setminus  (\gamma_{c,b^*}^{(1)} \cup \gamma_{a_2,c} \cup \gamma_{b,c})). \nonumber
\end{eqnarray}

We denote \index{Contours!Delta1starDelta2star@$\Delta_1^*$, $\Delta_2^*$} \index{Contours!Delta12@$\Delta_{1,2}$}
\begin{equation*}
 \widetilde{\Delta}_1 :=\gamma_{a_1,c}, \quad \widetilde{\Delta}_2 :=\gamma_{a_2,c}, \quad
    \Delta_{1,2} :=\gamma_{c,b},
\end{equation*}
and \index{Contours!Delta1Delta2@$\Delta_1$, $\Delta_2$} \index{Contours!Delta0@$\Delta_0$}
\begin{equation}  \label{eq:2.53}
  \Delta_1:= \widetilde{\Delta}_1 \cup \Delta_{1,2}, \quad
  \Delta_2:= \widetilde{\Delta}_2 \cup \Delta_{1,2}, \quad
  \Delta_0:= \Delta_1 \cup \Delta_2.
\end{equation}
We also denote \index{Contours!E0E1E2@$E_0$, $E_1$, $E_2$}
\begin{equation}  \label{eq:2.53'}
  E_2:= \gamma_1 \setminus \widetilde{\Delta}_1, \quad
  E_1:=  \gamma_2 \setminus \widetilde{\Delta}_2, \quad
  E_0:= \gamma_3 \setminus \Delta_{1,2}.
\end{equation}

\begin{proposition}  \label{prop:2.8}
The contour $\Gamma$ (see \eqref{eq:2.29}--\eqref{eq:2.30}) has the following structure:
\index{Contours!Gamma@$\Gamma$}\index{Contours!Gamma01Gamma02Gamma12@$\Gamma_{0,1}$, $\Gamma_{0,2}$, $\Gamma_{1,2}$}
\begin{eqnarray}
    \Gamma_{0,1} &:=& \{z: |\Phi_0(z)| = |\Phi_1(z)|\} = \widetilde{\Delta}_1, \nonumber \\
    \Gamma_{0,2} &:=& \{z: |\Phi_0(z)| = |\Phi_2(z)|\} = \Delta_2, \label{eq:2.54'} \\
    \Gamma_{1,2} &:=& \{z: |\Phi_1(z)| = |\Phi_2(z)|\} = \bigcup_{\ell=1}^3 \gamma_{b^*,c}^{(\ell)}
       = E_0 \cup E_1 \cup E_2 , \nonumber
\end{eqnarray}
and for the domains $\Omega_{j,k\ell}$ we have (see Figure~{\rm \ref{fig:14}}, third picture) \index{Domains!Omegajkl@$\Omega_{j,k,\ell}$}
\[   \Omega_{0,2,1} = G \cup \widetilde{G}, \]
where the components $G$ and $\widetilde{G}$ are defined by their
 boundaries \index{Domains!GGtilde@$G$, $\widetilde{G}$}
\begin{equation}  \label{eq:2.54}
   \partial G = E_1 \cup E_2, \quad \partial \widetilde{G} = E_0 \cup E_2,
\end{equation}
and
\[   \Omega_{0,1,2} = \overline{\mathbb{C}} \setminus \overline{\Omega}_{0,2,1}. \]
\end{proposition}

Now we introduce the measures.

\begin{theorem} \label{thm:2.5}
For $A \in \textup{IV}$ we have
\begin{enumerate}
\item The jump of $h_0$ on $\Delta_0$ produces a positive measure $\lambda$ of
total mass 2 \index{Measures!lambda@$\lambda$}
\[    \frac{1}{2\pi i} \Bigl( h_{0+}(\xi) - h_{0-}(\xi) \Bigr)\, d\xi =: d\lambda(\xi), \qquad \xi \in \Delta_0. \]
The measure $\lambda$ consists of two measures $\tilde{\lambda}_1$ and $\lambda_2$,
with \index{Measures!lambda1tildelambda2@$\tilde{\lambda}_1$, $\lambda_2$}
$|\tilde{\lambda}_1| + |\lambda_2|=2$, supported on $\widetilde{\Delta}_1$ and $\Delta_2$ respectively
\[   \lambda := \begin{cases}
            \tilde{\lambda}_1 & \textrm{on $\widetilde{\Delta}_1$}, \\
            \lambda_2 & \textrm{on $\Delta_2$},
      \end{cases} \]
and \index{Functions!m1m2@$m_1$, $m_2$}
\begin{eqnarray*}
   \tilde{\lambda}_1(\xi) &=& \frac{m_1(\xi)}{\sqrt{\xi-a_1}}, \qquad m_1 \in H(\widetilde{\Delta}_1), \\
   \lambda_2(\xi) &=& \frac{m_2(\xi)}{\sqrt{(\xi-a_2)(\xi-b)}}, \qquad m_2 \in H(\Delta_2).
\end{eqnarray*}
\item The jump of $h_2$ on $E_2$ produces a positive measure $\mu_2$ \index{Measures!mu2@$\mu_2$}
\[   \frac{1}{3\pi i} \Bigl( h_{2+}(\xi)-h_{2-}(\xi) \Bigr)\, d\xi =: d\mu_2(\xi), \qquad \xi \in E_2, \]
and \index{Functions!m3@$m_3$}
\[   \mu_2'(\xi) = \sqrt{\xi-b^*} \ m_3(\xi), \qquad m_3 \in H(E_2).  \]
\item There are connections among the total masses of these measures:
\[    |\tilde{\lambda}_1| + |\lambda_2|=2, \quad |\lambda_2|- |\mu_2|=1. \]
\end{enumerate}
\end{theorem}

We also denote
\index{Measures!lambda12lambda1lambda2tilde@$\lambda_{1,2}$, $\lambda_1$, $\tilde{\lambda}_2$}
\begin{equation}   \label{eq:2.55}
   \lambda_{1,2}:= \lambda\bigr|_{\Delta_{1,2}}, \quad \lambda_1 := \lambda\bigr|_{\Delta_1},
   \quad  \tilde{\lambda}_2:= \lambda\bigr|_{\Delta_2}.
\end{equation}

The proof of Proposition \ref{prop:2.8} and  Theorem \ref{thm:2.5} is
given in Subsection \ref{sec:3.1.5}.

\begin{remark}[\textbf{on critical cases}]
There are several cases of location of the branch points $A=\{a_1,b;a_2,b\}$ which are
the limits of the cases described above. For example, the case (considered by Kalyagin in \cite{16})
when \index{Point sets!A0@$A_0$}
\[   A_0 := \{ -1,0;1,0\}  \]
is the limit of the class I
\[  A_{\epsilon,\epsilon} := \{-1,-\epsilon;1,\epsilon\}, \qquad \textrm{as $\epsilon \to 0$}, \]
and at the same time it is the limit of the class $\textup{III}$
\[  A_{0,\epsilon} = \{-1,0;1,\epsilon\}, \qquad \textrm{as $\epsilon \to 0$}. \]
Also, the critical triple is both the limit of the strictly non-acceptable triple (case IV) and
of the case V. We call these cases \textit{the critical cases}. For the critical cases the
assignment of the Riemann surface and the definition of the global branches for the algebraic
functions $h$ and $\Phi$ is carried out by passing to the limit in the non-critical cases
described above. Even though one of the main ingredients of our analysis, namely the geometry of
the problem, is clear for the critical cases, we do not consider the \HP\ asymptotics for these
cases in this paper. (In what follows we exclude the critical triples from the class IV.)
The point is that the rigorous proof of the asymptotics for these cases requires a non-standard
local Riemann-Hilbert analysis.

Certain critical cases of $2\times 2$ matrix-valued \RH problems have been studied recently
and it was found that the local Riemann-Hilbert analysis can be constructed in an explicit way
with the isomonodromy characterization of certain Painlev\'e transcendents,
see \cite{69,70,71,72,73,75,74}.
Similar constructions can be expected to apply to some of the critical cases related to
this work, but there are other critical cases as well which are special for
$3 \times 3$ matrix-valued \RH problems and which do not occur for $2\times 2$ problems.

\end{remark}

\subsection{Weak asymptotics, convergence and vector potential problems}  \label{sec:2.5}
In this section we formulate corollaries from the strong asymptotics of the \HP\ polynomials
regarding the weak asymptotics of the poles of the \HP\ approximants and their convergence.
The weak results follow from the strong asymptotic formulas (Theorems~\ref{thm:2.11} and \ref{thm:2.12})
but they are stated here since to present these results we only need the basic functions
which we have defined during the description of the geometry of the problem.
The results will be different for each
geometrical case.

\subsubsection{Weak convergence} \label{sec:2.5.1}
We recall that we consider the \HP\ approximants for two functions (\ref{eq:2.18}),
see Definition~\ref{def:2.1} and convention (\ref{eq:2.16}): \index{Functions!f1f2@$f_1$, $f_2$}
\begin{equation}  \label{eq:2.59}
   f_j \in \mathcal{A}(a_j,\alpha_j;b_j,\beta_j;\Omega_j), \qquad j=1,2,
\end{equation}
where the domains $\Omega_j$ for the analytic continuation of the weight functions $w_{0,j}$, defining
$f_j$ (see (\ref{eq:2.5}) and (\ref{eq:2.2})), depend on the location of the branch points
\index{Parameters!a1a2@$a_1$, $a_2$}
\index{Parameters!b1b2@$b_1$, $b_2$}
\[   A = \{a_1,b_1;a_2,b_2\}. \]
We assume $\Omega_j$ is such that: \index{Domains!Omega1Omega2@$\Omega_1$, $\Omega_2$}
\begin{eqnarray} \label{eq:2.60}
   1) & & \Omega_j \supset \Delta_j, \quad j=1,2, \qquad \textrm{for $A \in \textup{I}$} \nonumber \\
   2) & & \Omega_j \supset \Delta_j \cup G, \quad j=1,2, \qquad \textrm{for $A \in \textup{II}$}, \\
   3) & & \begin{cases}
          \Omega_1 \supset \Delta_1^*, \\
          \Omega_2 \supset \Delta_2,
           \end{cases} \qquad \textrm{for $A \in \textup{III}$}, \nonumber \\
   4) & & \Omega_j \supset \Delta_j \cup G, \quad j=1,2, \qquad \textrm{for $A \in \textup{IV}$}.  \nonumber \\
  5) & & \begin{cases}
           \Omega_1 \supset \Delta_1^* \cup G_1, \\
           \Omega_2 \supset \Delta_2 \cup G_1,
          \end{cases}  \qquad \textrm{for $A \in \textup{V}$} \nonumber
 \end{eqnarray}
We recall that the definition of the classes $\textup{I}$, $\textup{II}$,
$\textup{III}$, IV, and $\textup{V}$ is given
in Definitions~\ref{def:2.2'}, \ref{def:2.5}, \ref{def:2.5*} and the corresponding arcs
$\Delta_1$, $\Delta_1^*$, $\Delta_2$ and sets $G_1$, $G$ are defined in (\ref{eq:2.35'}), (\ref{eq:2.375}),
(\ref{eq:2.373}), (\ref{eq:2.44*}), Proposition~\ref{prop:2.7}, (\ref{eq:2.53}), and (\ref{eq:2.54}).

The weak limit of the counting measures $\nu_{P_n}$, which have equal mass $1/(2n)$
at the poles of the \HP\ approximants (see (\ref{eq:1.4})), has a universal character.
\begin{theorem}  \label{thm:2.6}
Suppose that $A$ belongs to one of the cases I, II, III, IV, or V.
Then the poles of the \HP\ approximants \eqref{eq:1.2} for the functions
in \eqref{eq:2.59}--\eqref{eq:2.60} have a weak limit
\begin{equation}   \label{eq:2.61}
  \nu_{P_n} \stackrel{*}{\to} \lambda/2, \qquad n \to \infty,
\end{equation}
where the limiting measure $\lambda$ is defined in Theorems~\ref{thm:2.1}, \ref{thm:2.1B},
\ref{thm:2.2}, \ref{thm:2.5}, depending on the geometrical class.
\end{theorem}

Now we state a result about the finite zeros of the functions of the second
kind (\ref{eq:1.2}) \index{Functions!Pn@$P_n$}\index{Functions!Qnj@$Q_n^{(1)}$, $Q_n^{(2)}$}\index{Functions!Rnj@$R_n^{(1)}$, $R_n^{(2)}$}
\[    R_n^{(j)} := f_j P_n - Q_n^{(j)}, \qquad j=1,2. \]
These zeros represent the extra interpolation points. We use the notation
$\nu_{R_n^{(j)}}$ for the counting measures with equal mass $1/n$ at the finite
zeros of $R_n^{(j)}$. \index{Measures!nuRnj@$\nu_{R_n^{(1)}}$, $\nu_{R_n^{(2)}}$}
\begin{theorem}  \label{thm:2.7}
Consider $R_n^{(j)}$ for the functions (\ref{eq:2.59}) where
$\Omega_j$ is as in (\ref{eq:2.60}).
\begin{enumerate}
\item When $A \in \textup{I} \cup \textup{III}$ there are no extra
interpolation points for $n$ large enough.
\item When $A \in \textup{II} \cup \textup{IV}$ we have
\begin{equation}  \label{eq:2.62}
    \nu_{R_n^{(j)}} \stackrel{*}{\to} \mu_j, \qquad j=1,2.
\end{equation}
\item When $A \in \textup{V}$ the function $R_n^{(2)}$ has no finite zeros
for $n$ large enough and
\[   \nu_{R_n^{(1)}} \stackrel{*}{\to} \mu_1.  \]
\end{enumerate}
In 2 and 3 the limiting measures $\mu_1,\mu_2$ are defined in Theorems~\ref{thm:2.1B}, \ref{thm:2.5}
and in (\ref{eq:2.40B}), and (\ref{eq:2.48}), depending on
the geometrical class.
\end{theorem}

The next theorem describes the $n$th root asymptotics of the error term
\[      f_j - \pi_n^{(j)} = \frac{R_n^{(j)}}{P_n},
    \quad \pi_n^{(j)} = \frac{Q_n^{(j)}}{P_n}, \qquad j=1,2, \]
and the convergence of the \HP\ approximants (\ref{eq:1.2}). \index{Functions!pinj@$\pi_n^{(1)}$, $\pi_n^{(2)}$}

\begin{theorem} \label{thm:2.8}
Consider the \HP\ approximants for the functions (\ref{eq:2.59})
with $\Omega_j$  as in (\ref{eq:2.60}).
\begin{enumerate}
\item When $A \in \textup{I}$ we have
\[   \left| f_j - \pi_n^{(j)} \right|^{1/n} \to  \left| \frac{\Phi_0}{\Phi_j} \right|,  \]
uniformly on compact subsets of $\overline{\mathbb{C}} \setminus \Delta_j$, $j=1,2$, and
therefore for the \HP\ approximants $\pi_n^{(j)}$ we have
\[   \pi_n^{(j)} \to f_j, \qquad \textrm{on $ \Omega_{0,j}$}, \]
\[    | \pi_n^{(j)} - f_j | \to \infty, \qquad \textrm{on $
    \overline{\mathbb{C}} \setminus (\Delta_j \cup \overline{\Omega_{0,j}})$}, \quad j=1,2, \]
with geometric rate. Furthermore (see Figure~\ref{fig:11}) there exists an $A \in \textup{I}$ such
that  $\Omega_{1,0} \neq \emptyset$ or $\Omega_{2,0} \neq \emptyset$.
\item When $A \in \textup{III}$ we have
\[      \left| f_j - \pi_n^{(j)} \right|^{1/n} \to \left| \frac{\Phi_0}{\Phi_j} \right|  \]
uniformly on compact subsets of $\overline{\mathbb{C}} \setminus \Delta_1^*$ when $j=1$ and
$\overline{\mathbb{C}} \setminus \Delta_2$ when $j=2$, and therefore, with the
convention (\ref{eq:2.38}), we have
\[   \pi_n^{(1)} \to f_1 \qquad \textrm{on $\overline{\mathbb{C}} \setminus (\Delta_1^* \cup \overline{\Omega_{1,0}})$}, \]
\[  |\pi_n^{(1)}-f_1| \to \infty \qquad \textrm{on $\Omega_{1,0}$}, \]
\[   \pi_n^{(2)} \to f_2 \qquad \textrm{on $\overline{\mathbb{C}} \setminus \Delta_2$}. \]
Furthermore (see Figure~\ref{fig:13}) $\Omega_{1,0} \neq \emptyset$.
\item When $A \in \textup{V}$ we have (uniformly on compact subsets of the indicated sets)
\[  \left| f_1 - \pi_n^{(1)} \right|^{1/n} \to
    \begin{cases}
  \displaystyle \left| \frac{\Phi_0}{\Phi_1} \right| &  \textrm{on $(\overline{\mathbb{C}} \setminus \Delta_1^*)
    \setminus \overline{G}_1$} \\[12pt]
  \displaystyle \left| \frac{\Phi_0}{\Phi_2} \right| &  \textrm{on $G_1\setminus \Delta_2$}
    \end{cases} \]
\[   \left| f_2 - \pi_n^{(2)} \right|^{1/n} \to
       \left| \frac{\Phi_0}{\Phi_2} \right| \qquad  \textrm{on $\overline{\mathbb{C}} \setminus \Delta_2$} \]
and therefore (see Figure~\ref{fig:14})
\[  \pi_n^{(1)} \to f_1 \qquad \textrm{on $\overline{\mathbb{C}} \setminus (\Delta_1^* \cup \overline{\Omega_{1,0}})$} \]
\[   |\pi_n^{(1)}-f_1| \to \infty \qquad \textrm{on $\Omega_{1,0}$}, \]
where $\Omega_{1,0} \neq \emptyset$, and
\[   \pi_n^{(2)} \to f_2 \qquad \textrm{on $\overline{\mathbb{C}} \setminus \Delta_2$}. \]
\item When $A \in \textup{IV}$ we have (uniformly on compact subsets of the indicated sets)
\[   \left| f_1 - \pi_n^{(1)} \right|^{1/n} \to
    \begin{cases}
  \displaystyle   \left| \frac{\Phi_0}{\Phi_1} \right| &
        \textrm{on $(\overline{\mathbb{C}} \setminus \widetilde{\Delta}_1)
    \setminus G$} \\[12pt]
   \displaystyle  \left| \frac{\Phi_0}{\Phi_2} \right| &  \textrm{on $G\setminus \Delta_{1,2}$}
    \end{cases} \]
\[  \left| f_2 - \pi_n^{(2)} \right|^{1/n} \to
       \left| \frac{\Phi_0}{\Phi_2} \right| \qquad  \textrm{on $\overline{\mathbb{C}} \setminus
     ( \Delta_2\cup E_2 )$} \]
and therefore we have
\[   \pi_n^{(j)} \to f_j \qquad \textrm{on $\overline{\mathbb{C}} \setminus \Delta_j$}, \quad j=1,2. \]
\item When $A \in \textup{II}$ we have (uniformly on compact subsets of the indicated sets)
\[    \left| f_1 - \pi_n^{(1)} \right|^{1/n} \to
     \left| \frac{\Phi_0}{\Phi_1} \right| \qquad \textrm{on $\overline{\mathbb{C}} \setminus
     (\Delta_1 \cup E_1)$}, \]
\[ \left| f_2 - \pi_n^{(2)} \right|^{1/n} \to
    \begin{cases}
    \displaystyle  \left| \frac{\Phi_0}{\Phi_2} \right| &  \textrm{on $(\overline{\mathbb{C}} \setminus
   \Delta_2) \setminus G$}  \\[12pt]
     \displaystyle \left| \frac{\Phi_0}{\Phi_1} \right| &  \textrm{on $ G \setminus \Delta_{1,2}$}
    \end{cases} \]
and therefore we have
\[  \pi_n^{(j)} \to f_j \qquad \textrm{on $ \overline{\mathbb{C}} \setminus \Delta_j$}, \quad j=1,2. \]
\end{enumerate}
\end{theorem}

We recall that the definition of the branches of the algebraic function $\{ \Phi_0, \Phi_1, \Phi_2 \}$
and the domains of divergence $\Omega_{j,0}$, $j=1,2$, depending on the geometrical case is given
in (\ref{eq:2.37}), (\ref{eq:2.25})--(\ref{eq:2.26}), (\ref{eq:2.36'}), (\ref{eq:2.376}), (\ref{eq:2.46}), and
(\ref{eq:2.52}).

\subsubsection{Vector equilibrium problem} \label{sec:2.5.2}
Concluding this section we state a universal \textbf{vector equilibrium problem} for the logarithmic potentials
of the measures $\lambda=\lambda_1+\lambda_2$, $\mu_1$, $\mu_2$. These measures were introduced in
Theorems \ref{thm:2.1}, \ref{thm:2.1B}, \ref{thm:2.2}, and \ref{thm:2.5} for the description of
the weak limits of the poles and of the extra interpolation points (see Theorems \ref{thm:2.6} and \ref{thm:2.7}).

\begin{theorem}  \label{thm:2.9}
Suppose that
\[    A \in \textup{I} \cup \textup{III} \cup \textup{II} \cup \textup{V} \cup \textup{IV} , \]
then
\begin{description}
\item[I.] \begin{enumerate}
  \item[a)] There exist piecewise analytic arcs $\Delta_1$, $\Delta_2$ which make the functions
 $f_1,f_2$ in (\ref{eq:2.59}) holomorphic
 \[    f_j \in H(\overline{\mathbb{C}} \setminus \Delta_j), \qquad j=1,2, \]
 and a piecewise analytic contour $E$ containing the common part of $\Delta_1$ and $\Delta_2$
 \[   \Delta_{1,2}:= \Delta_1 \cap \Delta_2, \quad ( \Delta_{1,2} = \emptyset \Rightarrow E = \emptyset). \]
 \item[b)] There exists a triple of measures $(\lambda_1,\tilde{\lambda}_2,\mu_1)$ with supports
(we use the notation $S(\mu)=\supp(\mu)$)
 \[   S(\lambda_1) \subset \Delta_1, \quad
      S(\tilde{\lambda}_2) \subset \widetilde{\Delta}_2 := \Delta_2 \setminus \Delta_{1,2}, \quad
      S(\mu_1) \subset E, \]
 and with the relations on their total mass
 \[  \begin{cases}
      |\lambda_1| + |\tilde{\lambda}_2| = 2 \\
      |\lambda_1| - |\mu_1| = 1.
      \end{cases}  \]
 \item[c)] This triple of measures possesses the following equilibrium relations with some constants
  $\kappa_1$ and $\tilde{\kappa}_2$: \index{Parameters!kappa1kappa2tilde@$\kappa_1$, $\tilde{\kappa}_2$}
  \index{Functions!U1U2U3@$U_1$, $U_2$, $U_3$}
 \begin{eqnarray} \label{eq:2.70}
  U_1 &:=& 2V^{\lambda_1} + V^{\tilde{\lambda}_2} - V^{\mu_1}
   \begin{cases}
   =   \kappa_1, & \textrm{on $S(\lambda_1)$}, \\
   \geq \kappa_1, & \textrm{on $\Delta_1$},
   \end{cases}  \nonumber \\
   U_2 &:=& V^{\lambda_1} + 2V^{\tilde{\lambda}_2} + V^{\mu_1}
   \begin{cases}
   =   \tilde{\kappa}_2, & \textrm{on $S(\tilde{\lambda}_2)$}, \\
   \geq \tilde{\kappa}_2, & \textrm{on $\widetilde{\Delta}_2$},
   \end{cases}  \\
    U_3 &:=& -V^{\lambda_1} + V^{\tilde{\lambda}_2} + 2 V^{\mu_1}
   \begin{cases}
   =   \tilde{\kappa}_2-\kappa_1, & \textrm{on $S(\mu_1)$}, \\
   \geq \tilde{\kappa}_2-\kappa_1, & \textrm{on $E$}.
   \end{cases}  \nonumber
 \end{eqnarray}
 \item[d)] The supports of the measures possess the following symmetry relations
  \begin{equation}  \label{eq:2.71}
    \begin{cases}
   \displaystyle \frac{\partial U_1}{\partial n_+} = \frac{\partial U_1}{\partial n_-}, &
         \textrm{on $S(\lambda_1)$}, \\[10pt]
   \displaystyle \frac{\partial U_2}{\partial n_+} = \frac{\partial U_2}{\partial n_-}, &
         \textrm{on $S(\tilde{\lambda}_2)$}, \\[10pt]
   \displaystyle \frac{\partial U_3}{\partial n_+} = \frac{\partial U_3}{\partial n_-}, &
         \textrm{on $S(\mu_1)$},
     \end{cases}
  \end{equation}
  where $\partial/\partial n_{\pm}$ denotes the normal derivatives on the respective contours.
  \end{enumerate}
\item[II.] There is also a dual problem regarding the triple $(\lambda_2,\tilde{\lambda}_1,\mu_2)$ which
 can be obtained from the problem I.a)--I.d) by interchanging the indices 1 and 2.
\item[III.] The equilibrium measures $(\lambda_1,\tilde{\lambda}_2,\mu_1)$ and $(\lambda_2,\tilde{\lambda}_1,\mu_2)$
are related as follows
\[   \lambda:=\lambda_1 + \tilde{\lambda}_2 = \lambda_2 + \tilde{\lambda}_1, \quad
     S(\lambda) \subset \Delta_0 := \Delta_1 \cup \Delta_2, \]
\[ \mu:= \mu_1+\mu_2, \quad S(\mu_1) \cup S(\mu_2) = E, \]
\[ V^{\mu}\Bigr|_{E} = V^{\lambda_{1,2}}\Bigr|_E, \quad \lambda_{1,2} := \lambda \Bigr|_{\Delta_{1,2}}. \]
\item[IV.] The measure $\lambda/2$ is the weak limit (\ref{eq:2.61}) of the poles of the \HP\ approximants
of the functions (\ref{eq:2.59}), and the measures $\mu_1$ and $\mu_2$ are the weak limits (\ref{eq:2.62})
of the extra interpolation points.
\end{description}
\end{theorem}
The proof of Theorem \ref{thm:2.9} is given in Section \ref{sec:3.2}.

We now describe the basic notions of the equilibrium problems I and II of Theorem~\ref{thm:2.8}
\[  \begin{cases}
     \Delta_1, \widetilde{\Delta}_2,E \\
     \lambda_1, \tilde{\lambda}_2, \mu_1
    \end{cases}
    \qquad
    \begin{cases}
     \Delta_2, \widetilde{\Delta}_1, E \\
     \lambda_2, \tilde{\lambda}_1, \mu_2
    \end{cases}  \]
in conformity with the geometrical cases.

\begin{remark}  \label{rem:2.3}
\
\begin{enumerate}
\item In the case $\textup{I}$ we have $\Delta_{1,2} = \emptyset$ and therefore
\begin{equation}  \label{eq:2.73}
    \begin{cases}
      \Delta_j = \widetilde{\Delta}_j, \quad j=1,2, \\
       E =  \emptyset
    \end{cases}
\end{equation}
\[   \lambda_j = \tilde{\lambda}_j, \quad \mu_j = 0, \quad S(\lambda_j) = \Delta_j, \qquad j=1,2, \]
and in the systems of the equilibrium and symmetry relations (\ref{eq:2.70})--(\ref{eq:2.71}) only
the first two relations are needed.
\item In the case $\textup{III}$, we set \index{Contours!Delta1star@$\Delta_1^*$}\index{Contours!delta1star@$\delta_1^*$}
\begin{equation}  \label{eq:2.73'}
   \Delta_1 = \Delta_1^* \cup \delta_1^*, \quad \delta_1^* \subset \Omega_{1,0}
\end{equation}
where $\delta_1^*$ is an arbitrary rectifiable arc in $\Omega_{1,0,2}$ joining the points $b^*$ and $b_1$
(see (\ref{eq:2.37}) and Figure~\ref{fig:13}). For this case we have the same setting (\ref{eq:2.73})
as for the case I with the minor difference that $S(\lambda_1) = \Delta_1^* \subset \Delta_1$.
\item We note that for the cases $\textup{I}$ and $\textup{III}$ the equilibrium relations (\ref{eq:2.70}) reduce
to the equilibrium relations for an Angelesco system (\ref{eq:1.13}), (\ref{eq:1.16}) in the complex
plane.
\item In the case $\textup{V}$ we set (see Proposition~\ref{prop:2.7})
\index{Contours!Delta1Delta1tildeDelta2tilde@$\Delta_1$, $\widetilde{\Delta}_1$, $\widetilde{\Delta}_2$}
\[   \Delta_1:= \widetilde{\Delta}_1^* \cup \delta_1^* \cup \Delta_{1,2}, \quad
     \widetilde{\Delta}_j := \Delta_j \setminus \Delta_{1,2}, \quad j=1,2 , \]
where $\delta_1^*$ is an arbitrary rectifiable arc in $\Omega_{1,0,2}$ and joining the
corresponding end points of the arcs $\Delta_1^*$ and $\Delta_{1,2}$. Here we have \index{Contours!E@$E$}
\[    E = E_1 \neq \emptyset. \]
If we consider the dual problem (see II of Theorem~\ref{thm:2.9}) for the triple
\[   \begin{cases}
      \Delta_2, \widetilde{\Delta}_1, E  \\
       \lambda_2, \tilde{\lambda}_1, \mu_2
     \end{cases}  \]
then the solution of the Angelesco equilibrium problem for the case $\textup{III}$ \index{Parameters!kappa1tildekappa2@$\tilde{\kappa}_1$, $\kappa_2$}
\[ \begin{cases}
   2 V^{\lambda_2} + V^{\tilde{\lambda}_1} = \kappa_2, & \textrm{on $\Delta_2$}, \\
   V^{\lambda_2} + 2V^{\tilde{\lambda}_1} \geq \tilde{\kappa}_1, & \textrm{on $\Delta_1$}
  (=\kappa_1 \quad  \textrm{on $\Delta_1^*$})
   \end{cases}  \]
provides the solution of the dual equilibrium problem with $\mu_2=0$, because (see Figure~\ref{fig:14}, second picture)
\[   V^{\tilde{\lambda}_1} - \tilde{\kappa}_1 = V^{\lambda_2}-\kappa_2, \qquad \textrm{on $E$}. \]
In order to obtain the measure $\mu_1$ we have to take the balayage of the measure \index{Measures!lambda12@$\lambda_{1,2}$}
\[  \lambda_{1,2}:= \lambda_2 \Bigr|_{\Delta_{1,2}}  \]
on the contour $E$ or consider the problem (\ref{eq:2.70})--(\ref{eq:2.71}) for the triple
\[   \begin{cases}
      \Delta_1, \widetilde{\Delta}_2, E, \\
      \lambda_1, \tilde{\lambda}_2, \mu_1.
     \end{cases}      \]
\item In the cases IV and $\textup{II}$ the contours for the equilibrium problems are defined
in (\ref{eq:2.53})--(\ref{eq:2.53'}) and (\ref{eq:2.372})--(\ref{eq:2.375}). In these cases
(see Theorem~\ref{thm:2.1B} and \ref{thm:2.5}) there are no degeneracies of the components of the
equilibrium problem (\ref{eq:2.70}) and the dual problem.
\end{enumerate}
\end{remark}

\subsection{Szeg\H{o} functions for the \HP\ polynomials}  \label{sec:2.6}
The Szeg\H{o} function is an important ingredient of the strong asymptotic formulas
for orthogonal polynomials on an interval of the real line and on the unit circle.
In this section we introduce a generalization of the Szeg\H{o} function which we need
for the presentation of the strong asymptotics of the \HP\ polynomials. We shall define
the Szeg\H{o} function as a solution of a certain \BVP\ on the three-sheeted Riemann
surface
\[   \mathfrak{R} := \bigcup_{j=0}^2 \overline{\mathfrak{R}}_j.  \]
Before we state the \BVP\ we need some preparation. We fix the cuts $\delta_{j,k}$ (Jordan
arcs in the complex plane or a union of Jordan arcs) where the sheets $\mathfrak{R}_j$
and $\mathfrak{R}_k$ are glued together \index{Contours!deltalj@$\delta_{\ell,j}$}
\[  \mathfrak{R}_\ell =: \overline{\mathbb{C}} \setminus (\delta_{\ell,j} \cup \delta_{\ell,k}),
    \qquad \ell \notin \{k,j\}, \ k\neq j, \ \ell,k,j=0,1,2. \]
The two sides of the cuts $(\delta_{\ell,j} \cup \delta_{\ell,k})^{(\pm)}$ form a boundary
of the sheet $\mathfrak{R}_\ell$ (Jordan curves of the Riemann surface $\mathfrak{R}$)
\index{Contours!Rlpartial@$\partial \mathfrak{R}_{\ell}$} \index{Contours!Rlkpartial@$\partial \mathfrak{R}_{\ell,k}$}
\[  \partial \mathfrak{R}_\ell := \partial \mathfrak{R}_{\ell,j} \cup \partial \mathfrak{R}_{\ell,k}, \quad
    \partial \mathfrak{R}_{\ell,k} :=  \overline{\mathfrak{R}}_\ell \cap  \overline{\mathfrak{R}}_k
    = \delta_{\ell,k}^{(+)} \cup \delta_{\ell,k}^{(-)} . \]
Thus \index{Contours!deltalkpm@$\delta_{\ell,k}^{(+)}$, $\delta_{\ell,k}^{(-)}$} \index{Contours!Rpartial@$\partial \mathfrak{R}$}
\[  \mathfrak{R} =: \left( \bigcup_{j=0}^2 \mathfrak{R}_j \right) \cup \partial \mathfrak{R}, \quad
    \partial \mathfrak{R} := \partial \mathfrak{R}_{0,1} \cup \partial \mathfrak{R}_{0,2}
\cup \partial \mathfrak{R}_{1,2}. \]
We fix the orientation of the union of the Jordan curves $\partial \mathfrak{R}$ such
that \index{Contours!Rjkpartialpm@$\partial \mathfrak{R}_{j,k}^{(+)}$, $\partial \mathfrak{R}_{j,k}^{(-)}$}
\[   \partial \mathfrak{R}_{j,k}^{(+)} \subset \overline{\mathfrak{R}}_j, \quad
     \partial \mathfrak{R}_{j,k}^{(-)} \subset \overline{\mathfrak{R}}_k, \qquad j > k,\ j,k=0,1,2. \]
We also fix an orientation of the Jordan arcs $\delta_{j,k}$.
Our Riemann surfaces have four branch points at $\{a_1,b_1,a_2,b_2\}$, with a possible
replacement of $b_1$ by $b^*$ for the cases III, IV and V. We have
\begin{enumerate}
\item $\{a_j,b_j\} \subset \delta_{0,j}$, $j=1,2$ for $A \in \textup{I} \cup \textup{II}$,
\item $\{a_1,b^*\} \subset \delta_{0,1}$, $\{a_2,b\} \subset \delta_{0,2}$ for $A \in \textup{III}
\cup \textup{V}$,
\item $a_1 \in \delta_{0,1}$, $b^* \in \delta_{1,2}$, $\{a_2,b\} \subset \delta_{0,2}$ for $A \in \textup{IV}$.
\end{enumerate}
Let $\omega_j(z)$ be the branch of the function \index{Functions!omega1omega2@$\omega_1$, $\omega_2$}
\begin{equation}   \label{eq:2.74}
    \omega_j^2(z) = (z-a_j)(z-b_j), \quad \omega_j(z) = z + \cdots, \quad z \to \infty, \qquad j=1,2.
\end{equation}
For the cases III, IV and V we replace $b_1$ by $b^*$ in (\ref{eq:2.74}).
 From the statement of the problem (see Section~\ref{sec:2.1}) and from the definitions of the geometrical
cases (Sections \ref{sec:2.2}--\ref{sec:2.5}) it follows that the weight functions $w_1$ and $w_2$ are defined on the
arcs $\delta_{0,1}$ and $\delta_{0,2}$. We set \index{Functions!w1tildew2tilde@$\tilde{w}_1$, $\tilde{w}_2$}
\begin{equation}  \label{eq:2.75}
    \tilde{w}_j:= i w_j \omega_{j-}, \qquad \textrm{on $\delta_{0,j}$}, \quad j=1,2.
\end{equation}
We assume that the orientation of the arcs $\delta_{k,j}$ is fixed so that (see, for example,
Figure~\ref{fig:10B}) each $\tilde{w}_j$ has an analytic continuation on $\delta_{1,2}$ (when
$\delta_{1,2} \neq \emptyset$). Then we define \index{Functions!wtilde@$\tilde{w}$}
\begin{equation}  \label{eq:2.75'}
   \tilde{w} := \begin{cases}
       \tilde{w}_1, & \textrm{on $\delta_{0,1}$}, \\
       \tilde{w}_2, & \textrm{on $\delta_{0,2}$}, \\
       \tilde{w}_2/\tilde{w}_1, & \textrm{on $\delta_{1,2}$}.
                \end{cases}
\end{equation}
We note that, depending on the position of the points $A$, we have
\begin{eqnarray*}
   & & A \in \textup{I} \cup \textup{III} \cup \textup{V} \Rightarrow \delta_{1,2} = \emptyset,  \\
   & & A \in \textup{II} \cup \textup{IV} \Rightarrow \delta_{1,2} \neq \emptyset.
\end{eqnarray*}
Finally we duplicate the weight function $\tilde{w}$ on both sides of $\delta_{j,k}^{(\pm)}$, i.e., we
define $\tilde{w}$ on $\partial \mathfrak{R}$
\begin{equation}  \label{eq:2.76}
  \tilde{w}_+ = \tilde{w}_- \quad \textrm{on $\delta_{k,j}$} \Rightarrow \tilde{w} \textrm{ on $\partial \mathfrak{R}$}.
\end{equation}

Now we formulate the \BVP\ on $\mathfrak{R}$ whose solution defines the desired generalization of the
Szeg\H{o} function. We are looking for a piecewise holomorphic function $\mathcal{F}$ on $\mathfrak{R}$
such that \index{Functions!F@$\mathcal F$}
\begin{equation}  \label{eq:2.77}
  \begin{cases}
  1. & \mathcal{F} \in H(\mathfrak{R}\setminus \partial \mathfrak{R}), \\
  2a. &\exists \mathcal{F}_{\pm} \in C(\partial \mathfrak{R} \setminus \{a_1,b_1,a_2,b_2\}):\
  \mathcal{F}_+ = \mathcal{F}_- \tilde{w}  \textrm{ on $\partial \mathfrak{R}$} \\
  2b. & |\mathcal{F}^2(z) \tilde{w}(z)| = \O(1), \quad z \to A \\
  3. & \mathcal{F}(\infty^{(0)}) \mathcal{F}(\infty^{(1)}) \mathcal{F}(\infty^{(2)}) = 1.
  \end{cases}
\end{equation}

\begin{remark}  \label{rem:2.4}
The maximum principle implies that if $\mathcal{F}$ from (\ref{eq:2.77}) exists, then
\begin{enumerate}
\item[a)] For every $z \in \overline{\mathbb{C}}$
\begin{equation} \label{eq:2.78}
  \mathcal{F}(z^{(0)}) \mathcal{F}(z^{(1)}) \mathcal{F}(z^{(2)}) = 1.
\end{equation}
\item[b)] $\mathcal{F}$ is unique.
\end{enumerate}
Here $z^{(j)} := \pi^{-1}(z)$, where $\pi: \mathfrak{R} \to \overline{\mathbb{C}}$ is the projection
of the sheets of $\mathfrak{R}$ onto $\overline{\mathbb{C}}$. \index{Functions!pi@$\pi$}
\end{remark}
Due to (\ref{eq:2.76}) the index of the \BVP\ (\ref{eq:2.77}) is equal to zero. This gives the
existence of the solution of (\ref{eq:2.77}). In order to write an expression for the solution
of (\ref{eq:2.77}) we use the meromorphic (Cauchy) differential on $\mathfrak{R}$ \index{Functions!Mxi@$dM_{\xi}(z)$}
\[  dM_\xi(z) , \qquad z \in \mathfrak{R}, \xi \in \mathfrak{R}, \]
which has simple poles at the points $\xi \in \mathfrak{R}$ and at the points $\tilde{\xi} \in \mathfrak{R}$ where
$\tilde{\xi}$ has the same projection on $\overline{\mathbb{C}}$ as $\xi$, i.e., $\pi(\xi) = \pi(\tilde{\xi})$,
$\xi \neq \tilde{\xi}$:
\[   M'_\xi(z) \to \infty, \qquad z \to \{ \xi, \tilde{\xi}^{(k)}, \tilde{\xi}^{(\ell)} \}, \]
and the residues at the poles are
\[  \res_{z=\xi} M'_\xi = 2, \quad \res_{z = \tilde{\xi}} M'_\xi = -1, \quad \pi(\xi) = \pi(\tilde{\xi}), \
\xi \neq \tilde{\xi}. \]
Then Cauchy's theorem on $\mathfrak{R}$ and (\ref{eq:2.78}) give
\[  \frac{1}{2\pi i} \int_{\partial \mathfrak{R}} \log \tilde{w}(z)\, dM_\xi(z) =
    2 \log \mathcal{F}(\xi) - \log \mathcal{F}(\tilde{\xi}^{(k)}) - \log \mathcal{F}(\tilde{\xi}^{(\ell)})
   = 3 \log \mathcal{F}(\xi). \]
Summarising we have proved
\begin{theorem} \label{thm:2.10}
There exists a unique solution of the \BVP\ (\ref{eq:2.77}), which is given by
\[   \mathcal{F}(\xi) = \exp \left( \frac{1}{6\pi i} \int_{\partial \mathfrak{R}} \log \tilde{w}(z)\, dM_\xi(z) \right),
  \qquad \xi \in \mathfrak{R} \setminus \partial \mathfrak{R}. \]
\end{theorem}

Later on we shall use the notation for the branches of the multivalued function $\mathcal{F}$
\begin{equation}  \label{eq:2.79}
   F_\ell(z) := \mathcal{F}(z^{(\ell)}), \qquad \ell=0,1,2,
\end{equation}
where $z \in \overline{\mathbb{C}} \setminus (\delta_{\ell,j} \cup \delta_{\ell,k})$,
$\ell \notin \{k,j\}$, $k \neq j$, and $\ell,k,j=0,1,2$. \index{Functions!F0@$F_0$}  \index{Functions!F1F2@$F_1$, $F_2$}

\subsection{Strong asymptotics}  \label{sec:2.7}
We can now state the main results of this paper on the strong asymptotics of the \HP\ approximants.
The formulas for strong asymptotics contain the functions
$\Phi:=\{\Phi_0,\Phi_1,\Phi_2\}$ from (\ref{eq:1.28}) and the Szeg\H{o} functions
$\mathcal F:=\{F_0,F_1,F_2\}$ from (\ref{eq:2.77}) and (\ref{eq:2.79}). These functions are defined
by the Riemann surface $\mathfrak{R}$ of the corresponding geometrical cases (defined in Sections
\ref{sec:2.2}--\ref{sec:2.5}).

The strong asymptotics of the denominators $P_n$ of the \HP\ approximants have a universal character.
\begin{theorem}  \label{thm:2.11}
Suppose that
\[   A = \{ a_1,b_1;a_2,b_2\} \in \textup{I} \cup \textup{II} \cup \textup{III} \cup \textup{IV}
  \cup \textup{V}. \]
Then the denominators of the \HP\ approximants \eqref{eq:1.2} for the functions $f_j$, with
branch points $\{a_1,b_1,a_2,b_2\}$ (see \eqref{eq:2.59}--\eqref{eq:2.60}) have the following
asymptotic formulas  as $n \to \infty$, which hold uniformly on compact subsets of the indicated
sets:
\[  P_n(z) = \frac{F_0(\infty)}{F_0(z)} \bigl(C_0 \Phi_0(z) \bigr)^{-n} \bigl(1+\O(1/n)\bigr),
  \qquad z \in  \overline{\mathbb{C}} \setminus (\delta_{0,1} \cup \delta_{0,2}), \]
and
\begin{multline*}
  P_n(z) = \left[ \left( \frac{F_0(\infty)}{F_0(z)} \bigl(C_0 \Phi_0(z) \bigr)^{-n} \right)_+
           +  \left( \frac{F_0(\infty)}{F_0(z)} \bigl(C_0 \Phi_0(z) \bigr)^{-n} \right)_- \right]
 \bigl(1+o(1)\bigr), \\
  \qquad z \in  (\delta_{0,1} \cup \delta_{0,2}) \setminus A.
\end{multline*}
\end{theorem}
We recall that we have
\begin{equation}  \label{eq:2.82}
  \begin{cases}
   \delta_{0,j} = \Delta_j, \quad j=1,2, & \textrm{for $A \in \textup{I}$}, \\
   \delta_{0,1} = \Delta_1^*, \delta_{0,2} = \Delta_2, & \textrm{for $A \in \textup{III} \cup \textup{V}$}, \\
   \delta_{0,1} = \Delta_1, \delta_{0,2} = \widetilde{\Delta}_2, & \textrm{for $A \in \textup{II}$}, \\
   \delta_{0,1} = \widetilde{\Delta}_1, \delta_{0,2} = \Delta_2, & \textrm{for $A \in \textup{IV}$}.
  \end{cases}
\end{equation}
The strong asymptotics for the functions of the second kind $R_n^{(j)}$, $j=1,2$, have
different forms, depending on the geometrical case.

\begin{theorem}  \label{thm:2.12}
For the asymptotics, as $n \to \infty$, of the functions of the second kind $R_n^{(j)}$, given by \eqref{eq:2.7},
of the \HP\ approximants for the functions $f_j$ with branch points $\{a_1,b_1,a_2,b_2\}$ in
\eqref{eq:2.59}--\eqref{eq:2.60}, we have
\begin{enumerate}
\item When $A \in \textup{I} \cup \textup{III}$ we have
\begin{eqnarray*}
   R_n^{(j)}(z) &=& -\frac{F_0(\infty)}{\omega_j(z)F_j(z)} \bigl( C_0 \Phi_j(z) \bigr)^{-n} \bigl( 1+\O(1/n) \bigr),
    \qquad z \in  \overline{\mathbb{C}} \setminus \delta_{0,j}, \\
   R_{n\pm}^{(j)}(z) &=& \left( -\frac{F_0(\infty)}{\omega_j(z)F_j(z)} \bigl( C_0 \Phi_j(z) \bigr)^{-n}
\bigl( 1+\O(1/n) \bigr) \right)_{\pm},
    \qquad z \in  \delta_{0,j} \setminus A,
\end{eqnarray*}
for $j=1,2$. Here the $\delta_{0,j}$ are as in \eqref{eq:2.82}.
\item When $A \in \textup{V}$ we have
\begin{eqnarray*}
   R_n^{(2)}(z) &=& -\frac{F_0(\infty)}{\omega_2(z)F_2(z)} \bigl( C_0 \Phi_2(z) \bigr)^{-n} \bigl( 1+\O(1/n) \bigr),
    \qquad z \in  \overline{\mathbb{C}} \setminus \Delta_2, \\
   R_{n\pm}^{(2)}(z) &=& \left( -\frac{F_0(\infty)}{\omega_2(z)F_2(z)} \bigl( C_0 \Phi_j(z) \bigr)^{-n}
   \bigl( 1+\O(1/n) \bigr) \right)_{\pm},
    \qquad z \in  \Delta_2 \setminus A,
\end{eqnarray*}
and
\begin{eqnarray*}
   R_n^{(1)}(z) &=& -\frac{F_0(\infty)}{\omega_1(z)F_1(z)} \bigl( C_0 \Phi_1(z) \bigr)^{-n} \bigl( 1+\O(1/n) \bigr),
    \qquad z \in  (\overline{\mathbb{C}} \setminus \Delta_1^*)\setminus \overline{G}_1, \\
   R_n^{(1)}(z) &=&  -\frac{F_0(\infty)}{\omega_2(z)F_2(z)}  \bigl( C_0 \Phi_2(z) \bigr)^{-n} \frac{w_1(z)}{w_2(z)}
   \bigl( 1+\O(1/n) \bigr) ,     \qquad z \in  G_1 \setminus \Delta_{1,2}, \\
   R_{n\pm}^{(1)}(z) &= & \left( -\frac{F_0(\infty)}{\omega_1(z)F_1(z)} \bigl( C_0 \Phi_1(z) \bigr)^{-n}\right)_{\pm}
    \bigl( 1+\O(1/n) \bigr),     \qquad z \in  \Delta_1^* \setminus A, \\
   R_{n\pm}^{(1)}(z) &=& \left( -\frac{F_0(\infty)}{\omega_2(z)F_2(z)} \bigl( C_0 \Phi_2(z) \bigr)^{-n}
    \bigl( 1+\O(1/n) \bigr) \right)_{\pm}     \frac{w_1(z)}{w_2(z)} ,     \qquad z \in  \Delta_{1,2} \setminus A, \\
   R_{n}^{(1)}(z) & = & \frac{F_0(\infty)}{C_0^n} \left( -\frac{1}{\omega_1(z)F_1(z)\Phi_1^{n}(z)} -
   \frac{w_1(z)}{w_2(z)} \frac{1}{\omega_2(z)F_2(z)\Phi_2^n(z)} \right) \bigl( 1+\O(1/n) \bigr), \\
  && \qquad z \in \partial G_1 \setminus A.
\end{eqnarray*}
\item When $A \in \textup{IV}$ we have, with $\delta_{0,2}:=\Delta_2 \cup E_2$
\begin{eqnarray*}
   R_n^{(2)}(z) &=& -\frac{F_0(\infty)}{\omega_2(z)F_2(z)} \bigl( C_0 \Phi_2(z) \bigr)^{-n} \bigl( 1+\O(1/n) \bigr),
    \qquad z \in  \overline{\mathbb{C}} \setminus \delta_{0,2}, \\
   R_{n\pm}^{(2)} &= & \left( -\frac{F_0(\infty)}{\omega_2(z)F_2(z)} \bigl( C_0 \Phi_2(z) \bigr)^{-n}\right)_{\pm}
    \bigl( 1+\O(1/n) \bigr),     \qquad z \in  \Delta_2 \setminus A, \\
   R_{n}^{(2)}(z) & = & \frac{F_0(\infty)}{C_0^n} \left( -\frac{1}{(\omega_2(z)F_2(z)\Phi_2^{n}(z))_+} -
    \frac{1}{(\omega_2(z)F_2(z)\Phi_2^n(z))_-} \right) \bigl( 1+\O(1/n) \bigr), \\
    && \qquad z \in  E_2 \setminus \{b^*\},
\end{eqnarray*}
and with $G$ such that $\partial G = E_1 \cup E_2$
\begin{eqnarray*}
   R_n^{(1)}(z) &=& -\frac{F_0(\infty)}{\omega_1(z)F_1(z)} \bigl( C_0 \Phi_1(z) \bigr)^{-n} \bigl( 1+\O(1/n) \bigr),
    \qquad z \in  (\overline{\mathbb{C}} \setminus \widetilde{\Delta}_1)\setminus \overline{G}, \\
   R_n^{(1)}(z) &=&  -\frac{F_0(\infty)}{\omega_2(z)F_2(z)}  \bigl( C_0 \Phi_2(z) \bigr)^{-n} \frac{w_1(z)}{w_2(z)}
   \bigl( 1+\O(1/n) \bigr) ,     \qquad z \in  G \setminus \Delta_{1,2}, \\
   R_{n\pm}^{(1)} &= & \left( -\frac{F_0(\infty)}{\omega_1(z)F_1(z)} \bigl( C_0 \Phi_1(z) \bigr)^{-n}\right)_{\pm}
    \bigl( 1+\O(1/n) \bigr),     \qquad z \in  \widetilde{\Delta}_1, \\
   R_{n\pm}^{(1)}(z) &=& \left( -\frac{F_0(\infty)}{\omega_2(z)F_2(z)} \bigl( C_0 \Phi_2(z) \bigr)^{-n}
    \bigl( 1+\O(1/n) \bigr) \right)_{\pm}     \frac{w_1(z)}{w_2(z)} ,     \qquad z \in  \Delta_{1,2}, \\
   R_{n}^{(1)}(z) & = & \frac{F_0(\infty)}{C_0^n} \left( -\frac{1}{\omega_1(z)F_1(z)\Phi_1^{n}(z)} -
   \frac{w_1(z)}{w_2(z)} \frac{1}{\omega_2(z)F_2(z)\Phi_2^n(z)} \right) \bigl( 1+\O(1/n) \bigr), \\
   && \qquad z \in  E_1 \setminus \{b^*\}, \\
   R_n^{(1)}(z) &= & -\frac{F_0(\infty)}{C_0^n} \frac{1+\O(1/n)}{(\omega_1(z)F_1(z)\Phi_1^n(z))_+} \\
                &= & -\frac{F_0(\infty)}{C_0^n} \frac{1+\O(1/n)}{(\omega_2(z)F_2(z)\Phi_2^n(z))_-} \frac{w_1(z)}{w_2(z)},
    \qquad z \in  E_2 \setminus \{b^*\}.
\end{eqnarray*}
\item When $A \in \textup{II}$ we have the same asymptotic formulas as in 3) (the case $A \in \textup{IV}$) if we
interchange the indices 1 and 2 and replace $\{b^*\}$ by $\emptyset$.
\end{enumerate}
\end{theorem}

The proof of the strong asymptotics formulas is given in Section \ref{sec:4}.
It follows from a steepest descent analysis of the \RH problem (\ref{eq:2.13}).

\section{Proof of the geometrical results and equilibrium properties}  \label{sec:3}
\subsection{Proof of the geometrical results}  \label{sec:3.1}
Here we present the proofs of the propositions and the theorems stated in subsections
\ref{sec:2.2}--\ref{sec:2.4}. We start with some results of general character.

\subsubsection{Proof of Propositions \ref{prop:2.2}, \ref{prop:2.3} and \ref{prop:2.4}} \label{sec:3.1.1}
\begin{proof}
\
\begin{enumerate}
\item The parametrization of the algebraic curve (\ref{eq:2.20}) of order 3 and genus 0, given in
Proposition~\ref{prop:2.2}, can be checked directly. Indeed, substituting the polynomial
coefficients (\ref{eq:2.23'}) and (\ref{eq:2.24}) into the equation (\ref{eq:2.20}) for the
function $h$ and computing the discriminant
\[  \mathcal{D} := \frac{P_2^3 - \Pi_4 P_1^2}{\Pi_4^3}, \]
we obtain that the polynomial
\[  P_2^3 - \Pi_4 P_1^2 = \widetilde{\mathcal{D}} =
     (k P_1^2 + 3P_1 ps + s^3)^2 \]
is a complete square and therefore its zeros are not square root branch points of the algebraic function $h$ (the Puiseux
series at these points have integer exponents or cubic roots in some degenerate cases). Thus (in the non-degenerate cases) the only branch points of the function $h$
are the four zeros of the polynomial $\Pi_4$. The Riemann-Hurwitz formula (for a function of the third order
with four branch points) implies that $h$ has genus $0$.
\item To find the coefficients (\ref{eq:2.26}) for the equation of the algebraic function $\Phi$ in
Proposition~\ref{prop:2.3} we proceed as follows. We make a linear change of variable $z \to P_1:=z-c$
and substitute the series
\[   h = 1 + \sum_{k=1}^\infty \frac{c_k}{P_1^k}  \]
into the equation (\ref{eq:2.20}) to find two sets of coefficients $\{c_k^{(1)}\}$ and $\{c_k^{(2)}\}$.
Thus, we obtain the coefficients of the expansion at $P_1=\infty$ for the three branches of $h$
\[  1 + \sum_{k=1}^\infty \frac{c_k^{(1)}}{P_1^k}, \quad
    1 + \sum_{k=1}^\infty \frac{c_k^{(2)}}{P_1^k}, \quad
    -2 - \sum_{k=1}^\infty \frac{c_k^{(1)}+c_k^{(2)}}{P_1^k}. \]
Then, see (\ref{eq:2.25}), we integrate these three series, take the exponential function of the
results and again expand it into a power series around $P_1=\infty$
\[   \widetilde{\Phi}_1 = P_1 + \sum_{k=0}^\infty \frac{d_k^{(1)}}{P_1^k}, \quad
     \widetilde{\Phi}_2 = P_1 + \sum_{k=0}^\infty \frac{d_k^{(2)}}{P_1^k}, \quad
     \widetilde{\Phi}_0 =  \sum_{k=2}^\infty \frac{d_k^{(0)}}{P_1^k}. \]
Next, we find the constants $m_\ell$, $\ell=0,1,2$, such that the power series
\[   \sum_{\ell=0}^2 m_\ell \widetilde{\Phi}_\ell  \]
is a polynomial. In this way we obtain the series around $P_1=\infty$ for the branches of the function
$\Phi$
\begin{equation}  \label{eq:3.0}
    \Phi_\ell = m_\ell \widetilde{\Phi}_\ell.
\end{equation}
Finally, we use the Vi\`eta relations to obtain the coefficients (\ref{eq:2.26}) from the coefficients
of the series (\ref{eq:3.0}).
\item To obtain the algebraic parametrization of the curve $\Gamma$ in Proposition~\ref{prop:2.4} we note
that the coefficients of the polynomial $J(\nu)$
\begin{equation}  \label{eq:3.1*}
    J(\nu,z) := \left( \frac{\Phi_0(z)}{\Phi_1(z)} + \frac{\Phi_1(z)}{\Phi_0(z)} - \nu \right)
                \left( \frac{\Phi_0(z)}{\Phi_2(z)} + \frac{\Phi_2(z)}{\Phi_0(z)} - \nu \right)
                \left( \frac{\Phi_1(z)}{\Phi_2(z)} + \frac{\Phi_2(z)}{\Phi_1(z)} - \nu \right)
\end{equation}
are symmetric functions with respect to $\Phi_0, \Phi_1, \Phi_2$. Representing these symmetric
functions by means of the basic symmetric functions defined by the coefficients of the algebraic
equation for $\Phi$ we arrive at (\ref{eq:2.30})--(\ref{eq:2.31}).
\end{enumerate}
\end{proof}

\subsubsection{Proof of the geometric results for the case I} \label{sec:3.1.2}
Next we prove the results concerning the geometrical case I.

\begin{proof}[Proof of Theorem \ref{thm:2.1} and Proposition \ref{prop:2.4A}]
First we notice that (\ref{eq:2.370}) defines a charge of total mass 2. Indeed,
\[  \int_{\Delta_0} d\lambda = \frac{1}{2\pi i} \int_{\Delta_0} \bigl( h_{0+}(\xi) - h_{0-}(\xi) \bigr)\, d\xi
    = \frac{1}{2\pi i} \int_{\partial \mathfrak{R}_0} h(\xi)\, d\xi , \]
and by Cauchy's theorem and (\ref{eq:1.22}) we have
\[  \int_{\Delta_0} d\lambda = 2. \]
Analogously we have
\[  \int_{\Delta_j} d\lambda_j = 1, \qquad j=1,2. \]
Then we have to prove that the charges $\lambda_j$ are positive measures. From the definition (\ref{eq:2.25})
of $\Phi$ we have
\[   h = \bigl( \log |\Phi| + i \arg \Phi \bigr)', \quad
     h_0-h_j = \left( \log \left|\frac{\Phi_0}{\Phi_j}\right| + i \arg \frac{\Phi_0}{\Phi_j} \right)' . \]
Recalling (\ref{eq:2.35'}) and (\ref{eq:2.36'}) we have (see the notation in (\ref{eq:2.37}))
\[  \Delta_j = \Gamma_{0,j} \Leftrightarrow |\Phi_0| = |\Phi_j| \quad \textrm{on $\Delta_j=\gamma_j$}. \]
This relation clearly holds locally (in a neighborhood of the endpoints of $\Delta_j$) and since
Definition~\ref{def:2.2} and \ref{def:2.2'} of the case I and since the arcs $\gamma_j$ do not
intersect the cuts of $\mathfrak{R}$, it holds also globally on the whole of $\gamma_j$. Therefore
\begin{equation}  \label{eq:3.1}
  \lambda_j'(\xi) = \frac{1}{2\pi} \frac{\partial}{\partial \tau} \arg \frac{\Phi_0(\xi)}{\Phi_j(\xi)},
\qquad \xi \in \Delta_j (=\Gamma_{0,j}=\gamma_j).
\end{equation}
We now recall the parametrization of $\Gamma = \Gamma_{0,1} \cup \Gamma_{0,2} \cup \Gamma_{1,2}$
see (\ref{eq:2.29})
\[     \Gamma := \{ z : J(\nu,z) = 0, \nu \in [-2,2]\}. \]
Since the polynomial $J(\nu,z)$ has the form (\ref{eq:3.1*}), the parameter $\nu$ is related to the
argument of $\Phi_0/\Phi_j$ by
\[  \cos \left( \arg \frac{\Phi_0}{\Phi_j} \right) = \frac{\nu}2. \]
Hence, when $\xi$ goes along $\gamma_{a_j}$ from $a_j$ to $\beta_j=\gamma_{a_j}(-2)=\gamma_{b_j}(-2)$,
i.e., when $\nu$ changes from $2$ to $-2$, the argument of $\Phi_0/\Phi_j$ monotonically changes from
$0$ to $\pi$ and when $\xi$ continues to go from $\beta_j$ to $b_j$, the argument of $\Phi_0/\Phi_j$
monotonically changes from $\pi$ to $2\pi$ (a change of the argument here from $\pi$ to $0$ would
contradict the total variation of the argument, which follows from the argument principle applied to the function
$\Phi_0/\Phi_j$). Thus, the argument in (\ref{eq:3.1}) is a strictly monotone function and therefore,
with the proper orientation of $\gamma_j$,
\[  \lambda_j'(\xi) > 0, \qquad \xi \in \Delta_j, \quad j=1,2.  \]
The behavior of $\lambda_j'$ at the endpoints of $\Delta$ follows from the fact that the function
$h$ has poles (see (\ref{eq:1.23})) at the points $a_1,b_1,a_2,b2$. This proves Theorem~\ref{thm:2.1}.

The first part of Proposition~\ref{prop:2.4A} follows from the Cauchy-Riemann relations and (\ref{eq:3.1})
\[  - \frac{\partial}{\partial n} \log \left| \frac{\Phi_0}{\Phi_j} \right| =
      \frac{\partial}{\partial \tau} \arg \frac{\Phi_0}{\Phi_j} > 0. \]
The fact that for some $A \in \textup{I}$ the set $\Omega_{0,1}$ is nonemty can be verified by checking some particular
examples of $A$, see Figure \ref{fig:11}, second and third picture.
\end{proof}

\subsubsection{Proof of the geometric results for the case II} \label{sec:3.1.3}

The following statements (Proposition~\ref{prop:2.4B} and Theorem~\ref{thm:2.1B} are related to
the geometrical case II.

\begin{proof}[Proof of Proposition~\ref{prop:2.4B}]
For the geometrical case II (unlike for the case I) the analytic arcs $\gamma_1$ and
$\gamma_2$  (i.e., the trajectories forming the contour $\Gamma$ in (\ref{eq:2.29})) do not
coincide with the piecewise analytic arcs $\Delta_1$ and $\Delta_2$ (i.e., the cuts which make the
$f_1$ and $f_2$ holomorphic) and they do not coincide with the cuts $\delta_{0,1}$ and $\delta_{0,2}$
in (\ref{eq:1.28'}) of the sheets of $\mathfrak{R}$
\[  \gamma_j \neq \Delta_j, \quad \gamma_j \neq \delta_{0,j}, \qquad j=1,2. \]
Nevertheless, locally (in a neighborhood of the branch points of $\mathfrak{R}$) and therefore globally
(as long as the trajectories going from the branch points do not touch other cuts of $\mathfrak{R}$)
we have, with the notation (\ref{eq:2.372}) and (\ref{eq:2.37}),
\begin{equation}  \label{eq:3.2}
  \widetilde{\Delta}_j := \gamma_j \cap \Delta_j \subset \Gamma_{0,j} \Leftrightarrow
  |\Phi_0| = |\Phi_j| \quad \textrm{on $\widetilde{\Delta}_j$}, \quad j=1,2.
\end{equation}
Then the trajectory $\gamma_1$ passes through the cut $\delta_{0,2}$ and the branch $\Phi_0$
goes over to $\Phi_2$. Thus we have
\[  \gamma_1 \cap E_2 = E_2 \subset \Gamma_{1,2}.  \]
Analogously,
\[   \gamma_2 \cap E_1 = E_1 \subset \Gamma_{1,2}, \]
i.e.,
\begin{equation}  \label{eq:3.3}
    |\Phi_1| = |\Phi_2| \qquad \textrm{on $E_1 \cup E_2$}.
\end{equation}
It remains to determine which branches of $\Phi$ have the same modulus on $\Delta_{1,2}$.
The curve $\gamma_\alpha$ encloses one pair of the branch points and with the convention
(\ref{eq:2.38}) this is $\{a_2,b_2\}$. Indeed, it cannot enclose one or three branch points
because that would contradict the compactness of $\mathfrak{R}$, and it cannot enclose four branch points
because that would contradict the maximum principle. The fact that it encloses the pair
$\{a_j,b_j\}$ follows from the definition of the case I. Thus, we can join  infinity
and a point of $\widetilde{\Delta}_1$ by a path without crossing $\gamma_\alpha \setminus \Delta_{1,2}$,
and at the same time any path joining  infinity and a point of $\widetilde{\Delta}_2$ has to cross
$\gamma_\alpha \setminus \Delta_{1,2}$. Taking into account that in a neighborhood of infinity
the branch $\Phi_0$ has the smallest modulus and (\ref{eq:3.2}), we conclude that
\[   \gamma_\alpha \setminus \Delta_{1,2} \subset \Gamma_{1,2}, \]
and moreover we have that the neighborhood of infinity bounded by the parts of $\Gamma$ belongs to
$\Omega_{0,1,2}$. Also, the connected domain bounded by the parts of $\Gamma$ and containing
$\{a_2,b_2\}$ belongs to $\Omega_{0,2,1}$. Then the trajectory $\gamma_\alpha$ passes through the
cut $\delta_{0,2}$ and the branch $\Phi_2$ goes over to $\Phi_0$. Thus we have
\begin{equation}  \label{eq:3.4}
   \gamma_\alpha \subset \Delta_{1,2} = \Delta_{1,2} \subset \Gamma_{1,0}.
\end{equation}
To complete the proof of the proposition we notice that (\ref{eq:3.3}) and
(\ref{eq:3.4}) imply (see the notation (\ref{eq:2.373}))
\[   G \subset \Omega_{0,1,2}.  \]
\end{proof}

\begin{proof}[Proof of Theorem~\ref{thm:2.1B}] \
\begin{enumerate}
\item From Proposition~\ref{prop:2.4B} it follows that
\[  \overline{ \Omega_{0,1,2} \cup \Omega_{0,2,1}} = \overline{\mathbb{C}} . \]
Hence the modulus $|\Phi_0|$ is the minimal among all the moduli of all solutions
of the equation for $\Phi$. For a solution of the equation for $\Phi$ one has
\[  \log |\Phi| \in \textrm{Harm}(\mathbb{C} \setminus \{a_1,b_1,a_2,b_2\}) \]
hence the function $\log |\Phi_0|$ is a superharmonic function in $\mathbb{C} \setminus  \{a_1,b_1,a_2,b_2\}$,
and by the Riesz decomposition theorem, see e.g.\ \cite[Theorem II.3.1]{76},
\begin{equation}  \label{eq:3.5}
  \log |\Phi_0(z)| = V^{\lambda}(z) - \log |c_0|,
\end{equation}
see (\ref{eq:2.28}), where $\lambda$ is a positive measure. The total mass of this measure is
equal to $2$ due to (\ref{eq:2.28}) and this measure is supported on
$\Delta_0 = \Delta_1 \cup \widetilde{\Delta}_2$ because $\log |\Phi_0| \in \textrm{Harm}(\mathbb{C} \setminus
\Delta)$. The relation (\ref{eq:3.5}) can also be obtained by taking the real part of the
primitives of both sides of the Cauchy integral formula
\[  \frac{\Phi_0'(z)}{\Phi_0(z)} = \frac{1}{2\pi i} \int_{\Delta_0} \frac{d\xi}{\xi-z}
   \left[  \left( \frac{\Phi_0'(\xi)}{\Phi_0(\xi)} \right)_+
         - \left( \frac{\Phi_0'(\xi)}{\Phi_0(\xi)} \right)_- \right] . \]
Since $\Phi_0'/\Phi_0 = h$ and taking into account the pole type singularities of $h$ at the branch points,
we obtain part 1) of the theorem.
\item Since the jump of the function $h_0$ over $\Delta_1 \cup \widetilde{\Delta}_2$ produces the
positive measure $\lambda$, the jump of the function $h_1$ over $\Delta_1$ and the jump of the function
$h_2$ over $\widetilde{\Delta}_2$ produce measures with negative sign, i.e., $-\lambda_1$ and
$-\tilde{\lambda}_2$. Repeating the arguments of the proof of Theorem~\ref{thm:2.1} we see that the jump function $h_2$
over $\widetilde{\Delta}_2 \cup E_1$ produces a measure with negative sign
\[   -\mu_1 - \tilde{\lambda}_2, \quad \supp(\mu_1) = E_1, \quad \supp(\tilde{\lambda}_2) = \widetilde{\Delta}_2, \]
with total mass equal to $-1$
\[  -|\mu_1| - |\tilde{\lambda}_2| = -1. \]
The measure $\mu_1$ with positive sign can also be obtained by means of the jump of the function
$h_1$ over $E_1$. This proves part 2) of the theorem.
\item Since the total mass of the charge produced by the jump function $h_1$ over $\Delta_1 \cup E_1$
is equal to $-1$, see (\ref{eq:1.22})
\[    |\mu_1| - |\lambda_1| = -1, \]
we obtain the last part of the theorem
\end{enumerate}
\end{proof}

\subsubsection{Proof of the geometric results for the case III} \label{sec:3.1.4}

The following statements (Proposition \ref{prop:2.5}, \ref{prop:2.4II}, and Theorem \ref{thm:2.2}) are
related to the geometrical case III. Proposition~\ref{prop:2.5} is just a more explicit form
(relating to $\textup{III}$) of the general result from Proposition~\ref{prop:2.3}. The equation
(\ref{eq:2.43})--(\ref{eq:2.44}) in Proposition~\ref{prop:2.5} is well known (see \cite{20}; the
equivalent form of the algebraic equation for the Riemann surface of the case $\textup{III}$
goes back to Nuttall \cite{3,30}). The proofs of Theorem~\ref{thm:2.2} and Proposition~\ref{prop:2.4II}
are a repetition of the proofs of the corresponding results for the case I in Theorem~\ref{thm:2.1}
and Proposition~\ref{prop:2.4A}. A minor change in Theorem~\ref{thm:2.2} is the behavior of the weight
$\lambda_1'$ at the point $b^*$. This change is due to the regularity at $b^*$ of the function
$h$ given by (\ref{eq:2.40}), unlike in the case $A \in \textup{I}$. Another minor change in
Proposition~\ref{prop:2.4II} is that the domain $\Omega_{1,0}$ is always present for case $\textup{III}$;
for the case I it exists just for some $A \in \textup{I}$. To prove this fact we notice that all
trajectories $\gamma_{b^*}^{(j)}$, $j=1,2,3$ belong to $\Gamma_{0,1}$, see the notation in (\ref{eq:2.37}).
Indeed, the multiplicity of the root in (\ref{eq:2.44'}) and (\ref{eq:3.1*}) gives that all $\gamma_{b^*}^{(j)}$
belong to the same $\Gamma_{k,\ell}$, but $\gamma_{b^*}^{(1)} \subset \Gamma_{0,1}$, see Definition~\ref{def:2.4}-2.
Thus, from the outside the Jordan curve $J$
\[  J \subset \gamma_{b^*}^{(2)} \cup \gamma_{b^*}^{(3)} \subset \Gamma_{0,1}, \qquad
b^* \in J, \]
is a part of the boundary of $\Omega_{0,1,2}$ and inside it contains a tleast a part of $\gamma_2 \subset \Gamma_{0,2}$,
otherwise it would be in contradiction with the maximum principle, and therefore $\Omega_{1,0,2}$ is always
inside $J$ (see Figure~\ref{fig:13}).

\subsubsection{Proof of the geometric results for the cases IV and V} \label{sec:3.1.5}

As we already mentioned, Theorem~\ref{thm:2.4} gives a complete classification of the geometry for
the problem with three branch points

\begin{proof}[Proof of Theorem~\ref{thm:2.4} and Propositions \ref{prop:2.6}, \ref{prop:2.7}] \
\begin{enumerate}
\item According to the definition (\ref{eq:2.47'}) of the domains $D_j$, $j=1,2,3$, if a point $a$
belongs to the boundary between the domains $D_1$ and $D_2$, then the trajectory
$\{\gamma_{b^*}^{(i)}(\nu) : \nu \in (-2,2]\}$ has to pass through the point $b=0$ (see Figure~\ref{fig:13}).
To find this set of points $a$ we take the explicit form of the equation (\ref{eq:2.30}), substituting
in (\ref{eq:2.31}) the explicit expressions (\ref{eq:2.44}):
\[   J(\nu,z) = 0  \]
and set $z=0$. As a result we have
\[  J(\nu,z) = \frac{(\nu-2)(a^4+4(1-\nu)a^3+(22-8\nu)a^2+4(1-\nu)a+1)^2}{16a^2(a+1)^4} .  \]
If we set the numerator equal to zero, then we obtain the algebraic parametrization of $D_1$.
\item The boundary between $D_2$ and $D_3$ is formed by the critical triples (see Definition~\ref{def:2.4*}).
To find the set of points $a$ such that the triple $\{-1,0,a\}$ is critical, we can proceed in two ways.
The first way is just to substitute
\[   z = b^* = \frac{(a-1)^3}{9(a^2+a+1)} \]
from (\ref{eq:2.41}) in the explicit form of the equation (\ref{eq:2.30}). As a result
we get
\[   J(\nu,b^*) = \frac{(\nu-2) \widetilde{P}(a,\nu)^2}{11664a^4(a^2+a+1)^6(a+1)^4} = 0, \]
where $\widetilde{P}(a,\nu)$ is the same polynomial as in part 2) of Theorem~\ref{eq:2.4}. Thus we
obtained the characterization of the critical triples for which for some $\nu \in (-2,2]$ the
trajectories of $\Gamma$ pass through $b^*$.
The second way is to follow Definition~\ref{def:2.4*} and to find out when the algebraic curve $J(\nu,z)$
has branch points for $\nu \in (-2,2)$. For this purpose we compute the discriminant of $J(\nu,z)$ and as
a result we find
\[  \widetilde{\mathcal{D}}_J = C (\nu-2)^3(\nu+2)^3 P(a,\nu)^2 \widetilde{P}(a,\nu)^3, \]
where $C$ is a constant and $P(a,\nu)$ is a polynomial for which the zeros do not produce
branch points (because $P$ appears with a square). The zeros of the polynomial $\widetilde{P}(a,\nu)$
are the branch points and the local analysis $(z_j-z_i)^2 \simeq (\nu-\nu_0)^3$ shows that
the trajectories at these branch points behave like
\[   z = c_0 + c_2(\nu-\nu_0) \pm c_3 (\nu-\nu_0)^{3/2}, \]
i.e., at these branch points the trajectories meet each other at a zero angle. This proves
Theorem~\ref{thm:2.4}.
\item Now we complete the proof of Proposition~\ref{prop:2.7} (part A of this proposition has just been
proven). It remains to prove that for a strictly non-acceptable triple the analytic arcs
$\gamma_1$, $\gamma_2$, and $\gamma_3$ intersect at one point $c$. First we
notice that these arcs cannot join the point $b^*$ with the points $a_1,a_2,b$ without intersection,
because otherwise we get a contradiction with the irreducibility of $\Phi$. Indeed, $\Phi$ has
a quadratic branch point at $b^*$ and therefore along all three trajectories $\gamma_{b^*}^{(i)}$ ($i=1,2,3$)
we have the equality of the modulus of the same two branches and these trajectories therefore give
equality of the values of these branches to all other branch points $a_1,b_1,b$. The same reason
shows that $\gamma_1$, $\gamma_2$, and $\gamma_3$ cannot have only pairwise
intersections. Thus for this triple we must have a set of points $Z$ where the three arcs intersect
and where all three branches have the same modulus
\begin{equation} \label{eq:3.6}
   Z := \{ z : |\Phi_0(z)|=|\Phi_1(z)| =|\Phi_2(z)| \}.
\end{equation}
If we take the normalization (\ref{eq:2.27}) into account, then we have for these points the
prescribed values of $\Phi$
\[  Z := \{ z : \bigl\{\Phi_0(z),\Phi_1(z),\Phi_2(z)\} = \{1, e^{2\pi i/3}, e^{4\pi i/3}\} \bigr\}. \]
From the equation (\ref{eq:2.43}) for $\Phi$ it then follows that there exist at most two such points.
For an acceptable triple we may have all possibilities; see Figure~\ref{fig:13} where $Z=\emptyset$ in the
first picture, $Z$ has one point in the second picture, and $Z$ has two points in the third picture.
For a non-acceptable triple we have that $Z$ contains one point; as we already explained
$Z = \emptyset$ contradicts the irreducibility of $\Phi$ and two points in $Z$ leads to a contradiction
with the maximum principle. This proves Proposition~\ref{prop:2.7}.
\item Proposition~\ref{prop:2.6} follows immediately from Definition~\ref{def:2.5*} ($b \in \Omega_{0,1}$)
and the fact that $b  \in \Gamma_{0,2}$.
\end{enumerate}
\end{proof}
\medskip

The last set of statements (Proposition~\ref{prop:2.7} and Theorem~\ref{thm:2.5}) is related to the
geometrical case $\textup{IV}$. Their proofs go along the same lines as the proofs of the corresponding
statements for the geometrical case II. One minor change should be taken into account, namely
the regular behavior of the function $h$ at the branch point $b^*$. It causes a change of the local
behavior of the measures near this point.

\subsection{Proof of the vector equilibrium properties for the weak limits of the
\HP\ approximants} \label{sec:3.2}
The weak asymptotic formulas and the convergence theorems for the \HP\ approximants
(Theorems~\ref{thm:2.6}, \ref{thm:2.7}, and \ref{thm:2.8}) are direct consequences of the strong
asymptotic formulas (Theorems~\ref{thm:2.11} and \ref{thm:2.12}) which will be proven in the
next section. In this section we concentrate on the verification of the vector potential problem
related with the weak asymptotics.

\begin{proof}[Proof of Theorem~\ref{thm:2.9}]
In order to prove the theorem we have to verify the statements of the theorem for each
geometrical case of the position of the branch points $A := \{a_1,b_1;a_2,b_2\}$. We shall try
to do this in a unique way. The proof will be given in the following way.
\begin{enumerate}
\item First, we consider \index{Functions!g@$g$}
the real part of the Abelian integral (\ref{eq:1.19}), (\ref{eq:1.28'})
\begin{equation}  \label{eq:3.7}
   g(\xi) := \Re G(\xi), \qquad \xi \in \mathfrak{R},
\end{equation}
where $\mathfrak{R}$ is an arbitrary algebraic Riemann surface of order 3, and we introduce
universal global branches for $g:= \{g_0,g_1,g_2\}$. \index{Functions!g0g1g2@$g_0$, $g_1$, $g_2$}
These global branches of $g$ define a
sheet structure for $\mathfrak{R}$. This universal sheet structure for $\mathfrak{R}$ goes back to
Nuttall \cite{3} and is different from what we use for the definition of our geometrical cases.
\item Second, we define two universal measures $\Lambda$ and $M$ \index{Measures!Lambda@$\Lambda$} \index{Measures!M@$M$}
supported on the new cuts of
$\mathfrak{R}$. These measures have total mass $|\Lambda|=2$ and $|M|=1$ and they satisfy the
vector equilibrium property (\ref{eq:1.13}) with the matrix of interaction
\begin{equation}  \label{eq:3.8}
    \Bigl( d_{i,j} \Bigr)_{i,j=1,2} = \begin{pmatrix} 2 & -1 \\ -1 & 2 \end{pmatrix}.
\end{equation}
The vector potential problem (\ref{eq:1.13}) with (\ref{eq:3.8}) was introduced by Nikishin
\cite{24,4}.
\item Finally, we perform a \textit{reglueing} of the Riemann surface and make a correspondence
between the universal sheet structure of $\mathfrak{R}$ and the specific sheet structure
$\{\mathfrak{R}_0,\mathfrak{R}_1,\mathfrak{R}_2\}$ which we have defined for each geometrical case.
In practice this procedure gives us the relation (by means of some \textit{balayages}) between the
measures $(\Lambda,M)$ and $(\lambda,\mu)$ and transforms the universal vector potential problem
with the Nikishin matrix of interaction (\ref{eq:3.8}) to the potential problem (\ref{eq:2.70})
in the theorem.
\end{enumerate}
We emphasize that the potential problem of Theorem~\ref{eq:2.9} also has a universal character, i.e.,
it does not depend on the geometrical case under consideration. Although the problem (\ref{eq:2.70})
is more sophisticated in comparison with the potential problem (\ref{eq:1.13}) with (\ref{eq:3.8}),
this problem has an advantage because it is formulated in terms of the functions $f_1,f_2$
which we are approximating (i.e., in terms of the branch points and the cuts which make the functions
holomorphic) and not in terms of $\mathfrak{R}$. We also mention that the procedure to construct
the proper $\mathfrak{R}$ starting from the functions $f_1,f_2$ is the most difficult step in
proving the Nuttall conjectures (see \cite[Section 3.4]{3}) and this was in general not solved before.

We start the proof with an arbitrary algebraic $\mathfrak{R}$ of the third order and consider the
real part of the Abelian integral (\ref{eq:1.28'}), which is a single-valued (up to an additive constant)
function on $\mathfrak{R}$, with the only singularities at the points at infinity
\begin{equation}  \label{eq:3.9}
  \begin{cases}
    g_0(\xi) = - 2\log |\xi| + c_0 + \cdots, & \xi \to \infty^{(0)}, \\
    g_j(\xi) = \log |\xi| + c_j + \cdots, & \xi \to \infty^{(j)}, \ j=1,2,
   \end{cases}
\end{equation}
and  \index{Parameters!c0c1c2@$c_0$, $c_1$, $c_2$}
\[    g(\xi) \neq \infty, \qquad \xi \in \mathfrak{R} \setminus \pi^{-1}(\infty).  \]
This local definition of the branches of $g(\xi)$ can be continued on $\mathbb{C}$ as follows
\begin{equation}   \label{eq:3.10}
    \forall z \in \mathbb{C}: g_0(z) \leq g_1(z) \leq g_2(z).
\end{equation}
 From (\ref{eq:3.9}) and the maximum principle for harmonic functions we have
\begin{equation}  \label{eq:3.11}
   g_0(z) + g_1(z) + g_2(z)= c_0+c_1+c_2 = 0, \qquad \forall z \in \mathbb{C},
\end{equation}
and by (\ref{eq:3.11}) we have then fixed the additive constant for $g$. Thus we have
two piecewise analytic curves \index{Contours!GammaLambdaGammaM@$\Gamma_{\Lambda}$, $\Gamma_M$}
\begin{equation}  \label{eq:3.12}
    \begin{array}{l}
     \Gamma_\Lambda := \{ z \in \mathbb{C}: g_0(z)=g_1(z) \} \\[12pt]
      \Gamma_M := \{ z \in \overline{\mathbb{C}} : g_1(z) = g_2(z) \}
    \end{array}
\end{equation}
which define the Nuttall sheet structure of
$\mathfrak{R}$ \index{Surfaces!R0tildeR1tildeR2tilde@$\widetilde{\mathfrak{R}}_0$, $\widetilde{\mathfrak{R}}_1$, $\widetilde{\mathfrak{R}}_2$}
\begin{equation}  \label{eq:3.13}
   \widetilde{\mathfrak{R}}_0 := \overline{\mathbb{C}} \setminus \Gamma_\Lambda, \quad
   \widetilde{\mathfrak{R}}_1 := \overline{\mathbb{C}} \setminus (\Gamma_\Lambda \cup \Lambda_M), \quad
   \widetilde{\mathfrak{R}}_2 := \overline{\mathbb{C}} \setminus \Gamma_M .
\end{equation}
We shall denote the Riemann surface with the sheet structure (\ref{eq:3.13})
by $\widetilde{\mathfrak{R}}$:  \index{Surfaces!Rtilde@$\widetilde{\mathfrak{R}}$}
\[  \widetilde{\mathfrak{R}} := \overline{\bigcup_{j=0}^2 \mathfrak{R}_j} . \]
Then, like in the proof of Theorem~\ref{thm:2.1B}, we have by (\ref{eq:3.10}) that $g_0$ is a
superharmonic function in $\mathbb{C}$ and a harmonic function in $\mathbb{C} \setminus \Gamma_\Lambda$,
and $g_2$ is a subharmonic function in $\mathbb{C}$ and a harmonic function in $\overline{C} \setminus \Gamma_M$.
Therefore, by the Riesz decomposition theorem there exist positive measures $\Lambda$ and $M$ such that
\begin{equation}  \label{eq:3.14}
   \begin{cases}
      g_0(z) = V^\Lambda(z) + c_0, & \supp(\Lambda) = \Gamma_\Lambda, \\
      g_2(z) = - V^M(z) + c_2, & \supp(M) = \Gamma_M .
   \end{cases}
\end{equation}
 From (\ref{eq:3.9})  we have for the total mass of $\Lambda$ and $M$
\[   |\Lambda| = 2, \quad  |M| = 1, \]
and (\ref{eq:3.11})--(\ref{eq:3.12}) gives the Nikishin vector equilibrium relations
\begin{equation}  \label{eq:3.15}
   \begin{cases}
    2 V^\Lambda(z) - V^M(z) = -2c_0-c_2, & \textrm{on $\Gamma_\Lambda$}, \\
    -V^\Lambda(z) + 2 V^M(z) = 2c_2+c_0, & \textrm{on $\Gamma_M$}.
   \end{cases}
\end{equation}
We also note that, in the notation we used earlier, we have for $|z|$ large enough,
\begin{equation}   \label{eq:3.16}
    g_\ell(z) = \log |\Phi_\ell(z)|,
\end{equation}
however, the domains where the branches $\Phi_\ell$ are single-valued (which we defined earlier)
are different from the domains where the branches $g_\ell$ are single-valued. We define the second set of global
branches for the algebraic functions $h = \{\tilde{h}_0,\tilde{h}_1,\tilde{h}_2\}$ and
$\Phi = \{\widetilde{\Phi}_0,\widetilde{\Phi}_1,\widetilde{\Phi}_2\}$ in accordance with the Nuttall
sheet structure (\ref{eq:3.13}): \index{Functions!h0tildeh1tildeh2tilde@$\tilde{h}_0$, $\tilde{h}_1$, $\tilde{h}_2$}
\index{Functions!Phi0tildePhi1tildePhi2tilde@$\widetilde{\Phi}_0$, $\widetilde{\Phi}_1$, $\widetilde{\Phi}_2$}
\[  \tilde{h}_0, \widetilde{\Phi}_0 \in H(\overline{\mathbb{C}} \setminus \Gamma_\Lambda), \quad
    \tilde{h}_1, \widetilde{\Phi}_1 \in H(\mathbb{C} \setminus (\Gamma_\Lambda \cup \Gamma_M)), \quad
    \tilde{h}_2, \widetilde{\Phi}_2 \in H(\mathbb{C} \setminus \Gamma_M).  \]
Then we have
\begin{equation}  \label{eq:3.17}
  \begin{cases}
   d\Lambda(z) = \frac{1}{2\pi i} \Bigl( \tilde{h}_{0+}(z) - \tilde{h}_{0-}(z) \Bigr)\, dz, & z \in \Gamma_\Lambda, \\
   -dM(z)  = \frac{1}{2\pi i} \Bigl( \tilde{h}_{2+}(z) - \tilde{h}_{2-}(z) \Bigr)\, dz, & z \in \Gamma_M.
   \end{cases}:
\end{equation}
Now we perform the following \textit{reglueing} procedure to transform the Nuttall Riemann surface
$\widetilde{\mathfrak{R}}$ (\ref{eq:3.13}) into the Riemann surfaces with sheet structure defined
in Subsections \ref{sec:2.2}--\ref{sec:2.4}, for which we use the notation $\mathfrak{R}$. We note that
\[   \Gamma_\Lambda \cup \Gamma_M = \Gamma, \]
where $\Gamma$ is defined in (\ref{eq:2.29})--(\ref{eq:2.31}) and does not depend on the sheet
structure of $\mathfrak{R}$. The union of the cuts of the sheets $\{\widetilde{\mathfrak{R}}_0,
\widetilde{\mathfrak{R}}_1,\widetilde{\mathfrak{R}}_2\}$ coincides with $\Gamma$:
\[      \bigcup_{\ell,k=0}^2 \tilde{\delta}_{\ell,k} = \Gamma, \]
and at the same time the union of the cuts of $\{\mathfrak{R}_0,\mathfrak{R}_1,\mathfrak{R}_2\}$
occupies just a subset of $\Gamma$:
\[   \bigcup_{\ell,k=0}^2 \delta_{\ell,k} \subset \Gamma.  \]
To transform $\mathfrak{R}$ into $\widetilde{\mathfrak{R}}$ we consider the regions $\Omega_{j,k,\ell}$,
$j\neq k\neq \ell$ pairwise, $j,k,\ell=0,1,2$, see (\ref{eq:2.37'}). The contour $\Gamma$ divides
$\mathbb{C}$ into these regions. Then we take the components of $\Omega_{j,k,\ell}$ with indices
$(j,k,\ell) \neq (0,1,2)$, lift them to $\pi^{-1}(\Omega_{j,k,\ell})$ on $\mathfrak{R}$, cut $\mathfrak{R}$
along the boundaries $\partial (\pi^{-1}(\Omega_{j,k,\ell}))$, and interchange these pieces of
$\mathfrak{R}$ by reglueing them such that
\begin{eqnarray*}
     \pi_j^{-1}(\Omega_{j,k,\ell}) & \leftrightarrow & \pi_0^{-1}(\Omega_{j,k,\ell}) \\
      \pi_k^{-1}(\Omega_{j,k,\ell}) & \leftrightarrow & \pi_1^{-1}(\Omega_{j,k,\ell}) \\
     \pi_\ell^{-1}(\Omega_{j,k,\ell}) & \leftrightarrow & \pi_2^{-1}(\Omega_{j,k,\ell}).
\end{eqnarray*}
As a result we obtain $\widetilde{\mathfrak{R}}$ with Nuttall's structure of the sheets. During this
procedure some extra cuts have appeared. The transformation of $\widetilde{\mathfrak{R}}$ to
$\mathfrak{R}$ works in the reverse order and some cuts may disappear.
Since $\infty \in \Omega_{0,1,2}$ for all $\mathfrak{R}$ under consideration, this procedure
provides balayages of some parts of the measures $\lambda$ and $\mu$ to parts of the measures
$\Lambda$ and $M$. This allows to transform the potential problem (\ref{eq:3.15}) to the
potential problem (\ref{eq:2.70}).

\begin{figure}[ht]
\centering
\framebox{
\unitlength 2pt
\linethickness{0.5pt}
\begin{picture}(110,60)
\put(10,30){\circle*{2}}
\put(60,30){\circle*{2}}
\put(70,30){\circle*{2}}
\put(95,30){\circle*{2}}
\put(80,30){\circle*{2}}
\put(7,25){$a_1$}
\put(56,25){$b_1$}
\put(68,33){$a_2$}
\put(95,33){$b_2$}
\put(79,23){$c$}
\put(10,30){\line(1,0){50}}
\put(70,30){\line(1,0){25}}
\qbezier(65,30)(65,40)(72,40)
\qbezier(72,40)(77,40)(80,30)
\qbezier(65,30)(65,20)(72,20)
\qbezier(72,20)(77,20)(80,30)
\qbezier(80,30)(87,45)(95,45)
\qbezier(95,45)(105,45)(105,30)
\qbezier(80,30)(87,15)(95,15)
\qbezier(95,15)(105,15)(105,30)
\put(35,0){$\Omega_{0,1,2}$}
\put(40,40){$\Gamma_{0,1}$}
\put(42,38){\vector(0,-1){7}}
\put(50,40){\vector(2,-1){15}}
\put(70,10){$\Omega_{1,0,2}$}
\put(88,50){$\Omega_{0,2,1}$}
\put(103,13){\small $\Gamma_{1,2}$}
\put(72,15){\vector(0,1){10}}
\put(90,48){\vector(0,-1){10}}
\end{picture}}
\caption{The partition of $\mathbb{C}$ by $\Gamma$ (a generic case for $\textup{I}$ and $\textup{III}$)}
\label{fig:18}
\end{figure}
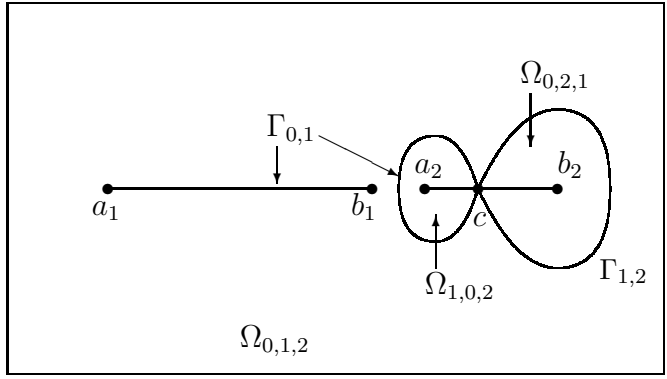

We now apply the transformation $\widetilde{\mathfrak{R}} \to \mathfrak{R}$ to prove the theorem.
We start with the cases I and $\textup{III}$. Of course, for these cases it is easy to prove the theorem
directly (see Remark~\ref{rem:2.3}, parts 1 and 3). However, for methodological reasons we use here the
more sophisticated general procedure. We consider a generic case of the classes I and $\textup{III}$
(see Figure~\ref{fig:11}, second picture, and Figure~\ref{fig:13}, second picture). For this case
all the regions $\Omega_{0,1,2}$, $\Omega_{0,2,1}$, and $\Omega_{1,0,2}$ are present. It follows
from Propositions~\ref{prop:2.4A} and \ref{prop:2.4II} that the contours $\Gamma_\Lambda$ and $\Gamma_M$
are (see Figure~\ref{fig:18})
\[   \Gamma_\Lambda := \Gamma_{0,1} \cup \gamma_{c,b_2} , \quad
     \Gamma_M := \Gamma_{1,2} \cup \gamma_{c,a_2}.  \]
Thus, the sheets $\{\widetilde{\mathfrak{R}}_0,\widetilde{\mathfrak{R}}_1,\widetilde{\mathfrak{R}}_2\}$
are as shown in Figure~\ref{fig:19}. If we re-glue the sheets $\widetilde{\mathfrak{R}}_0$ and
$\widetilde{\mathfrak{R}}_1$ along $\pi^{-1}(\Gamma_{0,1} \setminus \gamma_{a_1,b_1})$, interchange
the domains $\pi_0^{-1}(\Omega_{1,0,2})$ and $\pi_1^{-1}(\Omega_{1,0,2})$, \index{Functions!pi0pi1pi2@$\pi_0$, $\pi_1$, $\pi_2$}
and then glue
$\widetilde{\mathfrak{R}}_0$ and $\widetilde{\mathfrak{R}}_1$ along
$\pi^{-1}(\Gamma_{0,1} \setminus \gamma_{a_1,b_1})$,
and $\widetilde{\mathfrak{R}}_1$ and $\widetilde{\mathfrak{R}}_2$ along
$\pi^{-1}(\gamma_{a_2,c})$, \index{Contours!gammaa2c@$\gamma_{a_2,c}$}
we see that $\pi^{-1}(\gamma_{a_2,c})$ disappears from the first sheet and moves to the zero sheet.
If we do the same with the sheets $\widetilde{\mathfrak{R}}_1$ and $\widetilde{\mathfrak{R}}_2$ and
the domains $\pi_1^{-1}(\Omega_{0,2,1})$ and $\pi_2^{-1}(\Omega_{0,2,1})$, then we obtain the
Riemann surface $\mathfrak{R}$ from Figure~\ref{fig:10}.

\begin{figure}[ht]
\centering
\framebox{
\unitlength 2pt
\linethickness{0.5pt}
\begin{picture}(110,40)(0,10)
\put(0,45){$\widetilde{\mathfrak{R}}_0$}
\put(10,30){\circle*{2}}
\put(60,30){\circle*{2}}
\put(95,30){\circle*{2}}
\put(80,30){\circle*{2}}
\put(7,25){$a_1$}
\put(56,25){$b_1$}
\put(95,33){$b_2$}
\put(79,23){$c$}
\put(10,30){\line(1,0){50}}
\put(80,30){\line(1,0){15}}
\put(0,0){\bezier{10}(65,30)(65,40)(72,40)}
\put(0,0){\bezier{10}(72,40)(77,40)(80,30)}
\put(0,0){\bezier{10}(65,30)(65,20)(72,20)}
\put(0,0){\bezier{10}(72,20)(77,20)(80,30)}
\qbezier(70,25)(60,15)(60,5)
\put(69,23){\vector(2,3){2}}
\end{picture}}
\framebox{
\unitlength 2pt
\linethickness{0.5pt}
\begin{picture}(110,40)(0,10)
\put(0,45){$\widetilde{\mathfrak{R}}_1$}
\put(10,30){\circle*{2}}
\put(60,30){\circle*{2}}
\put(70,30){\circle*{2}}
\put(95,30){\circle*{2}}
\put(80,30){\circle*{2}}
\put(7,25){$a_1$}
\put(56,25){$b_1$}
\put(68,25){$a_2$}
\put(95,33){$b_2$}
\put(79,23){$c$}
\put(10,30){\line(1,0){50}}
\put(70,30){\line(1,0){25}}
\put(0,0){\bezier{10}(65,30)(65,40)(72,40)}
\put(0,0){\bezier{10}(72,40)(77,40)(80,30)}
\put(0,0){\bezier{10}(65,30)(65,20)(72,20)}
\put(0,0){\bezier{10}(72,20)(77,20)(80,30)}
\put(0,0){\bezier{15}(80,30)(87,45)(95,45)}
\put(0,0){\bezier{15}(95,45)(105,45)(105,30)}
\put(0,0){\bezier{15}(80,30)(87,15)(95,15)}
\put(0,0){\bezier{15}(95,15)(105,15)(105,30)}
\qbezier(70,35)(60,45)(60,55)
\put(69,37){\vector(2,-3){2}}
\qbezier(100,25)(110,15)(110,5)
\put(102,22){\vector(-2,3){2}}
\end{picture}}
\framebox{
\unitlength 2pt
\linethickness{0.5pt}
\begin{picture}(110,40)(0,10)
\put(0,45){$\widetilde{\mathfrak{R}}_2$}
\put(70,30){\circle*{2}}
\put(80,30){\circle*{2}}
\put(68,25){$a_2$}
\put(79,23){$c$}
\put(70,30){\line(1,0){10}}
\put(0,0){\bezier{15}(80,30)(87,45)(95,45)}
\put(0,0){\bezier{15}(95,45)(105,45)(105,30)}
\put(0,0){\bezier{15}(80,30)(87,15)(95,15)}
\put(0,0){\bezier{15}(95,15)(105,15)(105,30)}
\qbezier(100,35)(110,45)(110,55)
\put(102,38){\vector(-2,-3){2}}
\end{picture}}
\caption{The Nuttall sheet structure of $\widetilde{\mathfrak{R}}$ (the case $\textup{I}$)}
\label{fig:19}
\end{figure}

Decomposing the measure $\Lambda$ we have, due to (\ref{eq:3.17}) and the structure of
$\widetilde{\mathfrak{R}}$, \index{Measures!lambdaca2lambdabc2@$\lambda_{c,a_2}$, $\lambda_{c,b_2}$}
\begin{equation}  \label{eq:3.18}
  \Lambda = \lambda_1 + \lambda_{c,b_2} + \mu_{0,1},
\end{equation}
where \index{Measures!mu01mu12@$\mu_{0,1}$, $\mu_{1,2}$} \index{Contours!gammaca2gammacb2@$\gamma_{c,a_2}$, $\gamma_{c,b_2}$}
\[   \lambda_{c,b_2} := \lambda\Bigr|_{\gamma_{c,b_2}}, \quad
     \mu_{0,1} := \Lambda\Bigr|_{\Gamma_{0,1}\setminus \Delta_1}, \]
and $\lambda_1,\lambda_2$ are from Theorem~\ref{thm:2.1}. Analogously, we have
\begin{equation}  \label{eq:3.19}
   M = \lambda_{c,a_2} + \mu_{1,2},
\end{equation}
where
\[    \lambda_{c,a_2} := \lambda_2\Bigr|_{\gamma_{c,a_2}}, \quad
      \mu_{1,2} := M\Bigr|_{\Gamma_{1,2}} .  \]
The equality
\[    g_0 = \log |\Phi_0|, \]
see (\ref{eq:3.16}), is valid in the neighborhood of infinity bounded by $\Gamma$.
We have
\[ \log |\Phi_0| = V^{\lambda} - \log |C_0| \]
due to the definition of $\lambda$ in (\ref{eq:2.370}) and the relations between
$h$ and $\Phi$, see (\ref{eq:2.25}), (\ref{eq:2.27}), and (\ref{eq:2.28}).
Therefore in the neighbhorhood of infinity bounded by $\Gamma$, we have
\[   V^{\lambda_1} + V^{\lambda_{c,b_2}} + V^{\mu_{0,1}} + c_0 = V^{\lambda_1} + V^{\lambda_2} - \log |C_0|, \]
and we conclude that
\[   c_0 = - \log |C_0|,  \]
and the measure $\mu_{0,1}$ is the balayage of the measure $\lambda_{c,a_2}$ from
$\gamma_{c,a_2}$ to $\Gamma_{0,1}\setminus \Delta_1$, i.e.,
\begin{equation}   \label{eq:3.20}
     V^{\mu_{0,1}} = V^{\lambda_{c,a_2}}, \qquad \textrm{on $\Gamma_{0,1} \setminus \Delta_{1}$ and outside}.
\end{equation}
Analogously we have in the neighborhood of infinity bounded by $\Gamma_M$
\[   g_2 = \log |\Phi_2| \]
and
\[ \log |\Phi_2| = - V^{\lambda_2} - \log |C_2|. \]
Therefore in this domain we have
\[    V^{\lambda_{c,a_2}} + V^{\mu_{1,2}} = V^{\lambda_2}, \qquad c_2 = -\log |C_2|, \]
and we obtain another balayage relation
\begin{equation}  \label{eq:3.21}
   V^{\mu_{1,2}} = V^{\lambda_{c,b_2}}, \qquad \textrm{on $\Gamma_{1,2}$ and outside}.
\end{equation}
In addition to the balayage relations (\ref{eq:3.20})--(\ref{eq:3.21}), which are valid outside
$\Gamma_{0,1} \setminus \Delta_1$ and $\Gamma_{1,2}$, we need relations between the potentials
inside these curves. Let us consider the first relation of the Nikishin equilibrium problem (\ref{eq:3.15})
on the curve $\Gamma_{0,1} \setminus \Delta_1$. If we substitute the decompositions
(\ref{eq:3.18})--(\ref{eq:3.19}) into (\ref{eq:3.15}) and use the balayages (\ref{eq:3.20})--(\ref{eq:3.21})
we obtain
\begin{equation}  \label{eq:3.22}
   2 V^{\lambda_1} + V^{\mu_{0,1}} + V^{\lambda_{c,b_2}} = -2 c_0 - c_2, \qquad \textrm{on $\Gamma_{0,1} \setminus
 \Delta_1$}.
\end{equation}
The left hand side of (\ref{eq:3.22}) is a harmonic function inside the domain bounded by $\Gamma_{0,1} \setminus
\Delta_1$. Hence by the maximum principle for harmonic functions, the relation (\ref{eq:3.22}) is valid
in this domain. In the same way we obtain
from the second relation in (\ref{eq:3.15}) that
\begin{equation}  \label{eq:3.23}
     V^{\mu_{1,2}} + V^{\lambda_{a_2,c}} - V^{\lambda_1} = 2 c_2 + c_0, \qquad \textrm{on $\Gamma_{1,2}$ and
inside}.
\end{equation}
Now it is easy, for the cases under consideration, to derive the equilibrium relations of
Theorem~\ref{thm:2.9} from the equilibrium relations (\ref{eq:3.15}) using (\ref{eq:3.20})--(\ref{eq:3.23}).
The first relation of (\ref{eq:3.15}) considered on $\Delta_1$ gives us, using (\ref{eq:3.18})--(\ref{eq:3.19})
and (\ref{eq:3.20})--(\ref{eq:3.21}),
\[    2V^{\lambda_1} + V^{\lambda_2} = -2c_0 - c_2 =: \kappa_1,
   \qquad \textrm{on $\Delta_1$} . \]
The first relation of (\ref{eq:3.15}) on $\gamma_{c,b_2}$ (see Figure~\ref{fig:18}) gives us, using
also (\ref{eq:3.23}), \index{Parameters!kappa1kappa2@$\kappa_1$, $\kappa_2$}
\[    2V^{\lambda_2} + V^{\lambda_1} = -2c_0 - c_2- (2c_2+c_0) = -3(c_0+c_2) =: \kappa_2,
   \qquad \textrm{on $\gamma_{c,b_2}$} . \]
We get the same relation on $\gamma_{a_2,c}$ from the second relation of (\ref{eq:3.15}), from
(\ref{eq:3.21}) and with the help of (\ref{eq:3.22}), which is also valid on $\gamma_{a_2,c}$ as we have seen.
Thus, from the general Nikishin potential problem (\ref{eq:3.15})
we obtain the Angelesco potential problem (\ref{eq:1.13}), (\ref{eq:1.16})
\[   \begin{cases}
     2 V^{\lambda_1} + V^{\lambda_2} = \kappa_1, & \textrm{on $\Delta_1$}, \\
     V^{\lambda_1} + 2V^{\lambda_2} = \kappa_2, & \textrm{on $\Delta_2$}.
     \end{cases}  \]
Therefore for the cases I, $\textup{III}$, and $\textup{V}$ (see Remark~\ref{rem:2.3}, 1--4) the
equilibrium relations are verified. We emphasize the importance of the balayages (\ref{eq:3.20})--(\ref{eq:3.21})
for the description of the weak limits of the extra interpolation points in the case $\textup{V}$ (see
Remark~\ref{rem:2.3}-4).

Now we apply the same approach to verify the equilibrium relations (\ref{eq:2.7}) for the remaining
geometrical cases II and $\textup{IV}$. For the case II Proposition~\ref{prop:2.4B} gives
(see Figure~\ref{fig:11B} and (\ref{eq:2.372}), (\ref{eq:2.375}))
\begin{eqnarray*}
   \Gamma_\Lambda &=& \Delta_1 \cup \widetilde{\Delta}_2, \\
   \Gamma_M & = & E_1 \cup E_2 \cup E_0, \qquad E_0:= \gamma_\alpha \setminus \Delta_{1,2},
\end{eqnarray*}
and the corresponding measures (\ref{eq:3.14}) and (\ref{eq:3.17}) are given by \index{Measures!Lambda@$\Lambda$} \index{Measures!M@$M$}
\begin{eqnarray*}
  \Lambda & := & \lambda_1 + \tilde{\lambda}_2, \\
   M & := & \mu_1 + \mu_2 + \mu_0, \qquad \mu_2 := M \Bigr|_{E_2}, \ \mu_0:= M \Bigr|_{E_0},
\end{eqnarray*}
where the measures \index{Measures!mu0mu1mu2@$\mu_0$, $\mu_1$, $\mu_2$}
\[   \lambda_1 = \Lambda \Bigr|_{\Delta_1}, \quad
     \tilde{\lambda}_2 = \Lambda \Bigr|_{\widetilde{\Delta}_2}, \quad
     \mu_1 = M \Bigr|_{E_1}    \]
are defined in Theorem~\ref{thm:2.1B}. The reglueing of the Riemann surface $\widetilde{\mathfrak{R}}$
gives the balayage relation
\begin{equation}  \label{eq:3.24}
  V^{\tilde{\lambda}_2} = V^{\mu_2} + V^{\mu_0}, \qquad \textrm{on $E_2 \cup E_0$ and outside.}
\end{equation}
The second equilibrium relation of (\ref{eq:3.15}), the balayage (\ref{eq:3.24}) and the maximum
principle give
\begin{equation}  \label{eq:3.25}
  2V^{\mu_1} + V^{\mu_2} + V^{\mu_0} - V^{\lambda_1} = 2 c_2 + c_0, \qquad \textrm{on $E_2 \cup E_0$ and inside.}
\end{equation}
The first equilibrium relation in (\ref{eq:3.15}) on $\Delta_1$ and (\ref{eq:3.24}) give the
first equilibrium relation in (\ref{eq:2.70})
\[  2V^{\mu_1} + 2V^{\tilde{\lambda}_2} - V^{\mu_1} - V^{\mu_2} - V^{\mu_0}
   = 2 V^{\lambda_1} + V^{\tilde{\lambda}_2} - V^{\mu_1} = -2c_0 - c_2 := \kappa_1, \qquad \textrm{on $\Delta_1$}. \]
If we take the first relation in (\ref{eq:3.15}) on $\widetilde{\Delta}_2$ and use (\ref{eq:3.25}), then
we find the second relation in (\ref{eq:2.70})
\[  2V^{\mu_1} + 2V^{\tilde{\lambda}_2} - V^{\mu_1} - V^{\mu_2} - V^{\mu_0}
   =  V^{\lambda_1} + 2V^{\tilde{\lambda}_2} + V^{\mu_1} = \kappa_1 +2c_2 + c_0 := \tilde{\kappa}_2,
\qquad \textrm{on $\widetilde{\Delta}_2$}. \]
Finally, the second relation in (\ref{eq:3.15}) on $E = E_1 \cup E_2$ and the balayage (\ref{eq:3.24})
lead to the third relation of (\ref{eq:2.70})
\[  2V^{\mu_1} + 2V^{\mu_2} + 2V^{\mu_0} - V^{\lambda_1} - V^{\tilde{\lambda}_2}
   = 2 V^{\mu_1} + V^{\tilde{\lambda}_2} - V^{\lambda_1} = 2c_2 + c_0 = \tilde{\kappa}_2 - \kappa_1,
\qquad \textrm{on $E$}. \]
Thus the equilibrium relations (\ref{eq:2.70}) for the case II are verified.

If we verify the equilibrium relations for the last case $\textup{IV}$, then we will also get the dual
equilibrium problem (see II and III of Theorem~\ref{thm:2.9}). For the case $\textup{IV}$ Proposition~\ref{prop:2.8}
(see (\ref{eq:2.53}) and (\ref{eq:2.53'})) gives
\begin{eqnarray*}
   \Gamma_\Lambda & = & \Delta_1 \cup \widetilde{\Delta}_2 = \widetilde{\Delta}_1 \cup \Delta_{1,2}
     \cup \widetilde{\Delta}_2, \\
   \Gamma_M & = & E_0 \cup E_1 \cup E_2.
\end{eqnarray*}
The Nuttall sheet structure of $\widetilde{\mathfrak{R}}$ is shown in Figure~\ref{fig:20}. It defines
the global branches (\ref{eq:3.14}) for the real part of the Abelian integral (\ref{eq:3.9})
\[ \begin{cases}
     g_0 = V^{\tilde{\lambda}_1 + \tilde{\lambda}_2 + \lambda_{1,2}} + c_0, \\
     g_1 = - V^{\tilde{\lambda}_1 + \tilde{\lambda}_2 + \lambda_{1,2}} + V^{\mu_1 + \mu_2 + \mu_0} - c_0 - c_2, \\
     g_2 = -V^{\mu_1 + \mu_2 + \mu_0},
    \end{cases}  \]
where the measures
\[  \tilde{\lambda}_1 = \Lambda \Bigr|_{\widetilde{\Delta}_1}, \quad
    \tilde{\lambda}_2 = \Lambda \Bigr|_{\widetilde{\Delta}_2}, \quad
    \lambda_{1,2} = \Lambda \Bigr|_{\Delta_{1,2}}, \quad
    \mu_2 = M \Bigr|_{E_2}   \]
are the same as in Theorem~\ref{thm:2.5}, and
\[   \mu_1 := M \Bigr|_{E_1}, \quad \mu_0 := M\Bigr|_{E_0}. \]
Equating
\[ \begin{cases}
    g_0 = g_1 , & \textrm{on $\Gamma_\Lambda$}, \\
    g_1 = g_2,  & \textrm{on $\Gamma_M$},
   \end{cases}  \]
gives the Nikishin equilibrium (\ref{eq:3.5}). Observe that (\ref{eq:3.5}) has just two equilibrium
relations since the 0-sheet and the second sheet do not intersect.
\newpage 

\begin{figure}[ht]
\centering
\framebox{
\unitlength 1.8pt
\linethickness{0.4pt}
\begin{picture}(105,60)
\put(0,55){$\widetilde{\mathfrak{R}}_0$}
\put(10,10){\circle*{2}}
\put(75,45){\circle*{2}}
\put(70,25){\circle*{2}}
\put(75,13){\circle*{2}}
\qbezier(10,10)(35,25)(70,25)
\qbezier(75,13)(72.5,19.5)(70,25)
\qbezier(75,45)(72.5,35)(70,25)
\put(10,6){$a_1$}
\put(76,46){$a_2$}
\put(76,10){$b$}
\put(64,28){$c$}
\put(35,25){$\widetilde{\Delta}_1$}
\put(75,35){$\widetilde{\Delta}_2$}
\put(62,17){$\Delta_{1,2}$}
\end{picture}} \medskip

\framebox{
\unitlength 1.8pt
\linethickness{0.4pt}
\begin{picture}(105,60)
\put(0,55){$\widetilde{\mathfrak{R}}_1$}
\put(10,10){\circle*{2}}
\put(75,45){\circle*{2}}
\put(70,25){\circle*{2}}
\put(75,13){\circle*{2}}
\put(80,25){\circle*{2}}
\qbezier(10,10)(35,25)(70,25)
\qbezier(70,25)(75,25)(80,25)
\qbezier(75,13)(72.5,19.5)(70,25)
\qbezier(70,25)(60,57)(75,57)
\qbezier(75,57)(92,57)(80,25)
\qbezier(75,5)(90,5)(80,25)
\qbezier(70,25)(65,5)(75,5)
\qbezier(75,45)(72.5,35)(70,25)
\put(10,6){$a_1$}
\put(83,25){$b^*$}
\put(76,46){$a_2$}
\put(76,10){$b$}
\put(64,28){$c$}
\end{picture}} \medskip

\framebox{
\unitlength 1.8pt
\linethickness{0.4pt}
\begin{picture}(105,60)
\put(0,55){$\widetilde{\mathfrak{R}}_2$}
\put(70,25){\circle*{2}}
\put(80,25){\circle*{2}}
\qbezier(70,25)(75,25)(80,25)
\qbezier(70,25)(60,57)(75,57)
\qbezier(75,57)(92,57)(80,25)
\qbezier(75,5)(90,5)(80,25)
\qbezier(70,25)(65,5)(75,5)
\put(83,25){$b^*$}
\put(64,28){$c$}
\put(85,50){$E_0$}
\put(85,10){$E_1$}
\put(72,27){$E_2$}
\end{picture}}
\caption{The Nuttall sheet structure of $\widetilde{\mathfrak{R}}$ (the case $\textup{IV}$)}
\label{fig:20}
\end{figure}

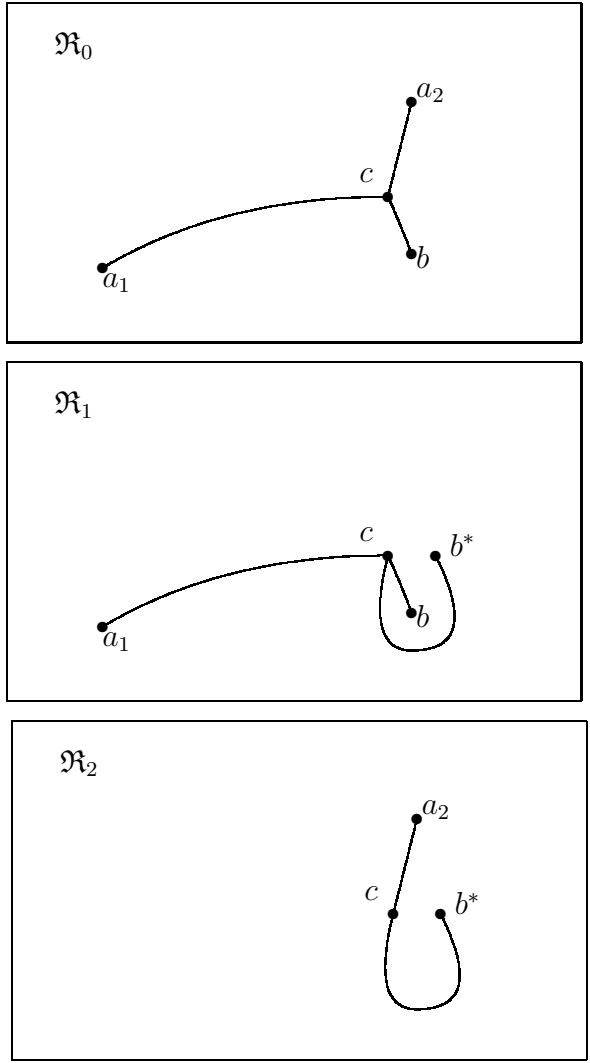
\begin{figure}[hbt]
\centering
\framebox{
\unitlength 1.8pt
\linethickness{0.4pt}
\begin{picture}(105,60)
\put(0,55){$\mathfrak{R}_0$}
\put(10,10){\circle*{2}}
\put(75,45){\circle*{2}}
\put(70,25){\circle*{2}}
\put(75,13){\circle*{2}}
\qbezier(10,10)(35,25)(70,25)
\qbezier(75,13)(72.5,19.5)(70,25)
\qbezier(75,45)(72.5,35)(70,25)
\put(10,6){$a_1$}
\put(76,46){$a_2$}
\put(76,10){$b$}
\put(64,28){$c$}
\end{picture}} \medskip

\framebox{
\unitlength 1.8pt
\linethickness{0.4pt}
\begin{picture}(105,60)
\put(0,55){$\mathfrak{R}_1$}
\put(10,10){\circle*{2}}
\put(70,25){\circle*{2}}
\put(75,13){\circle*{2}}
\put(80,25){\circle*{2}}
\qbezier(10,10)(35,25)(70,25)
\qbezier(75,13)(72.5,19.5)(70,25)
\qbezier(75,5)(90,5)(80,25)
\qbezier(70,25)(65,5)(75,5)
\put(10,6){$a_1$}
\put(83,25){$b^*$}
\put(76,10){$b$}
\put(64,28){$c$}
\end{picture}} \medskip

\framebox{
\unitlength 1.8pt
\linethickness{0.4pt}
\begin{picture}(105,60)
\put(0,55){$\mathfrak{R}_2$}
\put(70,25){\circle*{2}}
\put(80,25){\circle*{2}}
\put(75,45){\circle*{2}}
\qbezier(75,45)(72.5,35)(70,25)
\qbezier(75,5)(90,5)(80,25)
\qbezier(70,25)(65,5)(75,5)
\put(83,25){$b^*$}
\put(64,28){$c$}
\put(76,46){$a_2$}
\end{picture}}
\caption{Reglueing the Nuttall Riemann surface (the case $\textup{IV}$) along $E_0 \cup E_2$}
\label{fig:21}
\end{figure}

In order to obtain (\ref{eq:2.70}) we re-glue the Riemann surface shown in Figure~\ref{fig:20}.
We interchange the domains bounded by $E_0 \cup E_2$ on the sheets 1 and 2 and as a result
the cut along $E_0 \cup E_2$ between the sheets 1 and 2 disappears and we arrive at the
Riemann surface $\mathfrak{R}$ shown in Figure~\ref{fig:21}.

The new sheet structure defines new global branches for the real part of the \index{Functions!g@$g$}
Abelian integral (\ref{eq:3.9}) \index{Functions!g0tildeg1tildeg2tilde@$\tilde{g}_0$, $\tilde{g}_1$, $\tilde{g}_2$}
\[  g =: \begin{cases}
       \tilde{g}_0 = V^{\tilde{\lambda}_1+\tilde{\lambda}_2+\lambda_{1,2}} + c_0, \\
       \tilde{g}_1 = - V^{\tilde{\lambda}_1+\lambda_{1,2}} + V^{\mu_1} - c_0-c_2, \\
       \tilde{g}_2 = -V^{\tilde{\lambda}_2}-V^{\mu_1}+c_2.
         \end{cases} \]
Equating $\tilde{g}_0 = \tilde{g}_1$ on $\Delta_1$, where $\Delta_1$ is the cut joining the sheet
0 and 1, we have
\[  V^{\tilde{\lambda}_1} + V^{\tilde{\lambda}_2} + V^{\lambda_{1,2}} + c_0
   = -V^{\tilde{\lambda}_1} - V^{\lambda_{1,2}} + V^{\mu_1} - c_0 - c_2, \qquad \textrm{on $\Delta_1$}, \]
which gives the first relation in (\ref{eq:2.70}):
\[  2V^{\lambda_1} + V^{\tilde{\lambda}_2} - V^{\mu_1} = -2c_0 - c_2 := \kappa_1, \qquad \textrm{on $\Delta_1$}. \]
Equating $\tilde{g}_0=\tilde{g}_2$ on $\widetilde{\Delta}_2$ gives the second relation in (\ref{eq:2.70})
\[  V^{\lambda_1} + 2V^{\tilde{\lambda}_2} + V^{\mu_1} = c_2-c_0 =: \kappa_2, \qquad \textrm{on $\widetilde{\Delta}_2$}. \]
Finally, equating $\tilde{g}_1=\tilde{g}_2$ on $E_1$ gives the third relation in (\ref{eq:2.70}):
\[  -V^{\lambda_1} + V^{\tilde{\lambda}_2} + 2V^{\mu_1} = 2c_2+c_0 = \kappa_2-\kappa_1, \qquad \textrm{on $E_1$}. \]
Thus the equilibrium relations (\ref{eq:2.70}) are verified for the case $\textup{IV}$. \index{Parameters!kappa1kappa2@$\kappa_1$, $\kappa_2$}

To derive the dual equilibrium problem we re-glue the Riemann surface shown in Figure~\ref{fig:20}
by interchanging the domains bounded by $E_0 \cup E_1$ on the sheets 1 and 2, and as a result the
cut along $E_0 \cup E_1$ between the sheets 1 and 2 disappears and we arrive at the Riemann surface
shown in Figure~\ref{fig:17}. The new sheet structure defines new global branches for the
real part of the Abelian
integral (\ref{eq:3.9}): \index{Functions!g0hatg1hatg2hat@$\hat{g}_0$, $\hat{g}_1$, $\hat{g}_2$}
\[  g =: \begin{cases}
       \hat{g}_0 = V^{\tilde{\lambda}_1} + V^{\lambda_2} + c_0, \\
       \hat{g}_1 = - V^{\tilde{\lambda}_1} - V^{\mu_2} - c_0-c_2, \\
       \hat{g}_2 = -V^{\lambda_2}+V^{\mu_2}+c_2.
         \end{cases} \]
Equating $\hat{g}_0=\hat{g}_1$ on $\widetilde{\Delta}_1$ (see Figure~\ref{fig:17}) gives the first
relation for the dual equilibrium problem
\[  2V^{\tilde{\lambda}_1} + V^{\lambda_2} + V^{\mu_2} = \kappa_1, \qquad \textrm{on $\widetilde{\Delta}_1$} . \]
Equating $\hat{g}_0=\hat{g}_2$ on $\Delta_2$ gives the second relation
\[   V^{\tilde{\lambda}_1} + 2V^{\lambda_2} - V^{\mu_2} = \kappa_2, \qquad \textrm{on $\Delta_2$}, \]
and $\hat{g}_1=\hat{g}_2$ on $E_2$ gives
\[  V^{\tilde{\lambda}_1} - V^{\lambda_2} + 2V^{\mu_2} = \kappa_1 - \kappa_2, \qquad \textrm{on $E_2$}.  \]
Thus the equilibrium problem for the case $\textup{IV}$ (see II and III of Theorem~\ref{thm:2.9}) is verified.
The dual equilibrium problem can be verified in a similar way for the other geometrical cases.

To conclude the proof of Theorem~\ref{thm:2.9} we verify the symmetry property (\ref{eq:2.71}). We
show how to do this for the geometrical case I; one can use a similar reasoning for the other cases.
The equilibrium relations for the case I have the form
\begin{eqnarray*}
  U_1 & := & 2V^{\lambda_1} + V^{\lambda_2} = \kappa_1, \qquad \textrm{on $\Delta_1$}, \\
  U_2 & := & V^{\lambda_1} + 2 V^{\lambda_2} = \kappa_2, \qquad \textrm{on $\Delta_2$}.
\end{eqnarray*}
For the harmonic conjugate $\widetilde{U}_j$ one has
\[   \frac{\partial \widetilde{U}_j}{\partial n} = \frac{\partial U_j}{\partial \tau} = 0, \qquad
     \textrm{on $\Delta_j$}, \quad j=1,2, \]
hence one has to verify for $\mathcal{U}_j = U_j + i \widetilde{U}_j$ that
\[   \frac{\partial}{\partial n_+} \mathcal{U}_j(z) = \frac{\partial}{\partial n_-} \mathcal{U}_j(z),
    \qquad z \in \Delta_j \setminus \{a_j,b_j\}, \quad j=1,2. \]
This relation holds because the function $\mathcal{U}_j$ is analytic on $\Delta_j \setminus \{a_j,b_j\}$
(because $\mathcal{U}_j = \log \Phi_0/\Phi_j + \textrm{const}$), and therefore its derivative
does not depend on the direction. The theorem is now proved.
\end{proof}

\section{Asymptotic analysis of the matrix Riemann-Hilbert problem for the Hermite-Pad\'e polynomials} \label{sec:4}
In this section we prove the results about the strong asymptotics of the \HP\ polynomials, i.e.,
Theorem~\ref{thm:2.11} and Theorem~\ref{thm:2.12}. We demonstrate here the detailed proof for the case
$\textup{II}$. The proof for the case $\textup{I}$ is just a simplified version of the proof for the case II
and the most delightful pieces of the proof degenerate. The same is true for the cases $\textup{III}$ and $\textup{V}$
which are both a simplified version of the proof for the case $\textup{IV}$, and the proof for the case $\textup{IV}$
is just a repetition of the proof for the case II with minor differences in the local analysis
of the asymptotics in the neighborhood of the branch point $b^*$. For these reasons we present the proof
for the generic case II in detail and give a sketch of the proofs for the other cases.

For the proof we use the steepest descent method of Deift and Zhou (see \cite{35}). The method consists
of several consecutive transformations
\[ Y \mapsto Z \mapsto \widetilde{Z} \mapsto \widehat{Z} \mapsto \check{Z} \]
of the matrix Riemann-Hilbert problem (\ref{eq:2.13}) for the
matrix function $Y$ to the Riemann-Hilbert problem \index{Matrices!Zcheck@$\check{Z}$}
\begin{equation}  \label{eq:4.1}
   \begin{cases}
   \check{Z} \in H^{3\times 3}(\overline{\mathbb{C}} \setminus \check{\Sigma}), \\[8pt]
   \exists \check{Z}_{\pm} \in C(\check{\Sigma}): \check{Z}_+ = \check{Z}_- \bigl(I + \O(1/n)\bigr), &
   \textrm{uniformly on $\check{\Sigma}$ as $n\to \infty$}, \\[8pt]
   \check{Z}(z) = I + \O(1/z), & z \to \infty.
   \end{cases}
\end{equation}
It is known that for large $n$, the problem (\ref{eq:4.1}) has a solution which satisfies
\[   \check{Z} = I + \O(1/n), \qquad \textrm{uniformly in $\overline{\mathbb{C}} \setminus \check{\Sigma}$ as $n \to \infty$}. \]
Thus, after making the inverse transformations from $\check{Z}$ to $Y$ we obtain the existence
of the solution $Y$ of the Riemann-Hilbert problem (\ref{eq:2.13}) for large $n$ and asymptotics
for $n \to \infty$ for its components.

\subsection{Normalization of the \RH\ problem at infinity and decomposition of the jumps} \label{sec:4.1}
The goal of the first step is to transform the \RH\ problem (\ref{eq:2.13})--(\ref{eq:2.17}) such
that the solution of the new problem has the same normalization at infinity as the solution of
(\ref{eq:4.1}). For this purpose we use the algebraic function $\Phi = \{\Phi_0,\Phi_1,\Phi_2\}$ in
(\ref{eq:2.25})--(\ref{eq:2.28}). We denote \index{Matrices!S@$S$} \index{Matrices!C@$C$}
\begin{equation}  \label{eq:4.2}
   S := \begin{pmatrix}
        \Phi_0^n & 0 & 0 \\ 0 &  \Phi_1^n & 0 \\ 0 & 0 & \Phi_2^n
       \end{pmatrix}, \quad C := \begin{pmatrix} C_0^n & 0 & 0 \\ 0 & C_1^n & 0 \\ 0 & 0 & C_2^n \end{pmatrix},
\end{equation}
where $C_0,C_1,C_2$ are the inverses of the leading coefficients of the expansion of $\Phi_0, \Phi_1,\Phi_2$
near infinity, see (\ref{eq:2.28}). If we set \index{Matrices!Z@$Z$}
\begin{equation}  \label{eq:4.3}
    Z := C Y S,
\end{equation}
then, because of (\ref{eq:2.13}) and (\ref{eq:2.28})
\begin{equation}  \label{eq:4.4}
   Z(z) = I + \O\left( \frac{1}{z} \right), \qquad z \to \infty,
\end{equation}
and the piecewise analytic matrix function $Z$ satisfies the following jump condition: \index{Matrices!J@$J$}
\begin{equation}   \label{eq:4.5}
  Z_+ = Z_- J , \qquad \textrm{on } \widetilde{\Delta}_1 \cup \Delta_{1,2} \cup \widetilde{\Delta}_{2}
    \cup E_1 =: \Sigma,
\end{equation}
(see Figures \ref{fig:11B} and \ref{fig:21a}), and for the jump matrix $J$ we have from (\ref{eq:2.14})
that $J:= S_-^{-1} W S_+$ so that \index{Contours!Sigma@$\Sigma$}
\begin{equation}  \label{eq:4.6}
   J := \begin{cases}
         J(w_1,w_2), & \textrm{on $\Delta_{1,2}$}, \\
         J(w_1,0), & \textrm{on $\widetilde{\Delta}_1$}, \\
         J(0,w_2), & \textrm{on $\widetilde{\Delta}_2$}, \\
         J(0,0), & \textrm{on $E_1$},
        \end{cases}
\end{equation}
where we set \index{Matrices!Jw1w2@$J(w_1,w_2)$}
\begin{equation}  \label{eq:4.7}
   J(w_1,w_2) := \begin{pmatrix}
      \displaystyle \frac{\Phi_{0+}^n}{\Phi_{0-}^n} & \displaystyle \frac{\Phi_{1+}^n}{\Phi_{0-}^n}\ w_1
         & \displaystyle \frac{\Phi_{2+}^n}{\Phi_{0-}^n}\ w_2 \\[12pt]
      0 & \displaystyle \frac{\Phi_{1+}^n}{\Phi_{1-}^n} & 0 \\
      0 & 0 & \displaystyle \frac{\Phi_{2+}^n}{\Phi_{2-}^n}
     \end{pmatrix}.
\end{equation}
Note that, in comparison with $Y$, the RH problem for $Z$ has a new contour $E_1$ where $Z$ is discontinuous,
due to the jumps in the functions $\Phi_j$.

\begin{figure}[ht]
\centering
\framebox{
\unitlength 5pt
\linethickness{0.4pt}
\begin{picture}(60,50)(0,10)
\put(0,55){$\textup{II}$}
\put(18,52){\circle*{1}}
\put(14,18){\circle*{1}}
\put(40,48){\circle*{1}}
\put(32,18){\circle*{1}}
\put(28,42){\circle*{1}}
\put(24,26){\circle*{1}}
\put(30,43){\circle*{1}}
\put(26,24){\circle*{1}}
\put(17,54){$a_1$}
\put(12,16){$b_1$}
\put(41,49){$a_2$}
\put(33,17){$b_2$}
\put(24,42){$c_1$}
\put(20,26){$c_2$}
\put(32,42){$c_1^+$}
\put(28,23){$c_2^+$}
\qbezier(28,42)(24,34)(24,26)
\qbezier(28,42)(17,35)(24,26)
\put(0,0){\bezier{25}(28,42)(33,35)(24,26)}
\put(0,0){\bezier{30}(30,43)(36,35)(26,24)}
\qbezier(24,26)(27.5,21.5)(32,18)
\qbezier(24,26)(19.5,21.5)(14,18)
\qbezier(28,42)(33.5,45.5)(40,48)
\qbezier(28,42)(25.5,45.5)(18,52)
\qbezier(19.14,22.58)(18.89,21.83)(18.38,21.08)
\qbezier(18.38,21.08)(18.98,21.35)(19.84,21.43)
\qbezier(28.02,20.9)(28.59,20.75)(29.17,20.33)
\qbezier(29.17,20.33)(28.93,20.82)(28.59,21.48)
\qbezier(24.63,35.06)(24.72,34.63)(24.94,33.69)
\qbezier(24.94,33.69)(25.28,34.31)(25.75,34.69)
\qbezier(21.13,35.69)(21.34,35.28)(21.44,34.25)
\qbezier(21.44,34.25)(21.78,34.72)(22.25,35.31)
\qbezier(28.88,34.75)(29.13,34.16)(29.13,33.19)
\qbezier(29.13,33.19)(29.63,34.03)(30.13,34.5)
\qbezier(31.25,34.06)(31.5,33.16)(31.5,32.38)
\qbezier(31.5,32.38)(32.09,33.31)(32.56,33.75)
\qbezier(35.6,47.05)(35.3,46.49)(34.71,45.78)
\qbezier(34.71,45.78)(35.75,45.82)(36.2,45.71)
\qbezier(22,49.35)(22.15,48.5)(22.59,47.79)
\qbezier(22.59,47.79)(21.85,48.39)(21.11,48.39)
\put(21,51){$\widetilde{\Delta}_1$}
\put(35,48){$\widetilde{\Delta}_2$}
\put(18,36){\small $E_1$}
\put(28,38){\small $E_2$}
\put(33,36){\small $E_2^+$}
\put(23,30){$\Delta_{1,2}$}
\put(17,17){$\widetilde{\Delta}_1$}
\put(25,18){$\widetilde{\Delta}_2$}
\end{picture}}
\caption{Jump contours for the \RH\ problem for $Z$ and $\widetilde{Z}$ (normalization and the global lens)}
\label{fig:21a}
\end{figure}

Now we analyze the jump matrices (\ref{eq:4.6})--(\ref{eq:4.7}). We start with the usual decomposition
of the jumps on $\widetilde{\Delta}_1$ and $\widetilde{\Delta}_2$, which are commonly used in the
analysis of $2\times 2$ matrix-valued \RH\ problems. If we define \index{Matrices!D1D2@$D_1$, $D_2$} \index{Matrices!W1W2@$W_1$, $W_2$}
\begin{eqnarray}  \label{eq:4.8}
  D_1:= \begin{pmatrix}  1 & 0 & 0 \\ \displaystyle \frac{\Phi_0^n}{\Phi_1^n w_1} & 1 & 0 \\ 0 & 0 & 1 \end{pmatrix},
  \quad
  W_1:= \begin{pmatrix} 0 & w_1 & 0 \\ \displaystyle -\frac{1}{w_1} & 0 & 0 \\ 0 & 0 & 1 \end{pmatrix},  \nonumber \\
    && \\
  D_2:= \begin{pmatrix}  1 & 0 & 0 \\ 0 & 1 & 0 \\ \displaystyle \frac{\Phi_0^n}{\Phi_2^n w_2} & 0 & 1 \end{pmatrix},
  \quad
  W_2:= \begin{pmatrix} 0 & 0 & w_2 \\ 0 & 1 & 0 \\ \displaystyle -\frac{1}{w_2} & 0 & 0 \end{pmatrix},  \nonumber
\end{eqnarray}
then we have, since $\Phi_{0\pm} = \Phi_{1\mp}$ on $\widetilde{\Delta}_1$,
\begin{equation}  \label{eq:4.91}
  J = J(w_1,0) = \begin{pmatrix}
        \displaystyle \left(\frac{\Phi_{0}^n}{\Phi_{1}^n}\right)_+ &  w_1 & 0 \\
       0 & \displaystyle  \left(\frac{\Phi_{0}^n}{\Phi_{1}^n}\right)_- & 0 \\
     0 & 0 & 1 \end{pmatrix} = D_{1-} W_1 D_{1+} , \qquad \textrm{on $\widetilde{\Delta}_1$},
\end{equation}
and similarly, since $\Phi_{0\pm} = \Phi_{2\mp}$ on $\widetilde{\Delta}_2$,
\begin{equation}  \label{eq:4.92}
     J = J(0,w_2) = D_{2+} W_2 D_{2-}, \qquad \textrm{on $\widetilde{\Delta}_2$} .
\end{equation}
The following decompositions of the jumps on $\Delta_{1,2}$ and $E_1$ were not used before in the
analysis of \RH\ problems. If  we define
\[  D_{2,1} := \begin{pmatrix}
               1 & 0 & 0 \\ 0 & 1 & \displaystyle - \frac{\Phi_2^n w_2}{\Phi_1^n w_1} \\ 0 & 0 & 1
               \end{pmatrix} \]
then we observe that
\begin{equation}  \label{eq:4.10}
    J = J(w_1,w_2) = D_{2,1-} J(w_1,0) D_{2,1+}^{-1}, \qquad \textrm{on $\Delta_{1,2}$},
\end{equation}
i.e., after the decomposition (\ref{eq:4.10}) the jump has a block structure and $J(w_1,0)$ can then be
decomposed again in the usual way as in (\ref{eq:4.91}), since $\Phi_{0\pm} = \Phi_{1\mp}$
holds also on $\Delta_{1,2}$. If we define \index{Matrices!D12D21@$D_{1,2}$, $D_{2,1}$} \index{Matrices!W12@$W_{1,2}$}
\begin{equation}  \label{eq:4.10*}
 D_{1,2}:= \begin{pmatrix}  1 & 0 & 0 \\ 0 & 1 & 0 \\ 0 & \displaystyle -\frac{\Phi_1^n w_1}{\Phi_2^n w_2} & 1
\end{pmatrix},
  \quad
  W_{1,2}:= \begin{pmatrix}
             1 & 0 & 0 \\
             0 & 0 & \displaystyle - \frac{w_2}{w_1} \\
             0 & \displaystyle \frac{w_1}{w_2} & 0 \end{pmatrix},
 \end{equation}
then we have on $E_1$ (since $\Phi_{1\pm} = \Phi_{2\mp}$)
\[   J = J(0,0) = \begin{pmatrix}
              1 & 0 & 0 \\
              0 & \displaystyle \frac{\Phi_{1+}^n}{\Phi_{1-}^n} & 0 \\
              0 & 0 & \displaystyle \frac{\Phi_{2+}^n}{\Phi_{2-}^n}
              \end{pmatrix}
            = D_{1,2-}
             \begin{pmatrix}
             1 & 0 & 0 \\
             0 & 1 & \displaystyle \frac{\Phi_2^n w_2}{2\Phi_1^n w_1} \\
             0 & 0 & 1
             \end{pmatrix}_-
             W_{1,2}
             \begin{pmatrix}
             1 & 0 & 0 \\
             0 & 1 & 0 \\
             0 & \displaystyle \frac{\Phi_1^nw_1}{2\Phi_2^nw_2} & 1
            \end{pmatrix}_+^{-1}
            D_{2,1+}^{-1}.  \]
Since
\[   \begin{pmatrix} 1 & \displaystyle \frac{\Phi_2^n w_2}{2\Phi_1^n w_1} \\ 0 & 1 \end{pmatrix}_-
     \begin{pmatrix} 0 & \displaystyle - \frac{w_2}{w_1} \\ \displaystyle \frac{w_1}{w_2} & 0 \end{pmatrix} =
      \begin{pmatrix} 0 & \displaystyle - \frac{w_2}{w_1} \\ \displaystyle \frac{w_1}{w_2} & 0 \end{pmatrix}
      \begin{pmatrix} 1 & 0 \\ \displaystyle \frac{\Phi_1^nw_1}{2\Phi_2^nw_2} & 1 \end{pmatrix}_+^{-1}, \qquad
     \textrm{on $E_1$}, \]
we have
\[    J = J(0,0) = D_{1,2-} W_{1,2} \begin{pmatrix}
                                1 & 0 & 0 \\
                                0 & 1 & 0 \\
                                0 & \displaystyle \frac{\Phi_1^nw_1}{2\Phi_2^nw_2} & 1
                                \end{pmatrix}_+^{-2}
                D_{2,1+}^{-1}, \qquad \textrm{on $E_1$}, \]
and we obtain
\begin{equation}  \label{eq:4.11}
    J = J(0,0) = D_{1,2-} W_{1,2} D_{1,2+} D_{2,1+}^{-1}, \qquad \textrm{on $E_1$}.
\end{equation}
We also point out the commutation relations
\begin{equation}  \label{eq:4.12}
    D_1D_2 = D_2D_1, \quad D_1D_{1,2} = D_{1,2} D_1.
\end{equation}

\subsection{Opening a global lens (in preparation of opening local lenses)} \label{sex:4.2}
The goal of the second step is to transform the \RH\ problem (\ref{eq:4.4})--(\ref{eq:4.7}) for
the function $Z$ in (\ref{eq:4.3}) to a \RH\ problem with new jumps which can be decomposed
in the usual way by (\ref{eq:4.91})--(\ref{eq:4.92}) to jumps which do not depend on $n$ or to
jumps which tend to the identity matrix as $n \to \infty$.

We fix a neighborhood $O_{c_j}$ \index{Domains!Oc1Oc2@$O_{c_1}$, $O_{c_2}$}
of the point $c_j$ $(j=1,2)$ (see Figure~\ref{fig:21a}) and denote
\[   \widetilde{\Delta}_2^{(c_j)} := \widetilde{\Delta}_2 \cap O_{c_j}. \]
We join the endpoints \index{Contours!Delta2tildec1Delta2tildec2@$\widetilde{\Delta}_2^{(c_1)}$, $\widetilde{\Delta}_2^{(c_2)}$}
(different from $c_j$) of $\widetilde{\Delta}_2^{(c_j)}$, which we denote by $c_j^{+}$, by an arc
$E_2^+$ such that \index{Parameters!c1pc2p@$c_1^{+}$, $c_2^{+}$} \index{Contours!E2plus@$E_2^+$}
\begin{equation}  \label{eq:4.13}
    E_2^+ \subset \Omega_{2,1},
\end{equation}
where we used the notation (\ref{eq:2.37}).
We denote by $G^+$ the domain bounded by \index{Domains!Gplus@$G^+$}
\[  \partial G^+ := E_1 \cup \widetilde{\Delta}_2^{(c_1)} \cup \widetilde{\Delta}_2^{(c_2)} \cup E_2^+. \]
(see Figures~\ref{fig:12} and \ref{fig:21a}).
We define \index{Matrices!Ztilde@$\widetilde{Z}$}
\begin{equation}  \label{eq:4.14}
   \widetilde{Z} := \begin{cases}
                 Z D_{2,1}, & \textrm{in $G^+$}, \\
                 Z, & \textrm{in $\overline{\mathbb{C}} \setminus \overline{G^+}$}.
                    \end{cases}
\end{equation}
This piecewise analytic function $\widetilde{Z}$ has a jump \index{Matrices!Jtilde@$\widetilde{J}$} \index{Contours!Sigmatildeplus@$\widetilde{\Sigma}^+$}
\begin{equation}  \label{eq:4.15}
      \widetilde{Z}_+ = \widetilde{Z}_- \widetilde{J}, \qquad \textrm{on $\Sigma \cup E_2^+ =: \widetilde{\Sigma}^+$}
\end{equation}
where for $\widetilde{J}$ on $\widetilde{\Sigma}^+ \setminus ( \Delta_{1,2} \cup E_1 \cup E_2^+ \cup
\widetilde{\Delta}_2^{(c_1)} \cup \Delta_2^{(c_2)})$, see (\ref{eq:4.6}),
\begin{equation}  \label{eq:4.151}
  \widetilde{J} = J = \begin{cases}
   J(w_1,0) \\
   J(0,w_2)
   \end{cases}
  \qquad \textrm{on $\widetilde{\Delta}_1 \cup ( \widetilde{\Delta}_2 \setminus
  \bigcup_{j=1}^2 \widetilde{\Delta}_2^{(c_j)} )$}.
\end{equation}
On other parts of $\widetilde{\Sigma}$ the jump is changed to
\begin{equation}  \label{eq:4.152}
   \widetilde{J} = \begin{cases}
              J(0,w_2) D_{2,1}, & \textrm{on $\widetilde{\Delta}_2^{(c_1)} \cup \widetilde{\Delta}_2^{(c_2)}$}, \\
              D_{2,1}^{-1}, & \textrm{on $E_2^+$}.
                   \end{cases}
\end{equation}
On $E_1$ we obtain, using (\ref{eq:4.11}),
\[   \widetilde{Z}_+ = Z_+ D_{2,1+} = Z_- (D_{1,2-}W_{1,2} D_{1,2+} D_{2,1+}^{-1}) D_{2,1+} =
     \widetilde{Z}_- D_{1,2-} W_{1,2} D_{1,2+}, \]
and therefore we have
\begin{equation}  \label{eq:4.153}
    \widetilde{J} = D_{1,2-} W_{1,2} D_{1,2+}, \qquad \textrm{on $E_1$}.
\end{equation}
On $\Delta_{1,2}$ we have, using (\ref{eq:4.10}),
\[  \widetilde{Z}_+ = Z_+ D_{2,1+} = Z_- (D_{2,1-} J(w_1,0) D_{2,1+}^{-1}) D_{2,1+} =
    \widetilde{Z}_- J(w_1,0), \]
and therefore we have
\begin{equation}  \label{eq:4.154}
   \widetilde{J} = J(w_1,0), \qquad \textrm{on $\Delta_{1,2}$}.
\end{equation}
Thus the desired form of the jumps (\ref{eq:4.151})--(\ref{eq:4.154}) is indeed obtained
after the second transformation.

\subsection{Opening local lenses} \label{sec:4.3}
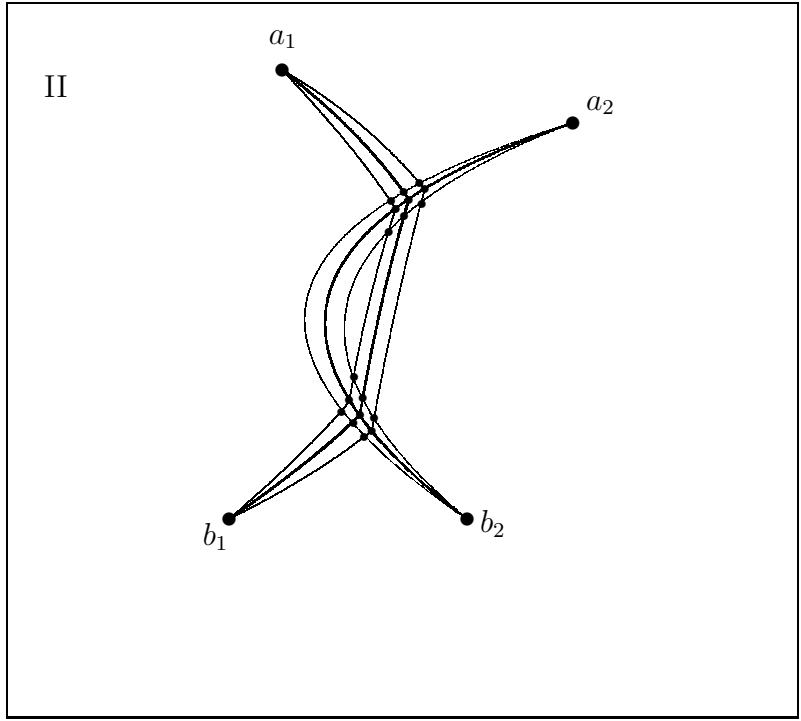
\begin{figure}[ht]
\centering
\unitlength 5.00pt 
\linethickness{0.6pt}
\framebox{
\begin{picture}(55,50)(0,5)
\put(0,50){$\textup{II}$}
\put(18,52){\circle*{1}}
\put(14,18){\circle*{1}}
\put(40,48){\circle*{1}}
\put(32,18){\circle*{1}}
\put(17,54){$a_1$}
\put(12,16){$b_1$}
\put(41,49){$a_2$}
\put(33,17){$b_2$}
\qbezier(40,48)(7.06,36.83)(32,18)
\qbezier(23.9,25.9)(20.48,22.49)(14,18)
\qbezier(18,52)(23.97,47.42)(27.6,42.2)
\qbezier(27.6,42.2)(25.36,34.73)(23.9,25.9)
\put(26.24,42.1){\circle*{.5}}
\put(27.2,42.8){\circle*{.5}}
\put(28.8,43){\circle*{.5}}
\put(27.6,42.2){\circle*{.5}}
\put(26.6,41.5){\circle*{.5}}
\put(23.9,25.9){\circle*{.5}}
\put(24.23,24.2){\circle*{.5}}
\put(22.52,26.13){\circle*{.5}}
\put(28.42,43.46){\circle*{.5}}
\put(23.1,27){\circle*{.5}}
\put(23.4,25.25){\circle*{.5}}
\linethickness{0.2pt}
\qbezier(32,18)(10.04,36.44)(40,48)
\qbezier(40,48)(3.94,36.98)(32,18)
\qbezier(18,52)(24.16,48.46)(28.8,43)
\qbezier(24.8,24.7)(20.66,21.33)(14,18)
\qbezier(14,18)(20.14,23.34)(23,26.8)
\qbezier(26.6,41.5)(23.34,46.42)(18,52)
\qbezier(26.6,41.5)(24.44,34.96)(23.1,27)
\qbezier(28.8,43)(26.44,34.15)(24.8,24.7)
\put(26.07,39.7){\circle*{.5}}
\put(27.22,40.9){\circle*{.5}}
\put(28.59,41.85){\circle*{.5}}
\put(23.44,28.76){\circle*{.5}}
\put(24.12,27.13){\circle*{.5}}
\put(24.96,25.66){\circle*{.5}}
\put(24.8,24.7){\circle*{.5}}
\end{picture}}
\caption{Jump contours of $\widehat{Z}$ (local lenses)}
\label{fig:22}
\end{figure}

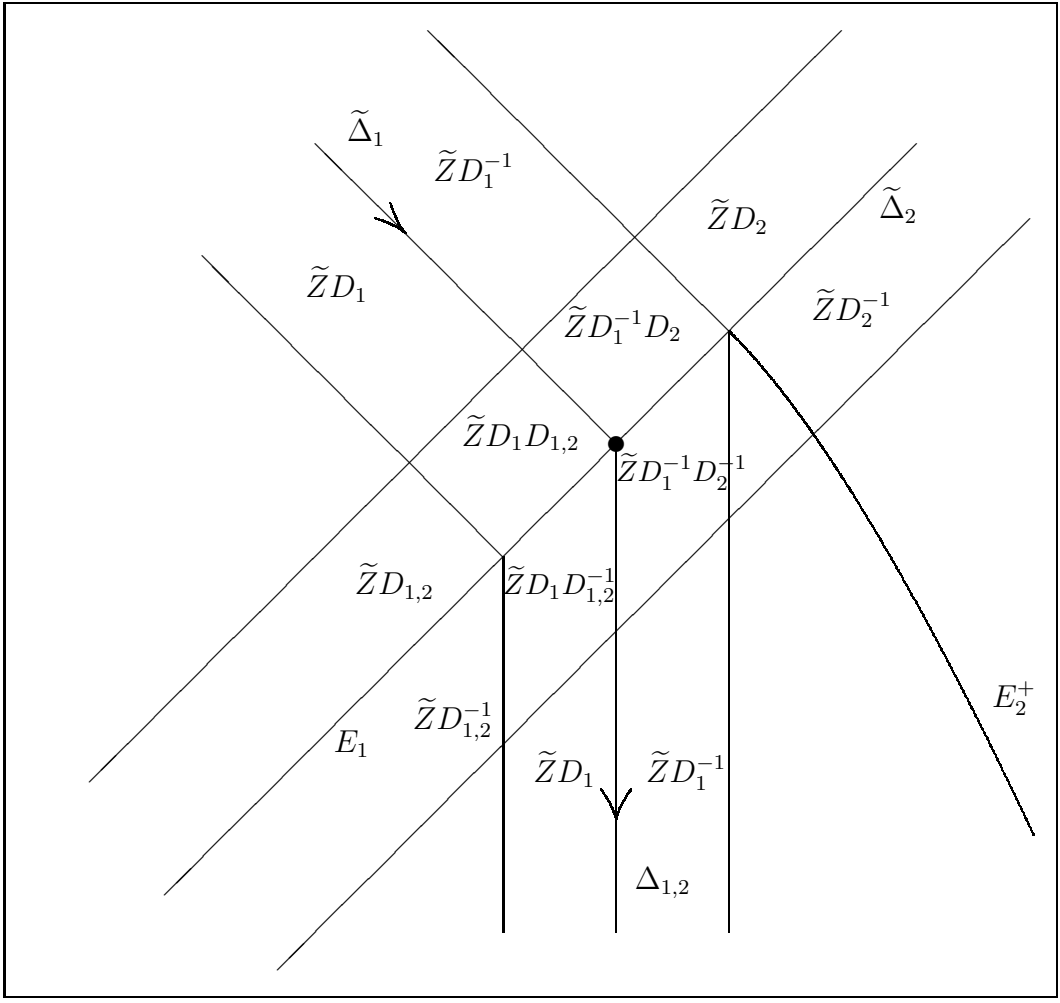
\begin{figure}[ht]
\centering
\framebox{
\unitlength 1mm 
\linethickness{0.4pt}
\begin{picture}(130,125)(0,0)
\put(75,70){\circle*{2}}
\put(105,125){\line(-1,-1){100}}
\put(115,110){\line(-1,-1){100}}
\put(130,100){\line(-1,-1){100}}
\put(75,70){\line(-1,1){40}}
\put(60,55){\line(-1,1){40}}
\put(90,85){\line(-1,1){40}}
\put(90,85){\line(0,-1){80}}
\put(75,70){\line(0,-1){65}}
\put(60,55){\line(0,-1){50}}
\qbezier(90,85)(106.75,68.5)(130.5,18)
\put(125,35){$E_2^+$}
\put(39.25,110.25){$\widetilde{\Delta}_1$}
\put(77.5,10.75){$\Delta_{1,2}$}
\qbezier(47,98)(46,99)(45,102)
\qbezier(47,98)(46,99)(43,100)
\qbezier(75,20)(75.5,22)(77,24)
\qbezier(75,20)(74.5,22)(73,24)
\put(37.5,29){$E_1$}
\put(110,100){$\widetilde{\Delta}_2$}
\put(40.25,50){$\widetilde{Z}D_{1,2}$}
\put(48,32){$\widetilde{Z}D_{1,2}^{-1}$}
\put(54.5,69.75){$\widetilde{Z}D_1D_{1,2}$}
\put(60,50){\small $\widetilde{Z}D_1D_{1,2}^{-1}$}
\put(68,84.25){$\widetilde{Z}D_1^{-1}D_2$}
\put(75,65){\small $\widetilde{Z}D_1^{-1}D_2^{-1}$}
\put(87,98.5){$\widetilde{Z}D_2$}
\put(101,86.5){$\widetilde{Z}D_2^{-1}$}
\put(64,25){$\widetilde{Z}D_1$}
\put(33.75,89.5){$\widetilde{Z}D_1$}
\put(50.75,105){$\widetilde{Z}D_1^{-1}$}
\put(79,25){$\widetilde{Z}D_1^{-1}$}
\end{picture}}
\caption{The function $\widehat{Z}$ (close-up around $c_1$)}
\label{fig:23}
\end{figure}

The goal of the third step is to transform the \RH\ problem (\ref{eq:4.15})--(\ref{eq:4.154}) for the
function $\widetilde{Z}$ from (\ref{eq:4.14}) to a \RH\ problem with jumps which do not depend on $n$
or which tend to the identity matrix as $n \to \infty$. We start with some notation. In addition to the
points
\[     c_j^+ := \partial O_{c_j} \cap \widetilde{\Delta}_2, \qquad j=1,2,  \]
(see Figures~\ref{fig:22}, \ref{fig:23}) we denote \index{Parameters!c1mc2m@$c_1^{-}$, $c_2^{-}$} \index{Parameters!c1tildeplusminusc2tildeplusminus@$\tilde{c}_1^{+}$, $\tilde{c}_2^{+}$, $\tilde{c}_1^{-}$, $\tilde{c}_2^{-}$ }
\[ c_j^- := \partial O_{c_j} \cap E_1, \quad \tilde{c}_j^+ := \partial O_{c_j} \cap \Delta_{1,2},
   \quad \tilde{c}_j^- := \partial O_{c_j} \cap \widetilde{\Delta}_1, \qquad j=1,2. \]
We define a local perturbation of the arcs $\widetilde{\Delta}_2$ joining the points $a_2$
with $\tilde{c}_1^{\pm}$ and $b_2$ with $\tilde{c}_2^{\pm}$ by Jordan arcs $\gamma_{a_2,\tilde{c}_1^{\pm}}$
and $\gamma_{b_2,\tilde{c}_2^{\pm}}$ \index{Contours!Delta1tildeplusminusDelta2tildeplusminus@$\widetilde{\Delta}_1^{+}$, $\widetilde{\Delta}_1^{-}$, $\widetilde{\Delta}_2^{+}$, $\widetilde{\Delta}_2^{-}$}
\[   \widetilde{\Delta}_2^{\pm} := \gamma_{a_2,\tilde{c}_1^{\pm}} \cup \gamma_{b_2,\tilde{c}_2^{\pm}}. \]
Analogously we define (see Figure~\ref{fig:22})
\[   \widetilde{\Delta}_1^{\pm} := \gamma_{a_1,c_1^{\pm}} \cup \gamma_{b_1,c_2^{\pm}}, \quad
     \Delta_{1,2}^{\pm} := \gamma_{c_1^{\pm},c_2^{\pm}}, \quad
     E_1^{\pm} := \gamma_{\tilde{c}_1^\pm,\tilde{c}_2^\pm} .  \]
\index{Contours!Delta12plusminus@$\Delta_{1,2}^{+}$, $\Delta_{1,2}^{-}$}
\index{Contours!E1plusminus@$E_{1}^{+}$, $E_{1}^{-}$}
\index{Contours!gammac1pmc2pm@$\gamma_{c_1^{+},c_2^{+}}$, $\gamma_{c_1^{-},c_2^{-}}$}
\index{Contours!gammac1tildepmc2tildepm@$\gamma_{\tilde{c}_1^{+},\tilde{c}_2^{+}}$, $\gamma_{\tilde{c}_1^{-},\tilde{c}_2^{-}}$}

Around the piecewise analytic arc $\widetilde{\Delta}_1 \cup \Delta_{1,2}$ we now define
the lens shaped domains $T_1^+$ and $T_1^-$ bounded by \index{Domains!T1T1plusT1minus@$T_1$, $T_1^+$, $T_1^-$}
\[   \partial T_1^\pm := \Delta_1 \cup (\widetilde{\Delta}_1^\pm \cup \Delta_{1,2}^\pm), \]
and
\[  T_1 := T_1^+ \cup T_1^- \cup \Delta_1 .  \]
Analogously, around the analytic arc $\widetilde{\Delta}_2$ we define domains $T_2^\pm$ bounded
by
\[  \partial T_2^\pm := \widetilde{\Delta}_2 \cup \widetilde{\Delta}_2^\pm, \]
and \index{Domains!T2T2plusT2minus@$T_2$, $T_2^+$, $T_2^-$}
\[   T_2 := T_2^+ \cup T_2^- \cup \widetilde{\Delta}_2 . \]
Also, around $E_1$ we define the domains $T_{E_1}^\pm$
\[  \partial T_{E_1}^\pm := E_1 \cup E_1^\pm, \]
and \index{Domains!TE1TE1plusTE1minus@$T_{E_1}$, $T_{E_1}^+$, $T_{E_1}^-$}
\[  T_{E_1} := T_{E_1}^+ \cup T_{E_1}^- \cup E_1.  \]
Now we can transform the \RH\ problem (\ref{eq:4.15})--(\ref{eq:4.154}). We define
(see Figure~\ref{fig:23}), \index{Matrices!Zhat@$\widehat{Z}$}
\begin{equation}   \label{eq:4.16}
\widehat{Z} := \begin{cases}
     \widetilde{Z} D_1^{-1}, & \textrm{in $T_1^+ \setminus T_2$}, \\
     \widetilde{Z} D_1, & \textrm{in $T_1^- \setminus T_{E_1}$}, \\
     \widetilde{Z} D_2^{-1}, & \textrm{in $T_2^+ \setminus T_1$}, \\
     \widetilde{Z} D_2, & \textrm{in $T_2^- \setminus T_1$}, \\
     \widetilde{Z} D_{1,2}^{-1}, & \textrm{in $T_E^+ \setminus T_1^-$}, \\
     \widetilde{Z} D_{1,2}, & \textrm{in $T_E^- \setminus T_1^-$}, \\
     \widetilde{Z} D_1^{-1}D_2, & \textrm{in $T_1^+ \cap T_2^-$}, \\
     \widetilde{Z} D_1^{-1}D_2^{-1}, & \textrm{in $T_1^+ \cap T_2^+$}, \\
     \widetilde{Z} D_1D_{1,2}, & \textrm{in $T_1^- \cap T_{E_1}^-$}, \\
     \widetilde{Z} D_1 D_{1,2}^{-1}, & \textrm{in $T_1^- \cap T_{E_1}^+$}, \\
     \widetilde{Z} , & \textrm{in $\overline{\mathbb{C}} \setminus (T_1 \cup T_2 \cup T_{E_1})$}.
                 \end{cases}
\end{equation}

The piecewise analytic matrix function $\widehat{Z}$ has a jump (see (\ref{eq:4.15})) \index{Matrices!Jhat@$\widehat{J}$}
\begin{equation}  \label{eq:4.17}
   \widehat{Z}_+ = \widehat{Z}_- \widehat{J}, \qquad \textrm{on $\widetilde{\Sigma}^+ \cup \partial T_1
                    \cup \partial T_2 \cup \partial T_{E_1} =: \widehat{\Sigma}$}.
\end{equation}
We now \index{Contours!Sigmahat@$\widehat{\Sigma}$}
describe the explicit form of the matrix function $\widehat{J}$ on the different parts
of $\widehat{\Sigma}$. From (\ref{eq:4.153}) we have on $E_1 \setminus T_1$
\[  \widehat{Z}_+ = \widetilde{Z}_+ D_{1,2+}^{-1} = \widetilde{Z}_- D_{1,2-} W_{1,2} D_{1,2+} D_{1,2+}^{-1}
    = \widehat{Z}_- W_{1,2}, \qquad \textrm{on $E_1 \setminus T_1$}. \]
The same relation holds on $E_1 \cap T_1$ since $D_{1+} = D_{1-}$ on $E_1$.
 From (\ref{eq:4.151}), (\ref{eq:4.154}) and (\ref{eq:4.91}) we have on $\Delta_1 \setminus (T_2 \cup T_{E_1})$
\[   \widehat{Z}_+ = \widetilde{Z} D_{1+}^{-1} = \widetilde{Z}_- D_{1-} W_1 D_{1+} D_{1+}^{-1}
     = \widehat{Z}_- W_1, \qquad \textrm{on $\Delta_1 \setminus (T_2 \cup T_{E_1})$}. \]
The remarkable fact is that this relation in fact holds on $\Delta_1 \cap (T_2 \cup T_{E_1})$ as well. Indeed,
on $\widetilde{\Delta}_1 \cap (T_2 \cup T_{E_1})$ we have (see Figure~\ref{fig:23})
\[  \widehat{Z}_+ = \widetilde{Z}_+ D_{1+}^{-1} D_{2+} = \widetilde{Z}_- (D_{1-}D_{1,2} D_{1,2}^{-1} D_{1-}^{-1})
    (D_{1-} W_1 D_{1+}^{-1}) D_{1+}^{-1} D_{2+} = \widehat{Z}_- D_{1,2}^{-1} W_1 D_{2+}. \]
Note that
\[     D_{1,2}^{-1} W_1 D_{2+} = W_1, \qquad \textrm{on $\widetilde{\Delta}_1$}, \]
and therefore
\[   \widehat{Z}_+ = \widehat{Z}_- W_1, \qquad \textrm{on $\widetilde{\Delta}_1 \cap (T_2 \cup T_{E_1})$}.  \]
Analogously, since
\[   D_{1,2} W_1 D_{2+}^{-1} = W_1, \qquad \textrm{on $\Delta_{1,2}$}, \]
we have
\begin{eqnarray*}
  \widehat{Z}_+ &=& \widetilde{Z}_+ D_{1+}^{-1}D_{2+}^{-1} = \widetilde{Z}_- D_{1-}W_1D_{1+}D_{1+}^{-1}D_{2+}^{-1} \\
                &=& \widetilde{Z}_- D_{1-}D_{1,2}^{-1} (D_{1,2} W_1 D_{2+}^{-1}) = \widehat{Z}_- W_1,
    \qquad \textrm{on $\widetilde{\Delta}_{1,2} \cap (T_2 \cup T_{E_1})$}.
\end{eqnarray*}
In the same way we have on $\widetilde{\Delta}_2 \setminus T_1$, see (\ref{eq:4.92});
\[  \widehat{Z}_+ = \widetilde{Z}_+ D_{2+}^{-1} = \widetilde{Z}_- D_{2-} W_2 D_{2+} D_{2-}^{-1}
                  = \widehat{Z}_- W_2, \qquad \textrm{on $\widetilde{\Delta}_2\setminus T_1$}, \]
and on $\widetilde{\Delta}_2 \cap T_1$, see (\ref{eq:4.151}) and (\ref{eq:4.92}),
\[  \widehat{Z}_+ = \widetilde{Z}_+ D_1^{-1} D_2^{-1} = \widetilde{Z}_- (D_2 W_2 D_{2+} D_{2,1+}) D_{1+}^{-1} D_2^{-1}. \]
Using the commutation relations (\ref{eq:4.12}) and the observation that
\[    D_1 W_2 D_{2+} D_{2,1+} D_{1+}^{-1} D_{2+}^{-1} = W_2, \qquad \textrm{on $\widetilde{\Delta}_2$}, \]
we arrive at
\[  \widehat{Z}_+ = \widetilde{Z}_- D_{1-}^{-1}D_{2-} W_2 = \widehat{Z}_- W_2, \qquad
    \textrm{on $\widetilde{\Delta}_2 \cap T_1$}. \]
Summarizing, we have for the jump matrix $\widehat{J}$ in (\ref{eq:4.17})
\begin{equation}  \label{eq:4.18}
  \widehat{J} = \begin{cases}
      D_1, & \textrm{on $\Delta_1^\pm$}, \\
      D_2, & \textrm{on $\widetilde{\Delta}_2^\pm$}, \\
      D_{1,2}, & \textrm{on $E1^\pm$}, \\
      D_{2,1}, & \textrm{on $E_2^+$}, \\
      W_1, & \textrm{on $\Delta_1$}, \\
      W_2, & \textrm{on $\widetilde{\Delta}_2$}, \\
      W_{1,2}, & \textrm{on $E_1$}.
               \end{cases}
\end{equation}
We also note that when we go around the point $c_j$ $(j=1,2)$, then the total jump of $\widehat{Z}$ is
\[   W_2 W_1^{-1} W_{1,2}^{-1} W_1 = I, \]
i.e., the intersection points $c_1$ and $c_2$ are not branch points of $\widehat{Z}$.
Note also that the matrices $W_1, W_2, W_{1,2}$ do not depend on $n$ and due to Proposition~\ref{prop:2.4B}
(see Figure~\ref{fig:11B}) and (\ref{eq:4.13}) the matrices $D_1, D_2, D_{1,2}, D_{2,1}$ tend to the
identity matrix $I$ on compact sets of the domains of their definition in (\ref{eq:4.18}). Therefore
we achieved the goal of the third transformation.

\subsection{Parametrix away from the branch points} \label{sec:4.4}
The goal of the fourth step is to construct a solution for the model \RH\ problem with jumps that do not
depend on $n$, so that later on we can eliminate these jumps in (\ref{eq:4.18}). We are looking
for a function $X$ for which \index{Matrices!X@$X$}
\begin{equation}  \label{eq:4.19}
   \begin{cases}
   X \in H^{3\times 3} (\mathbb{C}\setminus \Sigma), \\
   X_+ = X_- \widehat{W}, \qquad \textrm{ on } \stackrel{\circ}{\Sigma}, \\
   X(z) = I + \O(1/z), \qquad z \to \infty,
   \end{cases}
\end{equation}
where $\Sigma := \Delta_1 \cup \widetilde{\Delta}_2 \cup E_1$ (see Figure~\ref{fig:21a}),
and \index{Matrices!What@$\widehat{W}$}
\begin{equation}  \label{eq:4.20}
  \widehat{W} := \begin{cases}
            W_1, & \textrm{on $\Delta_1$}, \\
            W_2, & \textrm{on $\widetilde{\Delta}_2$}, \\
            W_{1,2}, & \textrm{on $E_1$},
                  \end{cases}
\end{equation}
see the notation in (\ref{eq:4.8}) and (\ref{eq:4.10*}). The function $X$ should also satisfy the
endpoint condition (\ref{eq:2.15})--(\ref{eq:2.17}).
We used the notation $\mbox{$\displaystyle \stackrel{\circ}{\Sigma} :=
\Sigma \setminus \{a_1,b_1,a_2,b_2\}$}$  in (\ref{eq:4.19}). \index{Contours!Sigmacirc@$\stackrel{\circ}{\Sigma}$}

To construct the solution of (\ref{eq:4.19})--(\ref{eq:4.20}) we use the Szeg\H{o} functions
(see subsection~\ref{sec:2.6}), i.e., the solution of the following system of
scalar \BVP s, see (\ref{eq:2.77}), (\ref{eq:2.79}),
\begin{eqnarray}  \label{eq:4.21}
 1) & & \begin{cases}
    F_0 \in H(\overline{\mathbb{C}} \setminus (\Delta_1 \cup \widetilde{\Delta}_2)), \\
    F_1 \in H(\overline{\mathbb{C}} \setminus (E_1 \cup \Delta_1)), \\
    F_2 \in H(\overline{\mathbb{C}} \setminus (E_1 \cup \widetilde{\Delta}_2)),
   \end{cases} \nonumber \\
 2) & &
    \begin{cases}
    F_{0\pm} = iF_{1\mp} \omega_{1-}w_1, & \textrm{on $\Delta_1 \setminus \{a_1,b_1\} =: \stackrel{\circ}{\Delta}_1$}, \\
  F_{0\pm} = iF_{2\mp}\omega_{2-}w_2, & \textrm{on $\widetilde{\Delta}_2 \setminus \{a_2,b_2\} =: \stackrel{\circ}{\widetilde{\Delta}}_2$}, \\
  F_{1\pm} = F_{2\mp} \frac{\omega_{2-}w_2}{\omega_{1-}w_1}, & \textrm{on $E_1$},
  \end{cases} \\
 3) & & \begin{cases}
      \textrm{normalization condition: } F_0F_1F_2 = 1 \textrm{ in } \overline{\mathbb{C}}, \\
      \textrm{local behavior (\ref{eq:2.77})-2b.}
\end{cases}  \nonumber
\end{eqnarray}
Here $\omega_j$ is the branch of $\omega_j^2(z)=(z-a_j)(z-b_j)$, see (\ref{eq:2.74}). The solution
of this system of \BVP s exists and is given in Theorem~\ref{thm:2.10}.
We define \index{Matrices!F@$F$} \index{Matrices!Finfinity@$F_{\infty}$}
\[  F(z) := \textrm{diag} (F_0(z),F_1(z),F_2(z)), \quad F_\infty := F(\infty), \]
and transform the \RH\ problem (\ref{eq:4.19})--(\ref{eq:4.20}) to the \RH\ problem
for the function \index{Matrices!Xtilde@$\widetilde{X}$}
\begin{equation}   \label{eq:4.22}
   \widetilde{X} := F_\infty^{-1} X F.
\end{equation}
We have
\[  \begin{cases}
    \widetilde{X} \in H^{(3\times 3)}(\mathbb{C} \setminus \stackrel{\circ}{\Sigma}), \\
    \widetilde{X}_+ = \widetilde{X}_- \widetilde{H}, \qquad \textrm{on $\stackrel{\circ}{\Sigma}$}, \\
    \widetilde{X}(z) = I + \O(1/z), \qquad z \to \infty,
    \end{cases}   \]
where the jump \index{Matrices!Htilde@$\widetilde{H}$}
\[  \widetilde{H} = F_-^{-1} \widehat{W}F_+ , \]
due to (\ref{eq:4.21})-2, is
\begin{equation}  \label{eq:4.23}
  \widetilde{H} := \begin{cases}
    \begin{pmatrix}
     0 & \displaystyle\frac{F_{1+}}{F_{0-}}\ w_1 & 0 \\
    \displaystyle \frac{-F_{0+}}{F_{1-}w_1} & 0 & 0 \\
    0 & 0 & 1
    \end{pmatrix} =
    \begin{pmatrix}
    0 & \displaystyle \frac{1}{i\omega_{1-}} & 0 \\
    -i\omega_{1-} & 0 & 0 \\
    0 & 0 & 1
    \end{pmatrix}, & \textrm{on $\stackrel{\circ}{\Delta}_1$}, \\
   \begin{pmatrix}
   0 & 0 & \displaystyle \frac{1}{i\omega_{2-}} \\
   0 & 1 & 0 \\
   -i\omega_{2-} & 0 & 0
   \end{pmatrix}, & \textrm{on $\stackrel{\circ}{\widetilde{\Delta}}_2$}, \\
   \begin{pmatrix}
   1 & 0 & 0 \\
   0 & 0 & \displaystyle -\frac{\omega_{1-}}{\omega_{2-}} \\
   0 & \displaystyle \frac{\omega_{2-}}{\omega_{1-}} & 0
   \end{pmatrix}, & \textrm{on $E_1$}.
  \end{cases}
\end{equation}
For the solution of the \RH\ problem (\ref{eq:4.22})--(\ref{eq:4.23}) we take the Riemann
surface $\mathfrak{R}$ from Definition~\ref{def:2.3B} and define on $\mathfrak{R}$
the rational function $\mathcal{X}_j$ with one pole (and one
zero) \index{Functions!XcaljXcalj0Xcalj1Xcalj2@$\mathcal{X}_j$, $\mathcal{X}_j^{(0)}$, $\mathcal{X}_j^{(1)}$, $\mathcal{X}_j^{(2)}$}
\begin{equation}  \label{eq:4.24*}
  \mathcal{X}_j := (\mathcal{X}_j^{(0)},\mathcal{X}_j^{(1)},\mathcal{X}_j^{(2)}) \in \mathfrak{M}(\mathfrak{R}):
  \begin{cases}
   \mathcal{X}_j(\xi) = \O(1/\xi), & \xi \to \infty^{(0)}, \\
   \mathcal{X}_j(\xi) = -\xi + \cdots , & \xi \to \infty^{(j)},
  \end{cases}, \qquad j=1,2.
\end{equation}
Then the function \index{Functions!xtildejk@$\tilde{x}_{j,k}$}
\begin{equation}   \label{eq:4.24}
   \widetilde{X} = \begin{pmatrix}
     \tilde{x}_{0,0} &  \tilde{x}_{0,1} & \tilde{x}_{0,2} \\
      \tilde{x}_{1,0} &  \tilde{x}_{1,1} & \tilde{x}_{1,2} \\
      \tilde{x}_{2,0} &  \tilde{x}_{2,1} & \tilde{x}_{2,2}
                    \end{pmatrix} =
    \begin{pmatrix}
    1 & -\frac{1}{i\omega_1} & -\frac{1}{i\omega_2} \\
    i\mathcal{X}_1^{(0)} & -\frac{1}{\omega_1} \mathcal{X}_1^{(1)} & - \frac{1}{\omega_2} \mathcal{X}_1^{(2)} \\
    i\mathcal{X}_2^{(0)} & -\frac{1}{\omega_1} \mathcal{X}_2^{(1)} & - \frac{1}{\omega_2} \mathcal{X}_2^{(2)}
   \end{pmatrix}
\end{equation}
is the solution of the \RH\ problem (\ref{eq:4.22})--(\ref{eq:4.23}). Indeed, the normalization at
infinity in (\ref{eq:4.22}) clearly holds because of (\ref{eq:4.24*}). To verify the jump condition
in (\ref{eq:4.22})--(\ref{eq:4.23}) we check the relation
\[  \widetilde{X}_+ = \widetilde{X}_- \widetilde{H}  \]
on the different parts of $\stackrel{\circ}{\Sigma}$. Substituting here (\ref{eq:4.23}) we
need to check that
\begin{equation}  \label{eq:4.25}
\begin{pmatrix}
 \tilde{x}_{0,0+} & \tilde{x}_{0,1+} & \tilde{x}_{0,2+} \\
  \tilde{x}_{1,0+} & \tilde{x}_{1,1+} & \tilde{x}_{1,2+} \\
  \tilde{x}_{2,0+} & \tilde{x}_{2,1+} & \tilde{x}_{2,2+}
 \end{pmatrix} = \begin{cases}
 \begin{pmatrix}
 -i\omega_{1-} \tilde{x}_{0,1-} & \frac{1}{i\omega_{1-}}\tilde{x}_{0,0-} & \tilde{x}_{0,2-} \\
  -i\omega_{1-}\tilde{x}_{1,1-} & \frac{1}{i\omega_{1-}}\tilde{x}_{1,0-} & \tilde{x}_{1,2-} \\
  -i\omega_{1-}\tilde{x}_{2,1-} & \frac{1}{i\omega_{1-}}\tilde{x}_{2,0-} & \tilde{x}_{2,2-}
 \end{pmatrix}, & \textrm{on $\stackrel{\circ}{\Delta}_1$}, \\
\begin{pmatrix}
 -i\omega_{2-} \tilde{x}_{0,2-} & \tilde{x}_{0,1-} & \frac{1}{i\omega_{2-}} \tilde{x}_{0,0-} \\
  -i\omega_{2-}\tilde{x}_{1,2-} & \tilde{x}_{1,1-} & \frac{1}{i\omega_{2-}}\tilde{x}_{1,0-} \\
  -i\omega_{2-}\tilde{x}_{2,2-} & \tilde{x}_{2,1-} & \frac{1}{i\omega_{2-}}\tilde{x}_{2,0-}
 \end{pmatrix}, & \textrm{on $\stackrel{\circ}{\widetilde{\Delta}}_2$}, \\
\begin{pmatrix}
 \tilde{x}_{0,0-} & \frac{\omega_{2-}}{\omega_{1-}}\tilde{x}_{0,2-} & -\frac{\omega_{1-}}{\omega_{2-}}\tilde{x}_{0,1-} \\
 \tilde{x}_{1,0-} & \frac{\omega_{2-}}{\omega_{1-}}\tilde{x}_{1,2-} & -\frac{\omega_{1-}}{\omega_{2-}}\tilde{x}_{1,1-} \\
 \tilde{x}_{2,0-} & \frac{\omega_{2-}}{\omega_{1-}}\tilde{x}_{2,2-} & -\frac{\omega_{1-}}{\omega_{2-}}\tilde{x}_{2,1-}
 \end{pmatrix}, & \textrm{on $E_1$}.
\end{cases}
\end{equation}
If we substitute the expression (\ref{eq:4.24}) for $\tilde{x}_{k,\ell}$ $(k,\ell=0,1,2)$
into (\ref{eq:4.25}), then we indeed see that (\ref{eq:4.25}) holds identically.
Thus the jump condition in (\ref{eq:4.22})--(\ref{eq:4.23}) holds and $\widetilde{X}$ in
(\ref{eq:4.24}) is the solution of the \RH\ problem (\ref{eq:4.22})--(\ref{eq:4.23}).
Finally we have from (\ref{eq:2.79}) and (\ref{eq:4.24*}) that \index{Matrices!X@$X$}
\begin{equation}  \label{eq:4.26}
  X := F_\infty \widetilde{X} F =
  \begin{pmatrix}
   \displaystyle \frac{F_0(\infty)}{F_0}  & \displaystyle -\frac{1}{i\omega_1} \frac{F_0(\infty)}{F_1} &
   \displaystyle - \frac{1}{i\omega_2} \frac{F_0(\infty)}{F_2} \\
   \displaystyle i \mathcal{X}_1^{(0)} \frac{F_1(\infty)}{F_0}  & \displaystyle -\frac{1}{\omega_1} \mathcal{X}_1^{(1)}\frac{F_1(\infty)}{F_1} &    \displaystyle - \frac{1}{\omega_2} \mathcal{X}_1^{(2)} \frac{F_1(\infty)}{F_2} \\
   \displaystyle i \mathcal{X}_2^{(0)} \frac{F_2(\infty)}{F_0}  & \displaystyle -\frac{1}{\omega_1} \mathcal{X}_2^{(1)}\frac{F_2(\infty)}{F_1} &    \displaystyle - \frac{1}{\omega_2} \mathcal{X}_2^{(2)} \frac{F_2(\infty)}{F_2}
   \end{pmatrix}
\end{equation}
is the desired solution of the model \RH\ problem (\ref{eq:4.19})--(\ref{eq:4.20}).

\begin{figure}[ht]
\centering
\unitlength 1.00mm 
\linethickness{0.4pt}
\ifx\plotpoint\undefined\newsavebox{\plotpoint}\fi 
\begin{picture}(119.13,88.13)(0,0)
\put(8,82){\circle*{1}}
\put(7.63,54.37){\circle*{1}}
\put(6,62){\circle*{1}}
\put(5.63,34.37){\circle*{1}}
\put(24,62){\circle*{1}}
\put(24,5.12){\circle*{1}}
\put(24,80){\circle*{1}}
\put(24,23.12){\circle*{1}}
\put(16,68){\circle*{1}}
\put(16,11.12){\circle*{1}}
\put(15.63,40.37){\circle*{1}}
\put(16,76){\circle*{1}}
\put(16,19.12){\circle*{1}}
\put(15.63,48.37){\circle*{1}}
\put(15,69){\circle*{1}}
\put(15,12.12){\circle*{1}}
\put(14.63,41.37){\circle*{1}}
\put(17,69){\circle*{1}}
\put(17,12.12){\circle*{1}}
\put(15,75){\circle*{1}}
\put(15,18.12){\circle*{1}}
\put(14.63,47.37){\circle*{1}}
\put(17,75){\circle*{1}}
\put(17,18.12){\circle*{1}}
\qbezier(17,75)(16.5,72)(17,69)
\qbezier(24,62)(22,60)(16,68)
\qbezier(15,69)(15.5,72)(15,75)
\qbezier(14.63,41.37)(15.13,44.37)(14.63,47.37)
\qbezier(8,81.88)(6.81,80.94)(9.88,79.25)
\qbezier(7.63,54.25)(6.44,53.31)(9.5,51.63)
\qbezier(9.88,79.25)(14.56,76.69)(15,74.88)
\qbezier(24,79.75)(24.5,78.63)(20,77.5)
\qbezier(24,22.88)(24.5,21.75)(20,20.63)
\qbezier(16,67.88)(6.5,60.69)(6,61.75)
\qbezier(15.63,40.25)(6.13,33.06)(5.63,34.13)
\qbezier(5.88,61.75)(4.5,63.31)(7.63,64.13)
\qbezier(5.5,34.13)(4.13,35.69)(7.25,36.5)
\qbezier(7.63,64.13)(13.69,65.81)(15,68.75)
\qbezier(7.88,81.88)(9.31,83.56)(16,76)
\qbezier(7.5,54.25)(8.94,55.94)(15.63,48.38)
\qbezier(16,75.88)(23.06,81.56)(23.88,80)
\qbezier(20.25,77.38)(17.38,76.38)(17,74.88)
\qbezier(17,68.88)(19.75,65.56)(22,64.5)
\qbezier(22,64.5)(25,63)(24,62)
\qbezier(22,7.63)(25,6.12)(24,5.12)
\qbezier(14.63,47.25)(9.88,43.94)(14.63,41.38)
\qbezier(15.63,48.25)(16.13,44.81)(15.63,40.13)
\qbezier(9.5,51.5)(12.38,49.75)(13.25,49)
\qbezier(7.25,36.63)(10.75,37.69)(12.75,40)
\qbezier(12.75,40.13)(7.75,44.44)(13.25,49)
\qbezier(20.13,20.63)(10.19,15.19)(22,7.5)
\qbezier(21.38,23.13)(6.88,16.5)(20.88,5.88)
\qbezier(21,5.63)(22.88,3.94)(23.75,5)
\qbezier(21.38,23.13)(24,24.44)(24.13,23)
\qbezier(49.13,42.5)(48.69,62.75)(77,63)
\qbezier(77,63)(111.81,63.31)(112.38,40.88)
\qbezier(112.38,40.88)(112.88,22.06)(79.88,21)
\qbezier(79.88,21)(49.5,20.06)(49.13,42.38)
\put(66.88,61.88){\circle*{1}}
\put(68,21.63){\circle*{1}}
\put(51.38,32.13){\circle*{1}}
\put(49.13,42.38){\circle*{1}}
\put(50.88,51.63){\circle*{1}}
\put(91.75,61.75){\circle*{1}}
\put(93,22.25){\circle*{1}}
\put(78.38,63){\circle*{1}}
\put(80.13,20.75){\circle*{1}}
\put(57,58.5){\circle*{1}}
\qbezier(66.88,62)(72.31,44.13)(68,21.75)
\qbezier(91.75,61.88)(86.75,43.69)(92.75,22.25)
\put(38,46.75){\makebox(0,0)[cc]{$\longrightarrow$}}
\put(38,50.25){\makebox(0,0)[cc]{$\zeta$}}
\put(0,1){\framebox(30,28)[cc]{}}
\put(0,30){\framebox(30,27)[cc]{}}
\put(0,58){\framebox(30,27)[cc]{}}
\put(43,15){\framebox(77,55)[cc]{}}
\put(110.63,65.63){\makebox(0,0)[cc]{$\zeta(\mathfrak{R}_0)$}}
\put(24.13,71.38){\makebox(0,0)[cc]{$\mathfrak{R}_0$}}
\put(23.13,43.38){\makebox(0,0)[cc]{$\mathfrak{R}_1$}}
\put(25.25,14.75){\makebox(0,0)[cc]{$\mathfrak{R}_2$}}
\put(25.38,80.63){\makebox(0,0)[cc]{$a_2$}}
\put(25.5,60.88){\makebox(0,0)[cc]{$b_2$}}
\put(3.63,61.13){\makebox(0,0)[cc]{$b_1$}}
\put(5.63,82.75){\makebox(0,0)[cc]{$a_1$}}
\put(3.63,33.13){\makebox(0,0)[cc]{$b_1$}}
\put(5.13,55){\makebox(0,0)[cc]{$a_1$}}
\put(25.63,23.63){\makebox(0,0)[cc]{$a_2$}}
\put(25.25,3.88){\makebox(0,0)[cc]{$b_2$}}
\put(78,65.63){\makebox(0,0)[cc]{$\zeta(a_2)$}}
\put(79.88,18.38){\makebox(0,0)[cc]{$\zeta(b_2)$}}
\put(48,30){\makebox(0,0)[cc]{$\zeta(b_1)$}}
\put(54.38,61){\makebox(0,0)[cc]{$\zeta(a_1)$}}
\put(79.75,42.75){\makebox(0,0)[cc]{$\zeta(\mathfrak{R}_2)$}}
\put(58,43){\makebox(0,0)[cc]{$\zeta(\mathfrak{R}_1)$}}
\put(98.25,42.5){\makebox(0,0)[cc]{$\zeta(\mathfrak{R}_1)$}}
\put(66.5,64.38){\makebox(0,0)[cc]{$\zeta(\tilde{c}_1)$}}
\put(16,78.25){\makebox(0,0)[cc]{$\tilde{c}_1$}}
\put(85,67){\vector(-1,-3){3}}
\put(85.5,66){\small $\log \tilde{w}_2$}
\put(60,65){\vector(1,-3){3}}
\put(51.5,65){\small $\log \tilde{w}_1$}
\put(73,53){\vector(-1,0){8}}
\put(69,54){\small $\log \tilde{w}_2-\log \tilde{w}_1$}
\end{picture}
\caption{Conformal mapping of $\mathfrak{R}$ to $\overline{\mathbb{C}}$ (case II)}
\label{fig:24}
\end{figure}
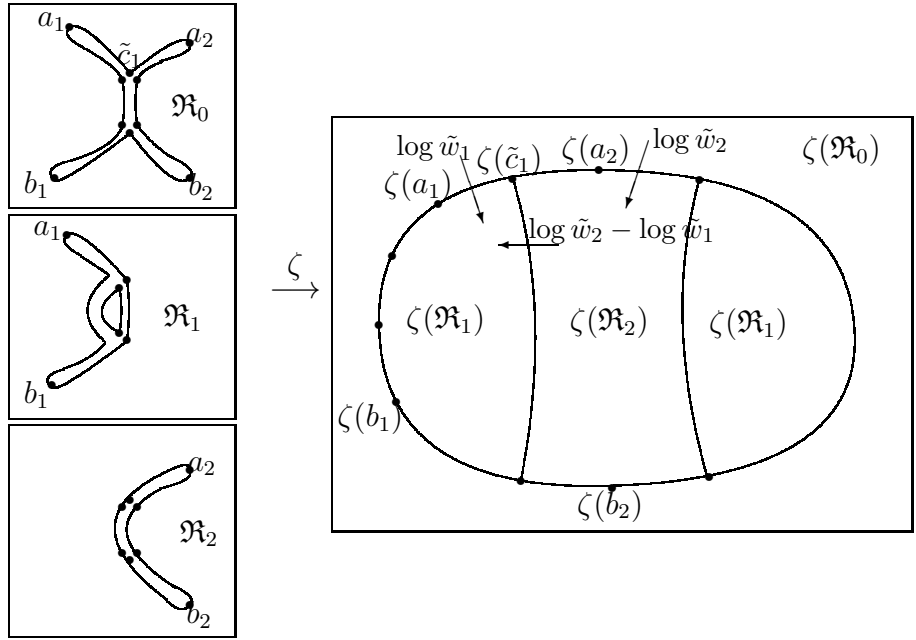

\begin{remark}
We recall that in Theorem~\ref{thm:2.10}  we presented the solution of the scalar \BVP\ (\ref{eq:2.77})
on the Riemann surface by means of the Cauchy integral with a meromorphic differential on $\mathfrak{R}$.
There is another way to solve the scalar \BVP\ (\ref{eq:2.77}); see \cite{48} for a similar approach.
The idea is the following. Since $\mathfrak{R}$ has genus zero, it is conformally equivalent to
$\overline{\mathbb{C}}$. Let $\zeta$ be a conformal mapping \index{Functions!zeta@$\zeta$}
\[  \zeta(\mathfrak{R}) = \overline{\mathbb{C}},  \quad \zeta(\infty^{(0)}) = \infty. \]
We consider $\mathfrak{R}$ from Definition~\ref{def:2.3B} (see Figure~\ref{fig:10B}). The image of
$\mathfrak{R}$ is given in Figure~\ref{fig:24}. Analyzing the jumps around the points $\zeta(c_j)$,
see the notation (\ref{eq:2.75}), we see that the special form of the jumps (\ref{eq:2.75'})
provides the existence of a continuous solution in the neighborhood of the points $c_j$ $(j=1,2)$.
\end{remark}

\subsection{Local parametrices}  \label{sec:4.5}
The function $X$ is not a good approximation near the branch points $A:=\{a_1,b_1,a_2,b_2\}$. We
need a local analysis near each of the branch points. The goal of the fifth step is to find a matrix
function $U$
in a neighborhood $O_e$ \index{Domains!Oe@$O_e$}
of each branch point $e \in A$ which asymptotically (as $n \to \infty$)
matches $X$ on the boundary $\partial O_e$ and which has the same jumps as $\widehat{Z}$. This
function is the solution of the following \RH\ problem: if $e_j$ denotes $a_j$ or $b_j$ $(j=1,2)$,
then \index{Parameters!e1e2@$e_1$, $e_2$}
we want to find a function $U_{e_j}$ such that \index{Matrices!Ue1Ue2@$U_{e_1}$, $U_{e_2}$}
\begin{equation} \label{eq:4.27}
 \begin{cases}
   1) & U_{e_j} \in H^{3\times 3}(O_{e_j} \setminus (\widetilde{\Delta}_j \cup \widetilde{\Delta}_j^{+} \cup
   \widetilde{\Delta}_j^{-})), \\[10pt]
   2) & U_{e_j+} = U_{e_j-} \widehat{J} \quad \textrm{on }
       \widehat{\Sigma}_{e_j} := O_{e_j} \cap
       (\widetilde{\Delta}_j \cup \widetilde{\Delta}_j^+ \cup
   \widetilde{\Delta}_j^-):   \\[10pt]
   3) & U_{e_j} = (I + \O(1/n)) X \qquad \textrm{uniformly on $\partial O_{e_j}$ as $n\to \infty$},
\end{cases}
\end{equation}
and \index{Contours!Sigmaejhat@$\widehat{\Sigma}_{e_1}$, $\widehat{\Sigma}_{e_2}$}
the elements of the matrix $U_{e_j}(z)$ have the same limiting behavior as the elements of
$\widehat{Z}(z)$ as $z \to e_j$.

We recall the explicit expressions of the jump matrix $\widehat{J}$, see (\ref{eq:4.18}) and (\ref{eq:4.8}),
\[ \widehat{J} = \begin{cases}
   \begin{pmatrix}
   1 & 0 & 0 \\
   \displaystyle \frac{\Phi_0^n}{\Phi_1^n w_1} & 1 & 0 \\
   0 & 0 & 1
   \end{pmatrix}, & \textrm{on $\widetilde{\Delta}_1^{(\pm)}$}, \\
  \begin{pmatrix}
   0 & w_1 & 0 \\
   \displaystyle -\frac{1}{w_1} & 0 & 0 \\
   0 & 0 & 1
   \end{pmatrix}, & \textrm{on $\widetilde{\Delta}_1$}, \\
   \begin{pmatrix}
   1 & 0 & 0 \\
   0 & 1 & 0 \\
   \displaystyle \frac{\Phi_0^n}{\Phi_2^n w_2} & 0 & 1
   \end{pmatrix}, & \textrm{on $\widetilde{\Delta}_2^{(\pm)}$}, \\
   \begin{pmatrix}
   0 & 0 & w_2 \\
   0 & 1 & 0 \\
   \displaystyle -\frac{1}{w_2} & 0 & 0
   \end{pmatrix}, & \textrm{on $\widetilde{\Delta}_2$}.
   \end{cases}  \]
We see that the jump matrices related to a branch point from $A$ have a $2\times 2$ block structure.
The solution of the \RH\ problem (\ref{eq:4.27}) in the $2\times 2$ case, with the local behavior
of the equilibrium measure as in (\ref{eq:2.377}), was obtained in \cite{41}. Thus, using the results in
\cite{41}, it is easy to write down the explicit form of the solution for the local \RH\ problem
(\ref{eq:4.27}). This solution is
\begin{equation}   \label{eq:4.28}
   U_{e_j} = E_{e_j} V_{e_j} A_{e_j}, \qquad j=1,2,
\end{equation}
where the matrices $A_{e_j}$ are \index{Matrices!Ae1Ae2@$A_{e_1}$, $A_{e_2}$}
\begin{eqnarray*}
  A_{e_1} & := & \textrm{diag} \left( \Bigl[ \left(\frac{\Phi_0}{\Phi_1} \right)^{n/2} \check{w}_{e_1} \Bigr]^{-1},
       \left( \frac{\Phi_0}{\Phi_1} \right)^{n/2} \check{w}_{e_1}, 1 \right), \\
  A_{e_2} & := & \textrm{diag} \left( \Bigl[ \left(\frac{\Phi_0}{\Phi_2} \right)^{n/2} \check{w}_{e_2} \Bigr]^{-1},
       1, \left( \frac{\Phi_0}{\Phi_2} \right)^{n/2} \check{w}_{e_2} \right),
\end{eqnarray*}
and, see (\ref{eq:2.5}), \index{Functions!wajcheckwbjcheck@$\check{w}_{a_j}$, $\check{w}_{b_j}$}
\begin{eqnarray*}
   \check{w}_{a_j}(\xi) &=& \left[ w_{0,j}(\xi)(b_j -\xi)^{\beta_j} (a_j -\xi)^{\alpha_j} \right]^{1/2}, \\
   \check{w}_{b_j}(\xi) &=& \left[ w_{0,j}(\xi)(\xi-b_j)^{\beta_j} (\xi-a_j)^{\alpha_j} \right]^{1/2}.
\end{eqnarray*}
The branch  of the square roots above is chosen such that
\[    \check{w}_{e_j+} \check{w}_{e_j-} = w_j, \qquad \textrm{on $\Delta_j \cap O_{e_j}$},\
     e_j \in \{a_j,b_j\},\ j=1,2. \]
The matrices $E_{e_j}$ in (\ref{eq:4.28}) are \index{Matrices!Ee1Ee2@$E_{e_1}$, $E_{e_2}$}
\index{Matrices!Ea1Ea2Eb1Eb2@$E_{a_1}$, $E_{a_2}$, $E_{b_1}$, $E_{b_2}$}
\begin{eqnarray*}
   E_{a_1} & := & \frac12 X \textrm{diag} \bigl(\check{w}_{a_1}, \check{w}_{a_1}^{-1}, 1 \bigr)
                   \begin{pmatrix} 1 & i & 0 \\ i & 1 & 0 \\ 0 & 0 & 1 \end{pmatrix}
                     \textrm{diag} \left( \sqrt{\pi n \varphi_1}, \frac{1}{\sqrt{\pi n \varphi_1}}, 1 \right), \\
   E_{b_1} & := & \frac12 X \textrm{diag} \bigl(\check{w}_{b_1}, \check{w}_{b_1}^{-1}, 1 \bigr)
                   \begin{pmatrix} 1 & -i & 0 \\ -i & 1 & 0 \\ 0 & 0 & 1 \end{pmatrix}
                     \textrm{diag} \left( \sqrt{\pi n \varphi_1}, \frac{1}{\sqrt{\pi n \varphi_1}}, 1 \right), \\
   E_{a_2} & := & \frac12 X \textrm{diag} \bigl(\check{w}_{a_2}, 1, \check{w}_{a_2}^{-1} \bigr)
                   \begin{pmatrix} 1 & i & 0 \\ i & 1 & 0 \\ 0 & 0 & 1 \end{pmatrix}
                     \textrm{diag} \left( \sqrt{\pi n \varphi_2}, 1, \frac{1}{\sqrt{\pi n \varphi_2}} \right), \\
   E_{b_2} & := & \frac12 X \textrm{diag} \bigl(\check{w}_{b_2}, 1, \check{w}_{b_2}^{-1} \bigr)
                   \begin{pmatrix} 1 & -i & 0 \\ -i & 1 & 0 \\ 0 & 0 & 1 \end{pmatrix}
                     \textrm{diag} \left( \sqrt{\pi n \varphi_2}, 1,  \frac{1}{\sqrt{\pi n \varphi_2}} \right).
\end{eqnarray*}
Here $X$ is given by (\ref{eq:4.26}) and we use the notation \index{Functions!phi1phi2@$\varphi_1$, $\varphi_2$}
\[     \varphi_j := \log \left( \frac{\Phi_0}{\Phi_j} \right).   \]
To present the explicit expression for the matrices $V_{e_j}$ in (\ref{eq:4.28}) we follow \cite{41}
and introduce the matrices \index{Matrices!Psia1Psia2Psib1Psib2@$\Psi_{a_1}$, $\Psi_{a_2}$, $\Psi_{b_1}$, $\Psi_{b_2}$}
\begin{eqnarray*}
  \Psi_{a_1} & := &
   \begin{pmatrix}
I_{\alpha_1} \left( -\frac{n\varphi_1}{2} \right) & -\frac{i}{\pi} K_{\alpha_1} \left( -\frac{n\varphi_1}{2} \right) & 0 \\
n\pi i \varphi_1 I_{\alpha_1}'\left( -\frac{n\varphi_1}{2} \right) &
  n \varphi_1 K_{\alpha_1}'\left( -\frac{n\varphi_1}{2} \right) & 0 \\
   0 & 0 & 1
   \end{pmatrix}, \\
  \Psi_{b_1} & := &
   \begin{pmatrix}
I_{\beta_1} \left( \frac{n\varphi_1}{2} \right) & \frac{i}{\pi} K_{\beta_1} \left( \frac{n\varphi_1}{2} \right) & 0 \\
n\pi i \varphi_1 I_{\beta_1}'\left( \frac{n\varphi_1}{2} \right) &
  -n \varphi_1 K_{\beta_1}'\left( \frac{n\varphi_1}{2} \right) & 0 \\
   0 & 0 & 1
   \end{pmatrix}, \\
\Psi_{a_2} & := &
   \begin{pmatrix}
I_{\alpha_2} \left( -\frac{n\varphi_2}{2} \right) & 0 & -\frac{i}{\pi} K_{\alpha_2} \left( -\frac{n\varphi_2}{2} \right) \\
   0 & 1 & 0 \\
n\pi i \varphi_1 I_{\alpha_2}'\left( -\frac{n\varphi_2}{2} \right) & 0 &
  n \varphi_2 K_{\alpha_2}'\left( -\frac{n\varphi_2}{2} \right)
    \end{pmatrix}, \\
 \Psi_{b_2} & := &
   \begin{pmatrix}
I_{\beta_2} \left( \frac{n\varphi_2}{2} \right) & 0 & \frac{i}{\pi} K_{\beta_2} \left( \frac{n\varphi_2}{2} \right) \\
   0 & 1 & 0 \\
  -n\pi i \varphi_2 I_{\beta_2}'\left( \frac{n\varphi_2}{2} \right) & 0 &
  n \varphi_2 K_{\beta_2}'\left( \frac{n\varphi_2}{2} \right)
    \end{pmatrix},
\end{eqnarray*}
where $I_{\alpha}$ and $K_\alpha$ are the modified Bessel functions of order $\alpha$ (see \cite{AS}).

We denote by $O_{e_j}^{(+)}$ \index{Domains!OejplusOejminus@$O_{e_j}^{(+)}$, $O_{e_j}^{(-)}$}
the sector domain bounded by
\[   \partial O_{e_j}^{(+)} := \partial O_{e_j} \cup (\widetilde{\Delta}_j \cap O_{e_j}) \cup
     ( \widetilde{\Delta}_j^+ \cap O_{e_j} ).  \]
Analogously, $O_{e_j}^{(-)}$ has the boundaries
\[   \partial O_{e_j}^{(-)} := \partial O_{e_j} \cup (\widetilde{\Delta}_j \cap O_{e_j}) \cup
     ( \widetilde{\Delta}_j^- \cap O_{e_j} ),  \]
and \index{Domains!Oejstar@$O_{e_j}^{(*)}$}
\[   O_{e_j}^{(*)} := O_{e_j} \setminus \overline{( O_{e_j}^{(+)} \cup O_{e_j}^{(-)})} . \]
The matrices $V_{e_j}$ in (\ref{eq:4.28}) are \index{Matrices!Va1Va2Vb1Vb2@$V_{a_1}$, $V_{a_2}$, $V_{b_1}$, $V_{b_2}$}
\[  V_{a_1} := \begin{cases}
     \Psi_{a_1} , & \textrm{in $O_{a_1}^{(*)}$}, \\
     \Psi_{a_1} \begin{pmatrix}
                1 & 0 & 0 \\ -e^{-\alpha_1\pi i} & 1 & 0 \\ 0 & 0 & 1 \end{pmatrix}, & \textrm{in $O_{a_1}^{(+)}$}, \\
     \Psi_{a_1} \begin{pmatrix}
                1 & 0 & 0 \\ e^{\alpha_1\pi i} & 1 & 0 \\ 0 & 0 & 1 \end{pmatrix}, & \textrm{in $O_{a_1}^{(-)}$},
               \end{cases} , \quad
  V_{b_1} := \begin{cases}
     \Psi_{b_1} , & \textrm{in $O_{b_1}^{(*)}$}, \\
     \Psi_{b_1} \begin{pmatrix}
                1 & 0 & 0 \\ e^{\beta_1\pi i} & 1 & 0 \\ 0 & 0 & 1 \end{pmatrix}, & \textrm{in $O_{b_1}^{(+)}$}, \\
     \Psi_{b_1} \begin{pmatrix}
                1 & 0 & 0 \\ -e^{-\beta_1\pi i} & 1 & 0 \\ 0 & 0 & 1 \end{pmatrix}, & \textrm{in $O_{b_1}^{(-)}$},
               \end{cases} \]
\[  V_{a_2} := \begin{cases}
     \Psi_{a_2} , & \textrm{in $O_{a_2}^{(*)}$}, \\
     \Psi_{a_2} \begin{pmatrix}
                1 & 0 & 0 \\ 0 & 1 & 0 \\  -e^{-\alpha_2\pi i} & 0 & 1 \end{pmatrix}, & \textrm{in $O_{a_2}^{(+)}$}, \\
     \Psi_{a_2} \begin{pmatrix}
                1 & 0 & 0 \\ 0 & 1 & 0 \\  e^{\alpha_2\pi i} & 0 & 1 \end{pmatrix}, & \textrm{in $O_{a_2}^{(-)}$},
               \end{cases} , \quad
  V_{b_2} := \begin{cases}
     \Psi_{b_2} , & \textrm{in $O_{b_2}^{(*)}$}, \\
     \Psi_{b_2} \begin{pmatrix}
                1 & 0 & 0 \\ 0 & 1 & 0 \\  e^{\beta_2\pi i} & 0 & 1 \end{pmatrix}, & \textrm{in $O_{b_2}^{(+)}$}, \\
     \Psi_{b_2} \begin{pmatrix}
                1 & 0 & 0 \\ 0 & 1 & 0 \\  -e^{-\beta_2\pi i} & 0 & 1 \end{pmatrix}, & \textrm{in $O_{b_2}^{(-)}$}.
               \end{cases} \]
This gives all the ingredients for the solution (\ref{eq:4.28}) of the local \RH\ problem (\ref{eq:4.27}).

\subsection{Final transformation. Asymptotic formulas} \label{sec:4.6}
\begin{figure}[ht]
\centering
\framebox{
\unitlength 4pt
\linethickness{0.2pt}
\ifx\plotpoint\undefined\newsavebox{\plotpoint}\fi 
\begin{picture}(50,55)(0,5)
\put(0,50){$\textup{II}$}
\put(18,52){\circle*{1}}
\put(14,18){\circle*{1}}
\put(40,48){\circle*{1}}
\put(32,18){\circle*{1}}
\put(17,54){$a_1$}
\put(12,16){$b_1$}
\put(41,49){$a_2$}
\put(33,17){$b_2$}
\put(27.6,42.2){\circle*{.5}}
\put(23.9,25.9){\circle*{.5}}
\qbezier(26.6,41.5)(24.44,34.96)(23.1,27)
\put(18,52){\circle{7}}
\put(40,48){\circle{7}}
\put(32,18){\circle{7}}
\put(14,18){\circle{7}}
\qbezier(36.38,47.38)(8.69,36.25)(28.75,19.63)
\qbezier(36.88,46.38)(12.13,36.06)(29.38,20.5)
\qbezier(17.25,19.38)(43.13,35.44)(21,50.25)
\qbezier(20.13,49.38)(24.63,45.81)(26.63,41.5)
\qbezier(16.5,20.38)(20.5,23.56)(23,27)
\qbezier(28.88,42.88)(26.5,35.81)(24.5,24.6)
\end{picture}}
\caption{The contour $\check{\Sigma}$ for the final transformation}
\label{fig:25}
\end{figure}
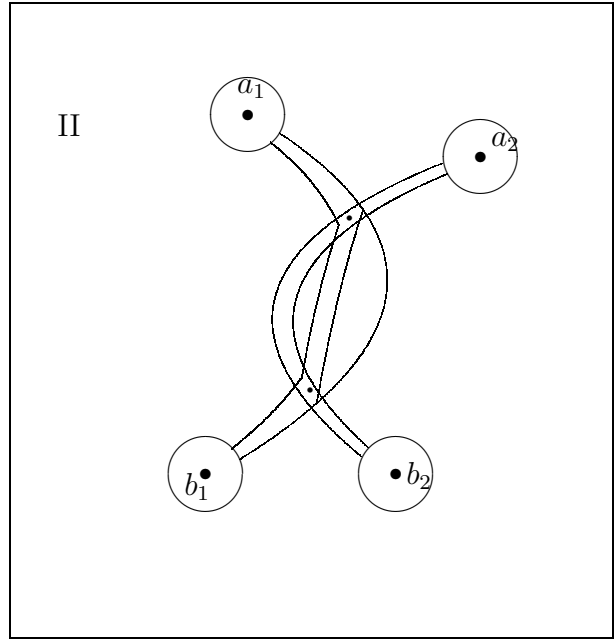

We now finish the transformation of the original \RH\ problem (\ref{eq:2.13})--(\ref{eq:2.17}) to
the \RH\ problem (\ref{eq:4.1}). The final transformation
is \index{Matrices!Zcheck@$\check{Z}$}
\begin{equation}  \label{eq:4.29}
 \check{Z} := \begin{cases}
    \widehat{Z} X^{-1}, & \textrm{in } \overline{\mathbb{C}} \setminus \bigcup_{j=1}^2 (O_{a_j} \cup O_{b_j}),\
    \textrm{(away from the branch points)}, \\
    \widehat{Z} U_{e_j}^{-1}, & \textrm{in $O_{e_j}$}, e_j \in \{a_j,b_j\}, j=1,2, \textrm{ (near the branch points)}.
    \end{cases}
\end{equation}
We have on $\Delta_1 \cup \widetilde{\Delta}_2 \cup E_1$
\[   \check{Z}_+ = \widehat{Z}_- \widehat{J} X^{-1} = \widehat{Z}_- X_-^{-1} ( X_- J X_+^{-1}) .  \]
 From (\ref{eq:4.19})--(\ref{eq:4.20}) we have
\[   X_- \widehat{J} X_+^{-1} = X_- \widehat{W} (X_- \widehat{W})^{-1}, \]
so that
\[   \check{Z}_+ = \check{Z}_- , \qquad \textrm{on } (\Delta_1 \cup \widetilde{\Delta}_2 \cup E_1) \setminus O_{e_1}. \]
Analogously, see (\ref{eq:4.27}), we have
\[   \check{Z}_+ = \check{Z}_- , \qquad \textrm{on } (\Delta_1 \cup \widetilde{\Delta}_2 \cup E_1) \setminus O_{e_2}. \]
Thus, denoting (see Figure~\ref{fig:25}) \index{Contours!Sigmacheck@$\check{\Sigma}$}
\[   \check{\Sigma} := \widetilde{\Delta}_1^{\pm} \cup \widetilde{\Delta}_2^{\pm} \cup
      \Delta_{1,2}^{\pm} \cup E_1^{\pm} \cup E_2^{+} \cup \bigcup_{e_j \in A} \partial O_{e_j}, \]
we have \index{Matrices!Jcheck@$\check{J}$}
\begin{equation}  \label{eq:4.30}
  \begin{cases}
   1) & \check{Z} \in H^{(3\times 3)}(\mathbb{C} \setminus \check{\Sigma}), \\[10pt]
   2) & \check{Z}_+ = \check{Z}_- \check{J}, \quad \textrm{on }
       \check{\Sigma}, \\[10pt]
   3) & \check{Z}(z) = I + \O(1/z), \quad z \to \infty.
  \end{cases}
\end{equation}

Since as $n \to \infty$
\[  \check{J} = I + \O(1/n), \qquad \textrm{uniformly on $\bigcup_{e_j \in A} \partial O_{e_j}$}, \]
and for some $C > 0$
\[  \check{J} = I + \O(e^{-Cn}), \qquad \textrm{uniformly on } \check{\Sigma} \setminus \bigcup_{e_j \in A} \partial O_{e_j}, \]
we can conclude that the solution of the problem (\ref{eq:4.30}) exists for $n$ large enough, and
\begin{equation}  \label{eq:4.301}
   \check{Z} = I + \O(1/n),
\end{equation}
uniformly for $z \in \overline{\mathbb{C}}$.

Now we can return through the sequence of transformations
(\ref{eq:4.29}), (\ref{eq:4.16}), (\ref{eq:4.14}), (\ref{eq:4.3}),
\begin{equation}  \label{eq:4.302}
    \check{Z} \to \widehat{Z} \to \widetilde{Z} \to Z \to Y
\end{equation}
and we can conclude that for large $n$ the solution of the \RH\ problem (\ref{eq:2.13}) exists, i.e., the
\HP\ approximants exists and the indices are normal, and we get
the asymptotic formulas of Theorems~\ref{thm:2.11} and \ref{thm:2.12}.
To obtain these formulas we note that in the transformation (\ref{eq:4.302}) it is sufficient to keep track of
the first row of the matrices. We start with the asymptotics on compact sets in $\overline{\mathbb{C}} \setminus
\widetilde{\Sigma}$, where $\widetilde{\Sigma} := \Sigma \cup E_2$,
see (\ref{eq:4.5}). \index{Contours!Sigmatilde@$\widetilde{\Sigma}$}
Comparing the inverse and direct transformation in (\ref{eq:4.302}), we have
\[  \widetilde{Z} = (I + \O(1/n)) X, \qquad \textrm{on $\overline{\mathbb{C}} \setminus
\widetilde{\Sigma}$}, \]
and
\[  \widetilde{Z} = Z  = \begin{pmatrix}
                  C_0^n \Phi_0^n P_n & C_0^n \Phi_1^n R_n^{(1)} & C_0^n \Phi_2^n R_n^{(2)} \\
                  * & * & * \\ * & * & * \end{pmatrix}, \qquad
    \textrm{in $(\overline{\mathbb{C}} \setminus \overline{G}) \setminus (\widetilde{\Delta}_1 \cup
      \widetilde{\Delta}_2)$}, \]
\[  \widetilde{Z} = Z D_{2,1} = \begin{pmatrix}
         C_0^n \Phi_0^n P_n & C_0^n \Phi_1^n R_n^{(1)} & C_0^n\Phi_2^nR_n^{(2)} -C_0^n \Phi_2^n R_n^{(1)}w_2/w_1 \\
                  * & * & * \\ * & * & * \end{pmatrix}, \qquad
    \textrm{in $G \setminus \Delta_{1,2}$}. \]
This gives
\begin{eqnarray}
   P_n & = & \frac{X_{0,0}}{C_0^n \Phi_0^n} (1 + \O(1/n)), \qquad \textrm{on $ \overline{\mathbb{C}} \setminus
\widetilde{\Sigma}$}, \label{eq:4.31} \\
   R_n^{(1)} & = & \frac{X_{0,1}}{C_0^n\Phi_1^n} (1+ \O(1/n)), \qquad \textrm{on $ \overline{\mathbb{C}} \setminus
\widetilde{\Sigma}$}, \label{eq:4.32} \\
  R_n^{(2)} & = & \begin{cases}
        \displaystyle \frac{X_{0,2}}{C_0^n\Phi_2^n} (1+\O(1/n)), & \textrm{on $ (\overline{\mathbb{C}} \setminus
\overline{G}) \setminus (\widetilde{\Delta}_1 \cup \widetilde{\Delta}_2)$}, \\[10pt]
 \displaystyle \left( \frac{X_{0,2}}{C_0^n\Phi_2^n} + \frac{X_{0,1}}{C_0^n\Phi_1^n} \frac{w_2}{w_1} \right), &
   \textrm{on $ G \setminus \Delta_{1,2}$}.
   \end{cases}  \label{eq:4.33}
\end{eqnarray}
The asymptotic formulas (\ref{eq:4.31})--(\ref{eq:4.33}) are valid uniformly on compact subsets of the indicated sets.
If we substitute (\ref{eq:4.26}) in these formulas, then we arrive at the formulas of Theorems~\ref{thm:2.11}
and \ref{thm:2.12} for $\overline{\mathbb{C}} \setminus \widetilde{\Sigma}$.

Next we consider the asymptotics on compact subsets of $\widetilde{\Sigma} \setminus \{ a_1,b_1,a_2,b_2\}$.
Since $X_{0,0}, \Phi_0 \in H(\mathbb{C} \setminus (\widetilde{\Delta}_1 \cup \widetilde{X}_2 \cup \Delta_{1,2}))$,
the asymptotic formula  (\ref{eq:4.31}) holds on compact subsets in $\overline{\mathbb{C}} \setminus
(\widetilde{\Delta}_1 \cup \widetilde{X}_2 \cup \Delta_{1,2}))$ and therefore in particular on $E_1 \cup E_2$.
For the same reason the asymptotic formula (\ref{eq:4.32}) holds on compact subsets of $\overline{\mathbb{C}}
\setminus \Sigma$, hence in particular on $E_2$, and the second asymptotic formula in (\ref{eq:4.33}) holds
on $E_2$. Thus the asymptotic behavior on $E_2$ is verified.

Now consider the asymptotics on $E_1$. On $E_1$ we have, see (\ref{eq:4.19})--(\ref{eq:4.20}),
\[   X_{0,1+} = X_{0,2-} \frac{w_1}{w_2}, \quad \Phi_{1+} =  \Phi_{2-}, \qquad \textrm{on $E_1$}, \]
hence we have from (\ref{eq:4.33})
\begin{equation}  \label{eq:4.34}
  R_n^{(2)} = \frac{X_{0,2-}}{C_0^n \Phi_{2-}^n} (1+\O(1/n)) = \frac{X_{0,1+}}{C_0^n \Phi_{1+}^n} \frac{w_2}{w_1}
       (1+\O(1/n)), \qquad \textrm{on $E_1$}.
\end{equation}
In $T_{E_1}$, see (\ref{eq:4.16}), we have
\[   \widehat{Z} = \begin{cases}
     Z D_{2,1} D_{1,2}^{-1}, & \textrm{in $T_{E_1}^{(+)}$}, \\
     Z D_{1,2}, & \textrm{in $T_{E_1}^{(-)}$},
     \end{cases}
      \]
and therefore
\[   \widehat{Z} = \begin{cases}
     \begin{pmatrix}
       C_0^n \Phi_0^n P_n & \displaystyle C_0^n \Phi_1^n \frac{w_1}{w_2} R_n^{(2)} &
      \displaystyle C_0^n\Phi_2^nR_n^{(2)} - C_0^n \Phi_2^n \frac{w_2}{w_1} R_n^{(1)} \\
      * & * & * \\ * & * & * \end{pmatrix}, & \textrm{on $T_{E_1}^{(+)}$}, \\
     \begin{pmatrix}
      C_0^n \Phi_0^n P_n & \displaystyle C_0^n \Phi_1^nR_n^{(1)} - C_0^n\Phi_1^n R_n^{(2)} \frac{w_1}{w_2} &
      C_0^n\Phi_2^nR_n^{(2)} \\
      * & * & * \\ * & * & * \end{pmatrix}, & \textrm{on $T_{E_1}^{(-)}$}.
     \end{cases} \]
On the other hand we have
\[ \widehat{Z} =  (I + \O(1/n)) X,
 \qquad \textrm{in $\overline{\mathbb{C}}$} \setminus \bigcup_{j=1}^2 (O_{a_j} \cup O_{b_j}), \]
hence using (\ref{eq:4.34}) we arrive at
\begin{eqnarray*}
  R_n^{(1)} &=& \left( \frac{X_{0,1+}}{C_0^n\Phi_{1+}^n} - \frac{X_{0,2+}}{C_0^n\Phi_{2+}^n} \frac{w_1}{w_2} \right)
    (1+\O(1/n)) \\
      &=& \left( \frac{X_{0,2-}}{C_0^n\Phi_{2-}^n}\frac{w_1}{w_2} + \frac{X_{0,1-}}{C_0^n\Phi_{2-}^n}  \right)
    (1+\O(1/n)), \qquad \textrm{on $E_1$}.
\end{eqnarray*}
Thus the asymptotics on $E_1$ are also verified.

Now consider the asymptotics on $\widetilde{\Delta}_j$ $(j=1,2)$. In $T_j$, away from $c_1$ and $c_2$, we have
\[ \widehat{Z} = \begin{cases}
      \widetilde{Z} D_j^{-1} = \begin{pmatrix}
      \displaystyle C_0^n \Phi_0^nP_n  - C_0^n \frac{\Phi_0^n}{w_j} R_n^{(j)}
      & C_0^n\Phi_1^n R_n^{(1)} & C_0^n \Phi_2^n R_n^{(2)} \\
      * & * & * \\ * & * & * \end{pmatrix}, & \textrm{on $T_j^{(+)}$}, \\
     \widetilde{Z} D_j = \begin{pmatrix}
      \displaystyle C_0^n \Phi_0^nP_n  + C_0^n \frac{\Phi_0^n}{w_j} R_n^{(j)}
       & C_0^n\Phi_1^n R_n^{(1)} & C_0^n \Phi_2^n R_n^{(2)} \\
      * & * & * \\ * & * & * \end{pmatrix}, & \textrm{on $T_j^{(-)}$}.
     \end{cases}  \]
Taking the limiting value we obtain
\begin{equation}  \label{eq:4.35}
   R_{n\pm}^{(1)} = \frac{X_{0,1\pm}}{C_0^n\Phi_{1\pm}^n} (1+\O(1/n)), \quad
   R_{n\pm}^{(2)} = \frac{X_{0,2\pm}}{C_0^n\Phi_{2\pm}^n} (1+\O(1/n)), \qquad \textrm{on $ \widetilde{\Delta}_j$},
  j=1,2.
\end{equation}
Using these asymptotics we arrive at
\begin{eqnarray}  \label{eq:4.36}
  P_n &=& \left( \frac{X_{0,0+}}{C_0^n \Phi_{0+}^n} + \frac{X_{0,j+}}{w_j C_0^n \Phi_{j+}^{n}} \right) (1+\O(1/n))
   \nonumber \\
      &=& \left( \frac{X_{0,0-}}{C_0^n\Phi_{0-}^n} - \frac{X_{0,j-}}{w_j C_0^n \Phi_{j-}^n} \right) (1+\O(1/n))
   \nonumber \\
      &=& \left[ \left( \frac{X_{0,0}}{C_0^n\Phi_0^n} \right)_+ + \left( \frac{X_{0,0}}{C_0^n\Phi_0^n} \right)_- \right]
      (1+\O(1/n)),
    \qquad \textrm{on $ \widetilde{\Delta}_k$}, j=1,2.
\end{eqnarray}

Finally we consider the asymptotics on $\Delta_{1,2}$. On $\Delta_{1,2}^{(\pm)}$ we have the limiting values
\[  \widehat{Z}_{\pm} = (\widetilde{Z} D_1^{-1})_{\pm} =
    \begin{pmatrix}
    C_0^n \Phi_0^n P_n - C_0^n \frac{\Phi_0^n}{w_1} R_n^{(1)} & C_0^n \Phi_1^n R_n^{(1)} &
    C_0^n \Phi_2^n R_n^{(2)} - C_0^n \frac{w_2}{w_1} \Phi_2^n R_n^{(1)} \\
    * & * & * \\ * & * & * \end{pmatrix}_{\pm}.  \]
From this we have that the asymptotics (\ref{eq:4.36}) for $P_n$ and (\ref{eq:4.35}) for $R_{n\pm}^{(1)}$
are also valid on $\Delta_{1,2}$. For $R_n^{(2)}$ we have
\[   C_0^n \Phi_2^n R_{n\pm}^{(2)} - C_0^n \Phi_2^n \frac{w_2}{w_1} \frac{X_{0,1\pm}}{C_0^n\Phi_1^n} =
     X_{0,2} (1+ \O(1/n)), \qquad \textrm{in $\Delta_{1,2}$}. \]
Taking into account that $|\Phi_0| = |\Phi_1| < |\Phi_2|$ in $\Delta_{1,2}$ (Proposition~\ref{prop:2.4B}) we
obtain
\begin{equation}  \label{eq:4.37}
    R_{n\pm}^{(2)} = \frac{w_2}{w_1} \frac{X_{0,1\pm}}{C_0^n \Phi_{1\pm}^n} (1+\O(1/n)), \qquad
    \textrm{in $\Delta_{1,2}$}.
\end{equation}
It remains to observe that the limiting values of the asymptotic formulas (\ref{eq:4.34})--(\ref{eq:4.37})
along $\widetilde{\Sigma}$ at the points $c_1$ and $c_2$ coincide.

This proves Theorems~\ref{thm:2.11} and \ref{thm:2.12} for the case II.

\begin{remark}  \label{rem:4.1}
We can get formulas for the local asymptotics for $P_n$ and $R_n^{(j)}$ $(j=1,2)$ from (\ref{eq:4.301}) by
using the transformations (\ref{eq:4.302}) in $O_{e_j}$, with $e_j \in \{a_j,b_j\}$. These formulas are
similar to the corresponding formulas from \cite{41} and contain Bessel functions.
\end{remark}

\subsection{Sketch of the proof for the other geometrical cases}  \label{sec:4.7}
\subsubsection*{Case I}
For this case $\Delta_1 \cap \Delta_2 = \emptyset$, hence we do not need the second transformation
(\ref{eq:4.11}), i.e., we do not need to open a global lens. All the other steps are as in the proof
for the case II. We just have to use for $\{\Phi_0,\Phi_1,\Phi_2\}$ the algebraic function (\ref{eq:2.25}),
(\ref{eq:2.36'}) and we have to drop some parts of the proof which involve the intersection $\Delta_1 \cap \Delta_2$.

\subsubsection*{Case IV}
This is also a generic case. All the steps of the proof above are present for this case. The proof now uses the
Riemann surface (\ref{eq:2.51}) and for the normalization we use the branches (\ref{eq:2.52}) of the
function $\Phi$:  \index{Matrices!Z@$Z$}
\[   Z := C Y S, \]
where we use the notation (\ref{eq:4.2}). The second step is to open the global lens around $\Delta_{1,2}$
in the domain $G^{+}$ with boundary $\partial G^{(+)} = E_2 \cup E_1^{+}$
(see Figure~\ref{fig:26} and \index{Contours!E1plus@$E_{1}^{+}$}
the notation (\ref{eq:2.53'})) and we obtain \index{Domains!Gplus@$G^+$}
\[  \widetilde{Z} := \begin{cases}
    Z D_{1,2}, & \textrm{in $G^{+}$}, \\
    Z, & \textrm{in $\overline{\mathbb{C}} \setminus G^{+}$}.
    \end{cases}  \]
Here we use the notation (\ref{eq:4.10*}). The function $\widetilde{Z}$ has the jump
\[  \widetilde{Z}_+ = \widetilde{Z}_- \widetilde{J}, \qquad \textrm{on $\Delta_1 \cup \widetilde{\Delta}_2
                        \cup E_2 \cup E_1^{+}$}, \]
where \index{Matrices!Jtilde@$\widetilde{J}$}
\[   \widetilde{J} := D_{1,2}, \qquad \textrm{on $E_1^{+}$}, \]
and, see (\ref{eq:4.8}),
\begin{equation}   \label{eq:4.38}
   \widetilde{J} := \begin{cases}
       D_{1-} W_1 D_{1+}, & \textrm{on $\widetilde{\Delta}_1$}, \\
       D_{2-} W_2 D_{2+}, & \textrm{on $\Delta_2 := \widetilde{\Delta}_2 \cup \Delta_{1,2}$}.
     \end{cases}
\end{equation}
If we use the analog of (\ref{eq:4.11})
\[   J(0,0) = D_{2,1-} W_{1,2} D_{2,1+} D_{1,2+}^{-1}, \qquad \textrm{on $E_2$}, \]
then we have
\begin{equation}  \label{eq:4.39}
   \widetilde{J} := D_{2,1-} W_{1,2} D_{2,1+}, \qquad \textrm{on $E_2$}.
\end{equation}
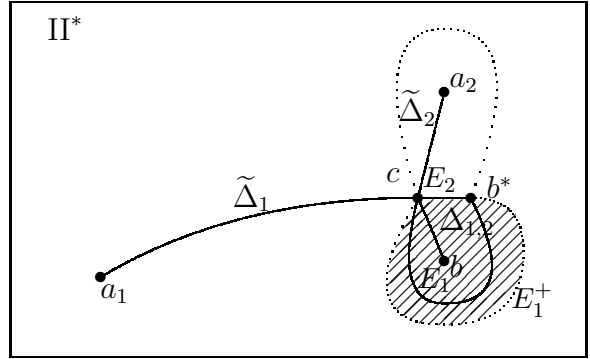
\begin{figure}[ht]
\centering
\framebox{\unitlength 2.00pt 
\linethickness{0.4pt}
\ifx\plotpoint\undefined\newsavebox{\plotpoint}\fi 
\begin{picture}(96.81,57)(0,0)
\put(0,55){$\textup{II}^*$}
\put(10,10){\circle*{2}}
\put(75,45){\circle*{2}}
\put(70,25){\circle*{2}}
\put(75,13){\circle*{2}}
\put(80,25){\circle*{2}}
\qbezier(10,10)(35,25)(70,25)
\qbezier(70,25)(75,25)(80,25)
\qbezier(75,13)(72.5,19.5)(70,25)
\put(0,0){\bezier{25}(70,25)(60,57)(75,57)}
\put(0,0){\bezier{25}(75,57)(92,57)(80,25)}
\qbezier(75,5)(90,5)(80,25)
\qbezier(70,25)(65,5)(75,5)
\qbezier(75,45)(72.5,35)(70,25)
\put(10,6){$a_1$}
\put(83,25){$b^*$}
\put(76,46){$a_2$}
\put(76,10){$b$}
\put(64,28){$c$}
\put(35,23.38){$\widetilde{\Delta}_1$}
\put(66.3,39.5){$\widetilde{\Delta}_2$}
\put(71,27){\small $E_2$}
\put(0,0){\bezier{25}(70,25)(56,1)(75,1)}
\put(0,0){\bezier{20}(75,1)(90,1)(90,15)}
\put(0,0){\bezier{20}(90,15)(90,25)(80,25)}
\put(87.75,4.13){\small $E_1^{+}$}
\put(74,19){$\Delta_{1,2}$}
\put(70,8){\small $E_1$}
\put(72.88,24.75){\line(-1,-1){5.7}}
\put(74.88,24.75){\line(-1,-1){9}}
\put(76.63,24.75){\line(-1,-1){12}}
\put(78.38,24.75){\line(-1,-1){13.8}}
\put(79.88,24){\line(-1,-1){15.5}}
\put(83,24.8){\line(-1,-1){18.2}}
\put(84.6,24.4){\line(-1,-1){19.3}}
\put(86,23.8){\line(-1,-1){19.9}}
\put(87.4,22.9){\line(-1,-1){20.1}}
 \put(88.5,21.8){\line(-1,-1){19.9}}
 \put(89.1,20.2){\line(-1,-1){18.8}}
 \put(89.6,18.7){\line(-1,-1){17.5}}
 \put(90,17){\line(-1,-1){16}}
 \put(90,15){\line(-1,-1){14}}
 \put(90,13){\line(-1,-1){12}}
 \put(89.5,10){\line(-1,-1){8}}
\end{picture}}
\caption{The global lens (case $\textup{IV}$)}
\label{fig:26}
\end{figure}

The third step is the opening of the local lenses around the arcs $\Delta_1, \widetilde{\Delta}_2, E_2$ (see
Section \ref{sec:4.3}). The decompositions (\ref{eq:4.38})--(\ref{eq:4.39}) give as a result a new function
$\widehat{Z}$ for which the jumps on these arcs do not depend on $n$
\[  \widehat{Z}_+ = \widehat{Z}_- \widehat{J}, \]
with \index{Matrices!Zhat@$\widehat{Z}$} \index{Matrices!Jhat@$\widehat{J}$}
\[  \widehat{J} = W := \begin{cases}
      W_1, & \textrm{on $\widetilde{\Delta}_1$}, \\
      W_2, & \textrm{on $\Delta_2$}, \\
      W_{1,2}, & \textrm{on $E_2$},
   \end{cases} \]
and on the lenses $\widetilde{\Delta}_1^{\pm}, \Delta_{1,2}^{\pm}, \widetilde{\Delta}_2^{\pm},
E_2^{\pm}$ \index{Contours!E2plusE2minus@$E_2^+$, $E_2^-$}
and $E_1^{+}$ (excluding the points $a_1,a_2,b,b^*$) the jump $\widehat{J}$ tends to the identity
matrix as $n \to \infty$ (due to Proposition~\ref{prop:2.8}).

The fourth step (the parametrix away from the branch points) is the solution of the \RH\ problem for
the function $X$ with the jump $W$. This is just a repetition of Section~\ref{sec:4.4} but exchanging
$E_1$ and $E_2$ and with another Riemann surface $\mathfrak{R}$, namely (\ref{eq:2.51}) (see Figure~\ref{fig:17}).

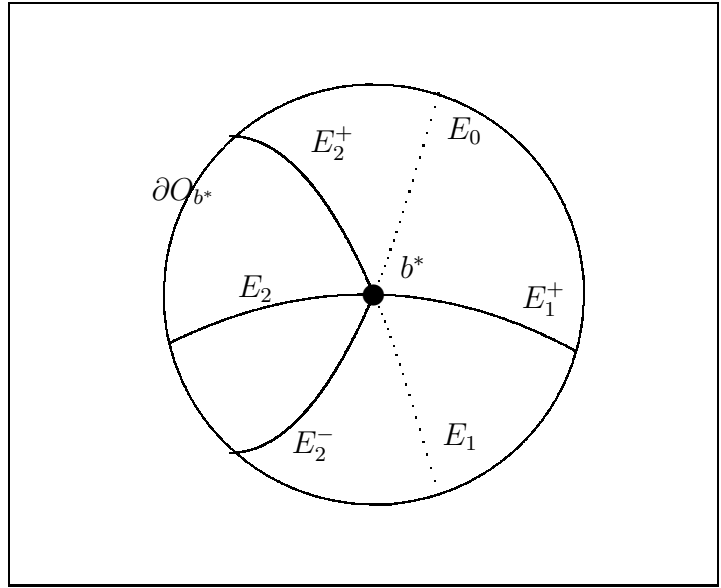
\begin{figure}[ht]
\centering
\framebox{
\unitlength 2pt
\linethickness{0.4pt}
\begin{picture}(120,100)(0,0)
\put(60,50){\circle*{4}}
\put(99.75,50){\line(0,1){1.423}}
\put(99.73,51.42){\line(0,1){1.421}}
\multiput(99.65,52.84)(-.0318,.3543){4}{\line(0,1){.3543}}
\multiput(99.52,54.26)(-.02964,.23532){6}{\line(0,1){.23532}}
\multiput(99.34,55.67)(-.03261,.20066){7}{\line(0,1){.20066}}
\multiput(99.12,57.08)(-.03093,.15506){9}{\line(0,1){.15506}}
\multiput(98.84,58.47)(-.03282,.13847){10}{\line(0,1){.13847}}
\multiput(98.51,59.86)(-.03146,.11434){12}{\line(0,1){.11434}}
\multiput(98.13,61.23)(-.0328,.104436){13}{\line(0,1){.104436}}
\multiput(97.71,62.59)(-.031648,.089436){15}{\line(0,1){.089436}}
\multiput(97.23,63.93)(-.032652,.08273){16}{\line(0,1){.08273}}
\multiput(96.71,65.25)(-.033499,.076714){17}{\line(0,1){.076714}}
\multiput(96.14,66.56)(-.03241,.067522){19}{\line(0,1){.067522}}
\multiput(95.52,67.84)(-.033066,.063003){20}{\line(0,1){.063003}}
\multiput(94.86,69.1)(-.033619,.058837){21}{\line(0,1){.058837}}
\multiput(94.16,70.34)(-.032599,.052588){23}{\line(0,1){.052588}}
\multiput(93.41,71.55)(-.033024,.049246){24}{\line(0,1){.049246}}
\multiput(92.61,72.73)(-.033375,.046111){25}{\line(0,1){.046111}}
\multiput(91.78,73.88)(-.033658,.04316){26}{\line(0,1){.04316}}
\multiput(90.9,75)(-.032668,.038933){28}{\line(0,1){.038933}}
\multiput(89.99,76.09)(-.032867,.036438){29}{\line(0,1){.036438}}
\multiput(89.04,77.15)(-.033012,.034063){30}{\line(0,1){.034063}}
\multiput(88.05,78.17)(-.03421,.03286){30}{\line(-1,0){.03421}}
\multiput(87.02,79.16)(-.036584,.032704){29}{\line(-1,0){.036584}}
\multiput(85.96,80.1)(-.040526,.033698){27}{\line(-1,0){.040526}}
\multiput(84.86,81.01)(-.04331,.033465){26}{\line(-1,0){.04331}}
\multiput(83.74,81.88)(-.046259,.033169){25}{\line(-1,0){.046259}}
\multiput(82.58,82.71)(-.049393,.032804){24}{\line(-1,0){.049393}}
\multiput(81.4,83.5)(-.052732,.032364){23}{\line(-1,0){.052732}}
\multiput(80.18,84.25)(-.058986,.033356){21}{\line(-1,0){.058986}}
\multiput(78.94,84.95)(-.06315,.032785){20}{\line(-1,0){.06315}}
\multiput(77.68,85.6)(-.067666,.032109){19}{\line(-1,0){.067666}}
\multiput(76.4,86.21)(-.076863,.033156){17}{\line(-1,0){.076863}}
\multiput(75.09,86.78)(-.082875,.032283){16}{\line(-1,0){.082875}}
\multiput(73.76,87.29)(-.095974,.033481){14}{\line(-1,0){.095974}}
\multiput(72.42,87.76)(-.104581,.032334){13}{\line(-1,0){.104581}}
\multiput(71.06,88.18)(-.11448,.03095){12}{\line(-1,0){.11448}}
\multiput(69.69,88.55)(-.13861,.0322){10}{\line(-1,0){.13861}}
\multiput(68.3,88.87)(-.1552,.03024){9}{\line(-1,0){.1552}}
\multiput(66.9,89.15)(-.2008,.03171){7}{\line(-1,0){.2008}}
\multiput(65.5,89.37)(-.23545,.02859){6}{\line(-1,0){.23545}}
\multiput(64.08,89.54)(-.3545,.0302){4}{\line(-1,0){.3545}}
\put(62.67,89.66){\line(-1,0){1.421}}
\put(61.25,89.73){\line(-1,0){1.423}}
\put(59.82,89.75){\line(-1,0){1.423}}
\put(58.4,89.72){\line(-1,0){1.421}}
\multiput(56.98,89.64)(-.3542,-.0334){4}{\line(-1,0){.3542}}
\multiput(55.56,89.5)(-.23518,-.03069){6}{\line(-1,0){.23518}}
\multiput(54.15,89.32)(-.20051,-.03351){7}{\line(-1,0){.20051}}
\multiput(52.75,89.08)(-.15492,-.03163){9}{\line(-1,0){.15492}}
\multiput(51.35,88.8)(-.13832,-.03343){10}{\line(-1,0){.13832}}
\multiput(49.97,88.46)(-.1142,-.03197){12}{\line(-1,0){.1142}}
\multiput(48.6,88.08)(-.104289,-.033265){13}{\line(-1,0){.104289}}
\multiput(47.24,87.65)(-.089294,-.032047){15}{\line(-1,0){.089294}}
\multiput(45.9,87.17)(-.082584,-.033021){16}{\line(-1,0){.082584}}
\multiput(44.58,86.64)(-.07231,-.031961){18}{\line(-1,0){.07231}}
\multiput(43.28,86.06)(-.067377,-.032711){19}{\line(-1,0){.067377}}
\multiput(42,85.44)(-.062855,-.033347){20}{\line(-1,0){.062855}}
\multiput(40.74,84.78)(-.056019,-.032341){22}{\line(-1,0){.056019}}
\multiput(39.51,84.06)(-.052442,-.032833){23}{\line(-1,0){.052442}}
\multiput(38.31,83.31)(-.049098,-.033244){24}{\line(-1,0){.049098}}
\multiput(37.13,82.51)(-.045962,-.033581){25}{\line(-1,0){.045962}}
\multiput(35.98,81.67)(-.041417,-.032597){27}{\line(-1,0){.041417}}
\multiput(34.86,80.79)(-.038787,-.032842){28}{\line(-1,0){.038787}}
\multiput(33.77,79.87)(-.036291,-.033029){29}{\line(-1,0){.036291}}
\multiput(32.72,78.91)(-.033916,-.033164){30}{\line(-1,0){.033916}}
\multiput(31.7,77.92)(-.032707,-.034356){30}{\line(0,-1){.034356}}
\multiput(30.72,76.89)(-.033703,-.038041){28}{\line(0,-1){.038041}}
\multiput(29.78,75.82)(-.033517,-.040676){27}{\line(0,-1){.040676}}
\multiput(28.87,74.73)(-.033272,-.043459){26}{\line(0,-1){.043459}}
\multiput(28.01,73.6)(-.032963,-.046407){25}{\line(0,-1){.046407}}
\multiput(27.19,72.44)(-.032584,-.049539){24}{\line(0,-1){.049539}}
\multiput(26.4,71.25)(-.033589,-.05528){22}{\line(0,-1){.05528}}
\multiput(25.66,70.03)(-.033093,-.059135){21}{\line(0,-1){.059135}}
\multiput(24.97,68.79)(-.032503,-.063295){20}{\line(0,-1){.063295}}
\multiput(24.32,67.52)(-.033574,-.071576){18}{\line(0,-1){.071576}}
\multiput(23.72,66.23)(-.032813,-.07701){17}{\line(0,-1){.07701}}
\multiput(23.16,64.92)(-.031913,-.083018){16}{\line(0,-1){.083018}}
\multiput(22.65,63.6)(-.033052,-.096123){14}{\line(0,-1){.096123}}
\multiput(22.18,62.25)(-.031867,-.104725){13}{\line(0,-1){.104725}}
\multiput(21.77,60.89)(-.03321,-.12503){11}{\line(0,-1){.12503}}
\multiput(21.4,59.51)(-.03158,-.13876){10}{\line(0,-1){.13876}}
\multiput(21.09,58.13)(-.03324,-.17475){8}{\line(0,-1){.17475}}
\multiput(20.82,56.73)(-.03082,-.20094){7}{\line(0,-1){.20094}}
\multiput(20.61,55.32)(-.03305,-.28269){5}{\line(0,-1){.28269}}
\put(20.44,53.91){\line(0,-1){1.418}}
\put(20.33,52.49){\line(0,-1){1.422}}
\put(20.26,51.07){\line(0,-1){1.423}}
\put(20.25,49.65){\line(0,-1){1.423}}
\put(20.29,48.22){\line(0,-1){1.42}}
\multiput(20.38,46.8)(.02797,-.28323){5}{\line(0,-1){.28323}}
\multiput(20.52,45.39)(.03174,-.23504){6}{\line(0,-1){.23504}}
\multiput(20.71,43.98)(.0301,-.17532){8}{\line(0,-1){.17532}}
\multiput(20.95,42.57)(.03232,-.15478){9}{\line(0,-1){.15478}}
\multiput(21.24,41.18)(.03096,-.12561){11}{\line(0,-1){.12561}}
\multiput(21.58,39.8)(.03248,-.11405){12}{\line(0,-1){.11405}}
\multiput(21.97,38.43)(.03373,-.104139){13}{\line(0,-1){.104139}}
\multiput(22.41,37.08)(.032445,-.08915){15}{\line(0,-1){.08915}}
\multiput(22.9,35.74)(.033389,-.082436){16}{\line(0,-1){.082436}}
\multiput(23.43,34.42)(.032283,-.072167){18}{\line(0,-1){.072167}}
\multiput(24.01,33.12)(.033011,-.06723){19}{\line(0,-1){.06723}}
\multiput(24.64,31.84)(.033627,-.062705){20}{\line(0,-1){.062705}}
\multiput(25.31,30.59)(.032591,-.055874){22}{\line(0,-1){.055874}}
\multiput(26.03,29.36)(.033067,-.052295){23}{\line(0,-1){.052295}}
\multiput(26.79,28.16)(.033462,-.048949){24}{\line(0,-1){.048949}}
\multiput(27.59,26.98)(.032486,-.044049){26}{\line(0,-1){.044049}}
\multiput(28.44,25.84)(.032781,-.041271){27}{\line(0,-1){.041271}}
\multiput(29.32,24.72)(.033014,-.03864){28}{\line(0,-1){.03864}}
\multiput(30.25,23.64)(.033191,-.036143){29}{\line(0,-1){.036143}}
\multiput(31.21,22.59)(.033315,-.033767){30}{\line(0,-1){.033767}}
\multiput(32.21,21.58)(.035692,-.033676){29}{\line(1,0){.035692}}
\multiput(33.24,20.6)(.038191,-.033533){28}{\line(1,0){.038191}}
\multiput(34.31,19.66)(.040825,-.033335){27}{\line(1,0){.040825}}
\multiput(35.41,18.76)(.043607,-.033077){26}{\line(1,0){.043607}}
\multiput(36.55,17.9)(.046554,-.032755){25}{\line(1,0){.046554}}
\multiput(37.71,17.09)(.049683,-.032362){24}{\line(1,0){.049683}}
\multiput(38.9,16.31)(.055429,-.033342){22}{\line(1,0){.055429}}
\multiput(40.12,15.58)(.059282,-.032829){21}{\line(1,0){.059282}}
\multiput(41.37,14.89)(.06344,-.03222){20}{\line(1,0){.06344}}
\multiput(42.64,14.24)(.071725,-.033254){18}{\line(1,0){.071725}}
\multiput(43.93,13.64)(.077155,-.032469){17}{\line(1,0){.077155}}
\multiput(45.24,13.09)(.088704,-.033645){15}{\line(1,0){.088704}}
\multiput(46.57,12.59)(.096269,-.032623){14}{\line(1,0){.096269}}
\multiput(47.92,12.13)(.104866,-.031399){13}{\line(1,0){.104866}}
\multiput(49.28,11.72)(.12518,-.03265){11}{\line(1,0){.12518}}
\multiput(50.66,11.36)(.1389,-.03096){10}{\line(1,0){.1389}}
\multiput(52.05,11.05)(.17489,-.03246){8}{\line(1,0){.17489}}
\multiput(53.45,10.79)(.20108,-.02992){7}{\line(1,0){.20108}}
\multiput(54.85,10.58)(.28283,-.03179){5}{\line(1,0){.28283}}
\put(56.27,10.42){\line(1,0){1.419}}
\put(57.69,10.32){\line(1,0){1.422}}
\put(59.11,10.26){\line(1,0){1.423}}
\put(60.53,10.25){\line(1,0){1.422}}
\put(61.95,10.3){\line(1,0){1.42}}
\multiput(63.37,10.39)(.28311,.02923){5}{\line(1,0){.28311}}
\multiput(64.79,10.54)(.2349,.03279){6}{\line(1,0){.2349}}
\multiput(66.2,10.74)(.17518,.03088){8}{\line(1,0){.17518}}
\multiput(67.6,10.98)(.15463,.03301){9}{\line(1,0){.15463}}
\multiput(68.99,11.28)(.12547,.03152){11}{\line(1,0){.12547}}
\multiput(70.37,11.63)(.11391,.03299){12}{\line(1,0){.11391}}
\multiput(71.74,12.02)(.09656,.031752){14}{\line(1,0){.09656}}
\multiput(73.09,12.47)(.089004,.032842){15}{\line(1,0){.089004}}
\multiput(74.43,12.96)(.077446,.031771){17}{\line(1,0){.077446}}
\multiput(75.74,13.5)(.072022,.032605){18}{\line(1,0){.072022}}
\multiput(77.04,14.09)(.067082,.033311){19}{\line(1,0){.067082}}
\multiput(78.31,14.72)(.059576,.032292){21}{\line(1,0){.059576}}
\multiput(79.56,15.4)(.055728,.03284){22}{\line(1,0){.055728}}
\multiput(80.79,16.12)(.052147,.0333){23}{\line(1,0){.052147}}
\multiput(81.99,16.89)(.0488,.03368){24}{\line(1,0){.0488}}
\multiput(83.16,17.69)(.043904,.032682){26}{\line(1,0){.043904}}
\multiput(84.3,18.54)(.041124,.032965){27}{\line(1,0){.041124}}
\multiput(85.41,19.43)(.038492,.033186){28}{\line(1,0){.038492}}
\multiput(86.49,20.36)(.035994,.033352){29}{\line(1,0){.035994}}
\multiput(87.53,21.33)(.033618,.033465){30}{\line(1,0){.033618}}
\multiput(88.54,22.33)(.033516,.035842){29}{\line(0,1){.035842}}
\multiput(89.52,23.37)(.033362,.03834){28}{\line(0,1){.03834}}
\multiput(90.45,24.45)(.033152,.040973){27}{\line(0,1){.040973}}
\multiput(91.34,25.55)(.032883,.043754){26}{\line(0,1){.043754}}
\multiput(92.2,26.69)(.032547,.046699){25}{\line(0,1){.046699}}
\multiput(93.01,27.86)(.033538,.051994){23}{\line(0,1){.051994}}
\multiput(93.78,29.05)(.033094,.055577){22}{\line(0,1){.055577}}
\multiput(94.51,30.28)(.032564,.059427){21}{\line(0,1){.059427}}
\multiput(95.2,31.53)(.033617,.066929){19}{\line(0,1){.066929}}
\multiput(95.84,32.8)(.032934,.071872){18}{\line(0,1){.071872}}
\multiput(96.43,34.09)(.032125,.0773){17}{\line(0,1){.0773}}
\multiput(96.97,35.4)(.033249,.088853){15}{\line(0,1){.088853}}
\multiput(97.47,36.74)(.032193,.096414){14}{\line(0,1){.096414}}
\multiput(97.92,38.09)(.03351,.11376){12}{\line(0,1){.11376}}
\multiput(98.33,39.45)(.03209,.12533){11}{\line(0,1){.12533}}
\multiput(98.68,40.83)(.03371,.15448){9}{\line(0,1){.15448}}
\multiput(98.98,42.22)(.03168,.17504){8}{\line(0,1){.17504}}
\multiput(99.24,43.62)(.02902,.20121){7}{\line(0,1){.20121}}
\multiput(99.44,45.03)(.03053,.28297){5}{\line(0,1){.28297}}
\put(99.59,46.44){\line(0,1){1.419}}
\put(99.69,47.86){\line(0,1){2.136}}
\qbezier(32.75,80)(47.13,80)(60,50)
\qbezier(32.75,20)(47.13,20)(60,50)
\qbezier(21.25,40.75)(60.25,60)(98.25,39.25)
\put(0,0){\bezier{20}(60,50)(64.88,38.88)(72.25,12)}
\put(0,0){\bezier{20}(60,50)(64.88,62.13)(72.25,88)}
\put(48,77){$E_2^{+}$}
\put(73.75,79.5){$E_0$}
\put(65,53.25){$b^*$}
\put(88,47.5){$E_1^{+}$}
\put(73,21.5){$E_1$}
\put(44.5,20){$E_2^{-}$}
\put(34.25,49.5){$E_2$}
\put(18,68){$\partial O_{b^*}$}
\end{picture}}
\caption{Local analysis around $b^*$}
\label{fig:27}
\end{figure}

The fifth step (local parametrices) has a new feature which is not present in the case II. The solutions
of the local matrix \RH\ problems (\ref{eq:4.27}) for the function $U_e$, with $e \in \{a_1,a_2,b\}$, are
the same as for the case II (see (\ref{eq:4.28})). However, around the branch point $b^*$ the
solution is different. Because of Theorem~\ref{thm:2.5}-2, the solution is represented by means of
Airy functions. The local \RH\ problem in the neighborhood of $b^*$ \index{Domains!Obstar@$O_{b^*}$}
is \index{Matrices!Ubstar@$U_{b^*}$}
\begin{equation}  \label{eq:4.40}
  \begin{cases}
   1) & U_{b^*} \in H^{(3\times 3)}(O_{b^*} \setminus \widehat{\Sigma}_{b^*}), \\[10pt]
   2) & U_{b^*+} = U_{b^*-} \widehat{J}, \quad \textrm{on $\widehat{\Sigma}_{b^*}$}, \\[10pt]
   3) & U_{b^*} = \bigl(I + \O(1/n)\bigr) X, \qquad \textrm{uniformly on $\partial O_{b^*}$ as $n \to \infty$},
 \end{cases}
\end{equation}
where (see Figure~\ref{fig:27}) \index{Contours!Sigmahatbstar@$\widehat{\Sigma}_{b^*}$}
\[   \widehat{\Sigma}_{b^*} := O_{b^*} \cap (E_2 \cup E_2^{+} \cup E_2^{-} \cup E_1^{+}), \]
and, see (\ref{eq:4.8}) and (\ref{eq:4.10*}) \index{Matrices!Jhat@$\widehat{J}$}
\begin{equation}  \label{eq:4.41}
   \widehat{J} := \begin{cases}
      D_{2,1}, & \textrm{on $E_2^{\pm} \cap O_{b^*}$}, \\
      W_{1,2}, & \textrm{on $E_2 \cap O_{b^*}$}, \\
      D_{1,2}, & \textrm{on $E_1^{+} \cap O_{b^*}$}.
   \end{cases}
\end{equation}

Observe the block structure of this matrix \RH\ problem, so that we can use the known solution
for the $2\times 2$ case. Following \cite{35,38} we have
\begin{equation}  \label{eq:4.42}
     U_{b^*} = E_{b^*} V_{b^*} A_{b^*},
\end{equation}
where \index{Matrices!Abstar@$A_{b^*}$} \index{Matrices!Ebstar@$E_{b^*}$}
\begin{eqnarray*}
  A_{b^*} &:=& \textrm{diag} \left( 1, \left( \frac{\Phi_1^nw_1}{\Phi_2^nw_2} \right)^{-1/2},
    \left( \frac{\Phi_1^nw_1}{\Phi_2^nw_2} \right)^{1/2} \right), \\
   E_{b^*} &:=& -\sqrt{\pi} e^{i\pi/6} X\ \textrm{diag} \left( 1, \left( \frac{w_1}{w_2} \right)^{-1/2},
     \left( \frac{w_1}{w_2} \right)^{1/2} \right) \\
   & & \times \begin{pmatrix} 2e^{-i\pi/6} & 0 & 0 \\ 0 & 1 & 1 \\ 0 & -i & i \end{pmatrix}
     \textrm{diag} \left( 1 , \left( \frac{3n}2 \varphi_{1,2} \right)^{-1/6},
          \left( \frac{3n}2 \varphi_{1,2} \right)^{1/6} \right),
\end{eqnarray*}
and \index{Functions!phi12@$\varphi_{1,2}$}
\[  \varphi_{1,2} := \log \frac{\Phi_1}{\Phi_2}. \]
To present $V_{b^*}$ we use the notation \index{Matrices!Psibstar@$\Psi_{b^*}$}
\[  \Psi_{b^*} := \begin{pmatrix}
                  1 & 0 & 0 \\
                  0 & \Ai\left( \left(\frac{3n}2 \varphi_{1,2}\right)^{2/3}\right) &
                      \Ai\left( \epsilon_3^2 \left(\frac{3n}2 \varphi_{1,2}\right)^{2/3}\right) \\
                 0 & \Ai'\left( \left(\frac{3n}2 \varphi_{1,2}\right)^{2/3}\right) &
                      \epsilon_3^2\Ai'\left( \epsilon_3^2 \left(\frac{3n}2 \varphi_{1,2}\right)^{2/3}\right)
                   \end{pmatrix} \tilde{\sigma}, \]
\[  \widetilde{\Psi}_{b^*} := \begin{pmatrix}
                  1 & 0 & 0 \\
                  0 & \Ai\left( \left(\frac{3n}2 \varphi_{1,2}\right)^{2/3}\right) &
                       -\epsilon_3^2\Ai\left( \epsilon_3^2 \left(\frac{3n}2 \varphi_{1,2}\right)^{2/3}\right) \\
                 0 & \Ai'\left( \left(\frac{3n}2 \varphi_{1,2}\right)^{2/3}\right) &
                      -\Ai'\left( \epsilon_3^2 \left(\frac{3n}2 \varphi_{1,2}\right)^{2/3}\right)
                   \end{pmatrix} \tilde{\sigma}, \]
where $\Ai$ is the usual Airy function and \index{Matrices!sigmatilde@$\tilde{\sigma}$}
\[  \tilde{\sigma} := \textrm{diag} \left( 1 , e^{i\pi/6} ,  e^{-i\pi/6} \right), \quad
    \epsilon_3 = e^{2\pi i/3}. \]
The matrix $V_{b^*}$ is then given by \index{Matrices!Vbstar@$V_{b^*}$}
\[  V_{b^*} := \begin{cases}
       \Psi_{b^*}, & \textrm{in $S(E_1^{+},E_2^{+})$}, \\
       \Psi_{b^*} \begin{pmatrix} 1 & 0 & 0 \\ 0 & 1 & -1 \\ 0 & 0 & 1 \end{pmatrix}, &
         \textrm{in $S(E_2^{+},E_2)$}, \\
       \widetilde{\Psi}_{b^*}   \begin{pmatrix} 1 & 0 & 0 \\ 0 & 1 & 1 \\ 0 & 0 & 1 \end{pmatrix}, &
         \textrm{in $S(E_2,E_2^{-})$}, \\
       \widetilde{\Psi}_{b^*}, & \textrm{in $S(E_2^{-},E_1^{+})$},
       \end{cases} \]
where $S(\gamma_1,\gamma_2)$ is the domain in $O_{b^*}$ bounded by $\gamma_1$ and $\gamma_2$ (see Figure~\ref{fig:27}).
\index{Domains!Sgamma1gamma2@$S(\gamma_1,\gamma_2)$}

\subsubsection*{Case V}
In this case (in comparison with the case IV) the branch point $b^*$ moves on the $0$-sheet of
$\mathfrak{R}$ (see Figure~\ref{fig:14} and (\ref{eq:2.45})). As a result the arc $E_2$ degenerates and
after the normalization (\ref{eq:4.3}) of the initial \RH\ problem (\ref{eq:2.13}) we have a new jump on
$\delta_1^*$ instead of a jump on $E_2$ (see Remark~\ref{rem:2.3}-4),
\[   Z_+ = Z_{-} D_1^*, \qquad \textrm{on $\delta_1^* \subset \Omega_{1,0,2}$}, \]
where \index{Matrices!D1star@$D_1^*$}
\[   D_1^* := \begin{pmatrix}
               1 & \displaystyle \left( \frac{\Phi_1}{\Phi_0} \right)^n w_1 & 0 \\
               0 & 1 & 0 \\ 0 & 0 & 1
               \end{pmatrix}, \]
which tends to the identity matrix as $n \to \infty$.
In this case we open a global lens around $\Delta_{1,2}$ in the
domain bounded by $E_1^{+}$ and then, in the usual way, we open local lenses around
$\Delta_2$ and $\widetilde{\Delta}_1$. The model \RH\ problem (the parametrix away from the branch points)
for the function $X$ for this case is identical with the corresponding problem of case I. Since
$b^*$ is now on the $0$-sheet of $\mathfrak{R}$, the local \RH\ problem around $b^*$ is modified. Instead
of the jump (\ref{eq:4.41}) in (\ref{eq:4.40}) we have the jump \index{Matrices!Jhat@$\widehat{J}$}
\begin{equation} \label{eq:4.43}
  \widehat{J} = \begin{cases}
       D_1, & \textrm{on $\Delta_1^{*(\pm)} \cap O_{b^*}$}, \\
       W_1 , & \textrm{on $\Delta_1^* \cap O_{b^*}$}, \\
       D_1^*, & \textrm{on $\delta_1^* \cap O_{b^*}$},
       \end{cases}, \quad
      \widehat{\Sigma}_{b^*} := O_{b^*} \cap (\Delta_1^{*(\pm)} \cup \Delta_1^* \cup \delta_1^*).
\end{equation}
Thus \index{Contours!Sigmahatbstar@$\widehat{\Sigma}_{b^*}$}
the Airy function solution given in (\ref{eq:4.42}) for the problem (\ref{eq:4.40}), (\ref{eq:4.43})
should be modified slightly: the non-trivial $2\times 2$ block now is in the left-upper corner
(instead of the right-lower corner in the case IV) and instead of $(\Phi_1,\Phi_2, w_1/w_2)$ we now
have $(\Phi_0,\Phi_1,w_1)$; compare (\ref{eq:4.41}) and (\ref{eq:4.43}). For the solution of an
identical  local $3\times 3$ \RH\ problem, see \cite{49,51,52,48,53,54}.

\subsubsection*{Case III}
The proof for this case is just a repetition of the proof for the case $\textup{V}$, with one
simplification. Since $\Delta_{1,2} = \emptyset$ and $\Delta_1 = \Delta_1^* \cup \delta_1$ for this case,
we do not need to open the global lens. Thus we are in the same situation as in the case I. The
only difference is the solution of the local \RH\ problem around $b^*$, which is the same as for the
case $\textup{V}$.

\section{Conclusion}  \label{sec:5}

In this conclusion we highlight the main results of this paper. They are
\begin{enumerate}
\item The \textbf{classification} of the sets $A := \{ a_1, b_1; a_2,b_2 \}$ such that the limiting
 counting measures for the poles and interpolation points (\ref{eq:1.4}) of the \HP\ approximants
(\ref{eq:1.2})
\[   \pi_n^{(j)} := \frac{Q_n^{(j)}}{P_n}, \quad  R_n^{(j)} := P_n f_j - Q_n^{(j)}, \]
for functions (\ref{eq:2.4})
\[   f_j \in \mathcal{A}(a_j,\alpha_j;b_j,\beta_j;\Omega), \]
 are described by an \textbf{algebraic function} $h$ of order $3$ and \textbf{genus $0$}, see (\ref{eq:1.20}),
 which give rise to measures $\lambda$, $\mu_1$, and $\mu_2$ such that
\begin{equation}  \label{eq:5.1}
  \nu_{P_n} \stackrel{*}{\to} \lambda/2, \quad  \nu_{R_n^{(j)}} \stackrel{*}{\to} \mu_j, \qquad j=1,2.
\end{equation}
\item A universal \textbf{vector-potential equilibrium problem} (\ref{eq:2.70}) for these limiting measures.
\item \textbf{Strong asymptotics} for the corresponding \HP\ approximants.
\end{enumerate}

More precisely
\begin{itemize}
\item In the definition of the \textbf{classes of sets $A$} we used a distinction in the formation of the system
of curves \index{Contours!Gamma@$\Gamma$}
\begin{equation}  \label{eq:5.2}
  \Gamma := \{ z \in \mathbb C : \Re \int h_j(z)\, dz = \Re \int h_k(z)\, dz, \ j \neq k, j,k=0,1,2\}
\end{equation}
with a certain normalization of the primitives. Since the genus of $h$ is $0$ the set $\Gamma$
admits an algebraic parametrization. It is formed by the trajectories
\[  \Gamma = \bigcup_{\ell=1}^6 \gamma_\ell, \quad
    \gamma_\ell := z_\ell(\eta), \qquad \ell=1,\ldots,6,\ \eta \in [-2,2], \]
given by the branches of an algebraic function $z(\eta)$ of order six, when $\eta$ runs from $2$ to $-2$.
The peculiarities of the behavior of these trajectories define different classes of $A$. Thus,
starting from the input data $A$ (i.e., the branch points $a_1,b_1,a_2,b_2$) we have a finite number
(because the genus of $h$ is $0$) of algebraic functions $h$ satisfying
\begin{equation}  \label{eq:5.3}
   h^3 - 3 \frac{P_2(z)}{\Pi_4(z)} h + 2 \frac{P_1(z)}{\Pi_4(z)} = 0,
\end{equation}
with $\Pi_4(z) = (z-a_1)(z-b_1)(z-a_2)(z-b_2)$. From the coefficients of the equation (\ref{eq:5.3}) we get
explicit expressions of the coefficients of the equation for the function $z(\eta)$. Then, observing the
behavior of the branches $z(\eta)$ when $\eta$ runs from $2$ to $-2$, we can conclude to which class $A$ belongs.
Depending on the class we define (globally in $\mathbb{C}$) the branches of the algebraic function
$h := \{h_0,h_1,h_2\}$ and the algebraic function
\[  \Phi := \exp \left( \int h(z)\, dz \right) = \{ \Phi_0,\Phi_1,\Phi_2\}, \]
i.e., we define the sheet structure of the Riemann surface $\mathfrak{R}$ of the function $h$.
Finally, the jumps of the function $h$ on certain parts (depending on the geometrical case) of the contour
$\Gamma$ given in (\ref{eq:5.2}), give the densities of the limiting measures $\lambda$ and
$\mu_1,\mu_2$ in (\ref{eq:5.1}). In particular
\[  d\lambda(\xi) = \frac{1}{2\pi i} \left( h_{0+}(\xi) - h_{0-}(\xi) \right)\, d\xi,
  \qquad \xi \in \Delta_0 \subset \Gamma, \]
where $\Delta_0$ is a union of cuts which form the boundary of the domain of analyticity of the
branch $h_0$, i.e., $\Delta_0$ are the cuts of the $0$-sheet $\mathfrak{R}_0$ of the Riemann surface
$\mathfrak{R}$.
\item The \textbf{vector-potential equilibrium problem} (Theorem~\ref{thm:2.9}) does not depend on the geometrical
class of $A$ (universality). The equilibrium relations of this problem are considered on cuts $\Delta_1$
and $\Delta_2$ joining the pairs $(a_j,b_j)$, $j=1,2$. In addition, if $\Delta_1 \cap \Delta_2
= \Delta_{1,2} \neq \emptyset$, an extra equilibrium relation is imposed on a curve $E$ containing $\Delta_{1,2}$,
see (\ref{eq:2.70}). This equilibrium problem reduces to the known
vector-potential problem for an Angelesco system when $\Delta_{1,2} = \emptyset$. The cut $\Delta_j$
makes the function $f_j$ holomorphic. Thus the sets $\Delta_1, \Delta_2, \Delta_{1,2}, E$ are natural
input data for the vector-potential problem and we expect that the limiting measures (\ref{eq:5.1}) of the
\HP\ approximants hold for a wider class of functions. We also recall that we impose in this paper an extra
analyticity condition on $(f_1,f_2)$ when $\Delta_{1,2} \neq \emptyset$: we require that the ratio
of the jumps of $f_1$ and $f_2$ on $\Delta_{1,2}$
\begin{equation} \label{eq:5.4}
   u(\xi) := \frac{f_{1+}-f_{1-}}{f_{2+}-f_{2-}}(\xi), \qquad \xi \in \Delta_{1,2}
\end{equation}
has a holomorphic (meromorphic) continuation from $\Delta_{1,2}$ to the domain $\overline{G}$
\begin{equation}  \label{eq:5.5}
     u \in H(\overline{G}), \quad \partial G = E.
\end{equation}
The analyticity condition (\ref{eq:5.4}) gives a link of our vector-potential problem
(\ref{eq:2.70}) with the known equilibrium problem for a Nikishin system.
\item For the derivation of the \textbf{strong asymptotics} for the multiple orthogonal
polynomials $P_n$ and for their functions of the second kind $R_n^{(j)}$ $(j=1,2)$, we use
a $3\times 3$ matrix-valued \RH\ problem. The increase of the order of the matrix functions
(in comparison with $2\times 2$ matrix-valued \RH\ problems for the usual orthogonal polynomials)
brings new features into the standard \RH\ technology. One of these features is a new
decomposition of the matrix jump on $\Delta_{1,2} \neq \emptyset$ which implies
the opening of a \textit{global lens} containing the domain $G$, see (\ref{eq:5.5}).
This procedure introduces a new effect (in comparison with usual orthogonal polynomials)
of oscillatory asymptotics on some curves
$\subset E$ in $\overline{\mathbb{C}} \setminus (\Delta_1 \cup \Delta_2)$ for the functions
of the second kind $R_n^{(j)}$ and, as a result, there is an accumulation of zeros of
$R_n^{(j)}$ on $E$.
\end{itemize}

\section*{Acknowledgments}
\addcontentsline{toc}{section}{Acknowledgments}
The authors would like to express their gratitude to B.~Beckermann, P.~Deift,
 A.~ Gonchar, V.~Lysov, H.~Stahl, and D.~Tulyakov
for helpful discussions and advise during the work on this paper.

This research was supported by
INTAS Research Network 03-51-6637. The first author was also supported by
Program No.~1 DMS RAS, grants RFBR-05-01-00522 and NSh-1551.2003.1.
The second and third author were also supported
by OT/04/21 of the Research Council of K.U.Leuven, by FWO-Flanders Project G.0455.04,
and by the Belgian Interuniversity Attraction Pole P06/02. The second author was
also supported by the European Science Foundation programme MISGAM and by a grant from
the Ministry of Education and Science of Spain, project code MTM2005-08648-C02-01.

\newpage
\addcontentsline{toc}{section}{Index of notation}
\printindex

\end{document}